\documentclass[a4paper, 10pt ]{article}

\usepackage[a4paper,left=1.25cm,right=1.25cm,top=2.5cm,bottom=2.5cm]{geometry}
\usepackage{hyperref}
\hypersetup{colorlinks,linkcolor={blue},citecolor={blue},urlcolor={black}}
\usepackage{subfig}
\usepackage{pgfplots}
\usepackage{pgfplotstable}
\usepackage{mathtools}
\usepackage{multicol}
\usepackage{comment}
\usepackage{booktabs}
\pgfplotsset{compat=1.5}
\usepackage{amssymb}
\usepackage{url}
\usepackage{bm}

\usepackage{pdfsync}
\usepackage{float}
\usepackage{tabularx}
\usepackage{enumerate}
\usepackage{array}
\usepackage{xspace}
\usepackage{tikz}
\usepackage{tikz-cd}
\usepackage{tikzsymbols}
\usetikzlibrary{calc,trees,positioning,arrows,chains,shapes.geometric,%
    decorations.pathreplacing,decorations.pathmorphing,shapes,%
    matrix,shapes.symbols, decorations.markings, patterns,fit}

\usepackage{siunitx}

\usepackage{authblk}

\usepackage[draft,inline,marginclue]{fixme}
\usepackage{mathrsfs}
\usepackage{color, colortbl}
\usepackage{multirow}

\usepackage{amsthm}
\usepackage{amsmath}
\usepackage{stmaryrd}
\usepackage[ruled,vlined,linesnumbered]{algorithm2e}
\usepackage{graphicx}
\usepackage{booktabs}
\usepackage{cleveref}
\usepackage{arydshln}
\usepackage[inkscapeformat=png]{svg}
\FXRegisterAuthor{dt}{adt}{\color{red}DT}
\FXRegisterAuthor{fr}{afr}{\color{blue}FR}

\theoremstyle{remark}

\newcommand{\x}{\ensuremath{\mathbf{x}}}

\DeclareMathOperator*{\argmin}{\arg\!\min}
\DeclareMathOperator*{\argmax}{\arg\!\max}

\newcolumntype{C}[1]{>{\centering\arraybackslash}m{#1}}

\definecolor{Gray}{gray}{0.9}

\begin{document}

\title{Explicable hyper-reduced order models on nonlinearly
approximated solution manifolds of compressible and incompressible Navier-Stokes equations}

\author[1]{Francesco~Romor\footnote{francesco.romor@sissa.it}}
\author[2]{Giovanni~Stabile\footnote{giovanni.stabile@uniurb.it}}
\author[1]{Gianluigi~Rozza\footnote{gianluigi.rozza@sissa.it}}

\affil[1] {Mathematics Area, mathLab, SISSA, via Bonomea 265, I-34136 Trieste,
  Italy
}
\affil[2] {Department of Pure and Applied Sciences, Informatics and Mathematics Section, University of Urbino Carlo Bo, Piazza della Repubblica 13, I-61029, Urbino, Italy}
\date{}
\maketitle

\begin{abstract}
A slow decaying Kolmogorov n-width of the solution manifold of a parametric partial differential equation precludes the realization of efficient linear projection-based reduced-order models. This is due to the high dimensionality of the reduced space needed to approximate with sufficient accuracy the solution manifold. To solve this problem, neural networks, in the form of different architectures, have been employed to build accurate nonlinear regressions of the solution manifolds. However, the majority of the implementations are non-intrusive black-box surrogate models, and only a part of them perform dimension reduction from the number of degrees of freedom of the discretized parametric models to a latent dimension. We present a new intrusive and explicable methodology for reduced-order modelling that employs neural networks for solution manifold approximation but that does not discard the physical and numerical models underneath in the predictive/online stage. We will focus on autoencoders used to compress further the dimensionality of linear approximants of solution manifolds, achieving in the end a nonlinear dimension reduction. After having obtained an accurate nonlinear approximant, we seek for the solutions on the latent manifold with the residual-based nonlinear least-squares Petrov-Galerkin method, opportunely hyper-reduced in order to be independent from the number of degrees of freedom. New adaptive hyper-reduction strategies are developed along with the employment of local nonlinear approximants. We test our methodology on two nonlinear time dependent parametric benchmarks involving a supersonic flow past a NACA airfoil with changing Mach number and an incompressible turbulent flow around the Ahmed body with changing slant angle.
\end{abstract}

\tableofcontents
\listoffixmes

\section{Introduction}
\label{sec:intro}
Real world numerical models, coming from systems of partial differential equations (PDEs), usually study a physical phenomenon under the influence of different parameters. For each parametric instance, a single numerical simulation could take from hours to weeks to complete. Such is the case for complex fluid dynamics models or large-scale geophysical simulations. Fortunately, in some cases, the outputs of these models show evident correlations among them, partially because they follow the same physical laws embedded in the same numerical models, and partially because the parameters' dependency affects the solutions only as relatively small perturbations. Reduced-order modelling (ROM) leverages these correlations among snapshots, i.e. single solutions corresponding to different parametric instances, to reduce the computational time. The most successful model order reduction (MOR) methods combine the knowledge from the physical and the numerical models with the information coming from a database of solutions. One of the most employed methods is the reduced basis method~\cite{hesthaven2016certified,rozza2022advanced}. As most numerical models search for the solutions on discrete finite dimensional vector spaces like the finite volumes method (FVM), the finite element method (FEM), the spectral element method (SEM) and the discontinuous Galerkin method (DGM), model order reduction exploits the prior information coming from a dataset of training snapshots to update these ansatz spaces. The results are very low-dimensional linear vector spaces for which seeking the solutions associated to new parametric instances is more efficient if these solutions are expected to be correlated with the training dataset. Fundamentally, ROMs amortize the cost of computing an initial training database of solutions and low-dimensional adapted ansatz spaces in the offline stage, through subsequent efficient evaluations of unseen solutions in the online stage. It is important to remark that the numerical models employed in the offline stage are still employed also in the online stage, so that the reduced solutions are discerned in the ansatz spaces through the satisfaction of the physical principles and mathematical constraints underneath the original numerical models.

Some difficulties arise when the solution manifold, that is the space of parameter dependent solutions, cannot be approximated with a satisfactory accuracy by linear low-dimensional spaces. If we consider a parameter space $\mathcal{P}\subset\mathbb{R}^p$, $p>0$, and a solution map $\mathbf{U}:\mathcal{P}\subset\mathbb{R}^p\rightarrow X_h\sim\mathbb{R}^d$ that associates to each parameter $\boldsymbol{\mu}\in\mathcal{P}$ the corresponding solution $\mathbf{U}(\boldsymbol{\mu})\in X_h\sim\mathbb{R}^d$ in the discretization space of choice $X_h$, where $d>0$ is the number of degrees of freedom, we can quantify the linear approximability of the solution manifold $X_h\supseteq\mathcal{M}=\mathbf{U}(\mathcal{P})$ with the Kolmogorov n-width (KnW):
\begin{equation}
  d_n(\mathcal{M}, X_h)=\inf_{\substack{W\subset\mathbb{R}^d\\\text{dim}W=n}}\sup_{\boldsymbol{\mu}\in\mathcal{P}}\inf_{\mathbf{V}\in X_h} \lVert \mathbf{V}-\mathbf{U}(\boldsymbol{\mu})\rVert_2.
\end{equation}
A slow decaying Kolmogorov n-width with respect to the dimension of the linear approximant precludes the realization of efficient ROMs. One of the most prominent defect of linear ROMs is that even simple physical models, like linear advection, suffer from a slow Kolmogorov n-width decay. These are cases for which the snapshots are poorly correlated, sometimes almost orthogonal in $X_h$.

Recently, with the diffusion of scientific machine learning, black-box surrogate models have tackled slow-decaying KnW solution manifolds thanks to nonlinear approximants represented by neural networks (NNs), in the form of different combined architectures. The majority of these surrogate models, being non-intrusive, do not even perform dimension reduction from the space of degrees of freedom $\mathbb{R}^d$ to a reduced or latent space $\mathbb{R}^r$, $r\ll d$. While a low-dimensional space is needed for linear projection-based ROMs to seek for the solutions efficiently in the online stage, surrogate models built with NNs rely on the fast evaluation of the nonlinear approximants for different inputs in the prediction phase. Apart from the imposition of additional inductive biases, the predicted solutions do not consider the physical and mathematical knowledge of the models under study: in fact, they are not obtained from the satisfaction of first principles like classical ROMs. Moreover, when dimension reduction is performed, it is essentially needed for features extraction rather than to increase the efficiency of the surrogate models. Nonetheless, when architectures like autoencoders (AE) are employed, the approximation error of the solution manifolds decays more rapidly with respect to the latent dimension when compared to linear subspaces. This can be quantified with an extension of the definition of KnW $\delta_n$ for continuous maps:
\begin{equation}
  \delta_{n}(\mathcal{M}, X_h) = \inf_{\substack{\psi\in\mathcal{C}(X_h, \mathbb{R}^{r})\\\phi\in\mathcal{C}(\mathbb{R}^{r}, X_h)}}\sup_{\boldsymbol{\mu}\in\mathcal{P}}\lVert \mathbf{U}(\boldsymbol{\mu})-(\phi\circ\psi)(\mathbf{U}(\boldsymbol{\mu}))\rVert_{2},
\end{equation}
where $\psi\in\mathcal{C}(X_h, \mathbb{R}^{r})$ and $\phi\in\mathcal{C}(\mathbb{R}^{r}, X_h)$ are continuous maps represented in our case by NNs. This enables the design of efficient intrusive ROMs even for models with a slow KnW decay.

The first employment of convolutional autoencoders for intrusive ROMs, namely Galerkin and least-squares Petrov-Galerkin nonlinear manifold methods appears in~\cite{lee2020model}. The major evident drawback is that both the architecture and the numerical schemes employed in the online predictive phase depend on the number of degrees of freedom (dofs), so the procedure itself is even slower than the full-order models. Typically, when performing MOR of nonlinear parametric PDEs hyper-reduction is employed to achieve independence with respect the number of dofs. In this case, another ingredient complicates the matter since the nonlinearity coming from the decoder map needs also to be treated and made independent on the number of dofs. One of the first approaches in this direction is introduced in~\cite{kim2022fast}. The architecture employed is a shallow masked autoencoder: the sparsity pattern imposed on the last decoder layers reduces the computational costs of the forward and Jacobian evaluations, while Gauss-Newton with approximated tensors (GNAT) is employed to hyper-reduce the residual. One problem that arises is that shallow autoencoders are sometimes not enough to accurately approximate complex solution manifolds and the methodology itself constraints the choice of architecture. The methodology was tested on the 2d Burgers equations solved with finite differences. Another strategy~\cite{romor2023non} uses teacher-student training to compress a generic architecture that performs dimension reduction, in this case a convolutional AE, onto a small feedforward NN. Possible combinations of hyper-reduction with reduced over-collocation~\cite{chen2021eim} only on the residual or for both the decoder and the residuals were taken into considerations. The methodology was tested on a 2d nonlinear conservation law test case and a 2d shallow water equations benchmark solved with \texttt{OpenFoam}~\cite{marcantoni2012high}. Afterwards, it is introduced a new implementation~\cite{barnett2022neural} that considers as nonlinear approximant of the solution manifold the sum of a linear subspace and a linear closure term whose coefficients are the output of a feedforward NN. The hyper-reduction method is the energy-conserving sampling and weighting method (ECSW) and it was tested on a 2d Burgers' equations model solved with finite differences. Another approach~\cite{cocola2023hyper}, directly employs a relatively small decoder from the latent space to the submesh identified by the reduced over-collocation hyper-reduction method as nonlinear approximant of the solution manifold. In this way, the training phase is more efficient and the solutions are finally reconstructed with the hyper-reduction linear basis from the collocation nodes. To increase the accuracy of shallow masked AE, in~\cite{diaz2023fast} they implement domain decomposition and build a local shallow masked AE for each subdomain. The procedure is tested on a 2d Burgers' equation model.

In this work, we introduce a new methodology and we test it on more challenging benchmarks than 2d Burgers' equations. For a moderately slow KnW decay, classical ROMs fail, but linear subspaces can still be employed as good approximants of the solution manifold. This is the rationale behind employing singular value decomposition modes (SVD)~\cite{fresca2022pod}, or other linear transforms or filters, as preprocessing step to dimension reduction with AE. In section~\ref{sec:results}, we will also show a case for which this assumption is not valid anymore and more deep NN architectures should be employed and reduced with teacher-student training following~\cite{romor2023non}. The novelties of our new approach are the following
\begin{itemize}
  \item a new collocated hyper-reduction procedure specific for our nonlinear manifold approximant that combines randomomized singular value decomposition modes with 1d convolutional autoencoders, in section~\ref{subsec:ahr}
  \item an adaptive gradient-based hyper-reduction strategy, in section~\ref{subsec:adaptive-hr}. Similar concepts of adaptation strategies are present in the literature~\cite{peherstorfer2020model, wen2023coupling}.
  \item the implementation of an efficient way to integrate local nonlinear manifolds with intrusive nonlinear least-squares Petrov-Galerkin through a local change of basis, in section~\ref{sec:loc}
  \item the validation of our methodology on challenging test cases with moderately slow Kolmogorov n-width, in section~\ref{sec:results}
\end{itemize}

In summary, in section~\ref{sec:nm} the nonlinear least-squares Petrov-Galerkin method is presented, our nonlinear manifold approximant is introduced and some specifications regarding the normalization of the datasets and the evaluation of the randomized singular value decomposition (rSVD) modes are made. Then, in section~\ref{sec:hr} the employed hyper-reduction methods are introduced. Particular attention is focused on the reduced over-collocation method and its new adaptive formulation in subsection~\ref{subsec:ahr}. A brief section~\ref{sec:loc} introduces a straight-forward way to include local nonlinear manifolds in the methodology through a linear change of rSVD basis. Finally, two benchmarks are introduced in section~\ref{sec:results}. A 2d nonlinear parametric time-dependent supersonic compressible Navier-Stokes equations model (CNS) is studied on a coarse and a finer mesh. A special focus is given to the comparison of different hyper-reduction techniques in subsection~\ref{subsec:coarseAirfoil} and to the implementation of local nonlinear manifolds in subsection~\ref{subsubsec:sudomains_finer_airfoil}. The last subsections involve the study of a 3d nonlinear time-dependent geometrically parametrized turbulent incompressible Navier-Stokes equations model (INS).

\section{Residual-based ROMs on nonlinear manifolds}
\label{sec:nm}
The starting point is a parametric time-dependent partial differential equation (PDE) on a computational domain $\Omega\subset\mathbb{R}^D$, $D=\{1,2,3\}$, with time interval $[0, T_{\boldsymbol{\mu}}]$ and parameter space $\boldsymbol{\mu}\in\mathcal{P}\subset\mathbb{R}^p$:
\begin{subequations}
    \begin{align}
        G(\mathbf{U}(\x, t), \partial_t\mathbf{U}(\x, t), \nabla \mathbf{U}(\x, t),\dots;\boldsymbol{\mu}) &= 0,\quad &(\x, t)\in\Omega_{\boldsymbol{\mu}}\times[0, T_{\boldsymbol{\mu}}],\\
        B(\mathbf{U}(\x, t); \boldsymbol{\mu}) &= 0,\quad &(\x, t)\in\partial\Omega_{\boldsymbol{\mu}}\times[0, T_{\boldsymbol{\mu}}],\\
        \mathbf{U}(\x, t) &= \mathbf{U}_{0, \boldsymbol{\mu}}(\x),\quad &(\x, t)\in\Omega_{\boldsymbol{\mu}}\times\{0\},
    \end{align}
\end{subequations}
where that state function $\mathbf{U}:\Omega_{\boldsymbol{\mu}}\times [0, T_{\boldsymbol{\mu}}]\subset\mathbb{R}^D\times\mathbb{R}\rightarrow\mathbb{R}^s$, $s\geq 1$ belongs to a Banach space $X(\mathbb{R}^D\times\mathbb{R}; \mathbb{R}^s)$ for all $\boldsymbol{\mu}\in\mathcal{P}^q$, $q\geq 1$, of vector-valued time-dependent functions. The function $\mathbf{U}_{0, \boldsymbol{\mu}}(\x)$ represents the possibly parametric dependent initial condition. The state function $\mathbf{U}$ is a synthetic notation that can include at the same time more than one physical field, like velocity, pressure, internal energy, and density for example. The function $G$ represents the PDE itself and has as arguments the state function and its partial derivatives with respect to time and space. We do not restrict only to first order, higher derivatives are omitted. The boundary conditions are expressed through the operator $B$ and are possibly parametric dependent. These definitions are introduced only to define the discretized systems we will work with. We have included also geometric parametrizations $\Omega_{\boldsymbol{\mu}}\subset\mathbb{R}^D$ through $\boldsymbol{\mu}$.

We will consider nonlinear time-dependent PDEs, but in general the framework we are going to introduce can be applied also to stationary PDEs and linear PDEs. The only requirement is the slow Kolmogorov n-width decay of the solution manifold, otherwise, it would be sufficient to apply the well-developed theory of linear projection-based ROMS. In fact, employing a nonlinear approximation of the solution manifold reduces the efficiency of linearly approximable solution manifolds, in general.

To be as general as possible, we will consider a generic discretization in space, with the constraints that it is supported on a computational mesh $\Omega_{\boldsymbol{\mu}, h}\subset\mathbb{R}^D$ and the discretized differential operators have local stencils in order to implement efficiently hyper-reduction schemes later~\cite{hesthaven2016certified}. So, in our framework, we include the Finite Volume Method (FVM), that we are employing, but also the Finite Element Method (FEM) and the Discontinuous Galerkin method (DGM), for example. Applying the method of lines, we discretize in space to obtain the following ordinary differential equation
\begin{subequations}
    \begin{align}
        G_h(\mathbf{U}_h(t), \partial_t\mathbf{U}_h(t), \nabla \mathbf{U}_h(t),\dots;\boldsymbol{\mu}) &= \mathbf{0},\quad t\in[0, T_{\boldsymbol{\mu}}],\\
        B_h(P_{\partial\Omega_{\boldsymbol{\mu}}}(\mathbf{U}_h(t)); \boldsymbol{\mu}) &= 0,\quad t\in[0, T_{\boldsymbol{\mu}}],\\
        \mathbf{U}_h(0) &= \mathbf{U}_{h,0},
    \end{align}
\end{subequations}
where in this case the discrete state function $\mathbf{U}_h(t)\in X_h(\Omega_{\boldsymbol{\mu}, h})$ belongs for all $t\in[0, T_{\boldsymbol{\mu}}]$ to a discretization space $X_h(\Omega_{\boldsymbol{\mu}, h})\sim\mathbb{R}^d$ where $d>0$ is the number of degrees of freedom, with its norm $\lVert\cdot\rVert_{X_h}$. The map $P_{\partial\Omega_{\boldsymbol{\mu}}}:X_h(\Omega_{\boldsymbol{\mu}, h})\rightarrow X_h(\partial\Omega_{\boldsymbol{\mu}, h})$ is the projection onto the discrete boundary $\partial\Omega_{\boldsymbol{\mu}, h}$ of the computational domain $\Omega_{\boldsymbol{\mu}, h}$.

Finally, we apply a discretization in time to obtain the discrete residual $G_{h,\delta t}:\mathcal{P}\times X_h\times X_h^{\mid I_t\mid}\rightarrow X_h$ at time $t\in\{t_0,\dots,t_{N_{\boldsymbol{\mu}}}\}=V_{\boldsymbol{\mu}}$
\begin{subequations}
\begin{align}
    G_{h, \delta t}(\boldsymbol{\mu}, \mathbf{U}_h^t, \{\mathbf{U}_h^s\}_{s\in I_t})&=\mathbf{0},\\
    B_{h, \delta t}(P_{\partial\Omega_{\boldsymbol{\mu}}}(\mathbf{U}_h^t); \boldsymbol{\mu}) &= 0,\\
        \mathbf{U}_h^0 &= \mathbf{U}_{h,0},
    \label{eq:PDEsystemDiscreteU}
\end{align}
\end{subequations}
where the time instances $t_i\in[0, T_{\boldsymbol{\mu}}],\ \forall i\in\{0,1,\dots, N_{\boldsymbol{\mu}}\}$ and $I_t$ is the set of previous time instances of the state variable $U_h^t\in X_h$ at time $t$, needed for the numerical time scheme of choice. For most cases of model order reduction, it is crucial that the discretization space $X_h\sim\mathbb{R}^d$ is not time or parametric dependent. With adaptive collocated hyper-reduction~\ref{subsec:adaptive-hr} these constraints can be relaxed.

\subsection{Nonlinear least-squares Petrov-Galerkin method}
\label{subsec:nlspg}
\begin{figure}
    \centering
    \includegraphics[width=0.49\textwidth]{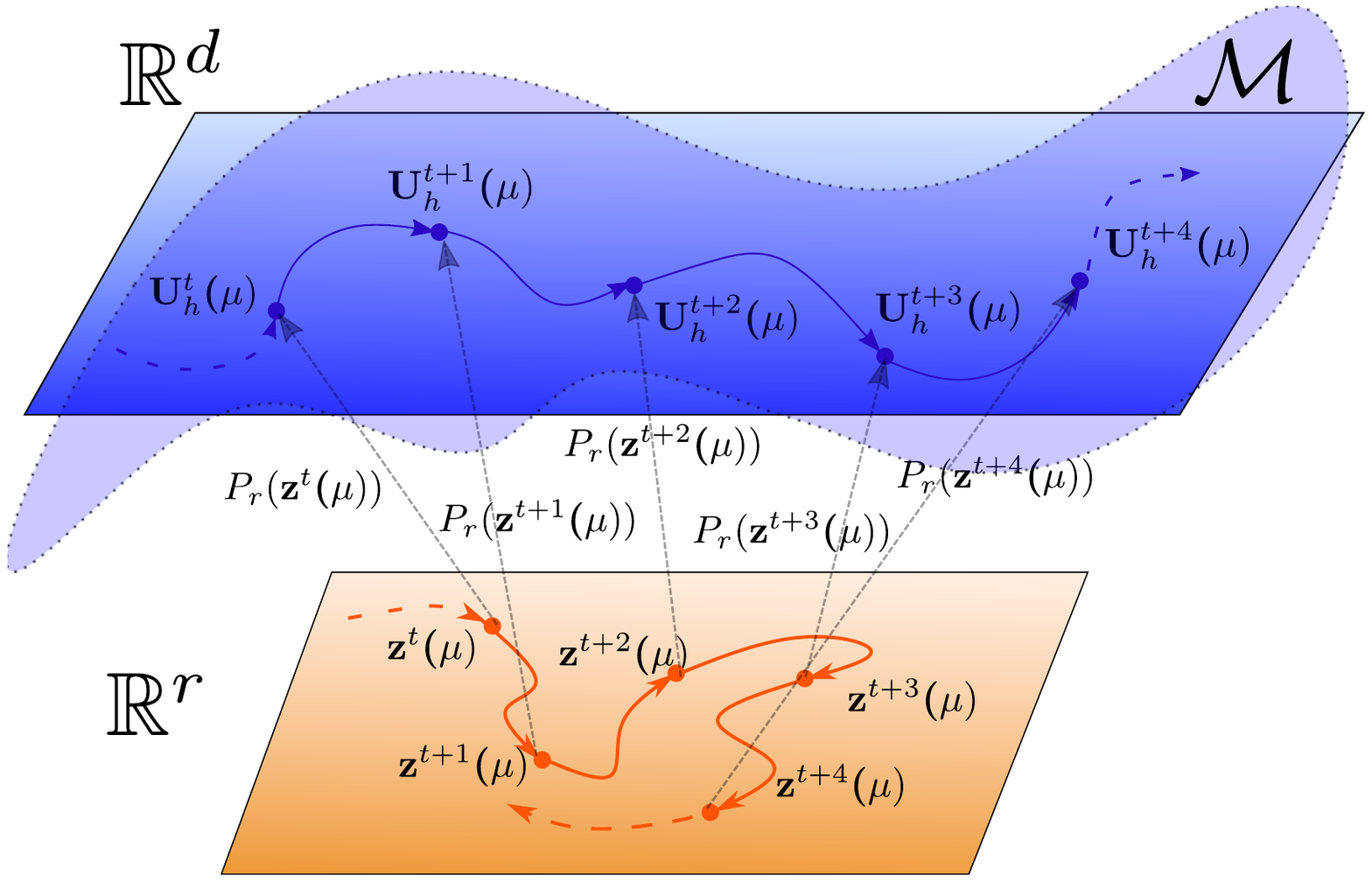}
    \includegraphics[width=0.49\textwidth]{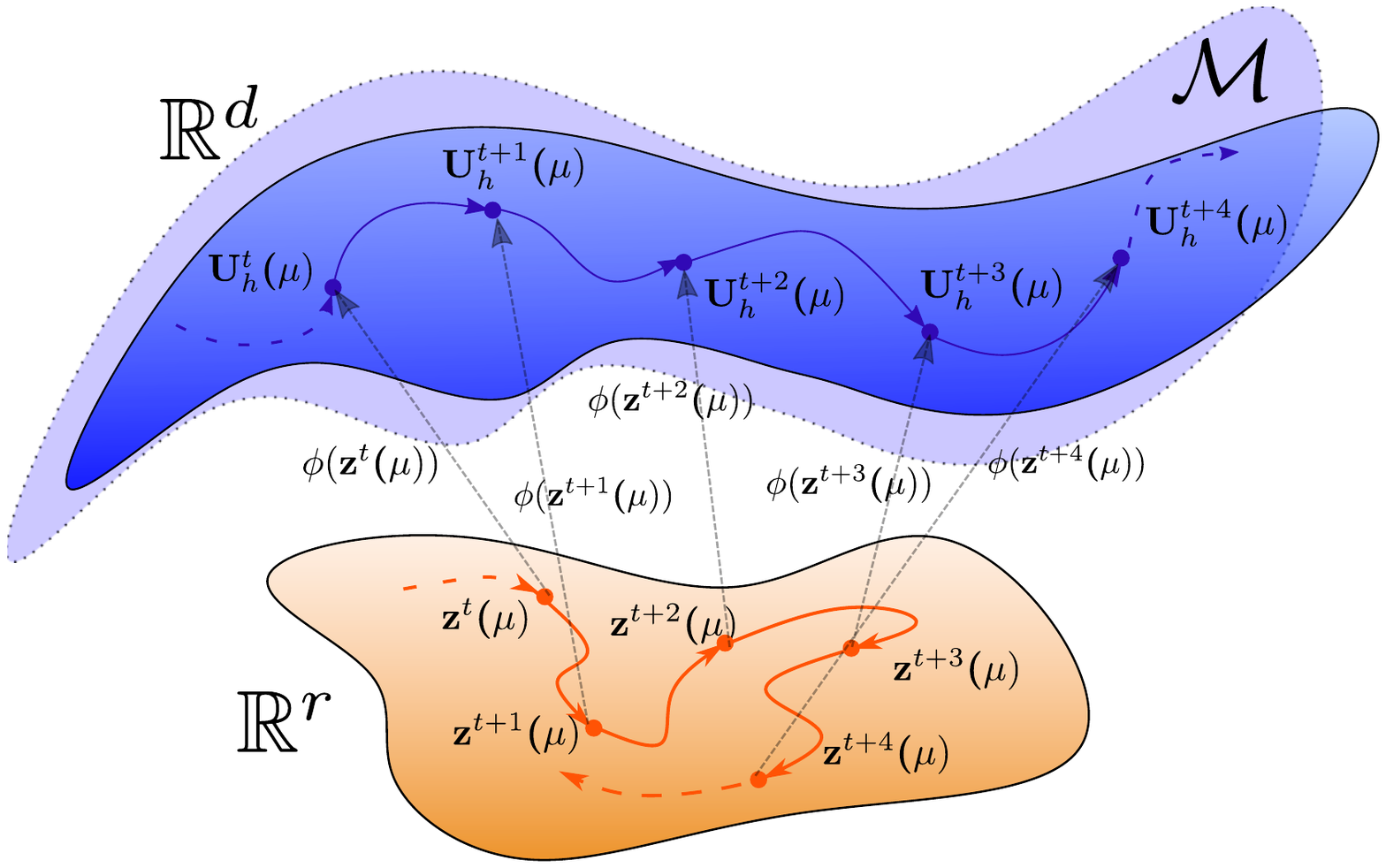}
    \caption{\textbf{Left:} evolution trajectories on the $r$-dimensional linear latent space and solution manifold embedded in the ambient space $\mathbb{R}^d$. The map $P_r$ is a linear projection. \textbf{Right:} evolution trajectories on the $r$-dimensional latent space and nonlinear solution manifold in the ambient space $\mathbb{R}^d$. The map $\phi$ is a single chart nonlinear parametrization of the approximating solution manifold.}
    \label{fig:linear_nonlinear_manifolds}
\end{figure}

We will introduce the nonlinear manifold least-squares Petrov-Galerkin method (NM-LSPG)~\cite{lee2020model} from its linear manifold version (LM-LSPG), see Figure~\ref{fig:linear_nonlinear_manifolds}.
To perform model order reduction with LM-LSPG, we need to define a linear projection map $P_r:\mathbb{R}^r\rightarrow X_h\sim\mathbb{R}^d$ from the reduced space $\mathbb{R}^r$, $r\ll d$ to the full-state space $X_h\sim\mathbb{R}^d$. We require that, fixed a tolerance $\epsilon \ll 1$, the relative reconstruction error of the linear solution manifold $\mathcal{M}=\{\mathbf{U}_h\in X_h \vert \exists\boldsymbol{\mu}\in\mathcal{P}:\ g(\boldsymbol{\mu}) = \mathbf{U}_h\}$ is small:
\begin{equation}
    \frac{\lVert \mathbf{U}_h-P_r P_r^T \mathbf{U}_h\rVert_2}{\lVert \mathbf{U}_h\rVert_2}<\epsilon,\quad \forall \mathbf{U}_h\in\mathcal{M}.
\end{equation}
Typically, this is achieved by sampling from the parameter space $\mathcal{P}$ a set of independent training parameters $\mathcal{P}_{\text{train}}\subset\mathcal{P}$. The corresponding training snapshots $\mathcal{U}_{\text{train}}=\{\mathbf{U}_{h, \boldsymbol{\mu}}\}_{\boldsymbol{\mu}\in\mathcal{P}_{\text{train}}}$ are employed to evaluate $P_r$. For every new parametric instance $\mu\in\mathcal{P}\backslash \mathcal{P}_{\text{train}}$, we can evaluate the corresponding solution solving the following nonlinear least-squares problem at each time step $t\in\{t_0,\dots,t_{N_{\boldsymbol{\mu}}}\}=V_{\boldsymbol{\mu}}$ in the reduced variables $\mathbf{z}\in \mathbb{R}^r$ and with initial condition $\mathbf{z}^0$
\begin{subequations}
    \begin{align}
        \label{eq:PDEsystemDiscreteLinear}
        \mathbf{z}^t &= \argmin_{\mathbf{z}\in \mathbb{R}^r}\ \lVert G_{h, \delta t}(\boldsymbol{\mu}, P_r(\mathbf{z}), \{P_r(\mathbf{z}^s)\}_{s\in I_t})\rVert^{2}_{X_h},\\
        \mathbf{z}^0 &= P_r^T(\mathbf{U}_{h, 0}),
    \end{align}
\end{subequations}
where $\mathbf{z}^s\in\mathbb{R}^r,\ \forall s\in I_t$ are the previous reduced coordinates needed at time $t$ by the numerical scheme. The nonlinear least-squares problem can be solved with optimization methods like Gauss-Newton with line-search~\cite{lee2020model}, Levenberg-Marquardt~\cite{romor2023non} and derivative-free Pounders~\cite{wildsolving} implemented in PETSc~\cite{balay2019petsc}, that we will employ.

If $G_{h, \delta t}$ is linear, then we can solve for \eqref{eq:PDEsystemDiscreteLinear} without reconstructing the solution onto the full-state space $\mathbb{R}^d$. If $G_{h, \delta t}$ is nonlinear, hyper-reduction techniques must be introduced in order to recover the independence from the number of degrees of freedom.

The evolution of the reduced trajectory in the latent space $\mathbb{R}^r$ in Figure~\ref{fig:linear_nonlinear_manifolds} is computed without reconstructing the full-states $U_h^t(\boldsymbol{\mu}) = P_r(\mathbf{z}(\boldsymbol{\mu}))$. The reconstruction from $\mathbb{R}^r$ to the ambient space $\mathbb{R}^d$ is performed only at the end with the projection map $P_r$.

The nonlinear counterpart of LM-LSPG, poses the approximability of the solution manifold $\mathcal{M}$ with a nonlinear manifold. We will define this approximating nonlinear manifold as the image of a nonlinear parametrization map $\phi:U\subset\mathbb{R}^r\rightarrow\mathbb{R}^d$, that is with a single chart, with $U$ an open subset of $\mathbb{R}^r$. For our purposes, mainly linked to the definition of the initial conditions, we also need an approximation of the right-inverse of $\phi$, that is $\psi:\phi(U)\rightarrow\mathbb{R}^r$ such that $\phi\circ\psi\approx I_d$.

There are many definitions that extend the notion of Kolmogorov n-width to nonlinear approximating spaces~\cite{devore1989optimal}. Our requirement is that the relative reconstruction error is below a fixed tolerance $\epsilon\ll 1$: 
\begin{equation}
    \label{eq:nonlinear_rec_err}
    \frac{\lVert \mathbf{U}_h-(\phi\circ\psi)(\mathbf{U}_h)\rVert_2}{\lVert \mathbf{U}_h\rVert_2}<\epsilon,\quad \forall \mathbf{U}_h\in\mathcal{M}.
\end{equation}

The nonliear least squares problem solved for each time instance $t\in\{t_0,\dots,t_{N_{\boldsymbol{\mu}}}\}=V_{\boldsymbol{\mu}}$ is similar to the linear case with the substitution of the linear projection map $P_r$ with $\phi$ and $\psi$:
\begin{subequations}
    \begin{align}
        \label{eq:PDEsystemDiscrete}
    \mathbf{z}^t &= \argmin_{\mathbf{z}\in \mathbb{R}^r}\ \lVert G_{h, \delta t}(\boldsymbol{\mu}, \phi(\mathbf{z}), \{\phi(\mathbf{z}^s)\}_{s\in I_t})\rVert^{2}_{X_h},\\
    \mathbf{z}^0 &= \psi(\mathbf{U}_{h, 0}).
    \end{align}
\end{subequations}

In this case, the sources of nonlinearity are the parametrization $\phi$ of the nonlinear approximation manifold and possibly also the residual $G_{h, \delta t}$. So, even if the residual $G_{h, \delta t}$ is linear we obtain a nonlinear least-squares problem to solve, due to the additional nonlinearity introduced with $\phi$. As mentioned in the introduction, this is often a necessary step to overcome the problem of a slow Kolmogorov n-width decay with a nonlinear approximating manifold that achieves a satisfactory accuracy with a lower latent/reduced dimension with respect to linear approximations.

Due to the nonlinearity, in general, solving~\eqref{eq:PDEsystemDiscrete} is inefficient since the dependence on the number of degrees of freedom cannot be overcome. There are two factors that contribute to making the formulation~\eqref{eq:PDEsystemDiscrete} not feasible as it was introduced in~\cite{lee2020model}. As for the linear case, the first is the nonlinearity of the residual $G_{h, \delta t}$, for which hyper-reduction techniques must be implemented. The second is the possibly expensive evaluation of $\phi$ and its dependence on the whole number of degrees of freedom $d$ since the image of $\phi$ is contained in $\mathbb{R}^d$. So, hyper-reduction or similar techniques must be implemented also for the map $\phi$, that in our case will be a neural network. See section~\ref{subsec:ahr} for more details on hyper-reduction and the next~\ref{subsec: inductive biases} for the definition of our nonlinear approximating solution manifold through the parametrization map $\phi$.

\subsection{Convolutional autoencoders: encodings and inductive biases}
\label{subsec: inductive biases}
Applying SVD or principal component analysis (PCA) in machine learning jargon, to extract small dimensional and meaningful features from data is a technique largely employed in the data science community. After PCA, the new features can be used to train a neural network architecture more efficiently or for other purposes like clustering. This reasoning is applied also to data representing physical fields in model order reduction. One of the first examples is introduced by Ghattas et al.~\cite{o2022derivative} in the context of inverse problems and model order reduction: the parameter-to-observable map is trained as a deep neural network (DNN) from the inputs reduced with active subspaces~\cite{constantine2015active} to the outputs reduced with proper orthogonal decomposition (POD). In the context of model order reduction with autoencoders this techinque is applied in~\cite{fresca2022pod}.

We remark that an autoencoder with linear activation functions can reproduce the accuracy of truncated SVD, and even the principal modes. So, employing the SVD to extract meaningful features instead of adding few linear neural networks layers does not make a difference in terms of accuracy. What is crucial, especially for problems with a huge number of degrees of freedom, is the efficiency of SVD and its randomized version (rSVD) compared to the training of neural network layers with linear activations.

Since we are considering physical fields supported on meshes we should not be limited to SVD: other compression algorithms possibly extracting meaningful features from spatial and temporal correlations are the Fourier and Wavelet transforms. In this context, model order reduction in frequency space is an active field of research~\cite{guglielmi2021pseudospectral,nobile2021non, hawkins2023model}. Recently, also the Radon-Cumulative-Distribution (RCD) transform was applied to advection-dominated problems in model order reduction~\cite{long2023novel}. We will represent such generic transforms with $f_{\text{filter}}:X_h\sim\mathbb{R}^d\rightarrow\mathbb{R}^p$ and their approximate or true left inverse $f^{-1}_{\text{filter}}:\mathbb{R}^p\rightarrow X_h\sim\mathbb{R}^d$, such that $f^{-1}_{\text{filter}}\circ f_{\text{filter}}\approx I_d$.

In this work we will only consider rSVD to define the filtering maps $f_{\text{filter}}$ and $f^{-1}_{\text{filter}}$, but the framework can be easily extended to other compression algorithms. The definitions of $f_{\text{filter}}$ and $f^{-1}_{\text{filter}}$ are reported in equation~\ref{eq:filtering maps} of the next section. We will show some testcases where the number of rSVD modes needed to achieve a satisfactory accuracy reaches $300$, underlying truly slow Kolmogorov n-width applications, while in the literature only a moderate number of modes has been employed. 

\begin{figure}
    \centering
    \includegraphics[width=0.6\textwidth]{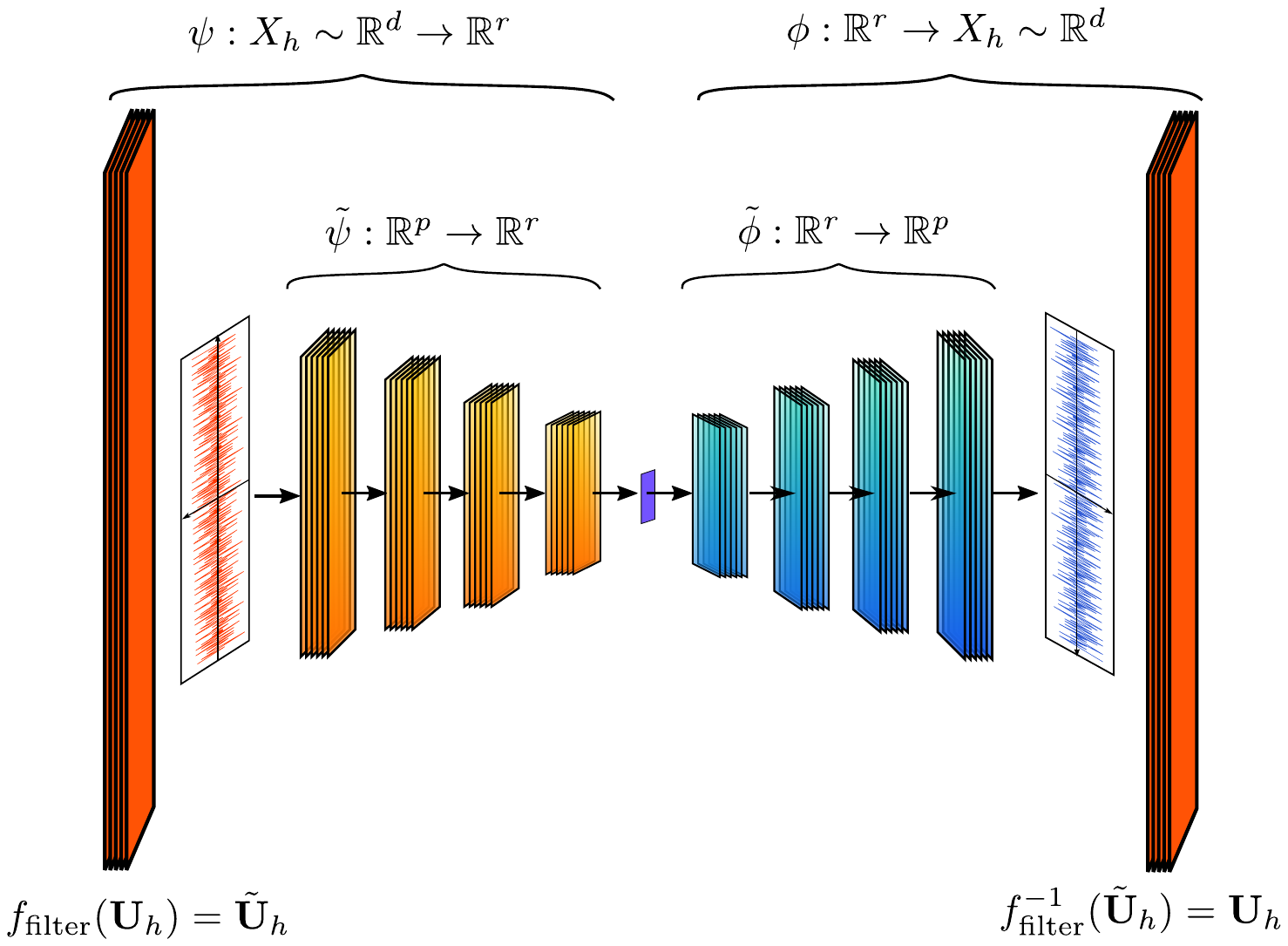}
    \caption{The image shows schematically our implementation of the nonlinear maps $\phi = f^{-1}_{\text{filter}}\circ \Tilde{\phi}$ and $\psi = \Tilde{\psi}\circ f_{\text{filter}}$ that define the approximate nonlinear solution manifold. The functions $\Tilde{\phi}$ and $\Tilde{\psi}$ represent neural networks composed by subsequent 1d-convolutional and transposed 1d-convolutional layers, their specifics are reported in Table~\ref{tab: cae 150} and ~\ref{tab: cae 300}. The maps $f_{\text{filter}}$ and $f^{-1}_{\text{filter}}$ are the linear projection onto the first $p$ rSVD modes and its transpose, their definition is reported in equation~\eqref{eq:filtering maps}.}
    \label{fig: POD-CNN}
\end{figure}

The neural network architectures we are going to use to define the maps $\phi:U\subset\mathbb{R}^r\rightarrow\mathbb{R}^d$ and $\psi:\phi(U)\rightarrow\mathbb{R}^r$ are reported in the Appendix in Table~\ref{tab: cae 150} and \ref{tab: cae 300} and shown in Figure~\ref{fig: POD-CNN}. They are composed by standard 1d-convolutional layers since the filtered states $\Tilde{U}_h = f_{\text{filter}}(\mathbf{U}_h)$ are not supported on a possible unstructured mesh anymore, but belong to the space of frequencies. So, in general, this approach is a viable alternative to graph neural networks~\cite{bronstein2017geometric} or other techniques to approximate physical fields supported on unstructured meshes. 

To separate the application of the filtering/transforms maps $f_{\text{filter}}$ and $f^{-1}_{\text{filter}}$ from the convolutional neural networks layers, we define $\Tilde{\psi}:\mathbb{R}^p\rightarrow\mathbb{R}^r$ and $\Tilde{\phi}:\mathbb{R}^r\rightarrow\mathbb{R}^p$ through the relations $\phi = f^{-1}_{\text{filter}}\circ \Tilde{\phi}$ and $\psi = \Tilde{\psi}\circ f_{\text{filter}}$. So, the CNN layers are encapsulated in $\Tilde{\phi}$ and $\Tilde{\psi}$.

If the state vector $\mathbf{U}_h\in X_h$ includes more than one physical field, we have decided not to extract the frequencies or SVD reduced variables $\Tilde{\mathbf{U}_h}$  for each field, but to do so altogether in a monolithic fashion. This results in single-channel input of the CNN $\Tilde{\psi}$ and a single-channel output of the CNN $\Tilde{\phi}$, instead of having as many channels as the number of physical fields included in $\mathbf{U}_h$.

We remark that our CAE is trained only in the space of frequencies as input-output spaces, achieving a relevant speedup thanks to this, as pointed out in~\cite{fresca2022pod}. The nonlinearity of the autoencoder is only exploited to further reduce the dimensionality from the frequency spaces and to directly approximate the solution manifold, which is effectively linear approximated in our case. For truly nonlinear solution manifold approximants see~\cite{romor2023non}.

\subsection{Parallel Randomized Singular Value Decomposition}
\label{subsec:rSVD}
We recap in brief the procedure of randomized singular value decomposition (rSVD)~\cite{Halko2011}, necessary to evaluate the modes when the snapshots matrix cannot be assembled altogether due to memory and computational constraints. An alternative is represented by the frequent directions algorithm~\cite{Ghashami2016}. We remark also that rSVD requires only matrix-vector evaluations and therefore can be applied also in a matrix-free fashion~\cite{Isaac2015}.

The only ingredient needed is the column-wise ordered matrix $A_{\text{train}}\in\mathbb{R}^{d\times n_{\text{train}}}$ of the training snapshots collection $\mathcal{U}_{\text{train}}:=\{\mathbf{U}_{\boldsymbol{\mu}, t}\}_{\boldsymbol{\mu}\in
\mathcal{P}_{\text{train}},\ t\in V_{\boldsymbol{\mu}}}$, with $V_{\boldsymbol{\mu}}=\{t_1,\dots,t_{N_{\boldsymbol{\mu}}}\}$,
\begin{equation}
    A_{\text{train}}=\left[
    \begin{pmatrix}
        \mid & \mid & \mid & \mid\\
        \mathbf{U}_{\boldsymbol{\mu}_1, t_1} & \mathbf{U}_{\boldsymbol{\mu}_1, t_2} & \dots & \mathbf{U}_{\boldsymbol{\mu}_1, t_{N_{\boldsymbol{\mu}_1}}} \\
        \mid & \mid & \mid & \mid
      \end{pmatrix},\dots, \begin{pmatrix}
        \mid & \mid & \mid & \mid\\
        \mathbf{U}_{\boldsymbol{\mu}_{|\mathcal{P}_{\text{train}}}|, t_1} & \mathbf{U}_{\boldsymbol{\mu}_{|\mathcal{P}_{\text{train}}}|, t_2} & \dots & \mathbf{U}_{\boldsymbol{\mu}_{|\mathcal{P}_{\text{train}}}|, t_{N_{\boldsymbol{\mu}_{|\mathcal{P}_{\text{train}}|}}}} \\
        \mid & \mid & \mid & \mid
      \end{pmatrix}\right]\in\mathbb{R}^{d\times n_{\text{train}}}
\end{equation}
with $|\mathcal{U}_{\text{train}}|=n_{\text{train}}$ and $\mathbf{U}_{\boldsymbol{\mu}, t}\in\mathbb{R}^d$ for all $\boldsymbol{\mu}\in
\mathcal{P}_{\text{train}},\ t\in V_{\boldsymbol{\mu}}$. In our case, we do not assemble $A_{\text{train}}$ since the computational domain can be partitioned and assigned to different processors during the evaluation of the full-order training solutions $\mathcal{U}_{\text{train}}$.

We represent with $M\in\mathbb{N}$ the number of cells $T\subset\Omega_h$ of the mesh $\{T_i\}_{i=1}^{M}\subset\mathcal{T}$ representing our discretized computational domain, and with $c\in\mathbb{N}$ the number of physical fields we are approximating: for the CNS test case $c=6$ for the INS $c=5$. The total number of degrees of freedom (dofs) is $d=c\cdot M$.

Since we will reduce with rSVD all the physical fields altogether in a monolithic fashion, we need to normalize the training snapshots with respect to the cell-wise measure $V\in\mathbb{R}^{d}$ (length, area or volume of each cell depending on the dimensionalty of the mesh) and the different order of magnitudes and unit of measurement of the different physical fields considered.

For example, for the CNS test case $c=6$ we consider velocity $\mathbf{u}_{\boldsymbol{\mu}, t}\in\mathbb{R}^M$, density $\rho_{\boldsymbol{\mu}, t}\in\mathbb{R}^M$, internal energy $e_{\boldsymbol{\mu}, t}\in\mathbb{R}^M$ and pressure $p_{\boldsymbol{\mu}, t}\in\mathbb{R}^M$, so that $$U_{\boldsymbol{\mu}, t}=(\mathbf{u}_{\boldsymbol{\mu}, t}, \rho_{\boldsymbol{\mu}, t}, e_{\boldsymbol{\mu}, t}, p_{\boldsymbol{\mu}, t})\in\mathbb{R}^{M\times c}=\mathbb{R}^d.$$ Similarly, for the INS test case $c=5$ we consider velocity $\mathbf{u}_{\boldsymbol{\mu}, t}\in\mathbb{R}^M$, pressure $p_{\boldsymbol{\mu}, t}\in\mathbb{R}^M$ and the turbulence viscosity $\nu_{\boldsymbol{\mu}, t}\in\mathbb{R}^M$, so that 
\begin{equation}
    U_{\boldsymbol{\mu}, t}=(\mathbf{u}_{\boldsymbol{\mu}, t}, p_{\boldsymbol{\mu}, t}, \nu_{\boldsymbol{\mu}, t})\in\mathbb{R}^{M\times c}=\mathbb{R}^d.
\end{equation}

The vector of cell-wise measures $\mathbf{V}\in\mathbb{R}^d$ is assembled from the vector $\mathbf{v}=\{\mathcal{L}_{\text{Lebesgue}}(T_i)\}_{i=1}^{M}$ where $\mathcal{L}_{\text{Lebesgue}}$ is the Lebesgue measure in $R^D\supset\Omega_h$. So that $\mathbf{V}=(\mathbf{v})_{i=1}^{c}\in\mathbb{R}^d$ is formed stacking $\mathbf{v}$ $c$-times. The pyhsical normalization field $\mathbf{N}\in\mathbb{R}^d$ is obtained from the maximum $L^2$-norm of each field: for the CNS test case we have
\begin{equation}
    u_{\text{max}}=\max_{\boldsymbol{\mu}, t} \lVert \mathbf{u}_{\boldsymbol{\mu}, t}\rVert_{2},\quad \rho_{\text{max}}=\max_{\boldsymbol{\mu}, t} \lVert \rho_{\boldsymbol{\mu}, t}\rVert_{2},\quad e_{\text{max}}=\max_{\boldsymbol{\mu}, t} \lVert e_{\boldsymbol{\mu}, t}\rVert_{2},\quad p_{\text{max}}=\max_{\boldsymbol{\mu}, t} \lVert p_{\boldsymbol{\mu}, t}\rVert_{2},
\end{equation}
so that $\mathbf{N}=(\mathbf{1}_{3M}\cdot u_{\text{max}}, \mathbf{1}_{M}\cdot\rho_{\text{max}}, \mathbf{1}_{M}\cdot e_{\text{max}}, \mathbf{1}_{M}\cdot p_{\text{max}})\in\mathbb{R}^d$, with $\mathbf{1}_M$ the $M$-dimensional vector of ones. Similarly, for the INS test case we consider 
\begin{equation}
    u_{\text{max}}=\max_{\boldsymbol{\mu}, t} \lVert \mathbf{u}_{\boldsymbol{\mu}, t}\rVert_{2},\quad \rho_{\text{max}}=\max_{\boldsymbol{\mu}, t} \lVert \rho_{\boldsymbol{\mu}, t}\rVert_{2},\quad e_{\text{max}}=\max_{\boldsymbol{\mu}, t} \lVert e_{\boldsymbol{\mu}, t}\rVert_{2},\quad p_{\text{max}}=\max_{\boldsymbol{\mu}, t} \lVert p_{\boldsymbol{\mu}, t}\rVert_{2},
\end{equation}
so that $\mathbf{N}=(\mathbf{1}_{3M}\cdot u_{\text{max}}, \mathbf{1}_{M}\cdot p_{\text{max}}, \mathbf{1}_{M}\cdot \nu_{\text{max}})\in\mathbb{R}^d$.

So the columns $\{A^{\text{train}}_{\boldsymbol{\mu}, t}\}_{\boldsymbol{\mu}\in
\mathcal{P}_{\text{train}},\ t\in V_{\boldsymbol{\mu}}}$ of $A_{\text{train}}$ are actually defined as
\begin{equation}
    \label{eq:normalization}
    \mathbb{R}^{d}\ni A^{\text{train}}_{\boldsymbol{\mu}, t} = \mathbf{U}_{\boldsymbol{\mu}, t}\oslash\mathbf{W},\qquad \mathbf{W}=\mathbf{N}\oslash\mathbf{V},
\end{equation}
where we have considered only element-wise divisions $\bullet\oslash\bullet$ between vectors in $\mathbb{R}^d$. Notice that in this way we obtain unit-less states $\mathbf{U}_{\boldsymbol{\mu}, t}/\mathbf{W}\in\mathbb{R}^d$. For the impact of physical normalization in model order reduction, see~\cite{parish2022impact}.

\IncMargin{1em}
\begin{algorithm}[htb]
    \SetKwInOut{Input}{input}\SetKwInOut{Output}{output}
    \Indm \Input{training column-wise ordered snapshots matrix $A_{\text{train}}\in\mathbb{R}^{d\times n_{\text{train}}}$,\\
    $r_{\text{rSVD}}$ reduced rSVD dimension,
    $p$ oversampling parameter}
    \Output{$U\in\mathbb{R}^{d\times r_{\text{rSVD}}}$ rSVD modes}
    \Indp \BlankLine 
    Define the sketch dimension $l=r_{\text{rSVD}}+p$.\\
    Draw the sketch matrix $\Omega\in\mathbb{R}^{n_{\text{train}}\times l}$ as a Gaussian random matrix.\\
    Assemble \textbf{in parallel} $Y=A_{\text{train}}\Omega$, with $Y\in\mathbb{R}^{d\times l}$.\\
    Evaluate the orthonormal basis $Q\in\mathcal{R}^{n_{\text{train}}}$ using $QR$-factorization $\mathbb{R}^{}Y=QR$, with $R\in\mathbb{R}^{n_{\text{train}}\times n_{\text{train}}}$.\\
    Project the snapshots into a random lower $l$-dimensional space \textbf{in parallel}: $S=Q^T A_{\text{train}}$, with $S\in\mathbb{R}^{l\times n_{\text{train}}}$.\\
    Compute the thin SVD $S=\Tilde{U}\Tilde{\Sigma} \Tilde{V}$, with $\Tilde{U}\in\mathbb{R}^{l\times r_{\text{rSVD}}},\ \Sigma\in\mathbb{R}^{r_{\text{rSVD}}\times r_{\text{rSVD}}},\ \Tilde{V}\in\mathbb{R}^{r_{\text{rSVD}}\times n_{\text{train}}}$.\\
    Evaluate the rSVD modes $U=Q\Tilde{U}$, with $U\in\mathbb{R}^{d\times r_{\text{rSVD}}}$.
    \DontPrintSemicolon \caption{Parallelized Randomized Singular Value Decomposition.}
    \label{alg:rsvd}
\end{algorithm}
\DecMargin{1em}

The reduced train and test rSVD coordinates are obtained with the following linear projection, employing the rSVD modes $U\in\mathbb{R}^{d\times r_{\text{rSVD}}}$ from Algorithm~\ref{alg:rsvd}:
\begin{subequations}
    \begin{align}
        \mathbb{R}^{r_{\text{rSVD}}\times n_{\text{train}}}\ni Y_{\text{train}}&=U^T A_{\text{train}},\qquad\qquad&\text{(reduced train rSVD coordinates)}\\
        \mathbb{R}^{r_{\text{rSVD}}\times n_{\text{test}}}\ni Y_{\text{test}}&=U^T A_{\text{test}},\qquad\qquad&\text{(reduced test rSVD coordinates)}
    \end{align}
\end{subequations}
and the reconstructed train and test fields are obtained employing the normalizing vector $\mathbf{W=\mathbf{N}\oslash\mathbf{V}}$:
\begin{subequations}
    \begin{align}
        \mathbb{R}^{d\times n_{\text{train}}}\ni A^{\text{rec}}_{\text{train}}&=\mathbf{W}\odot U Y_{\text{train}} = \mathbf{W}\odot UU^T A_{\text{train}},\qquad\qquad&\text{(reconstructed train snapshots)}\\
        \mathbb{R}^{d\times n_{\text{test}}}\ni A^{\text{rec}}_{\text{test}}&=\mathbf{W}\odot U Y_{\text{test}} = \mathbf{W}\odot UU^T A_{\text{test}},\qquad\qquad&\text{(reconstructed test snapshots)}
    \end{align}
\end{subequations}
where $\odot$ is the Hadamard columns-wise product, inverse operation of the normalization applied in~\eqref{eq:normalization}.

To decide if the number of rSVD modes is sufficient to achieve the desired accuracy for the problem at hand we consider the mean and max relative $L^2$ reconstruction error on the training and test sets:
\begin{subequations}
    \begin{align}
        \label{eq:rec_rsvd}
        \lVert A^{\text{rec}}_{\text{train}}\rVert_{2, mean} &= \frac{1}{n_{\text{train}}} \sum_{\boldsymbol{\mu}\in
        \mathcal{P}_{\text{train}},\ t\in V_{\boldsymbol{\mu}}}\frac{\lVert \mathbf{U}_{\boldsymbol{\mu}, t}^{\text{rec, train}}-\mathbf{U}^{\text{train}}_{\boldsymbol{\mu}, t}\rVert_2}{\lVert \mathbf{U}^{\text{train}}_{\boldsymbol{\mu},t}\rVert_2},\qquad &&\lVert A^{\text{rec}}_{\text{train}}\rVert_{2, max} = \max_{\boldsymbol{\mu}\in
        \mathcal{P}_{\text{train}},\ t\in V_{\boldsymbol{\mu}}}\frac{\lVert \mathbf{U}_{\boldsymbol{\mu}, t}^{\text{rec, train}}-\mathbf{U}^{\text{train}}_{\boldsymbol{\mu}, t}\rVert_2}{\lVert \mathbf{U}^{\text{train}}_{\boldsymbol{\mu},t}\rVert_2},\\
        \label{eq:rec_rsvd_test}
        \lVert A^{\text{rec}}_{\text{test}}\rVert_{2, mean} &= \frac{1}{n_{\text{test}}} \sum_{\boldsymbol{\mu}\in
        \mathcal{P}_{\text{test}},\ t\in V_{\boldsymbol{\mu}}}\frac{\lVert \mathbf{U}_{\boldsymbol{\mu}, t}^{\text{rec, test}}-\mathbf{U}^{\text{test}}_{\boldsymbol{\mu}, t}\rVert_2}{\lVert \mathbf{U}^{\text{test}}_{\boldsymbol{\mu},t}\rVert_2},\qquad &&\lVert A^{\text{rec}}_{\text{test}}\rVert_{2, max} = \max_{\boldsymbol{\mu}\in
        \mathcal{P}_{\text{test}},\ t\in V_{\boldsymbol{\mu}}}\frac{\lVert \mathbf{U}_{\boldsymbol{\mu}, t}^{\text{rec, test}}-\mathbf{U}^{\text{test}}_{\boldsymbol{\mu}, t}\rVert_2}{\lVert \mathbf{U}^{\text{test}}_{\boldsymbol{\mu},t}\rVert_2},
    \end{align}
\end{subequations}
where $A^{\text{rec}}_{\text{train}}=(\mathbf{U}_{\boldsymbol{\mu}, t}^{\text{rec, train}})_{\boldsymbol{\mu}\in
\mathcal{P}_{\text{train}},\ t\in V_{\boldsymbol{\mu}}}\in\mathbb{R}^{d\times n_{\text{train}}}$ and $A^{\text{rec}}_{\text{test}}=(\mathbf{U}_{\boldsymbol{\mu}, t}^{\text{rec, test}})_{\boldsymbol{\mu}\in
\mathcal{P}_{\text{test}},\ t\in V_{\boldsymbol{\mu}}}\in\mathbb{R}^{d\times n_{\text{test}}}$.

Finally, we want to explicitly define the filtering/transform map $f_{\text{filter}}:X_h\sim\mathbb{R}^d\rightarrow\mathbb{R}^p$ and its approximate left inverse $f^{-1}_{\text{filter}}:\mathbb{R}^p\rightarrow X_h\sim\mathbb{R}^d$ with $p=r_{\text{rSVD}}$:
\begin{equation}
    \label{eq:filtering maps}
    f_{\text{filter}}(\mathbf{U}_h) = U^T\left(\mathbf{U}_h\oslash\mathbf{W}\right),\quad f^{-1}_{\text{filter}}(\Tilde{\mathbf{U}_h}) = \left(\mathbf{W}\odot U\right)\Tilde{\mathbf{U}}_h.
\end{equation}

\section{Hyper-reduction}
\label{sec:hr}
\label{subsec:ahr}
As introduced in section~\ref{subsec:nlspg}, there are two main problematics that affect the efficient resolution of the nonlinear least squares problem in equation~\eqref{eq:PDEsystemDiscrete} at each time instance $t\in V_{\boldsymbol{\mu}}$ and for each intermediate optimization step $i\in\{1,\dots, N_{t, \boldsymbol{\mu}}\}$ required by the nonlinear least-squares method, and they are both linked to the evaluation of the residual 
\begin{equation}
  G_{h, \delta t}(\boldsymbol{\mu}, \phi(\mathbf{z}), \{\phi(\mathbf{z}^s)\}_{s\in I_t}).
\end{equation}
We recall that we employ the derivative-free Pounders solver~\cite{wildsolving} implemented in PETSc~\cite{balay2019petsc}. Also nonlinear least-squares problem optimizers that employ an approximation or the true Jacobian of the residual $G_{h, \delta t}$ can be employed.

The main problematics to efficiently evaluate $G_{h, \delta t}$ are the following:
\begin{enumerate}
  \item in general, the nonlinearity of $G_{h, \delta t}$ makes its evaluation dependent on the number of dofs $d$,
  \item the map $\phi:U\subset\mathbb{R}^r\rightarrow\mathbb{R}^d$ might be computationally heavy to evaluate for each $r$-dimensional input and depends on the number of dofs since its output is $d$-dimensional.
\end{enumerate}
The reason why we need this independence on the number of dofs is for our model order reduction procedure to be efficient even when $d$ increases. Our two test cases $CNS$ and $INS$ have approximately $M=30000$ and $M=200000$ cells, and $d=180000$ and $d=1000000$ dofs respectively, which are still a moderate number of dofs compared to real applications. If we want to extend the methodology to larger meshes, we have to guarantee the independence on the number of dofs of our procedure.

We will address first the nonlinearity coming from $G_{h, \delta t}$. Typically, in the case of LM-LSPG, if the residual has a nonlinear term directly coming from the parametric PDE model, a class of methods under the name of hyper-reduction can be applied to ameliorate the situation. The idea is to reconstruct the residual only from its evaluations on a subset of degrees of freedom, in general. To do so, from the physical fields of the model considered taken as inputs, the values of those fields on the stencil needed by the numerical discretization have to be computed. 

The simplest approach consists in collocating the residual on a subset of cells of the mesh, from now on called nodes or magic points. So we introduce two projection maps: the projection onto the magic points
\begin{equation}
  P_{r_h}:\mathbb{R}^d\rightarrow\mathbb{R}^{r_h},\quad P_{r_h}(\mathbf{U}_h) = \begin{pmatrix}
    \mid & \mid & \mid & \mid\\
    \mathbf{e}_{i_1} & \mathbf{e}_{i_2} & \dots & \mathbf{e}_{i_{r_h}} \\
    \mid & \mid & \mid & \mid
  \end{pmatrix}^{T}\in\mathbb{R}^{r_{h}\times d},\quad 0<r_h\ll d
\end{equation}
where $S_{r_h}=\{\mathbf{e}_{i_j}\}_{j=1}^{r_h}$ is a subset of the standard basis of $\mathbb{R}^d$ and the projection onto the submesh needed to evaluate the residual on the magic points
\begin{equation}
  P_{r_h}^s:\mathbb{R}^d\rightarrow\mathbb{R}^{s},\quad P^s_{r_h}(\mathbf{U}_h) = \begin{pmatrix}
    \mid & \mid & \mid & \mid\\
    \mathbf{e}_{i_1} & \mathbf{e}_{i_2} & \dots & \mathbf{e}_{i_{s}} \\
    \mid & \mid & \mid & \mid
  \end{pmatrix}^{T}\in\mathbb{R}^{s\times d},\quad 0<r_h < s\ll d
\end{equation}
where $S_{r_h}^s=\{\mathbf{e}_{i_j}\}_{j=1}^{s}$ is a subset of the standard basis of $\mathbb{R}^d$ containing $S_{r_h}\subset S_{r_h}^s$.

With these definitions the residual from equation~\eqref{eq:PDEsystemDiscrete} can be hyper-reduced as 
\begin{subequations}
  \begin{align}
    \mathbf{z}^t =& \argmin_{\mathbf{z}\in \mathbb{R}^r}\ \lVert P_{r_h}G_{h, \delta t}(\boldsymbol{\mu}, P_{r_h}^s(\phi(\mathbf{z})), \{P_{r_h}^s(\phi(\mathbf{z}^s))\}_{s\in I_t})\rVert^{2}_{\mathbb{R}^{r_h}}\\
    =&\argmin_{\mathbf{z}\in \mathbb{R}^r}\ \lVert P_{r_h}G_{h, \delta t}(\boldsymbol{\mu}, (P_{r_h}^s\circ f^{-1}_{\text{filter}})(\Tilde{\phi}(\mathbf{z})), \{(P_{r_h}^s\circ f^{-1}_{\text{filter}})(\Tilde{\phi}(\mathbf{z}^s))\}_{s\in I_t})\rVert^{2}_{\mathbb{R}^{r_h}}
    \label{eq:PDEsystemDiscreteHyperReduced2},
  \end{align}  
\end{subequations}
where we have employed the definitions of $\Tilde{\phi}:\mathbb{R}^r\rightarrow\mathbb{R}^p$ and $f^{-1}_{\text{filter}}:\mathbb{R}^p\rightarrow\mathbb{R}^d$ from section~\ref{subsec: inductive biases}. 

In this way we have addressed the problem coming from the nonlinearity of the residual $G_{h, \delta t}$. At the same time, thanks to the choice of $\phi=f^{-1}_{\text{filter}}\circ\Tilde{\phi}$ as composition of a linear projection depending on the dofs $f^{-1}_{\text{filter}}$ and a small nonlinear neural network $\Tilde{\phi}$ independent on the number of dofs, we have also tackled the second problem. In fact, also the parametrization map $(P_{r_h}^s\circ f^{-1}_{\text{filter}})(\Tilde{\phi})$ restricted to the submesh is now independent on the number of dofs.

\begin{figure}[ht!]
  \centering
  \includegraphics[width=0.7\textwidth]{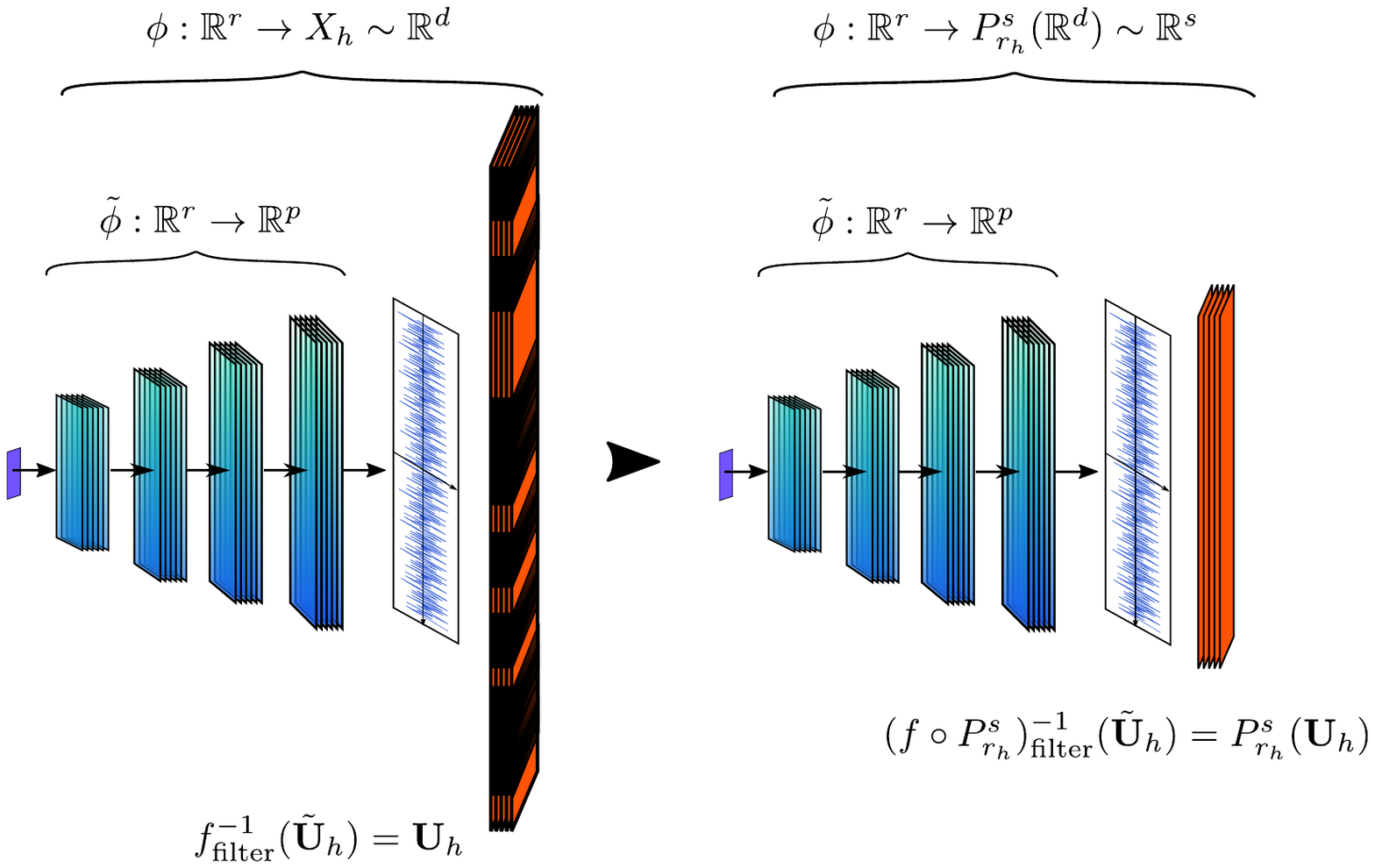}
  \caption{\textbf{Left:} the decoder map $\Tilde{\psi}:\mathbb{R}^r\rightarrow\mathbb{R}^p$ followed by the vector matrix multiplication with the rSVD modes $f^{-1}_{\text{filter}}(\Tilde{\mathbf{U}}_h):\mathbb{R}^p\rightarrow X_h\sim\mathbb{R}^d$. The latter are restricted to the submesh through the projection $P^s_{r_h}:X_h\sim\mathbb{R}^d\rightarrow\mathbb{R}^s$, that is the blackened dofs are discarded. \textbf{Right:} the map actually employed in the hyper-reduced nonlinear manifold least-squares Petrov-Galerkin method $\phi:\mathbb{R}^r\rightarrow P^s_{r_h}(X_h)\sim\mathbb{R}^s$. It is independent on the number of dofs and its evaluation is efficient thanks to the relatively small size of the decoder $\Tilde{\psi}$.}
  \label{fig:hr map}
\end{figure}

To be more specific, since we will be employing only rSVD as linear projections we have
\begin{equation}
  \label{eq:filtering maps hr}
  f^{-1}_{\text{filter}}(\Tilde{\mathbf{U}_h}) = P_{r_h}^s\circ\left(\left(\mathbf{W}\odot U\right)\Tilde{\mathbf{U}}_h\right) = \left(P_{r_h}^s(\mathbf{W})\odot P_{r_h}^s(U)\right)\Tilde{\mathbf{U}}_h.
\end{equation}
where we have employed definition~\ref{eq:filtering maps}. So the hyper-reduction affects only the rSVD modes $U\in\mathbb{R}^{d\times p}$ and the normalization vector $\mathbf{W}\in\mathbb{R}^d$. A schematic representation of the hyper-reduced approximate nonlinear manifold parametrization map is shown in Figure~\ref{fig:hr map}.

We want to remark that the hyper-reduction procedure presented is effective thanks to the choice of implementation of the parametrization map $\phi$ through a combination of neural networks and rSVD modes. However, in some cases the number of rSVD modes required $p=r_{\text{rSVD}}$ can become so large to guarantee a threshold accuracy in the relative $L^2$ reconstruction error that the methodology is no more efficient, since $p\gg 0$ and $\Tilde{\phi}:\mathbb{R}^r\rightarrow\mathbb{R}^{r_{\text{rSVD}}}$. This computational burden affects both the offline stage for the training of the NN and the online residual evaluation. In these cases, one may want to employ heavier and deep generic neural network architectures like CNNs for structured meshes or GNNs for unstructured meshes that recover a good approximation of the solution manifold. The employment of these deep NNs brings up the problem of how to hyper-reduce them. The methodology introduced in this section cannot be applied, but a solution is represented by a strategy called teacher-student training for which in a following training phase a smaller fast \textit{student} NN is tuned to replicate the results of the bigger slow \textit{teacher} NN. This alternative approach is presented and studied in~\cite{romor2023non}.

\subsection{Hyper-reduction methods and Magic Points Selection Algorithms}
We are left with the task of defining the set of magic points $S_{r_h}$, since the submesh $S_s\supset S_{r_h}$ is identified from the magic points and the choice of numerical scheme employed to discretized the parametric PDE. Until now, we have considered only the collocation of the residual on the magic points as hyper-reduction method. However, not only this is not the sole possible implementation but also not the most common one.

What we have actually described is part of the reduced over-collocation hyper-reduction method~\cite{chen2021eim}. In this section, we will introduce also the gappy discrete empirical interpolation method (DEIM)~\cite{chaturantabut2010nonlinear}, the energy-sampling and weighting method (ECSW)~\cite{farhat2015structure}, and DEIM with the quasi-optimal point selection algorithm (S-OPT) from~\cite{lauzon2022s}.

We have experimentally observed that for our test cases with slow Kolmogorov n-width and $r_{\text{rSVD}}=150$ and $r_{\text{rSVD}}=300$ rSVD modes, the reduced over-collocation method performs better. For more comments see section~\ref{sec:discussion}. Moreover, for test cases with a bigger number of dofs we have also developed a more successful adaptive magic points selection method introduced in section~\ref{subsec:adaptive-hr}.

So, what we will be particularly focused on are the magic points selection algorithms from DEIM, ECSW and S-OPT. The starting point for each one of them is the computation of a set of rSVD modes $U_{hr}\in\mathbb{R}^{d\times n_{hr}}$. Since we are hyper-reducing the residual $G_{h, \delta t}$ we need to collect a dataset of residual snapshots $\mathcal{G}_{\text{train}}=\{G_{h, \delta t}(\boldsymbol{\mu}, \mathbf{U}_h^t, \{\mathbf{U}_h^s\}_{s\in I_t})\}_{\boldsymbol{\mu}\in
\mathcal{P}_{\text{train}},\ t\in V_{\boldsymbol{\mu}}}=\{G_{\boldsymbol{\mu}, t}\}_{\boldsymbol{\mu}\in
\mathcal{P}_{\text{train}},\ t\in V_{\boldsymbol{\mu}}}$, with $V_{\boldsymbol{\mu}}=\{t_1,\dots,t_{N_{\boldsymbol{\mu}}}\}$, in the training phase and compress them with rSVD:
\begin{equation}
    G_{\text{train}}=\left[
    \begin{pmatrix}
        \mid & \mid & \mid & \mid\\
        G_{\boldsymbol{\mu}_1, t_1} & G_{\boldsymbol{\mu}_2, t_1} & \dots & G_{\boldsymbol{\mu}_1, t_{N_{\boldsymbol{\mu}_1}}} \\
        \mid & \mid & \mid & \mid
      \end{pmatrix},\dots, \begin{pmatrix}
        \mid & \mid & \mid & \mid\\
        G_{\boldsymbol{\mu}_{| \mathcal{P}_{\text{train}}|}, t_1} & G_{\boldsymbol{\mu}_{| \mathcal{P}_{\text{train}}|}, t_1} & \dots & G_{\boldsymbol{\mu}_{| \mathcal{P}_{\text{train}}|}, t_{N_{\boldsymbol{\mu}_{| \mathcal{P}_{\text{train}}|}}}} \\
        \mid & \mid & \mid & \mid
      \end{pmatrix}\right]\in\mathbb{R}^{d\times n_{\text{train}}},
\end{equation}
from this residual snapshots matrix $G_{\text{train}}$ the sRSVD modes are computed $U_G\in\mathbb{R}^{d\times r_{G}}$, $r_G>0$ and utilzed for magic points sampling after setting $U_{hr}=U_G$. In this case, $r_G$ might possibly be different than $r_{\text{rSVD}}$.

However, employing the physical fields' rSVD modes $U_{hr}=U\in\mathbb{R}^{d\times r_{\text{rSVD}}}$ also to perform the magic points' sampling is a fast alternative~\cite{choi2020sns}, since no additional residual snapshots need to be collected, apart from those used to define the map $\phi$, and the computational cost of an additional application of the rSVD algorithm, this time on the residual snapshots, is avoided.

\subsubsection{Gappy discrete empirical interpolation}
\label{subsubsec:deim}
In the gappy DEIM algorithm, after the computation of the hyper-reduction basis $U_{hr}\in\mathbb{R}^{d\times r_{hr}}$ from the physical fields $A_{\text{train}}\in\mathbb{R}^{d\times n_{\text{train}}}$ or residual snapshots $G_{\text{train}}\in\mathbb{R}^{d\times n_{\text{train}}}$, $U_{hr}$ is employed to find the magic points with a greedy algorithm: at each step, it is selected the cell of the mesh associated to the highest value of the hyper-reduction reconstruction error $\mathbf{r}\in\mathbb{R}^d$
\begin{equation}
  \mathbf{r}=A_{\text{train}}-U(P_{r_h}U)^{\dagger}(P_{r_h}A_{\text{train}}),\qquad\text{or}\qquad \mathbf{r}=G_{\text{train}}-U_{G}(P_{r_h}U_{G})^{\dagger}(P_{r_h}G_{\text{train}})
\end{equation}
where with the notation $(\bullet)^{\dagger}$ we represent the Moore-Penrose pseudo-inverse matrix and $P_{r_h}\in\mathbb{R}^{d\times\Tilde{r}}$ is the intermediate projection matrix that evaluates a vector on the full-state space $\mathbb{R}^d$ on the magic points  $S_{r_h}=\{\mathbf{e}_{i_j}\}_{j=1}^{\Tilde{r}}, \ 0<\Tilde{r}\leq r_h$, selected up to the considered step of the greedy algorithm. We remark that if vector-valued states $\mathbf{U}\in\mathbb{R}^d$ with $d=M\cdot c$, $c>1$ are considered, it is selected the cell of the mesh that maximizes the sum over each of the $c$ fields of the hyper-reduction reconstruction error.

For our implementation of the gappy DEIM greedy nodes selection algorithm, we have chosen the one studied in~\cite{carlberg2013gnat} and applied to the Gauss-Newton tensor approximation (GNAT) hyper-reduction method that is more general than DEIM.

Once the magic points set $S_{r_h}=\{\mathbf{e}_{i_j}\}_{j=1}^{r_h}$ has been evaluated, the magic points and submesh projections $P_{r_h}$, $P^s_{r_h}$ are computed and the following variants of nonlinear hyper-reduced least-squares problem~\eqref{eq:PDEsystemDiscrete} are solved for

\begin{subequations}
  \begin{align}
    \mathbf{z}^t &= \argmin_{\mathbf{z}\in \mathbb{R}^r}\ \lVert (P_{r_h}U)^{\dagger} \left(P_{r_h}G_{h, \delta t}(\boldsymbol{\mu}, P_{r_h}^s(\phi(\mathbf{z})), \{P_{r_h}^s(\phi(\mathbf{z}^s))\}_{s\in I_t})\oslash P_{r_h}\mathbf{W}\right)\rVert^{2}_{\mathbb{R}^{r_h}},\qquad&(\text{FB-DEIM})\label{eq:FB-DEIM}\\
    \mathbf{z}^t &= \argmin_{\mathbf{z}\in \mathbb{R}^r}\ \lVert (P_{r_h}U_{G})^{\dagger} \left(P_{r_h}G_{h, \delta t}(\boldsymbol{\mu}, P_{r_h}^s(\phi(\mathbf{z})), \{P_{r_h}^s(\phi(\mathbf{z}^s))\}_{s\in I_t})\oslash P_{r_h}\mathbf{W}_G\right)\rVert^{2}_{\mathbb{R}^{r_h}},\qquad&(\text{RB-DEIM})\label{eq:RB-DEIM}\\
    \mathbf{z}^t &= \argmin_{\mathbf{z}\in \mathbb{R}^r}\ \lVert P_{r_h}G_{h, \delta t}(\boldsymbol{\mu}, P_{r_h}^s(\phi(\mathbf{z})), \{P_{r_h}^s(\phi(\mathbf{z}^s))\}_{s\in I_t})\rVert^{2}_{\mathbb{R}^{r_h}},\qquad&(\text{C-DEIM})
    \label{eq:DEIM},
  \end{align}  
\end{subequations}
depending on the choice of hyper-reduction basis chosen $U_{hr}=U$ (FB-DEIM) or $U_{hr}=U_G$ (RB-DEIM) or if it is performed reduced over-collocation (C-DEIM). In the case of FB-DEIM and RB-DEIM, the residuals $G_{h, \delta t}$ are divided elment-wise by the normalization vectors $\mathbf{W}\in\mathbb{R}^d$ defined in~\eqref{eq:normalization} and $\mathbf{W}_G\in\mathbb{R}^d$ defined analogously from $G_{\text{train}}$.

\subsubsection{A quasi-optimal nodes sampling method}
\label{subsubsec:sopt}
As pointed out in~\cite{lauzon2022s}, the matrix $P_{r_h}U_{hr}\in\mathbb{R}^{r_{h}\times r_{hr}}$ from the DEIM algorithm loses the orthogonality of its columns with respect to $U_{hr}\in\mathbb{R}^{d\times r_{hr}}$ that is crucial for the numerical stability of DEIM and to minimize the error in the $L^2$-norm of the hyper-reduction interpolation. In fact, the hyper-reduction $L^2$ error can be decomposed~\cite{chaturantabut2010nonlinear} as the sum of the best approximation error on the linear subspace in $\mathbb{R}^d$ spanned by the columns of $U_{hr}\in\mathbb{R}^{d\times r_{hr}}$ and the distance from the projection onto it
\begin{subequations}
\begin{align}
  \lVert A_{\text{train}}-U_{hr}(P_{r_h}U_{r_h})^{\dagger}(P_{r_h}A_{\text{train}})\rVert^2_2 &= \lVert A_{\text{train}}-U_{hr} U_{hr}^TA_{\text{train}}\rVert^2_2 + \lVert  U_{hr} U_{r_h}^TA_{\text{train}}-U_{hr}(P_{r_h}U_{r_h})^{\dagger}(P_{r_h}A_{\text{train}})\rVert^2_2\\
  &= \lVert A_{\text{train}}-U_{hr} U_{hr}^TA_{\text{train}}\rVert^2_2 + \lVert U_{hr}(P_{r_h}U_{r_h})^{\dagger}P_{r_h}(U_{r_h}U_{r_h}^T-I_d)A_{\text{train}}\rVert^2_2\\
  &=\lVert A_{\text{train}}-U_{hr} U_{hr}^TA_{\text{train}}\rVert^2_2 + \lVert (P_{r_h}U_{r_h})^{\dagger}P_{r_h}(U_{r_h}U_{r_h}^T-I_d)A_{\text{train}}\rVert^2_2,
\end{align}
\end{subequations}
where in the last step we have used the fact that $U_{hr}\in\mathbb{R}^{d\times r_{rh}}$ has orthonormal columns. Only the second term depends on the magic points selection, so the optimal strategy would be to minimize the solution of the least-squares problem
\begin{equation}
  \label{eq:s-opt criterion}
  (P_{r_h}U_{r_h})^{\dagger}P_{r_h}(U_{r_h}U_{r_h}^T-I_d)A_{\text{train}}=\argmin_{X\in\mathbb{R}^{r_{hr}\times n_{\text{train}}}} \lVert P_{r_h}U_{r_h}X-P_{r_h}(U_{r_h}U_{r_h}^T-I_d)A_{\text{train}}\rVert^2_2,
\end{equation}
that is to minimize the hyper-reduction $L^2$ error of the remaing part of $A_{\text{train}}\in\mathbb{R}^{d\times n_{\text{train}}}$ after the difference with its projection on the subspace spanned by the hyper-reduction basis $U_{hr}\in\mathbb{R}^{d\times r_{hr}}$. 

A possible way, not necessarily optimal, to minimize~\eqref{eq:s-opt criterion} is to maximize the determinant of $P_{r_h}U_{r_h}\in\mathbb{R}^{r_h\times r_{hr}}$ and at the same time maximize its column orthogonality. In~\cite{shin2016nonadaptive}, they developed an efficient greedy algorithm to do so, maximizing
\begin{equation}
  P_{r_h}^{S}=\argmax_{P\in\mathbb{R}^{r_{h}\times r_{hr}}}\mathcal{S}(PU_{hr}),\qquad
  \mathcal{S}(PU_{hr}) = \left(\frac{\sqrt{\text{det}\left((PU_{hr})^TPU_{hr}\right)}}{\prod_i^{r_{hr}} \lVert (PU_{hr})_i\rVert_2}\right)^{\frac{1}{r_{hr}}}\in[0, 1],
\end{equation}
where $\{(PU_{hr})_i\}_{i=1}^{r_{rh}}$ are the columns of $PU_{hr}\in\mathbb{R}^{r_h\times r_{hr}}$, assuming $ \lVert (PU_{hr})_i\rVert_2\neq 0, \forall i\in\{1,\dots,r_{hr}\}$. In particular, in~\cite{shin2016nonadaptive} it is proved that $\mathcal{S}(PU_{hr})=1$ if and only if the columns of $PU_{hr}$ are mutually orthonormal. The same quasi-optimal criterion is employed for hyper-reduction in~\cite{lauzon2022s}, under the name of S-optimality (SOPT).

The DEIM algorithm with S-optimality magic points selection has the following formulation:
\begin{subequations}
  \begin{align}
    \mathbf{z}^t &= \argmin_{\mathbf{z}\in \mathbb{R}^r}\ \lVert (P^{S}_{r_h}U)^{\dagger} \left(P^{S}_{r_h}G_{h, \delta t}(\boldsymbol{\mu}, P^{S, s}_{r_h}(\phi(\mathbf{z})), \{P^{S, s}_{r_h}(\phi(\mathbf{z}^s))\}_{s\in I_t})\oslash P^{S}_{r_h}\mathbf{W}\right)\rVert^{2}_{\mathbb{R}^{r_h}},\qquad&(\text{FB-DEIM-SOPT})\label{eq:FB-DEIM-OPT}\\
    \mathbf{z}^t &= \argmin_{\mathbf{z}\in \mathbb{R}^r}\ \lVert (P^{S}_{r_h}U_{G})^{\dagger} \left(P^{S}_{r_h}G_{h, \delta t}(\boldsymbol{\mu}, P^{S, s}_{r_h}(\phi(\mathbf{z})), \{P^{S, s}_{r_h}(\phi(\mathbf{z}^s))\}_{s\in I_t})\oslash P^{S}_{r_h}\mathbf{W}_G\right)\rVert^{2}_{\mathbb{R}^{r_h}},\qquad&(\text{RB-DEIM-SOPT})\label{eq:RB-DEIM-OPT}\\
    \mathbf{z}^t &= \argmin_{\mathbf{z}\in \mathbb{R}^r}\ \lVert P^{S}_{r_h}G_{h, \delta t}(\boldsymbol{\mu}, P^{S, s}_{r_h}(\phi(\mathbf{z})), \{P^{S, s}_{r_h}(\phi(\mathbf{z}^s))\}_{s\in I_t})\rVert^{2}_{\mathbb{R}^{r_h}},\qquad&(\text{C-DEIM-SOPT})
    \label{eq:DEIM-SOPT},
  \end{align}  
\end{subequations}
where $P^{S, s}_{r_h}\in\mathbb{R}^{s\times d}$ is the submesh projection corresponding to $P^S_{r_h}\in\mathbb{R}^{r_h\times d}$. Also in this case, depending on the choice of hyper-reduction basis chosen $U_{hr}=U$ (FB-DEIM-SOPT) or $U_{hr}=U_G$ (RB-DEIM-SOPT) or if it is performed reduced over-collocation (C-DEIM-SOPT), there are three different formulations.

\subsubsection{Energy-conserving sampling and weighting method}
\label{subsubsec:ecsw}
Differently from DEIM, the ECSW hyper-reduction method finds a sparse integration formula to approximate the quantity of interest, if it is obtained through an integration on the computational domain $\Omega\subset\mathbb{R}^D$, like the residual $G_{h, \delta t}$ if it is calculated with the FVM, FEM or DGM.

The idea is to find a $\lVert\bullet\rVert_0$-sparse quadrature formula ($\lVert\mathbf{x}\rVert_0=\#\{i\in\mathbb{N}\vert\mathbf{x}_i\neq 0\},\ \forall\in\mathbb{R}^d$) such that the new weights $\mathbf{Q}^{hr}\in\mathbb{R}^d$  are sparse $\lVert\mathbf{Q}^{hr}\rVert_0\ll 1$ and approximate up to a tolerance $0<\tau\ll 1$ the sums of the training residual snapshots:
\begin{equation}
  \mathbf{Q}^{hr}=\argmin_{\mathbf{Q}\in\mathbb{R}^d_+}\lVert\mathbf{Q}\rVert_0\quad \text{s.t.}\quad\lVert \mathbf{Q}\odot U_{hr}-\mathbf{c}\rVert^2_2<\tau\lVert\mathbf{c}\rVert_2^2
\end{equation}
where $\mathbf{Q}\odot U_{hr}$ is the elment-wise multiplication of the quadrature weights vector $\mathbf{Q}\in\mathbb{R}^d$ with the columns of the training residual snapshots matrix $U_{hr}\in\mathbb{R}^{d\times r_{hr}}$, and $\mathbf{c}\in\mathbb{R}^{r_{hr}}$, $\mathbf{c}_j=\sum_{i=1}^{d} (U_{hr})_{ij},\ \forall j\in\{1,\dots, r_{hr}\}$ is the vector of integrals. With the notation $\mathbb{R}^d_+$ we consider the set of non-negative $d$-dimensional vectors. In practice, this NP-hard problem is relaxed to the following non-negative least-squares problem
\begin{equation}
  \mathbf{Q}^{rh}=\argmin_{\mathbf{Q}\in\mathbb{R}^d_+}\lVert \mathbf{Q}\odot U_{hr}-\mathbf{c}\rVert^2_2,
\end{equation}
we solve it with the non-negative least-squares algorithm based on~\cite{lawson1995solving} and implemented in the Eigen library~\cite{eigenweb}.

The nonlinear manifold least-squares problem~\eqref{eq:PDEsystemDiscreteLinear}, are hyper-reduced with the ECSW method in the following formulations:
\begin{subequations}
  \begin{align}
    \mathbf{z}^t &= \argmin_{\mathbf{z}\in \mathbb{R}^r}\ \lVert \Tilde{\mathbf{Q}}^{hr}\odot \left(P_{r_h}G_{h, \delta t}(\boldsymbol{\mu}, P_{r_h}^s(\phi(\mathbf{z})), \{P_{r_h}^s(\phi(\mathbf{z}^s))\}_{s\in I_t})\oslash P_{r_h}\mathbf{W}\right)\rVert^{2}_{\mathbb{R}^{r_h}},\qquad&(\text{FB-ECSW})\label{eq:FB-ECSW}\\
    \mathbf{z}^t &= \argmin_{\mathbf{z}\in \mathbb{R}^r}\ \lVert \Tilde{\mathbf{Q}}^{hr}\odot \left(P_{r_h}G_{h, \delta t}(\boldsymbol{\mu}, P_{r_h}^s(\phi(\mathbf{z})), \{P_{r_h}^s(\phi(\mathbf{z}^s))\}_{s\in I_t})\oslash P_{r_h}\mathbf{W}_G\right)\rVert^{2}_{\mathbb{R}^{r_h}},\qquad&(\text{RB-ECSW})\label{eq:RB-ECSW}\\
    \mathbf{z}^t &= \argmin_{\mathbf{z}\in \mathbb{R}^r}\ \lVert P_{r_h}G_{h, \delta t}(\boldsymbol{\mu}, P_{r_h}^s(\phi(\mathbf{z})), \{P_{r_h}^s(\phi(\mathbf{z}^s))\}_{s\in I_t})\rVert^{2}_{\mathbb{R}^{r_h}},\qquad&(\text{C-ECSW})
    \label{eq:ECSW},
  \end{align}  
\end{subequations}
where $\Tilde{\mathbf{Q}}^{hr}\in\mathbb{R}^{r_h}$ is the quadrature weights vector obtained from the restriction of $\mathbf{Q}^{hr}\in\mathbb{R}^d$ to its non-zero entries. We also define $P_{r_h}\in\mathbb{R}^{d\times r_h}$ as the boolean matrix that selects the non-zero entries of $\mathbf{Q}^{hr}\in\mathbb{R}^d$ so that $\Tilde{\mathbf{Q}}^{hr}=P_{r_h}\mathbf{Q}^{hr}\in\mathbb{R}^d$, and as a consequence the projection onto the submesh $P_{r_h}^s\in\mathbb{R}^{d\times s}$.

Also for this case, depending on the choice of hyper-reduction basis chosen $U_{hr}=U$ (FB-ECSW) or $U_{hr}=U_G$ (RB-ECSW) or if it is performed reduced over-collocation (C-ECSW), there are three different formulations.

\subsection{Gradient-based adaptive hyper-reduction}
\label{subsec:adaptive-hr}

The most successful hyper-reduction strategy for our implementation of the NM-LSPG method is the adaptive reduced over-collocation method (C-UP) that we introduce now. For comments on the results see section~\ref{sec:discussion}. The method employs the standard formulation of the reduced over-collocation hyper-reduction method
\begin{equation}
  \mathbf{z}^t = \argmin_{\mathbf{z}\in \mathbb{R}^r}\ \lVert P_{r_h}G_{h, \delta t}(\boldsymbol{\mu}, P_{r_h}^s(\phi(\mathbf{z})), \{P_{r_h}^s(\phi(\mathbf{z}^s))\}_{s\in I_t})\rVert^{2}_{\mathbb{R}^{r_h}},\qquad(\text{C-UP})
\end{equation}
with the difference that the magic points are sampled adaptively during the time-evolution of the NM-LSPG trajectories. Its cost is amortized over the successive NM-LSPG time evaluations.

A heuristic approach relies on the positioning of the magic points where the sensitivities of the residual at time $t$ have greater components in $L^2$-norm. If we define the residual map at time $t$
\begin{equation}
  \mathbb{R}^r\ni\mathbf{z}\mapsto G_{h, \delta t}(\boldsymbol{\mu}, P_{r_h}^s(\phi(\mathbf{z})), \{P_{r_h}^s(\phi(\mathbf{z}^s))\}_{s\in I_t})=G_{h, \delta t}(\phi(\mathbf{z}))\in\mathbb{R}^d,\qquad
  G_{h, \delta t}:\mathbb{R}^r\rightarrow\mathbb{R}^d,
\end{equation}
losing for brevity the dependencies on the previous time steps,
its sensitivities with respect to the latent coordinate $\mathbf{z}^t$ at time $t$ are for all $i\in\{1,\dots,r\}$,
\begin{subequations}
  \begin{align}
    \mathbb{R}^d\ni\mathbf{J}_i &= \left.\frac{\partial G_{h, \delta t}(\phi(\mathbf{z}))}{\partial\mathbf{z}_i}\right|_{\mathbf{z}=\mathbf{z}^t}=\left[\left.\frac{\partial G_{h, \delta t}(\phi(\mathbf{z}))}{\partial\boldsymbol{\phi}}\right|_{\phi(\mathbf{z})=\phi(\mathbf{z}^t)}\right]_{d\times s}\circ\left[\left.\frac{\partial \phi(\mathbf{z})}{\partial\mathbf{z}_i}\right|_{\mathbf{z}=\mathbf{z}^t}\right]_{s\times 1}\\
    &=\left[\left.\frac{\partial G_{h, \delta t}(\phi(\mathbf{z}))}{\partial\boldsymbol{\phi}}\right|_{\phi(\mathbf{z})=\phi(\mathbf{z}^t)}\right]_{d\times s}\circ\left[P_{r_h}^s\circ f^{-1}_{\text{filter}}\right]_{s\times r_{\text{rSVD}}}\circ\left[\left.\frac{\partial \Tilde{\phi}(\mathbf{z})}{\partial\mathbf{z}_i}\right|_{\mathbf{z}=\mathbf{z}^t}\right]_{r_{\text{rSVD}}\times 1}\\
    &=\left[\mathbf{G}\right]_{d\times s}\left[\mathbf{F}\right]_{s\times r_{\text{rSVD}}}\left[\mathbf{\Phi}_i\right]_{r_{\text{rSVD}}\times 1}
  \end{align}
\end{subequations}  
with
\begin{equation}
  \mathbf{G}=\left.\frac{\partial G_{h, \delta t}(\phi(\mathbf{z}))}{\partial\boldsymbol{\phi}}\right|_{\phi(\mathbf{z})=\phi(\mathbf{z}^t)}\in\mathbb{R}^{d\times s},\quad\mathbf{F}=P_{r_h}^s\circ f^{-1}_{\text{filter}}\in\mathbb{R}^{s\times r_{\text{rSVD}}},\quad \mathbf{\Phi}=\left.\frac{\partial \Tilde{\phi}(\mathbf{z})}{\partial\mathbf{z}}\right|_{\mathbf{z}=\mathbf{z}^t}\in\mathbb{R}^{r_{\text{rSVD}}\times r}
\end{equation}
where with the chain rule we have highlighted the terms that compose the evaluation of the Jacobian matrix of the residual.

The exact evaluation of the Jacobian $\mathbf{\Phi}\in\mathbb{R}^{r_{\text{rSVD}}\times r}$ of the parametrization map $\Tilde{\phi}:\mathbb{R}^r\rightarrow\mathbb{R}^{r_{\text{rSVD}}}$ is efficient enough to be employed by our adaptive hyper-reduction procedure. However, the matrix multiplication $\mathbf{G}\mathbf{F}\mathbf{\Phi}$ is inefficient since $\mathbf{G}\in\mathbb{R}^{d\times s}$ depends on the number of degrees of freedom and the submesh size $s$. If we wanted to apply DEIM to recover the residual on the full-space $\mathbb{R}^d$ from $P_{r_h}G_{h,\delta t}$, so considering the sensitivities of $U_{hr}(P_{r_h}U_{rh})^{\dagger}P_{r_h}G_{h,\delta t}$, we would need to compute the pseudo-inverse $(P_{r_h}U_{rh})^{\dagger}$ which could become a heavy task if repeated multiple times during the NM-LSPG method.

Our solution consists instead in evaluating $O:U\subset\mathbb{R}^r\times \mathcal{P}\rightarrow\mathbb{R}^{d\times r}$,
\begin{equation}
  O(\mathbf{z}^t, \boldsymbol{\mu})=f^{-1}_{\text{filter}}\circ\left(\left.\frac{\partial \Tilde{\phi}(\mathbf{z})}{\partial\mathbf{z}}\right|_{\mathbf{z}=\mathbf{z}^t}-\left.\frac{\partial \Tilde{\phi}(\mathbf{z})}{\partial\mathbf{z}}\right|_{\mathbf{z}=\mathbf{z}^{t-1}}\right)=\left[\mathbf{W}\odot U\right]_{d\times r_{\text{rSVD}}}\left[\left.\frac{\partial \Tilde{\phi}(\mathbf{z})}{\partial\mathbf{z}}\right|_{\mathbf{z}=\mathbf{z}^t}-\left.\frac{\partial \Tilde{\phi}(\mathbf{z})}{\partial\mathbf{z}}\right|_{\mathbf{z}=\mathbf{z}^{t-1}}\right]_{r_{\text{rSVD}}\times r},
\end{equation}
considering the current time-step $t$ and the previous one $t-1$. In this way, we are employing the sensitivities of the neural network itself, rather then including also the information coming from the residual through the Jacobian $G\in\mathbb{R}^{d\times s}$.

We select only the first $r_h$ cells of the computational domain that maximize
\begin{equation}
  \label{eq:ahrEq}
  \argmax_{k=\{1,\dots,M\}}\left(\argmax_{j=1,\dots,r} \left(\sum_{i=1}^c O_i(\mathbf{z}^t, \boldsymbol{\mu})^2\right)_{kj}\right)
\end{equation}
where $O_i(\mathbf{z}^t, \boldsymbol{\mu})\in\mathbb{R}^M$ for $i\in\{1,\dots,c\}$ are the degrees of freedom corresponding to the $i$-th physical field composing the full-state $\mathbf{U}_h\in X_h\sim\mathbb{R}^d$. When we update the magic points and consequently the submesh every $n$ time steps, we use the notation C-Un, for example, C-U50 for an adaption of the magic points every $50$ time steps.

\begin{figure}[ht!]
  \centering
  \includegraphics[width=0.32\textwidth, trim={0 0 500 0}, clip]{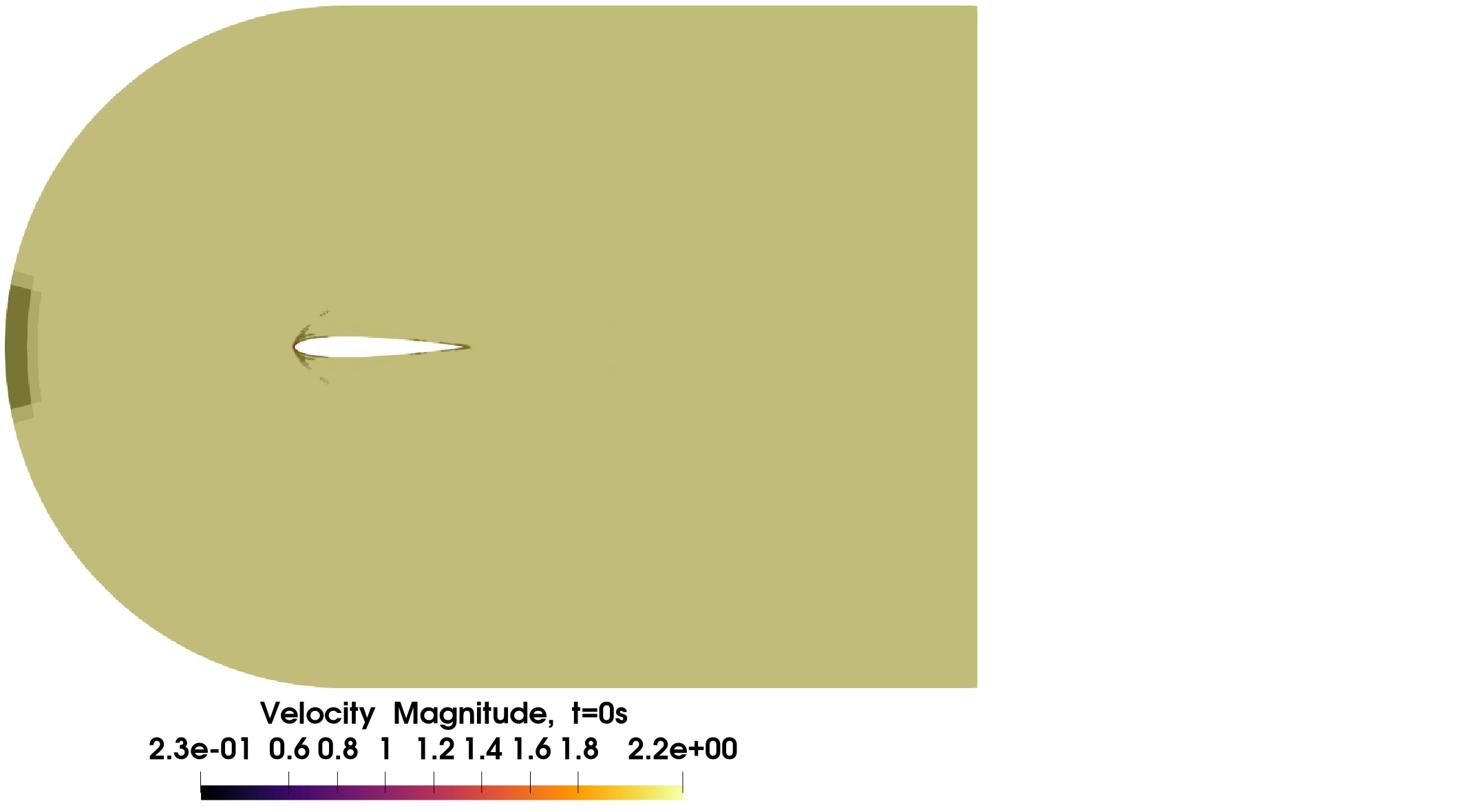}
  \includegraphics[width=0.32\textwidth, trim={0 0 500 0}, clip]{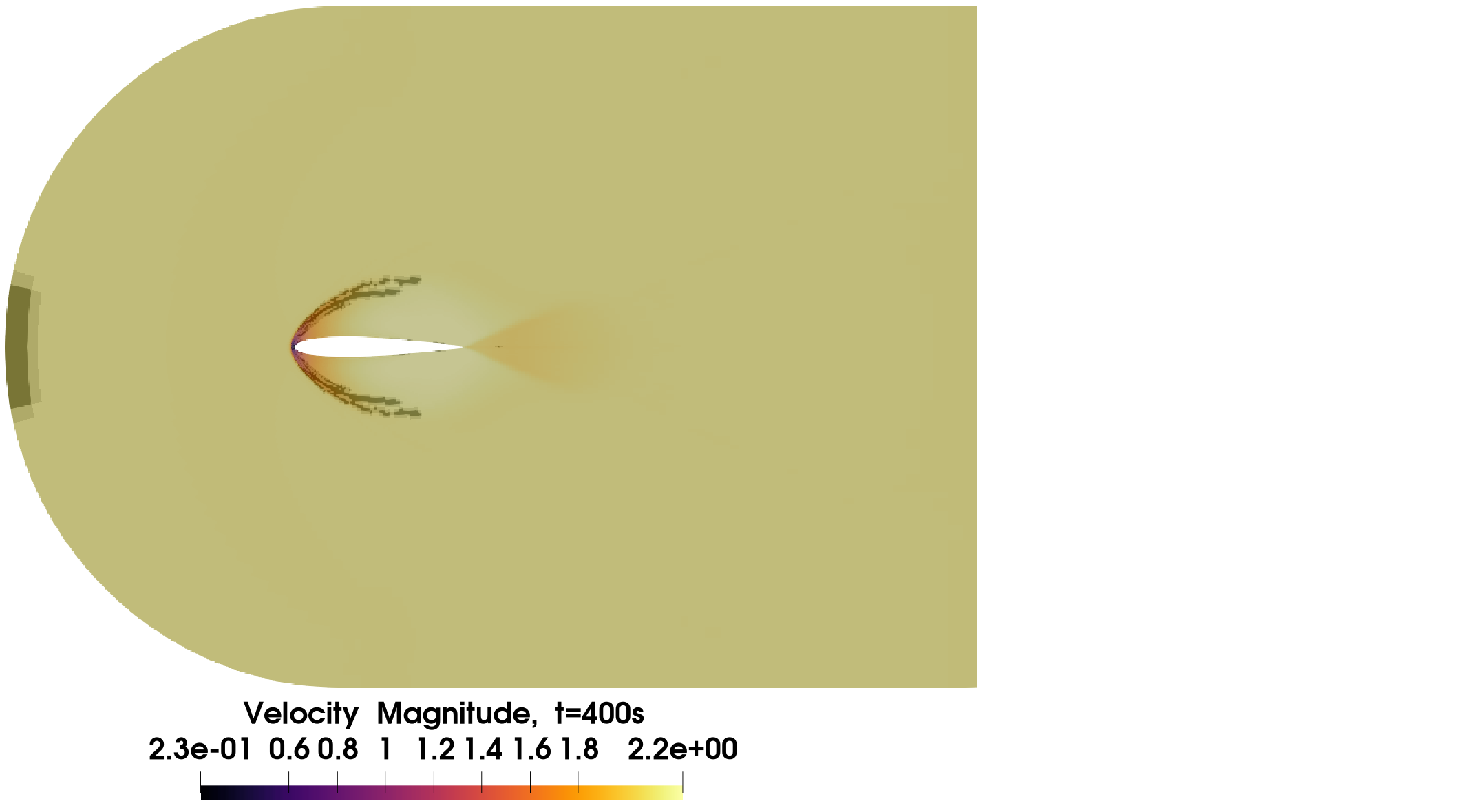}
  \includegraphics[width=0.32\textwidth, trim={0 0 500 0}, clip]{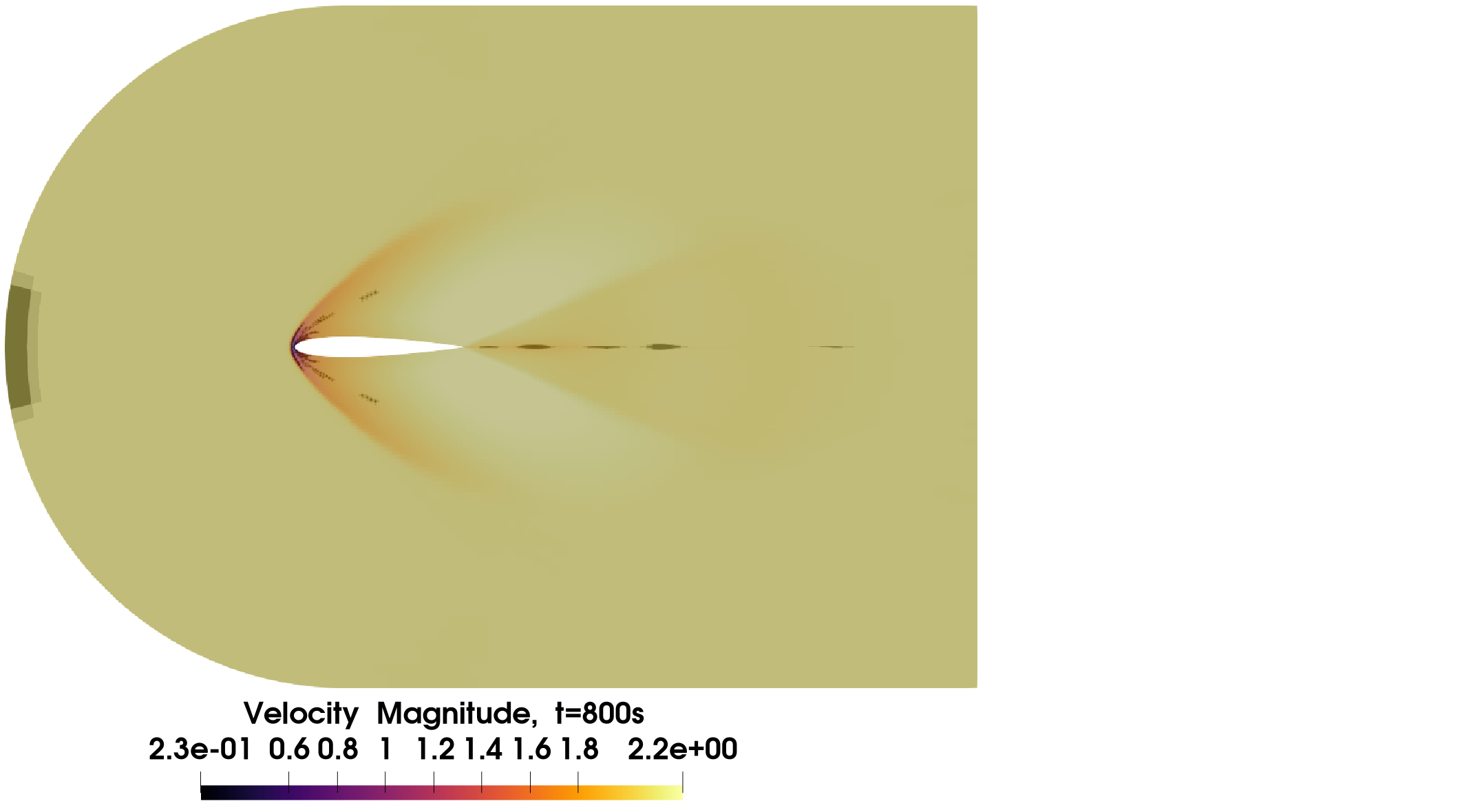}\\
  \includegraphics[width=0.32\textwidth, trim={0 0 500 0}, clip]{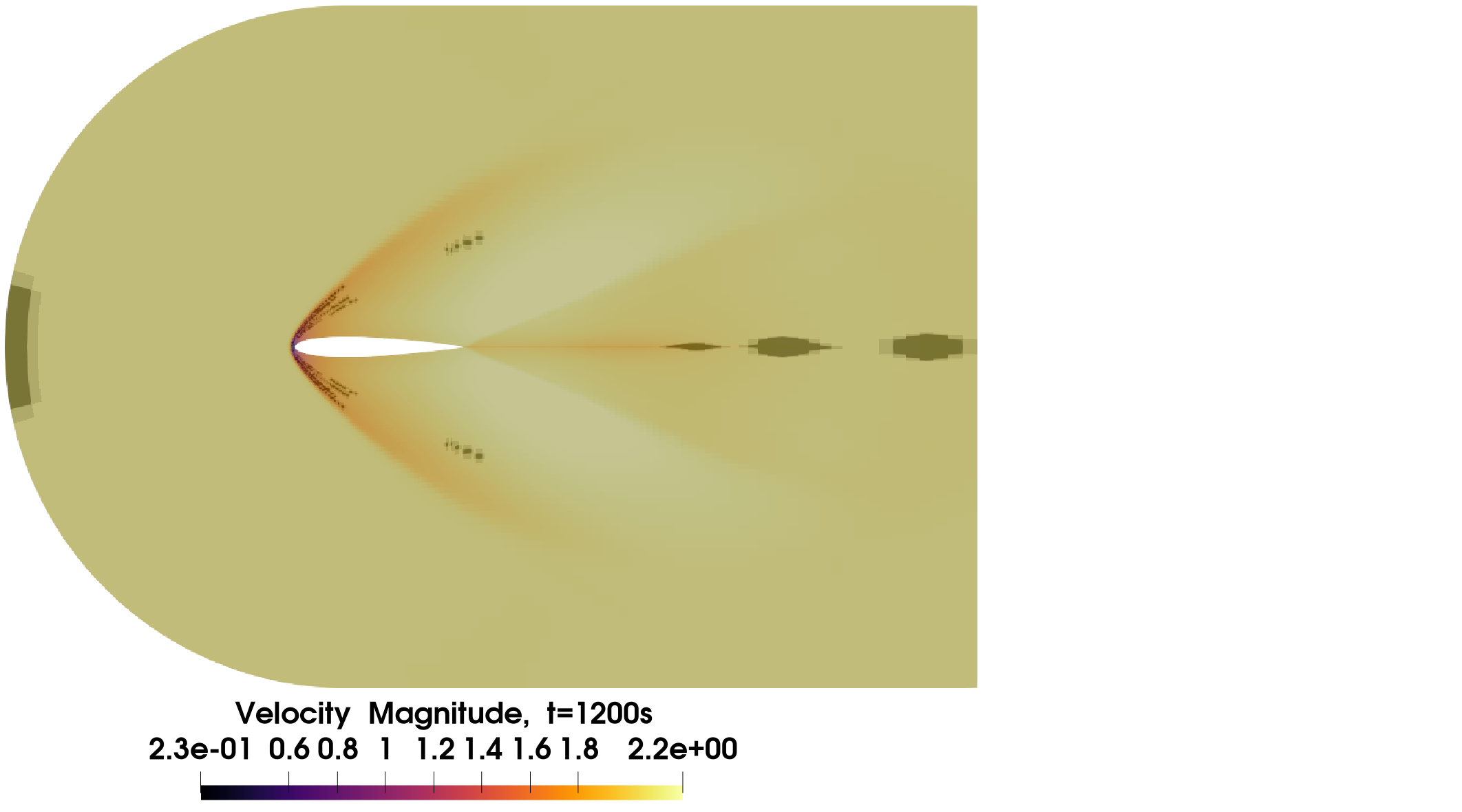}
  \includegraphics[width=0.32\textwidth, trim={0 0 500 0}, clip]{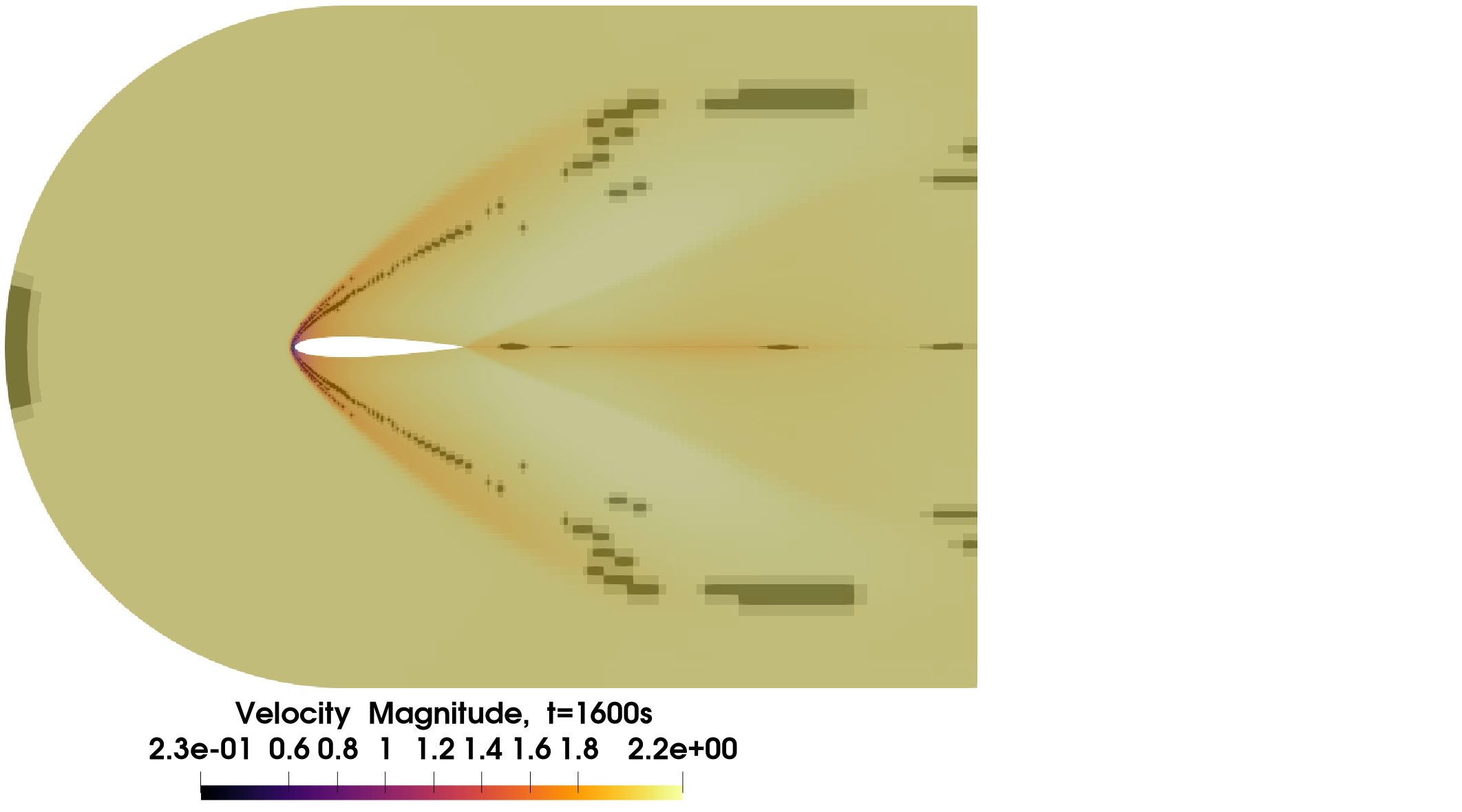}
  \includegraphics[width=0.32\textwidth, trim={0 0 500 0}, clip]{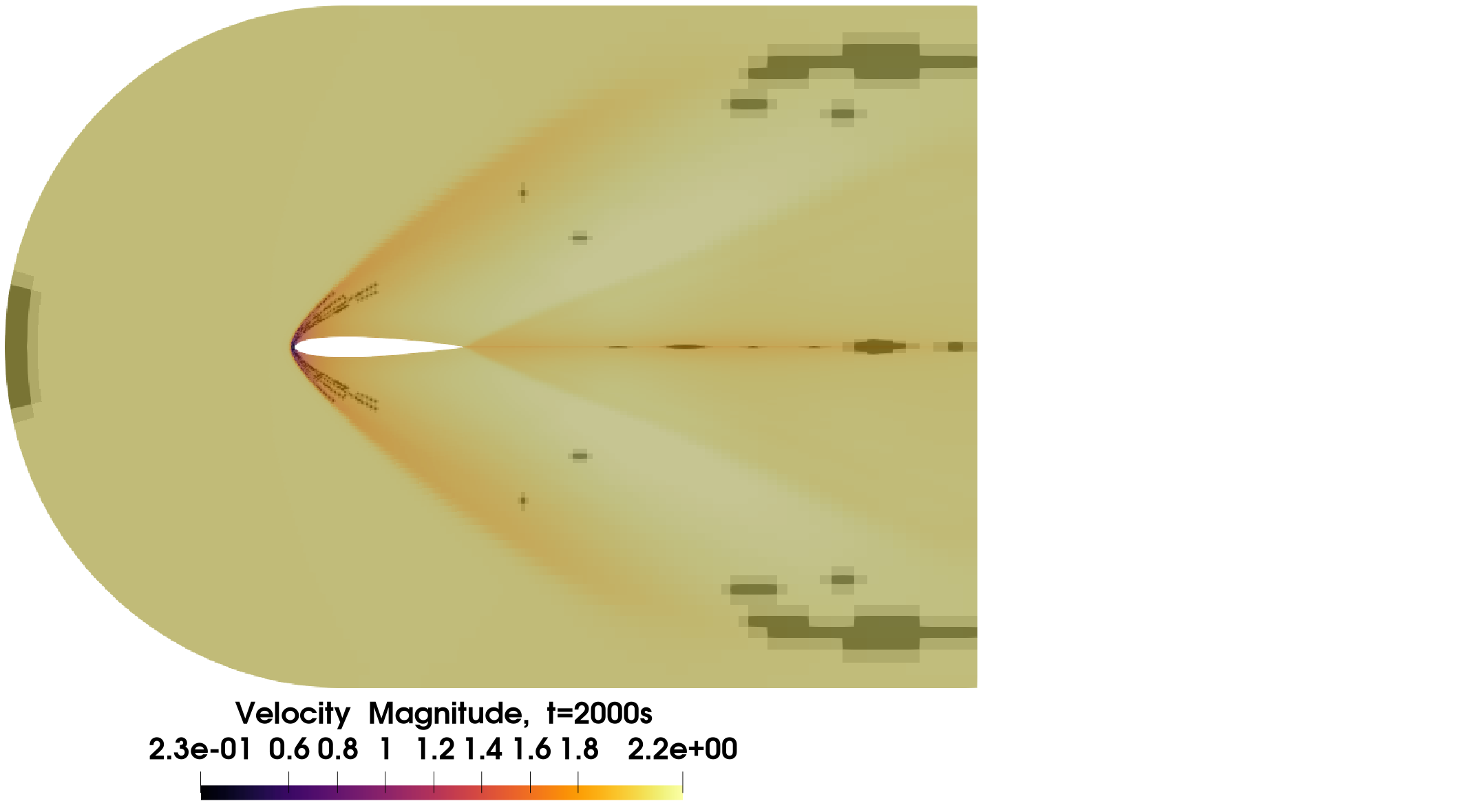}
  \caption{Visual example of adaptive gradient-based reduced over-collocation. Predicted velocity magnitude time snapshots, for $t\in\{0, 400, 800, 1200, 1600, 2000\}$, of the 2d compressible Navier-Stokes equations' model described in subsection~\ref{subsubsec:sudomains_finer_airfoil} for a value of the Mach number $\text{Ma}=2.12422656$, the fourth of five test parameters. It can be seen that the $700$ collocation nodes in black and the corresponding submesh in grey, updated with gradient-based adaptive hyper-reduction through equation~\eqref{eq:ahrEq} every $50$ time steps, follow the transient dynamics. The time instants reported refer to the reference time interval $t\in[0, 2500]$, the real time instants are obtained applying the scaling reported in equation~\eqref{eq:timesAirfoil}}
  \label{fig:adaptive}
\end{figure}

The methodology presented derives from heuristic considerations and the necessity to keep the computational costs as low as possible. Experimentally it is shown that it is able to track the main moving features of our transient numerical simulations, thus adapting the magic points' position close to the most informative regions of the computational domain. See Figure~\ref{fig:adaptive}. Additional collocation nodes are imposed on the boundaries to force the satisfaction of boundary conditions as can be noticed in Figure~\ref{fig:adaptive}, at the inflow left boundary.

\section{Local nonlinear manifold}
\label{sec:loc}
Sometimes the $L^2$-norm relative reconstruction error~\eqref{eq:rec_rsvd} cannot approximate with sufficient accuracy the nonlinear solution manifold in terms of reconstruction error, or the autoencoder architecture $\phi\circ\psi:\mathbb{R}^d\rightarrow\mathbb{R}^d$ or $\Tilde{\phi}\circ\Tilde{\psi}:\mathbb{R}^{r_{\text{rSVD}}}\rightarrow\mathbb{R}^{r_{\text{rSVD}}}$ does not meet the tolerance requirement~\eqref{eq:nonlinear_rec_err}. This happens when there are regions in the parameter space that correspond to less correlated solutions in the solution manifold.

A possible solution is to partition the parameter space in subdomains where the approximation properties of the rSVD modes and the autoencoder are satisfactory for the problem at hand. There are many implementations of local ROMs, often under the name of dictionary-based ROMs~\cite{daniel2020model}. Generally, they also can be efficiently applied to our framework, thanks to the definition of our autoencoder through linear projections. 

The setting is introduced only for two communicating local solution manifolds. For $i=1, 2$, with the notations 
\begin{subequations}
  \begin{align}
    &\mathcal{P}^i_{n_{\text{train}, i}}\subset\mathcal{P}_{n_{\text{train}}}, \qquad&\text{(local parameter subset)}\\
    &\phi^i:\mathbb{R}^{r}\rightarrow X_h\sim\mathbb{R}^d,\quad\Tilde{\phi}^i:\mathbb{R}^{r}\rightarrow\mathbb{R}^p,\quad\psi^i:X_h\sim\mathbb{R}^d\rightarrow \mathbb{R}^r,\quad\Tilde{\psi}^i:\mathbb{R}^p\rightarrow\mathbb{R}^r,\qquad&\text{(local autoencoder)}\\
    &f^i_{\text{filter}}:X_h\sim\mathbb{R}^d\rightarrow\mathbb{R}^p,\quad(f^i_{\text{filter}})^{-1}:\mathbb{R}^p\rightarrow X_h\sim\mathbb{R}^d,\qquad&\text{(local linear projections)}\\
    &A^i_{n_{\text{train}, i}}\in\mathbb{R}^{d\times n_{\text{train}, i}},\quad \mathbf{W}^i\in\mathbb{R}^d,\quad U^i\in\mathbb{R}^{d\times r_{\text{rSVD}}}\qquad&\text{(local rSVD quantities)}
  \end{align}
\end{subequations}
we denote the corresponding parametric subsets, the decoder maps, the encoder maps, the linear filter/transform maps, the training snapshots matrices, the normalizing vectors and the local rSVD basis of the two local solution manifolds. In principle, the local latent dimensions $r=r^1=r^2$ with our notation, can be different $r^1\neq r^2$ and the same is valid for $p=r_{\text{rSVD}}=p^1=p^2$, that is $p^1=r_{\text{rSVD}}^1\neq r_{\text{rSVD}}^2=p^2$. In fact, it is possible to adapt the latent and linear filter dimensions of the nonlinear approximating manifold parametrized by $\phi^i$ to the local Kolmogorov n-width decay of the subset of the parameter space under consideration. In a sense, the rationale is similar to the one behind hp-FEM or hp-DG methods.

Gluing two local solution manifolds requires care especially because it may not be guaranteed that the corresponding full-states are close $((f^2_{\text{filter}})^{-1}\circ f^2_{\text{filter}}\circ (f^1_{\text{filter}})^{-1})(\Tilde{\mathbf{U}_h^1})=(f^1_{\text{filter}})^{-1}(\Tilde{\mathbf{U}_h^1})$ to each other with sufficient accuracy. Many techniques have been developed to tackle this problem, but, for the moment, we will study only the most simple one. In fact, a possible way to avoid this consists in just overlapping the training datasets $A_{n_{\text{train}, 1}}\in\mathbb{R}^{d\times n_{\text{train}, 1}}$ and $A_{n_{\text{train}, 2}}\in\mathbb{R}^{d\times n_{\text{train}, 2}}$, that is considering a bigger intersection of their corresponding parameters subsets $\mathcal{P}^1_{n_{\text{train}, 1}}\subset\mathcal{P}_{n_{\text{train}}}$ and $\mathcal{P}^2_{n_{\text{train}, 2}}\subset\mathcal{P}_{n_{\text{train}}}$, $\mathcal{P}^1_{n_{\text{train}, 1}}\cap\mathcal{P}^2_{n_{\text{train}, 2}}\neq\varnothing$.

\begin{figure}[ht!]
  \centering
  \includegraphics[width=0.8\textwidth, trim={0 0 180 0}, clip]{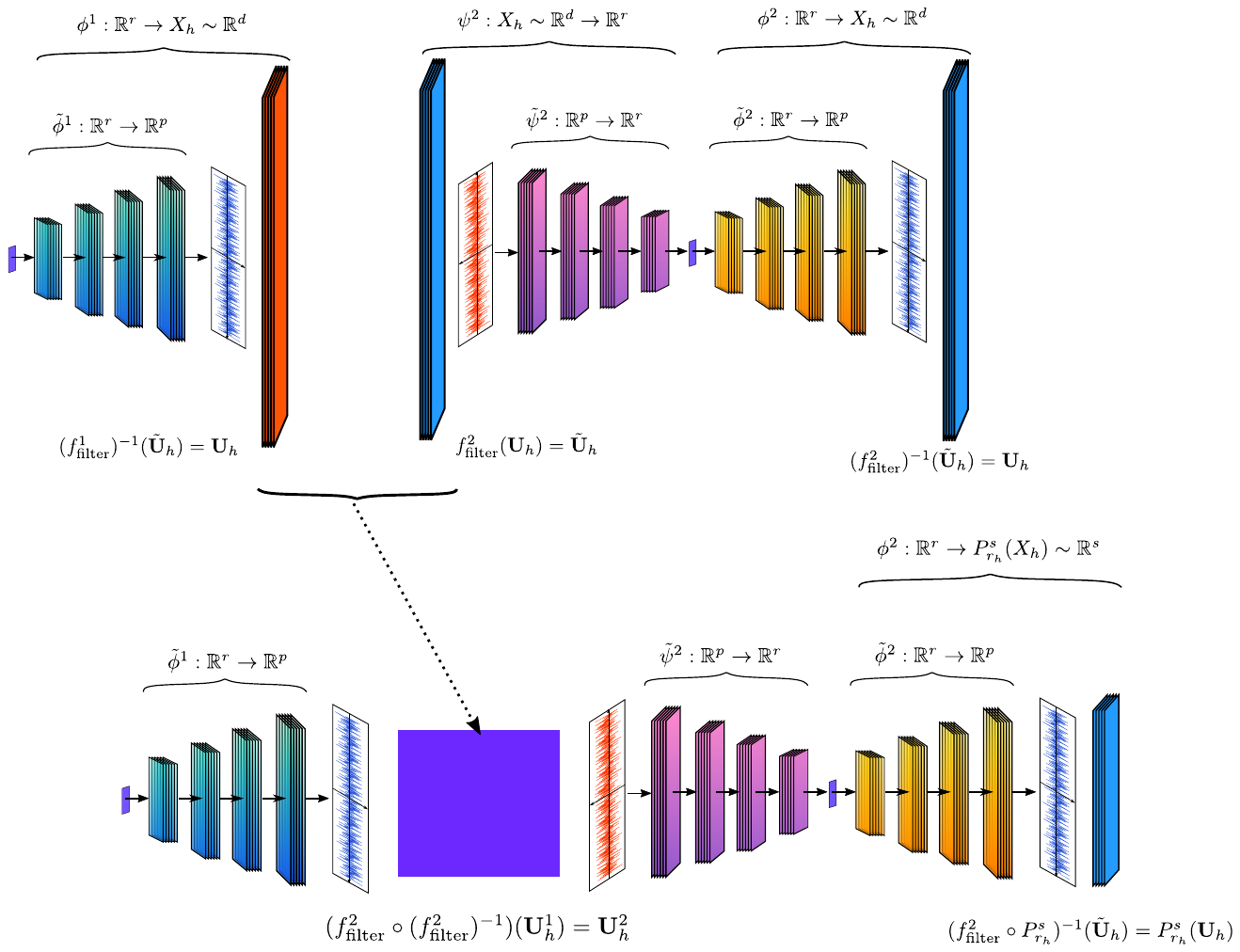}
  \caption{\textbf{Above:} the decoder $\phi^1$ of subdomain $1$ on the left, and the encoder $\psi^2$ followed by the decoder $\phi^2$ of subdomain $2$. To pass from one subdomain whose solution manifold is spanned by the rSVD modes associated to $(f^1)^{-1}_{\text{filter}}$ to the one spanned by the rSVD modes associated to $f^2_{\text{filter}}$ and $(f^2)^{-1}_{\text{filter}}$, an efficient change of basis needs to be performed during the online stage. The dependency on the number of dofs is avoided by multiplying in the offline stage the rSVD basis associated to subdomain $1$ with the ones associated to subdomain $2$. \textbf{Below: } the same maps as above, only now an efficient change of basis between subdomains $\Tilde{\mathbf{U}}_h^2=f^{12}_{\text{filter}}(\Tilde{\mathbf{U}}_h^1)=(f^2_{\text{filter}}\circ (f^1_{\text{filter}})^{-1})(\Tilde{\mathbf{U}}_h^1)$ can be used online, since it is precomputed offline as a change of basis matrix of sizes $p^2\times p^1$.}
  \label{fig:change_of_basis}
\end{figure}

In our case, we have that the change of basis linear map $f^{12}_{\text{filter}}=f^2_{\text{filter}}\circ (f^1_{\text{filter}})^{-1}:\mathbb{R}^{r_{\text{rSVD}}}\rightarrow\mathbb{R}^{r_{\text{rSVD}}}$ is computed offline as
\begin{equation}
  \label{eq:change of basis}
  \Tilde{\mathbf{U}}_h^2=f^{12}_{\text{filter}}(\Tilde{\mathbf{U}}_h^1)=(f^2_{\text{filter}}\circ (f^1_{\text{filter}})^{-1})(\Tilde{\mathbf{U}}_h^1) = (U^2)^T\left(\left(\left(\mathbf{W}^1\odot U^1\right)\Tilde{\mathbf{U}}_h^1\right)\oslash\mathbf{W}^2\right).
\end{equation}
This change of basis between communicating local solution manifolds is represented schematically in Figure~\ref{fig:change_of_basis}.

We will consider only two local solution manifolds such that only the time interval is partitioned, see section~\ref{subsubsec:sudomains_finer_airfoil}.

\section{Numerical experiments}
\label{sec:results}
We will test the presented methodology on two challenging benchmarks from the point of view of model order reduction, as they evidently suffer from a slow Kolmogorov n-width decay. Employing classical linear projection-based ROMs would be unfeasible as they would need more than hundreds of rSVD modes. For the comparison of classical methods with intrusive ones exploiting nonlinear approximants see~\cite{kim2022fast,barnett2022neural}.

The test cases we will present are developed with the open source CFD software library \texttt{OpenFoam}~\cite{Weller1998}. Being rather small test cases, discretized with the finite volumes method on meshes with \textbf{4500}, \textbf{32160} and \textbf{198633} cells, the speedups obtained are relatively small. However, since our ROMs are independent on the number of dofs, as long as increasing the resolution does not bring to a slower Kolmogorov n-width decay, more evident speedups can be achieved. Compared to the finite differences method employed in~\cite{kim2022fast,barnett2022neural}, the FVM implementation in \texttt{OpenFoam} is highly optimized and it is therefore more challenging to achieve a speedup for small test cases. The new methdology is implemented on top of the open-source software library for model order reduction \texttt{ITHACA-FV}~\cite{Stabile2017CAIM}.

The first test case we consider involves the compressible Navier-Stokes equations (CNS) in the supersonic regime. We will consider as parameters time and the Mach number $\text{Ma}\in[2, 5]$ imposed at the inflow boundary. In order to show, that the method scales well increasing the number of dofs, we employ two meshes: a coarse one with \textbf{4500} cells and a finer one with \textbf{32160} cells, see Figure~\ref{fig:domain_airfoil}. On the first mesh in subsection~\ref{subsec:coarseAirfoil}, we test different hyper-reduction methods. On the second mesh in subsection~\ref{subsubsec:sudomains_finer_airfoil} we test the use of local nonlinear manifold approximants and the straight-forward change of basis we have presented in section~\ref{sec:loc}. The \texttt{OpenFoam} solver we employ is the sonicFoam~\cite{marcantoni2012high} solver.

The second test case we consider is the incompressible turbulent flow around the Ahmed body (INS). The parameters are time and the slant angle of the Ahmed body. The test case presented is introduced and studied in~\cite{zancanaro2021hybrid}. Steady state solutions are obtained with the solver SIMPLE~\cite{issa1991solution}(Semi-Implicit Method for Pressure-Linked Equations), however we employ PISO~\cite{patankar1983calculation} (Pressure Implicit with Splitting of Operator) since we want to reduce the transient dynamics that has a slow Kolmogorov n-width decay instead. Adding time as a parameter brings to a much more complex solution manifold to approximate.

The specifics of the convolutional autoencoders used are reported in the Appendix~\ref{sec:arch} along with the training costs and other hyper-parameters. For all the CAE trainings we employed the ADAM~\cite{kingma2014adam} optimizer with initial learning rate of $0.001$ and a scheduler that halves it every $200$ epochs if the validation loss has not improved. Every architecture is trained for $3000$ epochs on a single GPU NVIDIA Quadro RTX 4000. The wallclock time expended on training is between $1$ hour and $1$ hour and half for the test cases in~\ref{subsec:coarseAirfoil},~\ref{subsubsec:sudomains_finer_airfoil} and~\ref{subsubec:smallAhmed}. Substantial computational savings depend on the choice of the architecure, especially on the fact that the training is independent with respect to the number of dofs since the inputs and outputs of the autoencoder belong to the lower dimensional space generated by the rSVD basis.

In order not to interrupt the presentation of the numerical results, we postponed to Appendix~\ref{sec:snaps} the showcase of some predicted solutions of each test case: Figure~\ref{fig:snapsCoarseAirfoil} refers to the \textbf{CNS1} test case in~\cref{subsec:coarseAirfoil}, Figure~\ref{fig:snapsFinerAirfoil} refers to the \textbf{CNS2} test case in~\cref{subsubsec:sudomains_finer_airfoil} and finally Figures~\ref{fig:ahmedSmallSnaps} and~\ref{fig:adaptiveSmall} refers to the test case \textbf{INS1} in~\cref{subsubec:smallAhmed}.

\subsection{Supersonic flow past a NACA airfoil}
\label{subsec:compr}

\begin{figure}[ht!]
  \centering
  \includegraphics[width=0.37\textwidth]{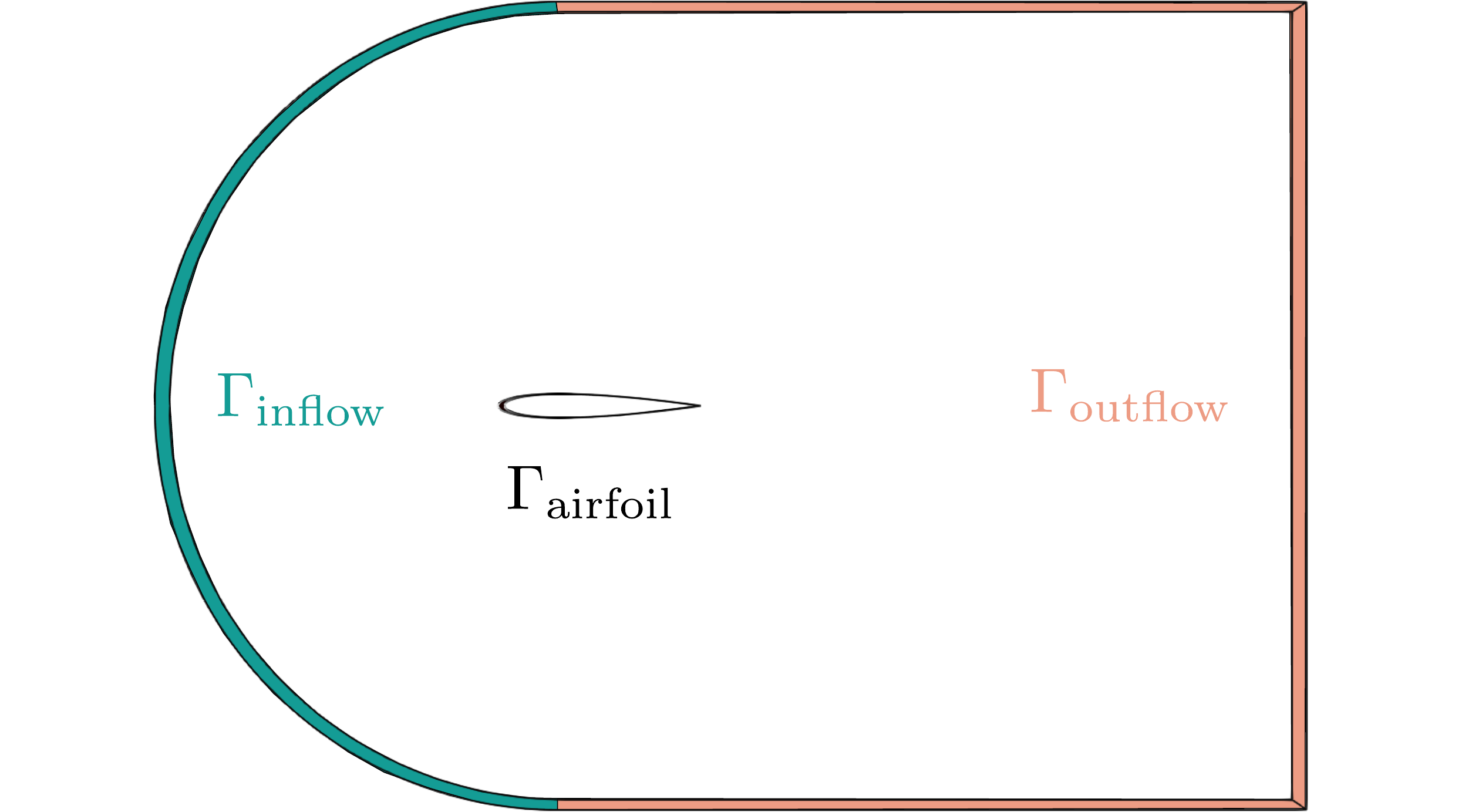}
  \includegraphics[width=0.3\textwidth, trim={160 0 100 0}, clip]{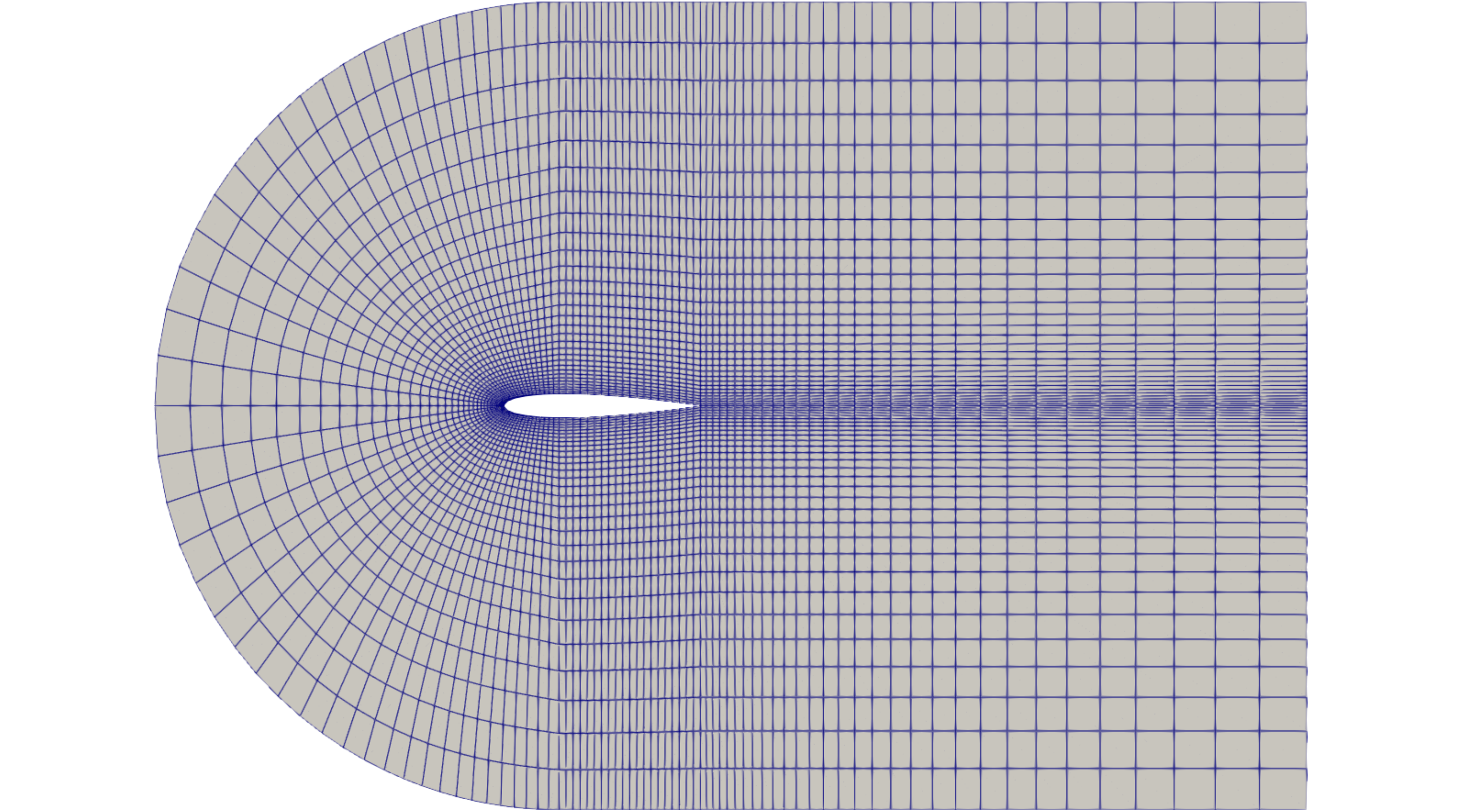}
  \includegraphics[width=0.3\textwidth, trim={160 0 100 0}, clip]{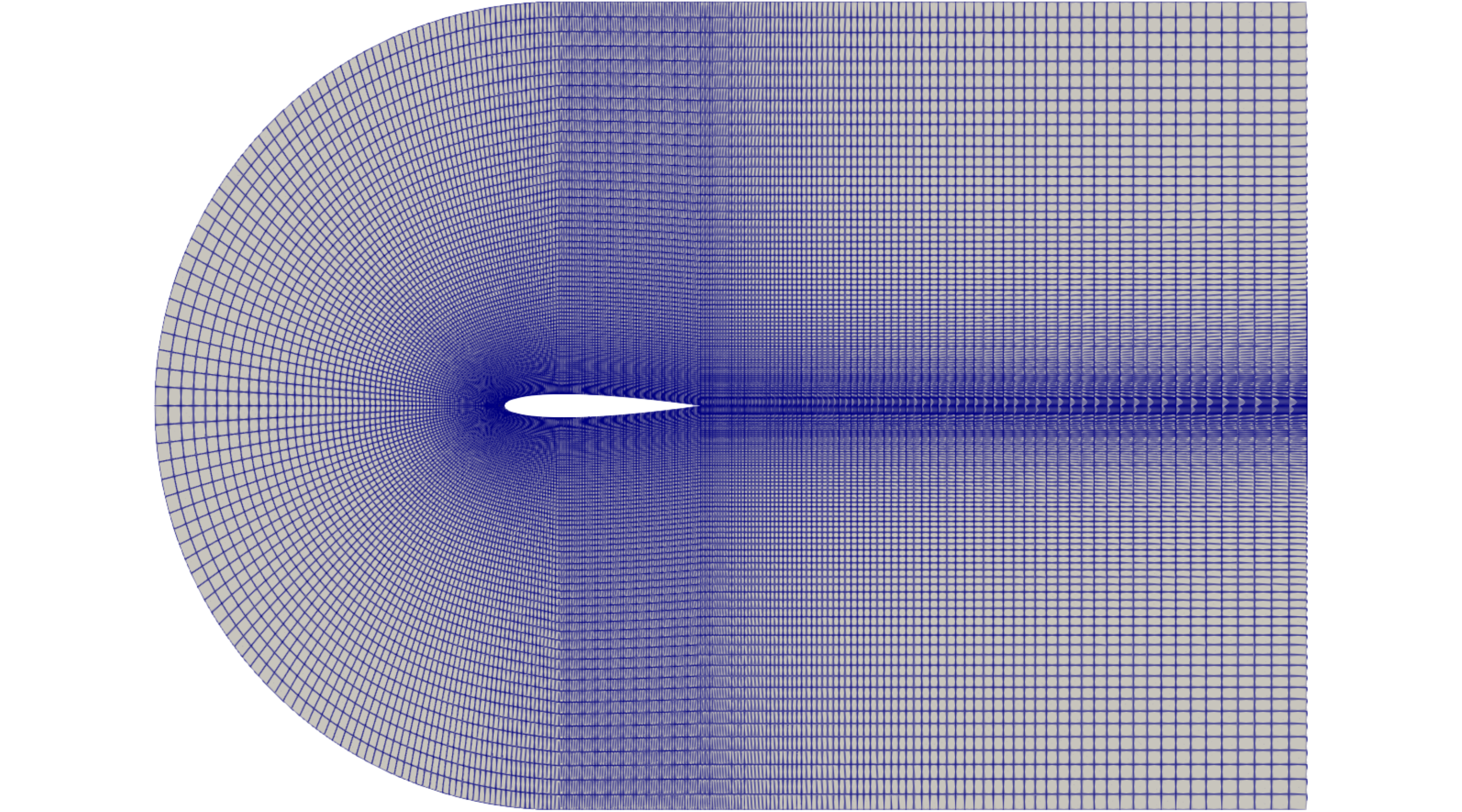}
  \caption{\textbf{Left:} computational domain of the compressible flow past a NACA airfoil test case. \textbf{Center:} coarse mesh with $4500$ cells. \textbf{Right:} finer mesh with $32160$ cells.}
  \label{fig:domain_airfoil}
\end{figure}

The first nonlinear time-dependent parametric PDE model we consider are the compressible Navier-Stokes equations for a perfect diatomic gas with ratio of specific heat $\gamma=7/5$. A speed of sound of $c=\sqrt{\gamma \frac{R}{M} T }=1 \text{m/s}$ is imposed through the choice of molar mass $M$ at a temperature of $T=1 K$, where $R$ is the universal gas constant. The following system of PDEs is solved on the 2d computational domain $\Omega_h\subset\mathbb{R}^2$ shown in Figure~\ref{fig:domain_airfoil}:
\begin{subequations}
  \begin{align}
    \partial_t \rho + \nabla\cdot (\rho\mathbf{u})=0&,\quad &\text{(mass conservation)}\label{eq:mass}\\
    \partial_t (\rho\mathbf{u}) + \nabla\cdot (\rho\mathbf{u}\otimes\mathbf{u})+\nabla p - \nabla\cdot \left(\nu\left(\nabla\mathbf{u}+\nabla\mathbf{u}^T+\tfrac{1}{3}(\nabla\cdot \mathbf{u})I_d\right)\right)= 0&,\quad &\text{(momentum conservation)}\label{eq:momentum}\\
    \partial_t \rho(e+K) + \nabla\cdot ((\rho(e+K)+p)\mathbf{u})- \nabla\cdot \left(\nu\left(\nabla\mathbf{u}+\nabla\mathbf{u}^T+\tfrac{1}{3}(\nabla\cdot \mathbf{u})I_d\right)\mathbf{u}\right)= 0&,\quad &\text{(energy conservation)}\label{eq:energy}\\
    e=Tc_v,\quad pc_v = \rho Re,\quad K=\tfrac{1}{2}\rho\mathbf{u}^2&,\quad &\text{(state equations)},\label{eq:state}
  \end{align}
\end{subequations}
the parameters we consider are time and the inlet Mach number $\text{Ma}=|\mathbf{u}|/c=|\mathbf{u}|$, $\mu=\text{Ma}\in[2, 5]$, $t\in V_{\mu}=\{1,\dots, N_{\mu}\}$. The viscosity is fixed at $\nu=1e-5$. The boundary conditions are imposed at the inflow $\Gamma_{\text{inflow}}$, outflow $\Gamma_{\text{outflow}}$ and airfoil $\Gamma_{\text{airfoil}}$ boundaries, see Figure~\ref{fig:domain_airfoil}. The initial and boundary conditions for the velocity $\mathbf{u}$, pressure $p$ and temperature $T$ fields are:
\begin{align*}
  \begin{cases}
    \mathbf{u}(\x, t) = \text{Ma},\quad & (\x, t)\in(\mathring{\Omega}_h\times \{t=0\})\cup(\Gamma_{\text{inflow}}\times {[0, T_{\boldsymbol{\mu}}]})\\
    p(\x, t) = 1 ,\quad & (\x, t)\in(\mathring{\Omega}_h\times \{t=0\})\cup(\Gamma_{\text{inflow}}\times {[0, T_{\boldsymbol{\mu}}]})\\
    T(\x, t) = 1,\quad & (\x, t)\in(\mathring{\Omega}_h\times \{t=0\})\cup(\Gamma_{\text{inflow}}\times {[0, T_{\boldsymbol{\mu}}]})
  \end{cases},\\
  \begin{cases}
    \mathbf{n}\cdot\nabla\mathbf{u}(\x, t) = 0,\quad &\x\in\Gamma_{\text{outflow}}\\
    \text{(non-reflective)} ,\quad &\x\in\Gamma_{\text{outflow}}\\
    \mathbf{n}\cdot\nabla T(\x, t) = 0,\quad &\x\in\Gamma_{\text{outflow}}
  \end{cases},\qquad
  \begin{cases}
    \mathbf{u}(\x, t) = 0,\quad &\x\in\Gamma_{\text{airfoil}}\\
    \mathbf{n}\cdot\nabla p(\x, t) = 0 ,\quad &\x\in\Gamma_{\text{airfoil}}\\
    \mathbf{n}\cdot\nabla T(\x, t) = 0,\quad &\x\in\Gamma_{\text{airfoil}}
  \end{cases},
\end{align*}
where non-reflective boundary conditions are imposed on the pressure field at the outflow boundaries. 

The Mach number training and test instances are sampled from the parameter space $\mathcal{P}=[2, 5]\times V_{\mu}$, such that the Mach angle $\alpha$ is sampled uniformly, and the time step $\Delta t$ and final time $T_{\mu}$ are chosen depending on the Mach number from the reference time step $(\Delta t)_{\text{ref}}=0.001$ and final time $(T_{\mu})_{\text{ref}}=2.5$:
\begin{equation}
  \label{eq:timesAirfoil}
  \alpha = \arcsin{\frac{1}{\text{Ma}}},\quad \Delta t=(\Delta t)_{\text{ref}}\frac{\text{Ma}_{\text{ref}}}{\text{Ma}},\quad T_{\mu}=(T_{\mu})_{\text{ref}}\frac{\text{Ma}_{\text{ref}}}{\text{Ma}}
\end{equation}
in this way the training and test time series have the same length even if the final times are different.

We take into account two different meshes a coarse one with \textbf{4500} cells and \textbf{2700} dofs and a finer one with \textbf{32160} cells and \textbf{192960} dofs. We will refer to these two test cases with the notation CNS-1 for the coarse mesh and CNS-2 for the finer mesh. The employment of two meshes permits us to show that our methodologies achieve more significant speedups when the number of dofs is increased since they are independent on the number of dofs: the timings and relative speedups can be observed from Tables~\ref{tab: decoder costs coarseAirfoil} and~\ref{tab: decoder costs finer CNS}.

 The solver we will employ is \texttt{OpenFoam}'s~\cite{Weller1998} pressure-based transonic/supersonic solver for compressible gases sonicFoam~\cite{marcantoni2012high}. SonicFoam algorithm~\ref{alg:sonic} solves for the solution at the $n$-th time instant follows PIMPLE predictor-corrector scheme, a combination of PISO~\cite{patankar1983calculation} and SIMPLE~\cite{issa1991solution}. In algorithm~\ref{alg:sonic}, we highlight the predictor-corrector scheme and the main steps. We employ the Euler scheme in time. Starting from the previous fields $\mathbf{u}^{n}$ velocity, $\rho^{n}$ density, $e^{n}$ internal energy, and $p^{n}$ pressure, the solutions at time step $n+1$ are evaluated. After an intermediate density evaluation $\rho^*$ in line $2$ corresponding to the continuity equation~\eqref{eq:mass}, begins the PIMPLE corrector loop from line $3$: this outer loop comes from SIMPLE and relaxes the intermediate solution fields after every iteration. Then, at line $4$, the intermediate velocity field $\mathbf{u}^*$ is evaluated implicitly solving the momentum equation~\eqref{eq:momentum}: the diagonal $A[\rho^*, \Delta t]$ and over-diagonal $H[\mathbf{u}^*, \rho^*, \Delta t]$ parts of the system are highlighted along with the pressure contribution $\nabla p^*$ since they will be employed later in the pressure-Poisson equation at line $10$. The energy predictor step at line $5$ is evaluated afterwards corresponding to the energy conservation~\eqref{eq:energy}. Then, the thermodynamics properties corresponding to the state equations~\eqref{eq:state} are corrected and the PISO pressure corrector loop begins at line $7$. Inside the non-orthogonal corrector loop in case of non-orthogonal meshes, the pressure-Poisson equation is repeatedly solved and, afterwards, the velocity is corrected to satisfy the continuity equation and the density is updated with the new pressure through the equations of state. Many steps have not been reported for simplicity, for a more detailed analysis see~\cite{marcantoni2012high}. What we want to show is that in comparison, the nonlinear least-squares Petrov-Galerkin scheme for the compressible Navier-Stokes (CNS) equations (NM-LSPG-CNS) is simplified, at each $n$-th time step and $i$-th intermediate optimization step in algorithm~\ref{alg:nmsonic}: having a solution manifold, trained on a previous database, as ansatz space, the corrector loops can be avoided and only the residual evaluation of the mass~\eqref{eq:mass}, momentum~\eqref{eq:momentum}, energy~\eqref{eq:energy} and pressure-Poisson equations are needed. For the same reasoning the orthogonal corrector loops can be avoided as the solutions searched for on the approximant nonlinear manifold should be already corrected.

 Along the lines of the previous considerations, we can employ larger time steps. In fact, it will be shown that using four times the reference time step, that is $(\Delta t)_{\text{ref}}=0.004$, brings a speedup also in the case of a coarse mesh, see Table~\ref{tab: decoder costs coarseAirfoil}.

\IncMargin{1em}
\vspace{0.5cm}
\begin{figure}[!t]
  \caption{\textbf{INS}. Comparison between the $n$-th time instant iterations of the full-order model numerical scheme sonicFoam (\textbf{Left}) and the nonlinear manifold least-squares Petrov Galerkin (NM-LSPG-CNS) method (\textbf{Right}).}
\begin{minipage}{0.47\linewidth}
\begin{algorithm}[H]
  \caption{sonicFoam $n$-th iteration}\label{alg:sonic}
  Start with an initial velocity field $\mathbf{u}^{n}$, density field $\rho^{n}$, internal energy field $e^{n}$, and pressure field $p^{n}$ at the $n$-th time step.\\
  Density prediction step:
  $$\rho^*-\rho^n + \Delta t \nabla\cdot (\rho^*\mathbf{u}^n)=0$$\\
  \While{Pressure-velocity PIMPLE corrector loop}{
    Momentum predictor step:
    $$ A[\rho^{*}]\mathbf{u}^{*} = H[\mathbf{u}^{*}, \rho^{*}]-\nabla p^{*}.$$\\
    Energy predictor step:
    \begin{align*}
      \partial_t (\rho^* (e^{*} + K^{*})) &+ \nabla\cdot(\rho^*\mathbf{u}(e^{*}+K^{*})+\\
      &\mathbf{u}^*p^*)= 0
    \end{align*}\\
    Correct thermodynamics properties.\\
  \While{Pressure corrector loop}{
    \While{Non-orthogonal corrector loop}{
      Evaluate the pressure-corrector $p^{*}$:
      \begin{align*}
        \nabla\cdot(A[\rho^{*}]^{-1}\nabla p^{*}) = \nabla\cdot (A[\rho^{*}]^{-1}H[\mathbf{u}^{n, i}, \nu_t^{n, i}])
      \end{align*}\\
    }
  Correct density and velocity
  \begin{align*}
    \mathbf{u^{*}} \leftarrow  \mathbf{u}^{*}-A^{-1}\nabla p^{*}
  \end{align*}\\
  }
  }
\end{algorithm}
\end{minipage}
\hspace{1cm}
\begin{minipage}{0.47\linewidth}
  \begin{algorithm}[H]
    \caption{NM-LSPG-CNS $(i,n)$-th step}\label{alg:nmsonic}
    Start with an initial velocity field $\mathbf{u}^{n}$, density field $\rho^{n}$, internal energy field $e^{n}$, and pressure field $p^{n}$ at the $n$-th time step.\\
    Density residual evaluation:
    $$r_{\rho}=\rho^*-\rho^n + \Delta t \nabla\cdot (\rho^*\mathbf{u}^n)$$
    Momentum residual evaluation:
    $$ r_{\mathbf{u}}=A[\rho^{*}]\mathbf{u}^{n, i} - H[\mathbf{u}^{n, i}, \nu_t^{n, i}]+\nabla p^{n, i}.$$\\
    Energy residual evaluation:
    \begin{align*}
      r_{e}=\partial_t (\rho^* (e^{*} + K^{*})) &+ \nabla\cdot(\rho^*\mathbf{u}(e^{*}+K^{*})+\mathbf{u}^*p^*)
    \end{align*}\\
    Pressure-Poisson residual evaluation:
    $$ r_{p}=\nabla\cdot(A[\rho^{*}]^{-1}\nabla p^{n, i}) - \nabla\cdot (A[\rho^{*}]^{-1}H[\mathbf{u}^{n, i}, \nu_t^{n, i}]).$$\\
    Normalization of the residuals:
    \begin{align*}
      r_{\rho}\leftarrow \frac{r_{\rho}}{\max_{i} \rho_i},\quad &r_{\mathbf{u}}\leftarrow \frac{r_{\mathbf{u}}}{\max_{i} \mathbf{u}_i},\quad r_{e}\leftarrow \frac{r_{e}}{\max_{i} e_i},\\
      &r_{p}\leftarrow \frac{r_{p}}{\max_{i} p_i}
    \end{align*}
    \vspace{1.8cm}
  \end{algorithm}
  \end{minipage}
\DecMargin{1em}
\end{figure}
\vspace{0.5cm}

\begin{figure}[htpb!]
  \centering
  \includegraphics[width=0.24\textwidth, trim={0 0 600 0}, clip]{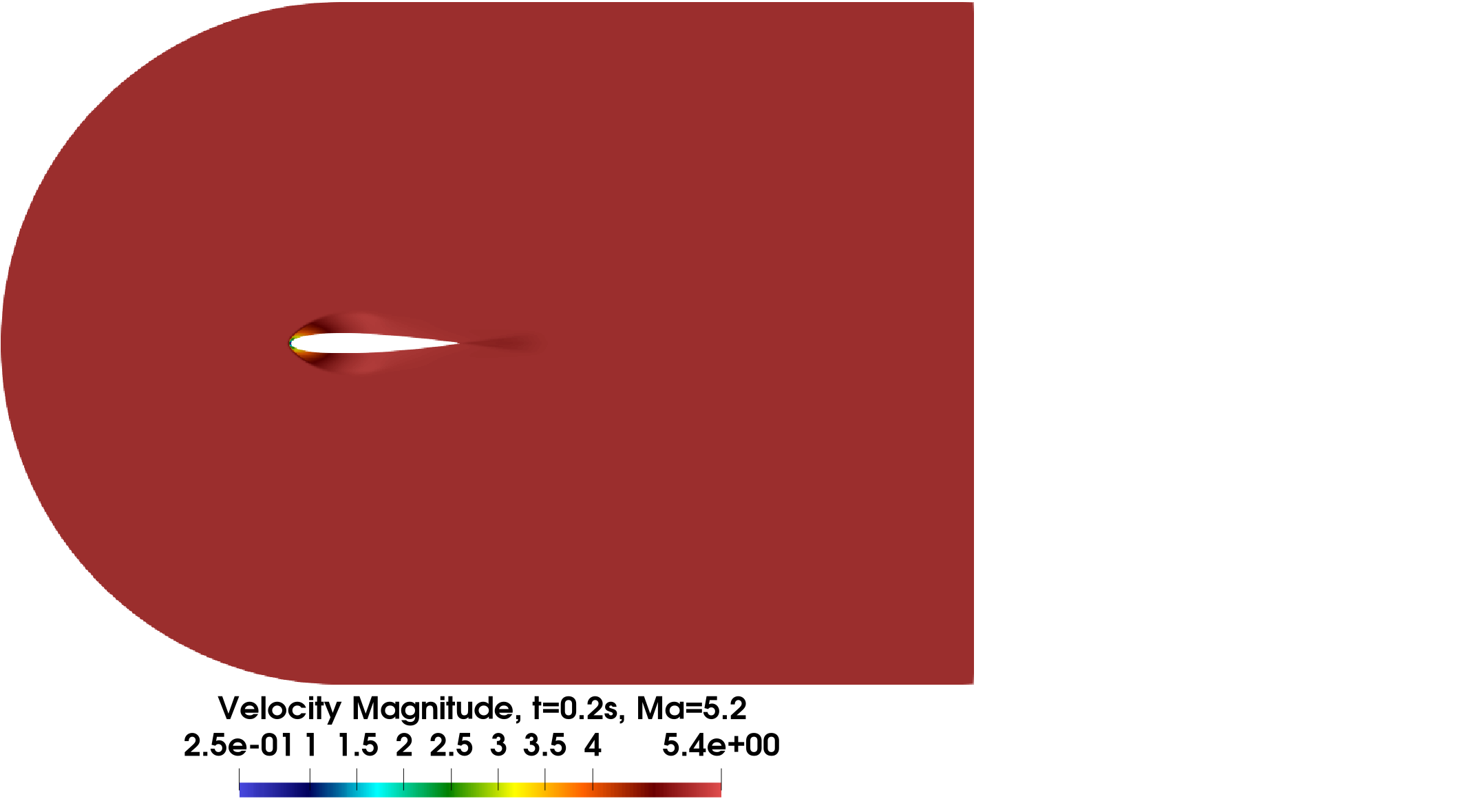}
  \includegraphics[width=0.24\textwidth, trim={0 0 600 0}, clip]{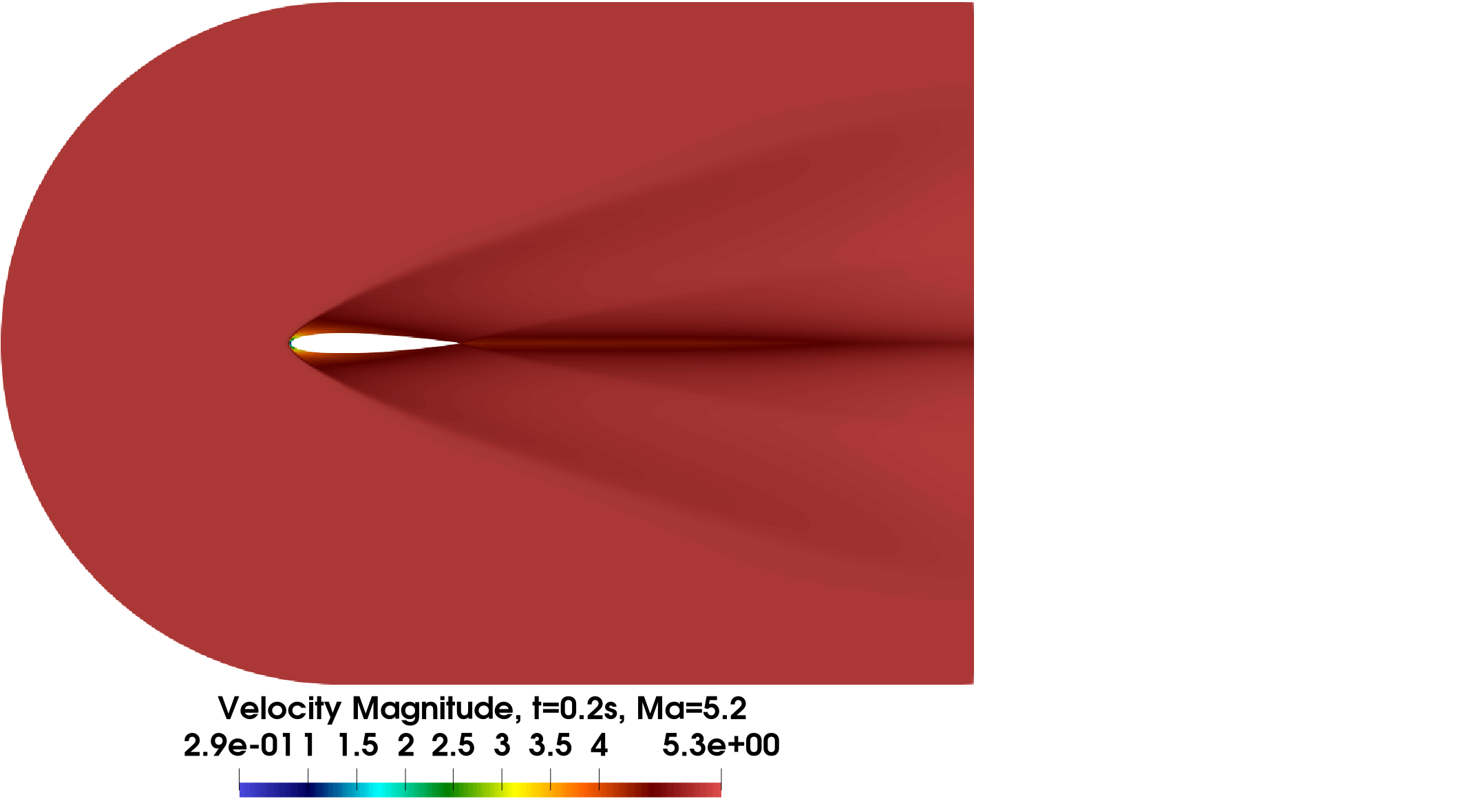}
  \includegraphics[width=0.24\textwidth, trim={0 0 600 0}, clip]{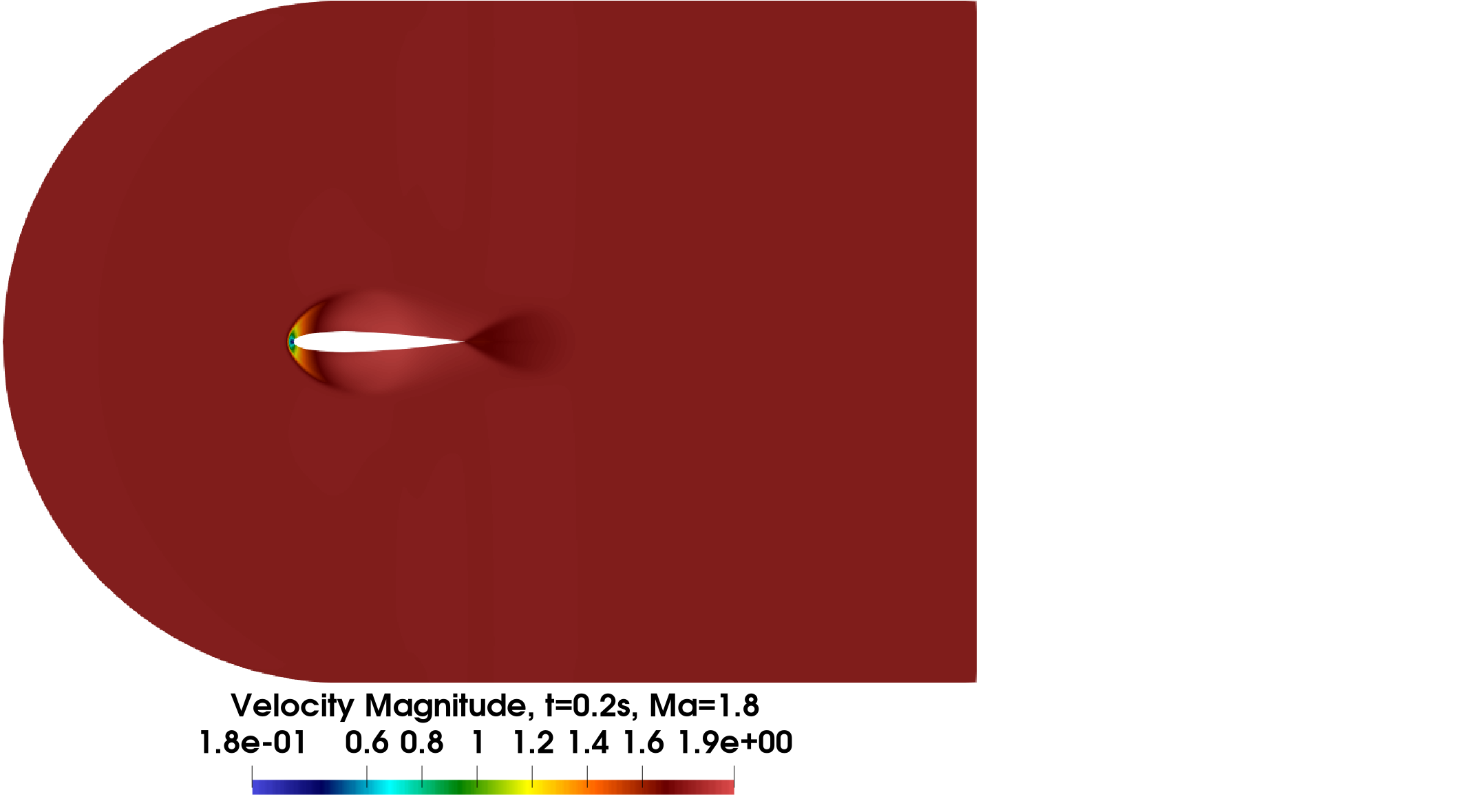}
  \includegraphics[width=0.24\textwidth, trim={0 0 600 0}, clip]{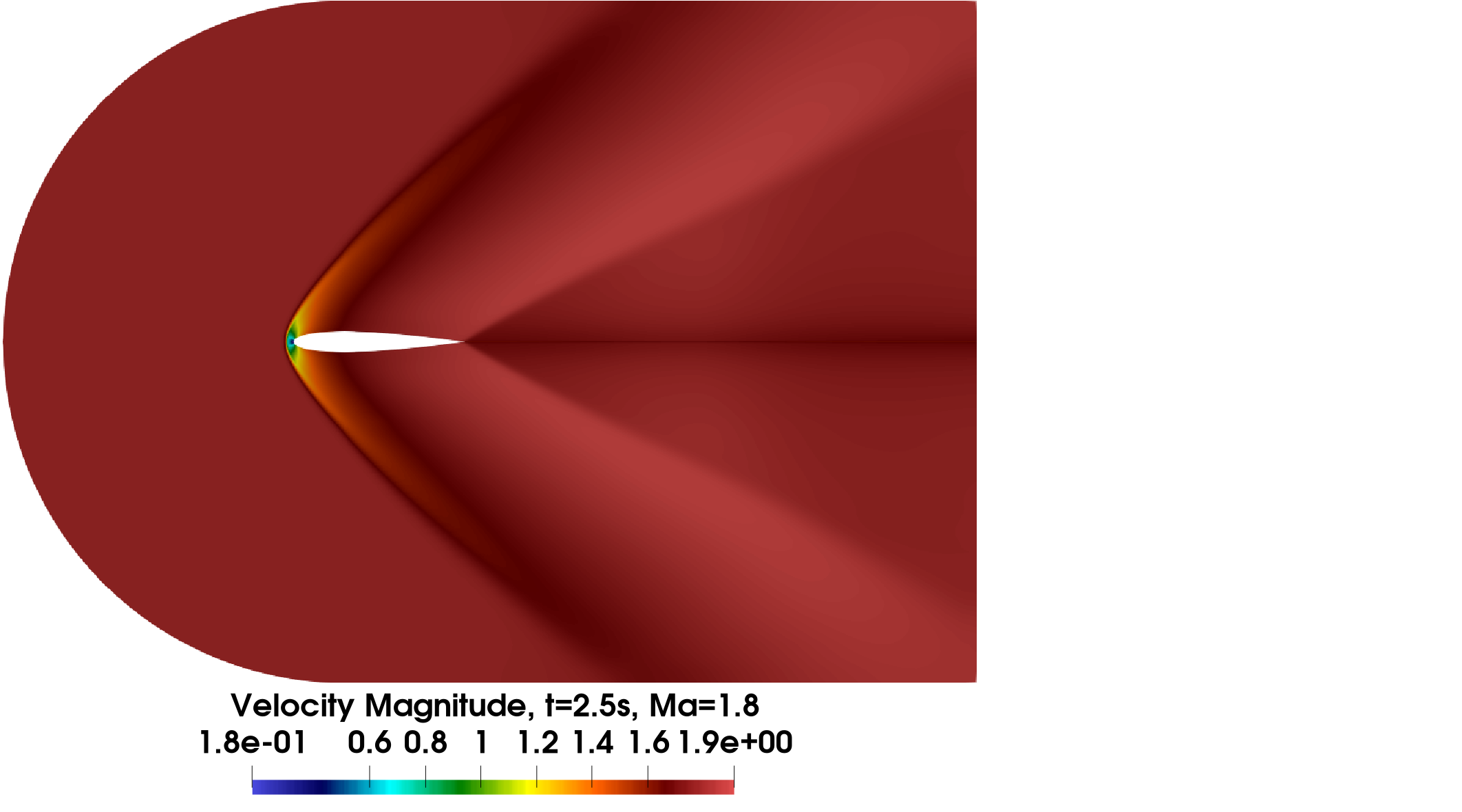}
  \caption{\textbf{CNS-1 and CNS-2}. \textbf{From left to right}: test full-order snapshots representing the velocity field magnitude corresponding to the parametric values $(t=0.2s,\ \text{Ma}=5.2)$, $(t=2.5s,\ \text{Ma}=5.2)$, $(t=0.2s,\ \text{Ma}=1.8)$, and $(t=2.5s,\ \text{Ma}=1.8)$. The whole extension of the solution manifold includes the transient dynamics from the time instants $t=0s$ to $t=2.5s$ opportunely rescaled through equation~\eqref{eq:timesAirfoil} and the different test Mach number values $\text{Ma}\in[1.8, 5.2]$. The presence of moving discontinuities at different Mach angles makes this test case difficult to reduced with classical linear projection based ROMs.}
  \label{fig:solutionManifoldCNS}
\end{figure}

To get a grasp of the extension of the solution manifold for the test cases CNS-1 and CNS-2, we show $4$ snapshots corresponding to the time instants and Mach numbers $t=0.2s,\ \text{Ma}=5.2$, $t=2.5s,\ \text{Ma}=5.2$, $t=0.2s,\ \text{Ma}=1.8$, and $t=2.5s,\ \text{Ma}=1.8$ in Figure~\ref{fig:solutionManifoldCNS}.

The value of the $L^2$ relative error is low for the internal energy field because its absolute value is higher than the absolute errors as can be seen in Appendix~\ref{sec:snaps} for the test cases CNS1 and CNS2.

\subsubsection{Interpolation and quadrature based hyper-reduction on a coarse mesh (CNS-1)}
\label{subsec:coarseAirfoil}
As anticipated, in this subsection we will consider a coarse mesh of \textbf{4500} cells, reported in Figure~\ref{fig:domain_airfoil}, for a total of $d=27000$ dofs. The training interval is $\text{Ma}\in[2, 5]$ and $t\in\{0, \Delta t,\dots,2500 \cdot\Delta t\}=V_{\boldsymbol{\mu}}$ with $\text{Ma}_{\text{ref}}=2$ and $(T_{\mu})_{\text{ref}}=2$, and the time steps opportunely scaled with respect the current Mach number through equation~\eqref{eq:timesAirfoil}. We consider the following training and test Mach numbers 
\begin{subequations}
  \begin{align}
    \mathcal{P}_{\text{train}}&=\{5, 4.256, 3.709, 3.291, 2.960,
    2.693, 2.473, 2.290, 2.134, 2\}\times V_{\boldsymbol{\mu}},\quad\vert\mathcal{P}_{\text{train}}\vert=10\cdot 2500,\\
    \mathcal{P}_{\text{test}}&=\{\mathbf{5.2}       , 4.042, 3.314 , 2.816, 2.455,
    2.182, \mathbf{1.970}, \mathbf{1.8}       \}\times V_{\boldsymbol{\mu}},\quad\vert\mathcal{P}_{\text{test}}\vert=8\cdot 2500,\label{eq:testCNSCoarse}
  \end{align}
\end{subequations}
it can be observed that the first and last two test parameters are in the extrapolation regime as they don't belong to the interval $[2, 5]$. From now on, the test parameters will be numbered with the order in which they appear in equation~\eqref{eq:testCNSCoarse} and refer only to the Mach number: from test parameter $1$ with $\text{Ma}=5.2$ to test parameter $8$ with $\text{Ma}=1.8$, we will use this notation also in the following figures. A grasp of the solution manifold extension can be observed in Figure~\ref{fig:solutionManifoldCNS}. So, the training dataset $A_{\text{train}}\in\mathbb{R}^{d\times n_{\text{train}}}$ is represented by $n_{\text{train}}=\vert\mathcal{P}_{\text{train}}\vert\cdot 1250=12500$ training snapshots, as only one every two time instants is saved. The test dataset $A_{\text{test}}\in\mathbb{R}^{d\times n_{\text{test}}}$ is composed by $n_{\text{test}}=\vert\mathcal{P}_{\text{test}}\vert\cdot 500=5000$ snapshots as only one every four time instants is saved. The number of rSVD modes considered is $r_{\text{rSVD}}=150$ evaluated from the training dataset $A_{\text{train}}$. The CAE architecture employed is reported in Table~\ref{tab: cae 150}. 

\begin{figure}[ht!]
  \centering
  \includegraphics[width=0.8\textwidth, trim={0 0 0 0}, clip]{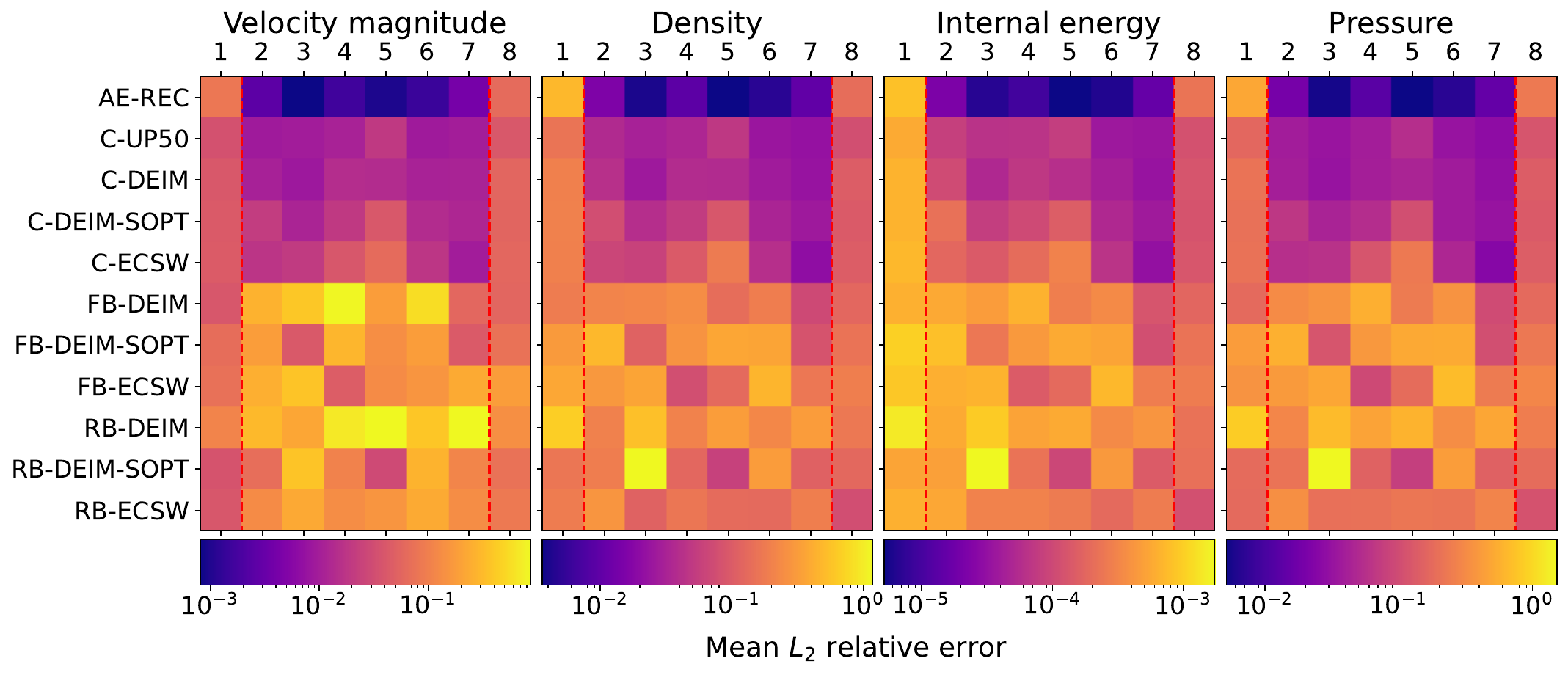}
  \caption{\textbf{CNS-1}. Comparison of hyper-reduction methods based on the mean $L^2$ relative error evaluated on each physical field of interest for the CNS test case. The acronyms correspond to: \textbf{AE-REC} is the autoencoder reconstruction error from~\eqref{eq:nonlinear_rec_err}, \textbf{C-UP50} is the reduced collocation method with the gradient based-sampling strategy from subsection~\ref{subsec:adaptive-hr} applied every $50$ time steps over $2000$ time istants; the other notations are introduced in subsections~\ref{subsubsec:deim}, \ref{subsubsec:sopt}, and \ref{subsubsec:ecsw}. For all the methods, it is employed a fixed number of $r_h=\mathbf{500}$ collocation nodes and $r_{\text{rSVD}}=\mathbf{150}$ rSVD modes used for both the definition of the nonlinear approximant map introduced in section~\ref{subsec: inductive biases} and the hyper-reduction basis. With the current monolithic hyper-reduction implementation, the best methods are the collocated ones, for further comments see section~\ref{sec:discussion} on this matter.}
  \label{fig:coarseHrConvergence}
\end{figure}

As a first study, we show in Figure~\ref{fig:coarseHrConvergence} the accuracy of the different hyper-reduction methods introduced in section~\ref{sec:hr}. For all the methods, it is employed a fixed number of $r_h=\mathbf{500}$ collocation nodes and $r_{\text{rSVD}}=\mathbf{150}$ rSVD modes used for both the definition of the nonlinear approximant map introduced in section~\ref{subsec: inductive biases} and the hyper-reduction basis. When residual basis \textbf{RB} are employed they are evaluated separately with respect to the ones used to define the nonlinear approximant $\phi:\mathbb{R}^r\rightarrow\mathbb{R}^d$. It can be seen that the reconstruction error~\eqref{eq:nonlinear_rec_err} of the autoencoder represents the baseline accuracy that we want to reach in terms of mean $L^2$ relative error. It is also clear that with our implementation of hyper-reduction the most accurate but also performing methods are the collocated ones. The lower accuracy of the other methodologies may be attributed to our choice of considering the physical fields of interest altogether in a monolithic fashion, both for the evaluation of the rSVD basis employed in the hyper-reduction offline stage and the computation of the normalization vectors from equation~\eqref{eq:normalization}. Moreover, collocation methods are more efficient as the collocated residuals do not need to be multiplied further by a pseudo-inverse or a vector of weights as in DEIM and ECSW methods. Further comments are added in section~\ref{sec:discussion}.

\begin{figure}[ht!]
  \centering
  \includegraphics[width=1\textwidth, trim={0 0 0 0}, clip]{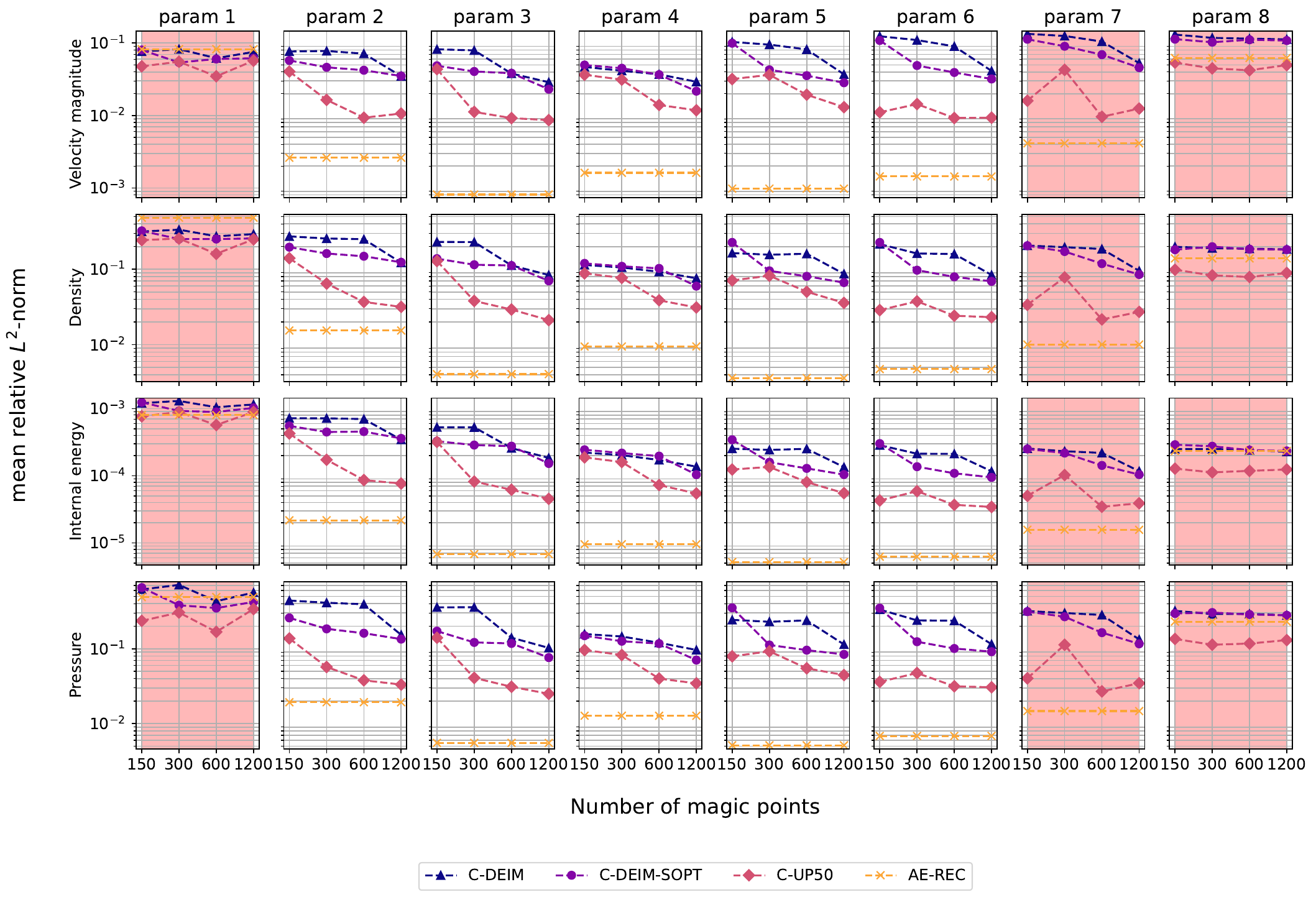}
  \caption{\textbf{CNS-1}. The decay of the mean $L^2$ relative error is assessed for the \textbf{C-DEIM}, \textbf{C-DEIM-SOPT}, and \textbf{C-UP50} reduced over-collocation methods with sampling strategies corresponding to a greedy one, the quasi-optimal (SOPT) one and the gradient-based adaptive one every $50$ time steps. In this case we employ the same reference time step as the full-order model $(\Delta t)_{\text{ref}}=0.001$. The baseline represented by the autoencoder mean $L^2$ reconstruction error \textbf{AE-REC} is reported. The extrapolation parameters have a shaded red background.}
  \label{fig:hrConvergenceCoarseCNS}
\end{figure}

Since we observed that collocated hyper-reduction reached a better accuracy, we show in the next study the decay of the mean $L^2$ relative error associated to the C-DEIM, C-DEIM-SOPT, and C-UP50 methods, varying the number of collocation nodes/magic points from $150$ to $1200$. The results are shown in Figure~\ref{fig:hrConvergenceCoarseCNS}. The most performing method is C-UP50, the gradient-based adaptive one, that pays the additional cost of a submesh update every $50$ time steps.

\begin{figure}[ht!]
  \centering
  \includegraphics[width=1\textwidth, trim={0 0 0 0}, clip]{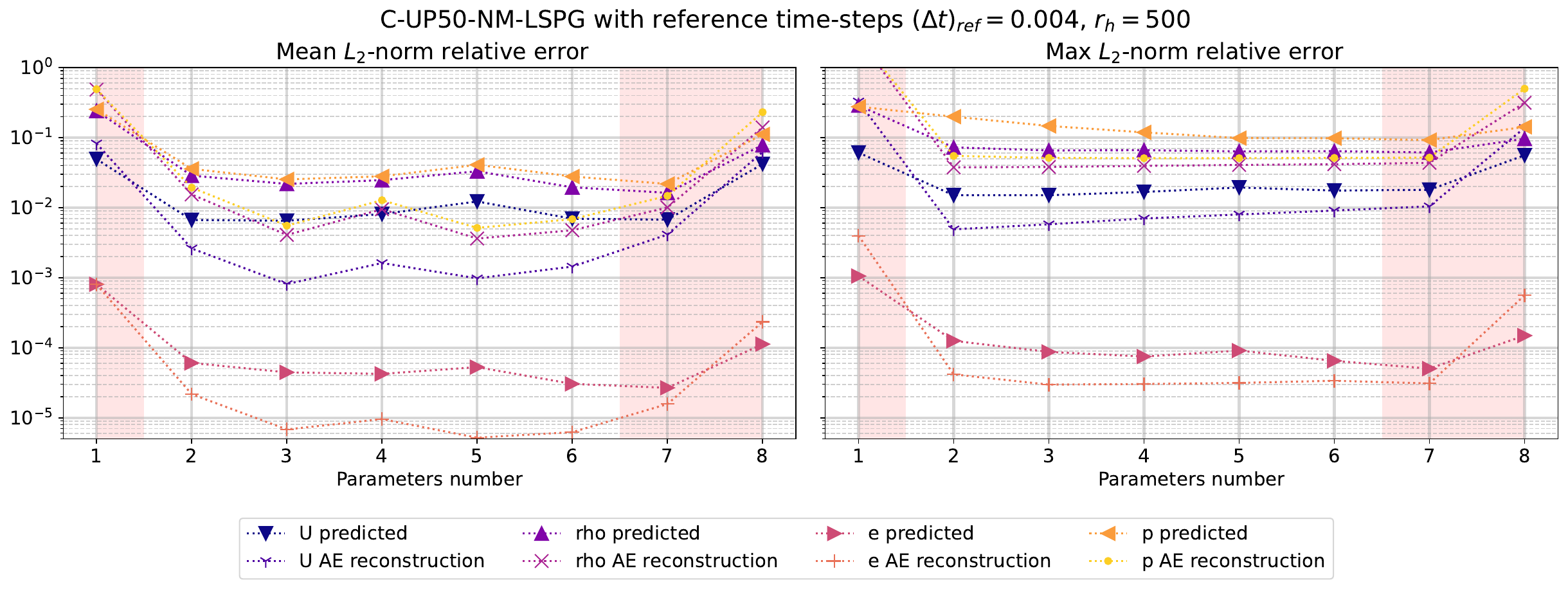}
  \caption{\textbf{CNS-1}. Accuracy of the nonlinear manifold least-squares Petrov-Galerkin (NM-LSPG) method with adaptive gradient-based reduced over-collocation (C-UP50) every $50$ time steps. The results are reported in terms of the mean $L^2$ relative error for the physical fields of interest with respect to the $8$ test paramters. The extrapolation regimes are delimited by a red shaded background. The number of magic points is $h_r=500$. It is important to notice that the reference time step is $(\Delta t)_{\text{ref}}=0.004$ and not $(\Delta t)_{\text{ref}}=0.001$ equal to the full-order model one. This permits the methodology to reach a speedup even for this small test cases with only \textbf{4500} cells.}
  \label{fig:bigStepsCNSCoarse}
\end{figure}

The advantage of having a continuous nonlinear approximant for the solution manifold, enables the possibility to choose a bigger reference time step in the online prediction stage. The results in terms of mean and max $L^2$ relative error are shown in Figure~\ref{fig:bigStepsCNSCoarse}. Thanks to this choice of reference time step $(\Delta t)_{\text{ref}}=0.004$ four times bigger than the full-order one $(\Delta t)_{\text{ref}}=0.001$ and twice as the sampling step used to select the training snapshots $A_{\text{train}}\in\mathbb{R}^{d\times n_{\text{train}}}$, a speedup can be achieved also for this small test case.

The average computational times of the NM-LSPG method with gradient-based adaptive hyper-reduction C-UP50 are shown in Table~\ref{tab: decoder costs coarseAirfoil}. The average is performed considering all $8$ test parameters for different Mach numbers. The average total time includes the cost of submesh updates introduced by the C-UP50 hyper-reduction. There is no evident speedup with respect to the full-order model. However, when a four time bigger reference time step $(\Delta t)_{\text{ref}}=0.004$ is imposed a speedup of almost $2$ is achieved also for this small test case.


\begin{table}[hbtp!]
  \centering
  \caption{\textbf{CNS-1}. Timings of the CNS on a coarse mesh, subsection~\ref{subsec:coarseAirfoil}. The timings of the reduced over-collocation method \textbf{C-UP50} relative to the choice of $r_h=\{150, 300, 600, 1200\}$ collocation nodes is reported along the full-order model (\textbf{FOM}) timings of the sonicFoam solver. The same reference time step as the \textbf{FOM} one is employed $(\Delta t)_{\text{ref}}=0.001$. Before the \textbf{FOM} results, the timings of \textbf{C-UP50} but with a larger reference time step of $(\Delta t)_{\text{ref}}=0.004$ are reported in bold for a number of magic points $r_h=500$. For the timings of the CNS on the finer mesh, see Table~\ref{tab: decoder costs finer CNS}.}
  \footnotesize
  \begin{tabular}{ c | c | c | c | c | c | c}
      \hline
      \hline
      collocation nodes ($r_h$)  & mean time-step  & mean update every 50 & average total time\\
      \hline
      \hline
      $r_h=150$, $(\Delta t)_{\text{ref}}=0.001$  & 11.937 [ms]  & 56.830 [ms]   & 29.842 [s]     \\
      \hline
      $r_h=300$, $(\Delta t)_{\text{ref}}=0.001$  & 25.294 [ms]  & 97.628 [ms]   & 63.235 [s]     \\
      \hline
      $r_h=600$, $(\Delta t)_{\text{ref}}=0.001$ & 36.894 [ms]  & 110.110 [ms]  & 92.234 [s]      \\
      \hline
      $r_h=1200$, $(\Delta t)_{\text{ref}}=0.001$ & 56.277 [ms]   & 119.707 [ms] & 140.691 [s]      \\
      \hline
      \hline 
      $r_h=\mathbf{500}$, $(\Delta t)_{\text{ref}}=0.004$ & \textbf{28.736} [ms]  & \textbf{82.337} [ms]   & \textbf{17.960} [s]     \\
      \hline
      \hline
      FOM, $(\Delta t)_{\text{ref}}=0.001$ & 13.440 [ms]  & -   & 33.614 [s]    \\
      \hline
      \hline
  \end{tabular}
  \label{tab: decoder costs coarseAirfoil}
\end{table}

\subsubsection{Collocated hyper-reduction on a finer mesh (CNS-2)}
\label{subsubsec:sudomains_finer_airfoil}

\begin{figure}[ht!]
  \centering
  \includegraphics[width=1\textwidth, trim={0 0 0 0}, clip]{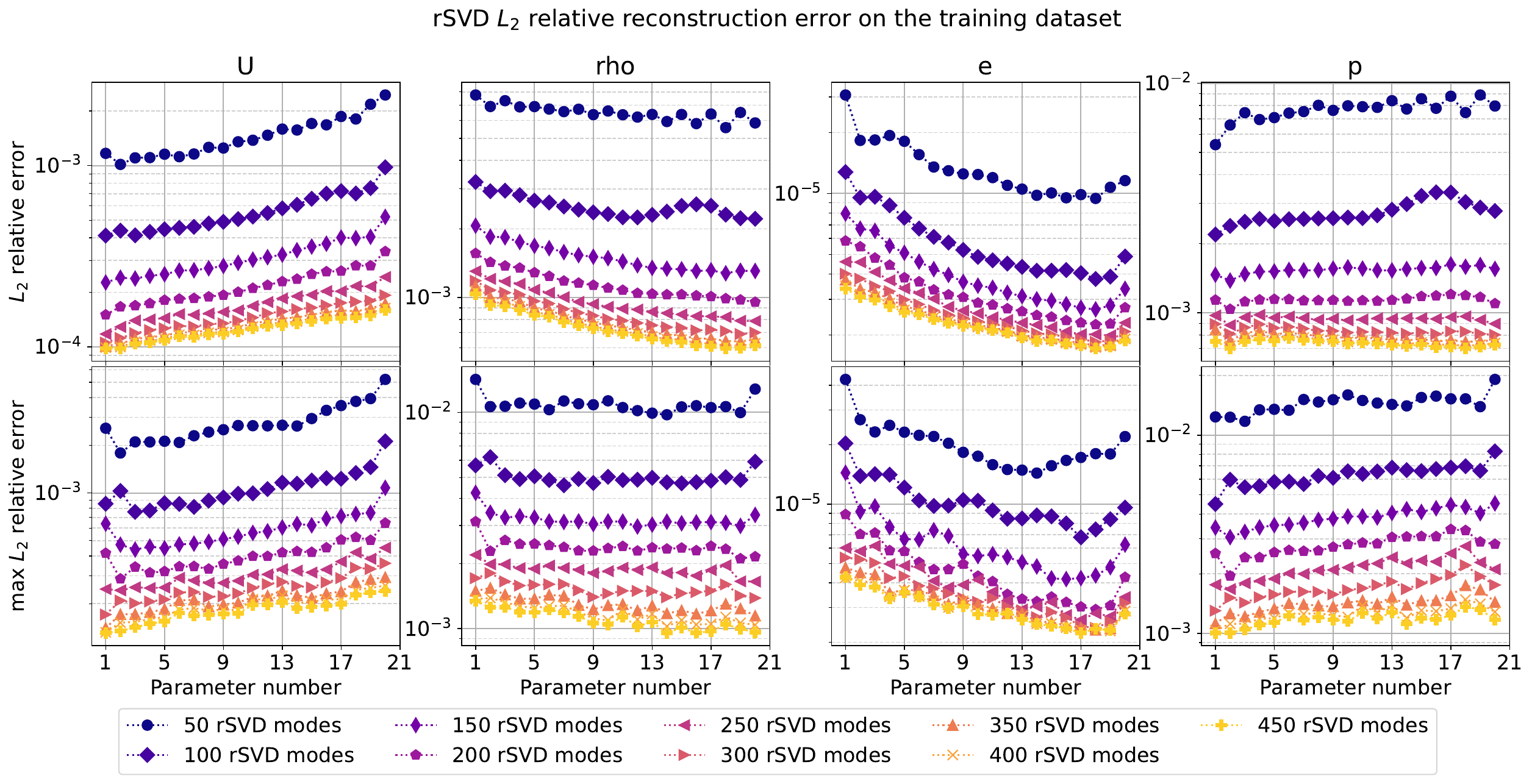}
  \caption{\textbf{CNS-2}. The decay of the mean and max $L^2$ reconstruction errors defined trough equations~\eqref{eq:rec_rsvd} for the $20$ train parameters introduced in equation~\eqref{eq:trainCNSFine}, is shown. The number of rSVD modes initially chosen on the whole parameter space for this test case is $r_{\text{rSVD}}=\mathbf{300}$. Even though the solution manifold includes moving discontinuities at different Mach angles as shown in Figure~\ref{fig:solutionManifoldCNS}, a moderately high number of rSVD modes can still be employed to obtain a good reconstruction error but at the same time precludes the implementation of efficient linear projection-based ROMs. A heuristic understanding of the moderately slow Kolmogorov n-width decay for this test case can be observed in the degradation of the reconstruction error on the test set~\ref{fig:decayCNS2test}.}
  \label{fig:decayCNS2train}
\end{figure}

\begin{figure}[ht!]
  \centering
  \includegraphics[width=1\textwidth, trim={0 0 0 0}, clip]{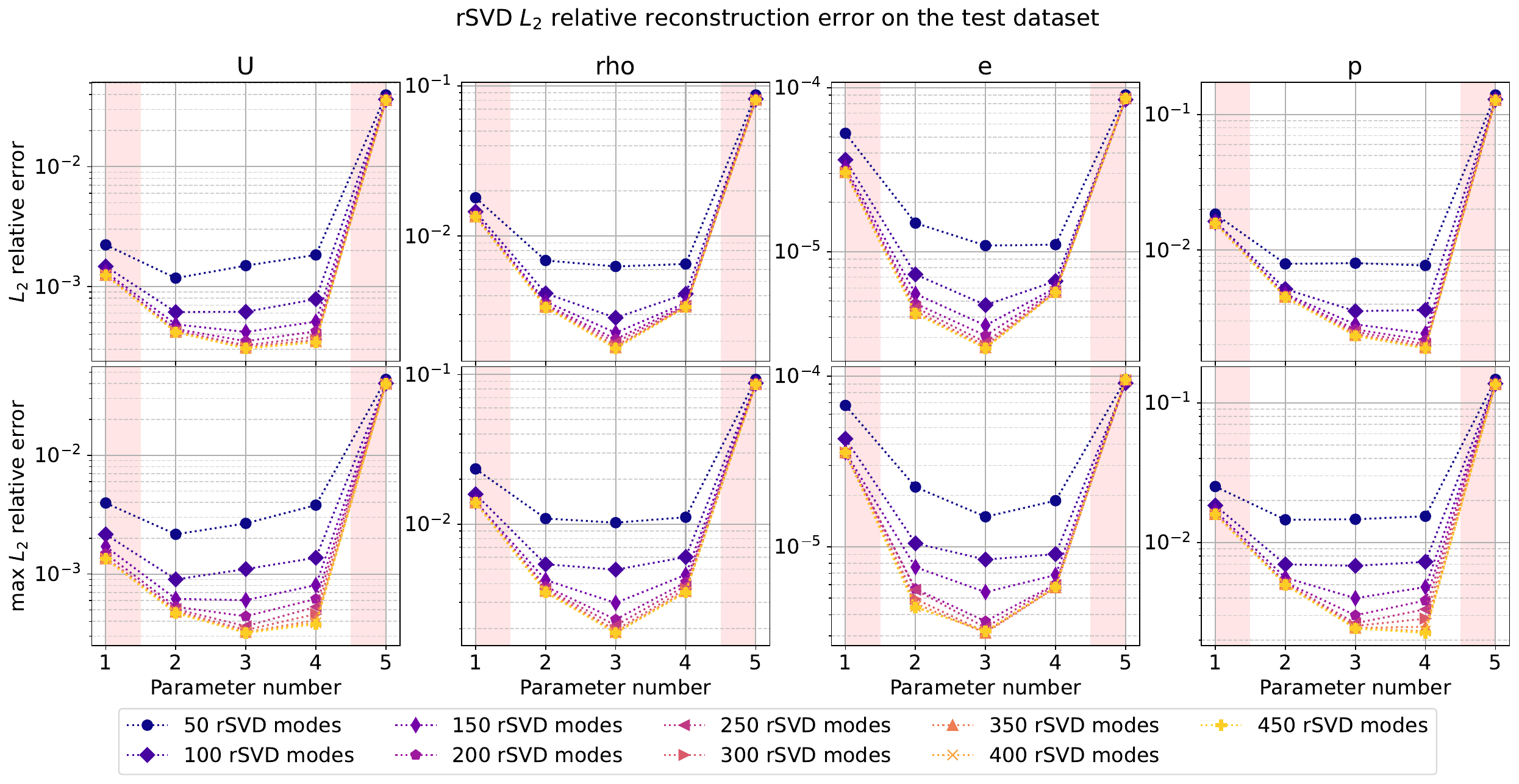}
  \caption{\textbf{CNS-2}. The decay of the mean and max $L^2$ reconstruction errors defined trough equations~\eqref{eq:rec_rsvd} for the $5$ test parameters introduced in equation~\eqref{eq:rec_rsvd_test}, is shown. The number of rSVD modes initially chosen on the whole parameter space for this test case is $r_{\text{rSVD}}=300$. The presence of moving discontinuities at different Mach angles as shown in Figure~\ref{fig:solutionManifoldCNS} causes an evident degradation of the reconstruction error with respect to the training reconstruction error in Figure~\ref{fig:decayCNS2train}. A shaded red background identifies the extrapolation regime.}
  \label{fig:decayCNS2test}
\end{figure}

To be sure to obtain a good approximation with rSVD modes, we increase the number of training snapshots for this CNS-2 test case on the finer mesh of \textbf{32160} cells shown in Figure~\ref{fig:domain_airfoil}. So we consider the following $20$ training and $5$ test parameters:
\begin{subequations}
  \begin{align}
    \mathcal{P}_{\text{train}}=\{&5.        , 4.617, 4.290, 4.007, 3.760,
    3.543, 3.350, 3.178, 3.024, 2.885,\nonumber\\
    &2.758, 2.644, 2.539, 2.442, 2.354,
    2.272, 2.196, 2.126, 2.061, 2.        \}\times V_{\boldsymbol{\mu}},\quad\vert\mathcal{P}_{\text{train}}\vert=20\cdot 2500,\label{eq:trainCNSFine}\\
    \mathcal{P}_{\text{test}}=\{&\textbf{5.2}       , 3.469, 2.622, 2.124, \textbf{1.8}       \}\times V_{\boldsymbol{\mu}},\quad\vert\mathcal{P}_{\text{test}}\vert=5\cdot 2500,\label{eq:testCNSFine}
  \end{align}
\end{subequations}
the first and last test parameters in bold correspond to the extrapolation regime. We remark that we did not optimize the number of training snapshots: possibly, a smaller number of them is needed to obtain an accurate regression of the solution manifold. As for the previous test case we will refer to the test Mach numbers from $\text{Ma}=5.2$ to $\text{Ma}=1.8$ with the numbers from $1$ to $5$ in this order. A grasp of the extension of the solution manifold we want to approximate is shown in Figure~\ref{fig:solutionManifoldCNS}. The training dataset $A_{\text{train}}\in\mathbb{R}^{d\times n_{\text{train}}}$ is represented by $n_{\text{train}}=\vert\mathcal{P}_{\text{train}}\vert\cdot 625=20\cdot 625 = 12500$ training snapshots, since only one every four time instants is saved. The test dataset $A_{\text{test}}\in\mathbb{R}^{d\times n_{\text{test}}}$ is composed by $n_{\text{test}}=\vert\mathcal{P}_{\text{test}}\vert\cdot 625=5\cdot 625=3125$ snapshots as only one every four time instants is saved. At first, the number of rSVD modes considered on the whole parameter space is $r_{\text{rSVD}}=300$ evaluated from the training dataset $A_{\text{train}}$. Later, we will extract rSVD on two separate time intervals to study the performance of local nonlinear solution manifolds. The decay of the mean and max $L^2$ reconstruction errors is shown for the fields of interest in Figures~\ref{fig:decayCNS2train} and ~\ref{fig:decayCNS2test} for the training and test datasets. The rSVD basis evaluated from the training dataset as explained in section~\ref{subsec:rSVD} are the same used to evaluate the test reconstruction error in Figure~\ref{fig:decayCNS2test}. The presence of moving discontinuities coming from the transient dynamics and the different Mach angles causes an evident degradation of the test reconstruction error with respect to the training reconstruction error. 

For this test case, a moderately high number of rSVD modes equal to $r_{\text{rSVD}}=300$ can be employed to approximate with sufficient accuracy the solution manifold. However, it can be observed that more complex parameter dependencies can exacerbate this behavior and make the approximation of the solution manifold with a linear rSVD basis unfeasible. In those cases, fully nonlinear NN architectures directly supported on the dofs of the physical fields of interest can be employed and hyper-reduced with the variant of the nonlinear manifold LSPG method introduced in~\cite{romor2023non}.

We want to study the implementation of local nonlinear manifolds approximants and how to efficiently change from one subdomain to the other in the online stage. We consider only two parametric subdomains determined by the splitting of the reference time interval $[0s, 2.5s]$ into to two subintervals $[0s, 1.2s]\supset V_{\boldsymbol{\mu}}^1$ and $[0.6s, 2.5s]\supset V_{\boldsymbol{\mu}}^2$, with $\vert V_{\boldsymbol{\mu}}^1\vert = 1.2/(\Delta t)_{\text{ref}} = 300$ and $\vert V_{\boldsymbol{\mu}}^2\vert = (2.5-0.6)/(\Delta t)_{\text{ref}} = 625-150$ since one every $4$ time instants in $[0, 2.5]$ is saved with a reference time step of $(\Delta t)_{\text{ref}}=0.001$. Notice that they overlap in order to achieve a good accuracy when the change of basis is performed. The change of basis between the domains is performed at the reference time instant $t=\mathbf{0.8s}$. More sophisticated techniques can be implemented~\cite{zimmermann2021manifold}. The notation we employ to distinguish between the subdomains is introduced in section~\ref{sec:loc}. So the parameter spaces we consider are:
\begin{subequations}
  \begin{align}
    \mathcal{P}^{1}_{n_{\text{train}, 1}}=\{&5.        , 4.617, 4.290, 4.007, 3.760,
    3.543, 3.350, 3.178, 3.024, 2.885,\nonumber\\
    &2.758, 2.644, 2.539, 2.442, 2.354,
    2.272, 2.196, 2.126, 2.061, 2.        \}\times V_{\boldsymbol{\mu}}^1,\quad\vert\mathcal{P}^1_{n_{\text{train}, 1}}\vert=20\cdot 300,\\
    \mathcal{P}^{2}_{n_{\text{train}, 2}}=\{&5.        , 4.617, 4.290, 4.007, 3.760,
    3.543, 3.350, 3.178, 3.024, 2.885,\nonumber\\
    &2.758, 2.644, 2.539, 2.442, 2.354,
    2.272, 2.196, 2.126, 2.061, 2.        \}\times V_{\boldsymbol{\mu}}^2,\quad\vert\mathcal{P}^2_{n_{\text{train}, 2}}\vert=20\cdot (625-150),
  \end{align}
\end{subequations}
and consequently the training snapshots matrices $A^i_{n_{\text{train}, i}}\in\mathbb{R}^{d\times n_{\text{train}, i}}$ $i=1, 2$ are assembled and the two rSVD basis $U^i\in\mathbb{R}^{d\times r_{\text{rSVD}}}$ $i=1, 2$ evaluated. Each one has $r_{\text{rSVD}}=300$ modes, but in principle a different number can be used. The number of training snapshots are $n_{\text{train}, 1}=20\cdot \vert V_{\boldsymbol{\mu}}^1\vert = 20\cdot 300$ and $n_{\text{train}, 2}=20\cdot \vert V_{\boldsymbol{\mu}}^2\vert=20\cdot (625-150)$. The same CAE architecture for the two nonlinear approximants $\phi^i:\mathbb{R}^{r}\rightarrow X_h\sim\mathbb{R}^d$ $i=1, 2$ is employed and reported in Table~\ref{tab: cae 300}. The latent dimension is $4$ in both cases. We remind that the change of basis matrix from equation~\eqref{eq:change of basis} is computed offline so the methdology maintains the independence with respect to the number of dofs.

\begin{figure}[ht!]
  \centering
  \includegraphics[width=1\textwidth, trim={0 0 0 0}, clip]{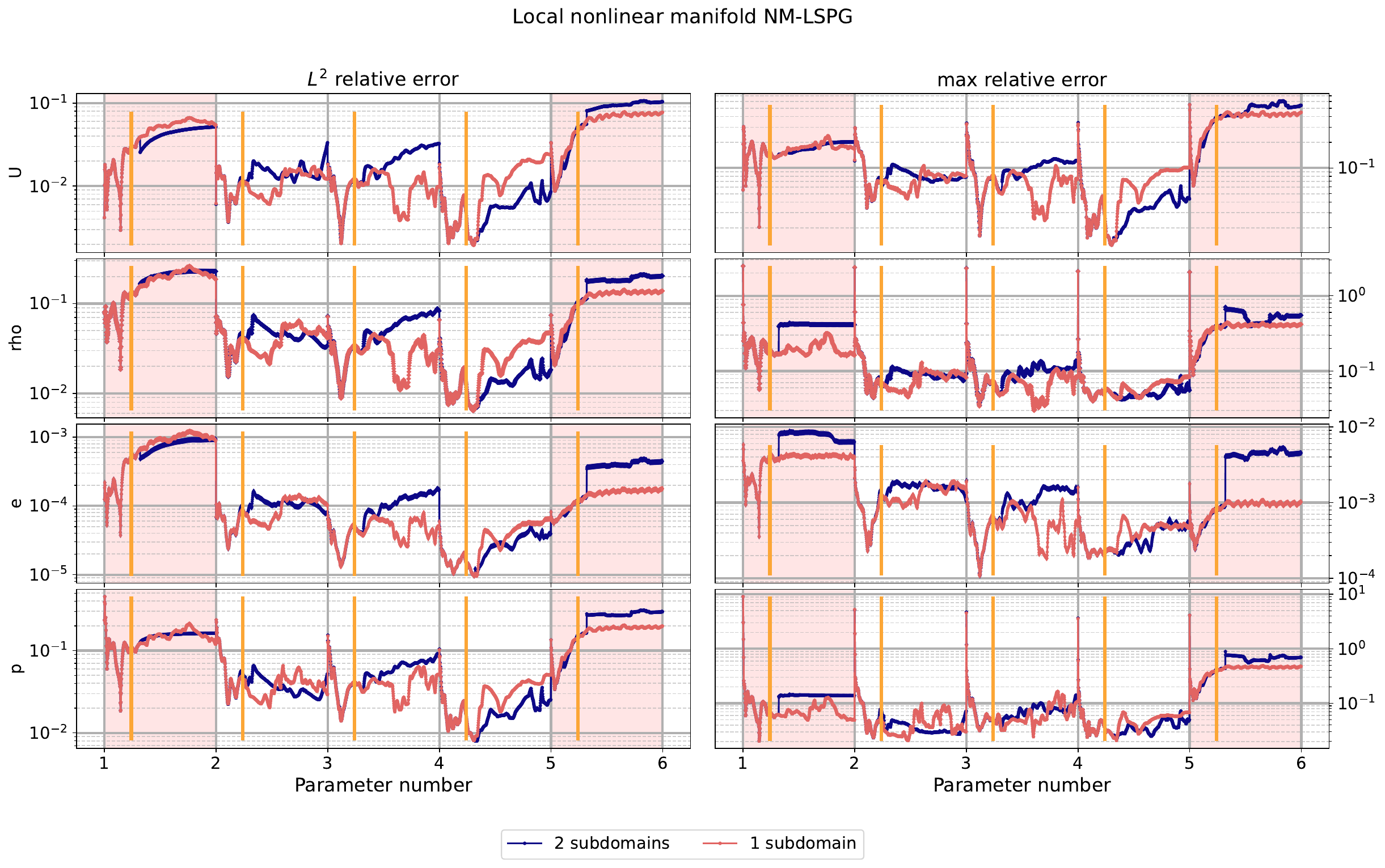}
  \caption{\textbf{CNS-2}. Results corresponding to the employment of two local nonlinear manifolds corresponding to the reference time intervals $[0s, 1.2s]=V_{\boldsymbol{\mu}}^1$ and $[0.6s, 2.5s]=V_{\boldsymbol{\mu}}^2$, as introduced in section~\ref{sec:loc}. The mean $L^2$ relative error and relative errors in max norm are shown for each of the \textbf{625} time instances associated with each one of the \textbf{5} test parameters corresponding to the test Mach numbers in equation~\eqref{eq:testCNSFine}. The reference time instant in which the change of basis is performed is $t=0.8s$ and it is highlighted by an orange vertical line for each of the $5$ test parameters. The extrapolation parameters $1$ and $5$ have a shaded red background.}
  \label{fig:local}
\end{figure}

The results for the choice of $r_h=800$ collocation nodes, $(\Delta t)_{\text{ref}}=0.001$ reference time step and $C-UP50$ hyper-reduction method, are reported in Figure~\ref{fig:local}. There, the whole trajectories corresponding to the orderly numbered $5$ test parameters are shown. In particular, the extrapolation regimes involving parameters $1$ and $5$ are shown, and the instant $t=0.8s$ in reference time scale, is highlighted for each test trajectory by an orange vertical line. It must be observed that not always the employment of local nonlinear submanifolds permits to reach smaller prediction errors: the only test parameter affected by an improvement is test parameter $4$, that nonetheless corresponds to Mach number $\text{Ma}=2.124$, and therefore it is more difficult to approximate since the Mach angle is wider. For this test case, we don't notice discontinuities from the passage of the solution from one local solution manifold to the other, as it can be seen also from the continuity of the errors in Figure~\ref{fig:local}.

\begin{figure}[htpb!]
  \centering
  \includegraphics[width=0.9\textwidth, trim={0 0 0 0}, clip]{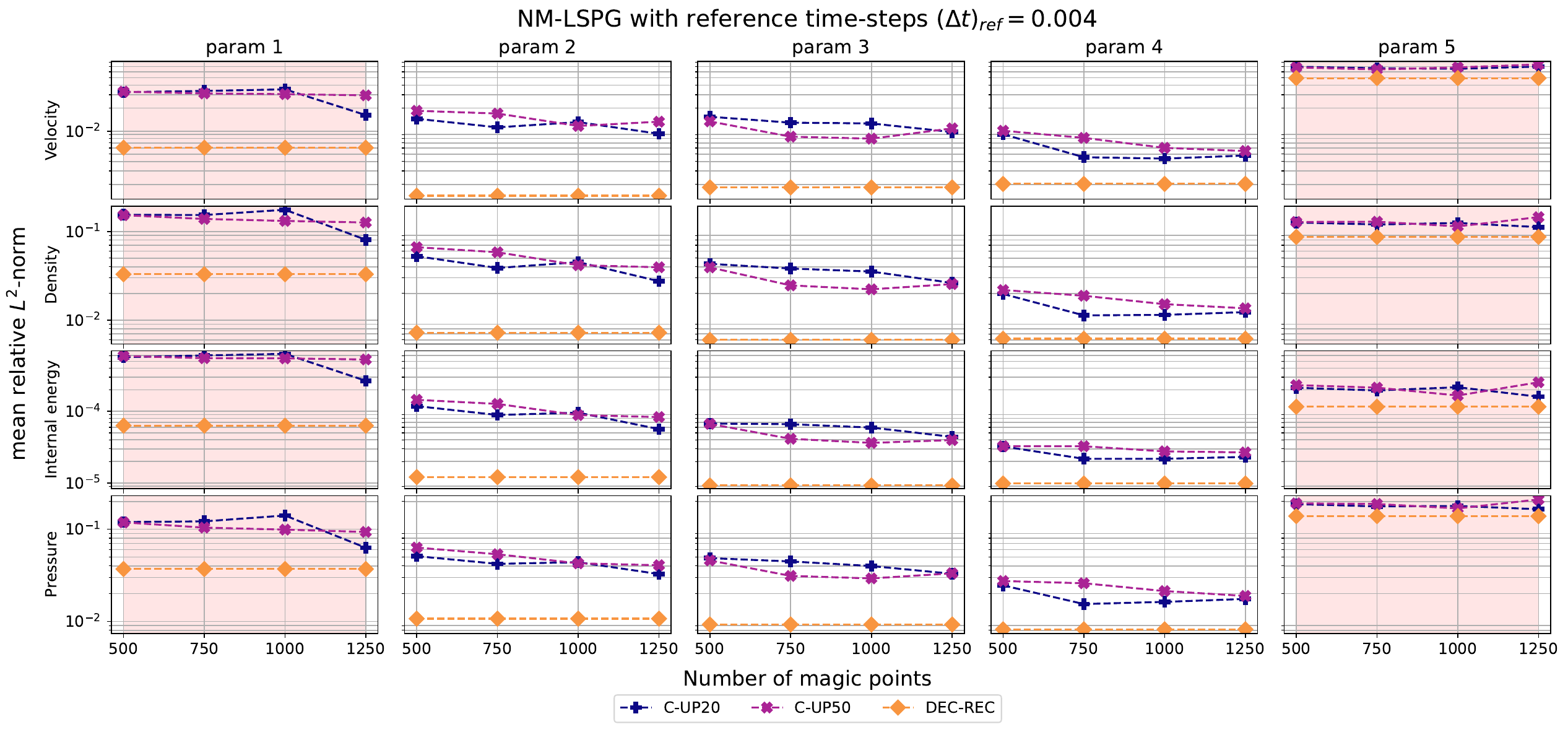}
  \caption{\textbf{CNS-2}. The decay of the mean $L^2$ relative error is assessed for the \textbf{C-UP20}, and \textbf{C-UP50} reduced over-collocation methods with sampling strategies corresponding to the gradient-based adaptive one every $50$ and $20$ time steps, respectively. In this case we employ a reference time step of four times $(\Delta t)_{\text{ref}}=0.004s$ the full-order model one $(\Delta t)_{\text{ref}}=0.001s$. We always consider two local nonlinear manifolds. The baseline represented by the autoencoder mean $L^2$ reconstruction error \textbf{AE-REC} is reported. The extrapolation parameters $1$ and $5$ have a shaded red background.}
  \label{fig:finerConvergence}
\end{figure}

Similarly to the previous test case CNS-1, the decay of the mean $L^2$ relative error is studied in Figure~\ref{fig:finerConvergence} for the C-UP20 and C-UP50 hyper-reduction methods with reference time step $(\Delta t)_{\text{ref}}=0.004$, that is four times the reference time step of the FOM $(\Delta t)_{\text{ref}}=0.001$. In this case always two local subdomains are considered. Associated with these studies, the computational costs for mean time instant evaluation and mean total time expended for the whole trajectory evaluation with reference time step $(\Delta t)_{\text{ref}}=0.004$, are shown in Table~\ref{tab: decoder costs finer CNS}. It can be observed that the computational costs increase from the same number of collocation nodes (MP) in Table~\ref{tab: decoder costs finer CNS} with respect to Table~\ref{tab: decoder costs coarseAirfoil}: this is due mostly to the number of rSVD basis changed from 150 to 300.

\begin{table}[htpb!]
  \centering
  \caption{\textbf{CNS-2}. Timings of the CNS on a finer mesh, subsection~\ref{subsubsec:sudomains_finer_airfoil}. The timings of the reduced over-collocation method \textbf{C-UP50} relative to the choice of $r_h=\{150, 300, 600, 1200\}$ collocation nodes is reported along the full-order model (\textbf{FOM}) timings of the sonicFoam solver. The reference time step of \textbf{C-UP50} is $(\Delta t)_{\text{ref}}=0.004$, four times bigger than the \textbf{FOM}. The different reference time steps affect the average total time and not the mean time steps computational costs. For the timings of the CNS on the coarse mesh, see Table~\ref{tab: decoder costs coarseAirfoil}. The change of basis between the local nonlinear manifolds is irrelevant but nonetheless included in the average total time.}
  \footnotesize
  \begin{tabular}{ l | c | c | c | c | c | c}
      \hline
      \hline
        & mean time-step  & mean update every 20 & average total time\\
      \hline
      \hline
      150  & 38.615 [ms]  & 286.103 [ms]   & 34.215 [s]     \\
      \hline
      300  & 52.209 [ms]  & 325.225 [ms]   & 42.774 [s]     \\
      \hline
      600 & 62.888 [ms]  & 330.269 [ms]  & 49.642 [s]      \\
      \hline
      1200 & 74.759 [ms]   & 362.719 [ms] & 58.302 [s]      \\
      \hline
      \hline
      FOM & 110.334 [ms]  & -   & 275.944 [s]    \\
      \hline
      \hline
  \end{tabular}
  \hspace{-0.05cm}
  \begin{tabular}{ l | c | c | c | c | c | c}
    \hline
    \hline
     mean time-step  & mean update every 50 & average total time\\
    \hline
    \hline 
    35.297 [ms]  & 302.612 [ms]   & 25.890 [s]     \\
    \hline
    48.940[ms]  & 316.992 [ms]   & 34.717 [s]     \\
    \hline
    61.803 [ms]  & 358.956 [ms]  & 42.887 [s]      \\
    \hline
    75.128 [ms]   & 362.236 [ms] & 51.510 [s]      \\
    \hline
    \hline 
    110.334 [ms]  & -   & 275.944 [s]    \\
    \hline
    \hline
\end{tabular}
  \label{tab: decoder costs finer CNS}
\end{table}
\subsection{Incompressible turbulent flow around the Ahmed body}
\label{subsec:incompr}

\begin{figure}[!thpb]
  \centering
  \includegraphics[width=0.49\textwidth]{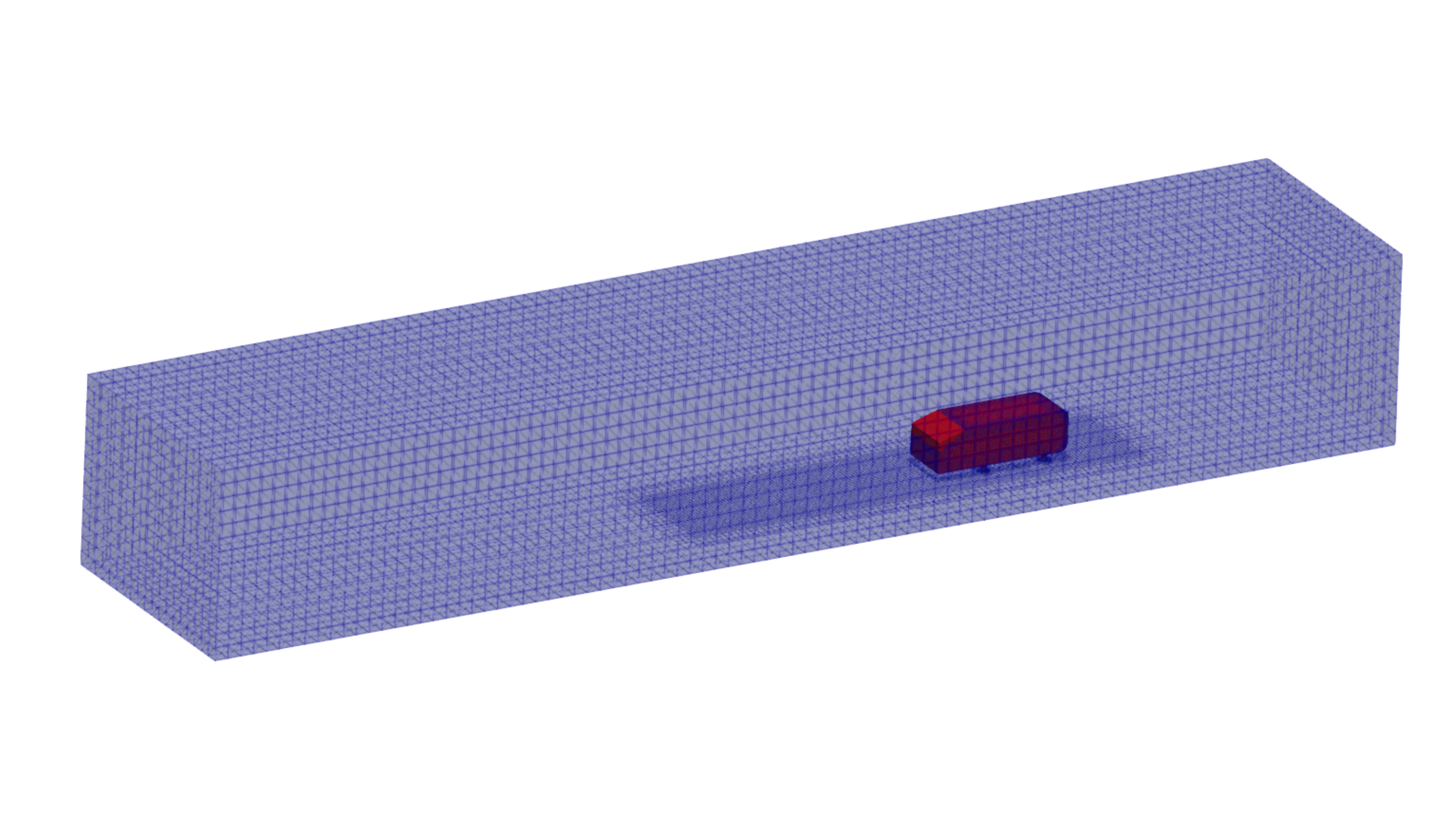}
  \includegraphics[width=0.49\textwidth]{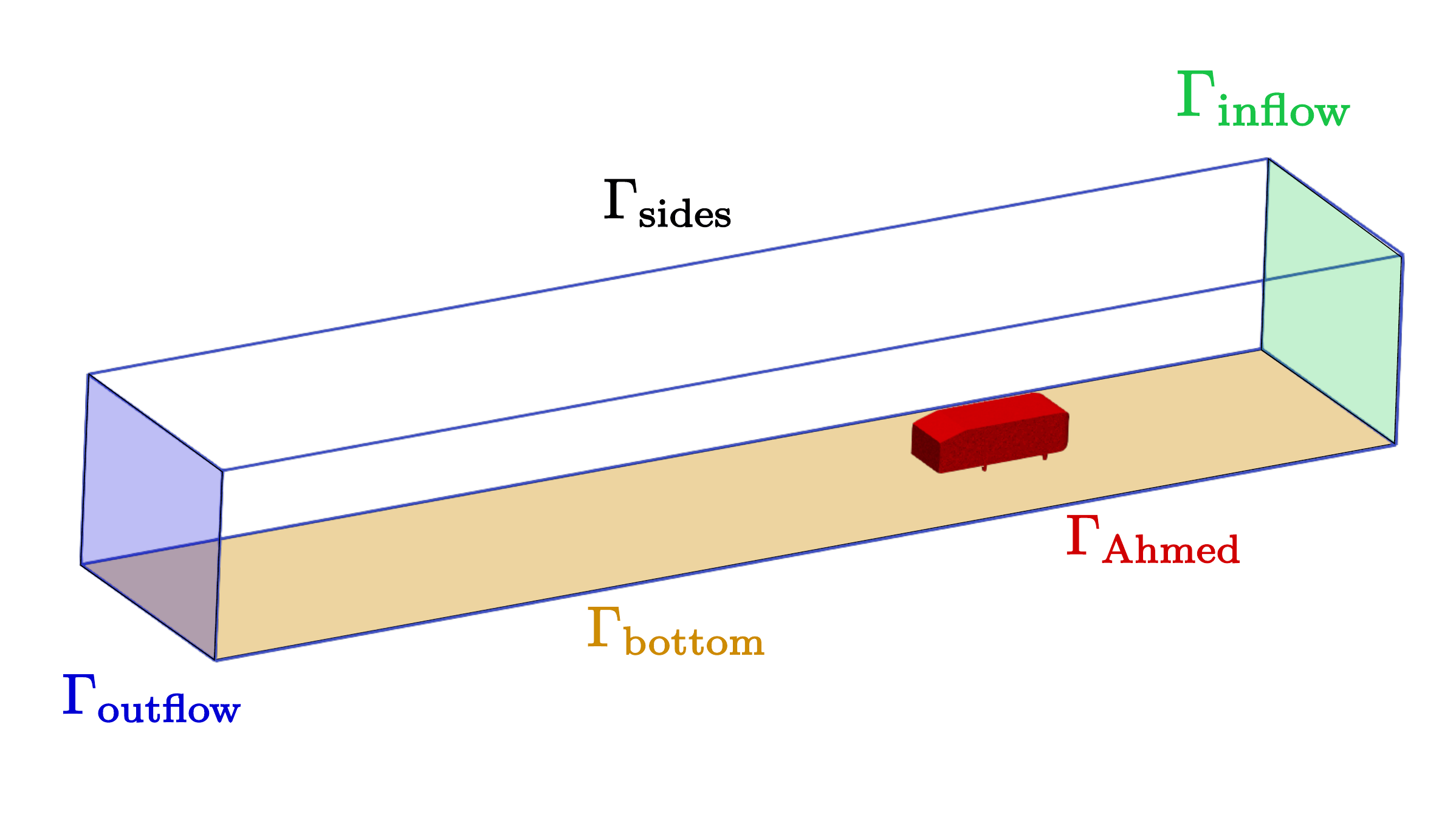}
  \caption{\textbf{INS.} \textbf{Left:} computational domain of the incompressible turbulent flow past the Ahmed body test case. The mesh has \textbf{198633} cells, for a total of $d=993165$ dofs for each geometrical deformation corresponding to a different choice of the slant angle $\theta$. \textbf{Right:} description of the boundaries of the computational domain: a fixed velocity $40\text{ms}^{-1}$ and pressure $p=0$ fields are imposed at the inflow boundary $\Gamma_{\text{inflow}}$ for each geometrical deformation.}
  \label{fig:domain_ahmed}
\end{figure}

The other benchmark (INS) we introduce to test our methodology involves the Reynolds-averaged Navier-Stokes equations (RANS) used to model an incompressible flow around the Ahmed body:
\begin{subequations}
  \begin{align}
    \partial_t \bar{\mathbf{u}} + \nabla\cdot (\bar{\mathbf{u}}\otimes\bar{\mathbf{u}})+\nabla \bar{p} - \nabla\cdot \left((\nu+\nu_t)\left(\nabla\bar{\mathbf{u}}+\nabla\bar{\mathbf{u}}^T\right)\right)= 0&,\quad &\text{(momentum conservation)}\label{eq:momINS}\\
    \nabla\cdot\bar{\mathbf{u}}=0&,\quad &\text{(mass conservation)}\label{eq:contINS}
  \end{align}
\end{subequations}
where $\bar{\mathbf{u}}$ and $\bar{p}$ are the time averaged velocity and the kinematic pressure fields, and $\nu_t$ is the kinematic Eddy turbulent viscosity. The turbulence is modelled with the Spalart-Allmaras one-equation model~\cite{spalart1992one, wilcox1998turbulence}. The Reynolds' number is in the order of $\text{Re}=2.8\cdot 10^6$. The parameters we consider are the slant angle $\theta$ of the Ahmed body and time $\mathcal{P}=[\theta_{\text{min}}, \theta_{\text{max}}]\times V_{\boldsymbol{\mu}}$. We will consider two parameter ranges for the slant angle: one associated to small geometrical deformations in subsection~\ref{subsubec:smallAhmed} (INS-1) with $\theta_{\text{min}}=15^\circ$ and $\theta_{\text{max}}=18.265^\circ$ and one associated to large geometrical deformations in subsection~\ref{subsubec:largeAhmed} (INS-2) with $\theta_{\text{min}}=15^\circ$ and $\theta_{\text{max}}=35^\circ$.The computational domain is shown in Figure~\ref{fig:domain_ahmed}. The geometry of the Ahmed body and the mesh are deformed with radial basis function interpolation as described in~\cite{zancanaro2021hybrid}.

The mesh, geometries, initial and boundary conditions are taken from the studies performed in~\cite{zancanaro2021hybrid} where a classical linear projection-based method is applied to reduce the SIMPLE~\cite{issa1991solution} numerical scheme with the employment of a neural network to approximate the Eddy viscosity. In that case, steady-state solutions are predicted while we focus on the transient dynamics since it is more challenging from the point of view of solution manifold approximation. In fact, we employ the PISO~\cite{patankar1983calculation} numerical scheme to achieve, from the initial conditions, convergence towards periodic cycles rather than stationary solutions.

The time interval is not dependent on the slant angle $\theta$, that is the final time is $T_{\boldsymbol{\mu}}= T = 0.1s$, the time step is $(\Delta t)_{\boldsymbol{\mu}} = \Delta t =0.0001s$ and the collection of time instants is $V_{\boldsymbol{\mu}} = V = \{0s, 0.0001s,\dots, 0.1s\}$. The initial and boundary conditions for the velocity $\bar{\mathbf{u}}$ and pressure $\bar{p}$ averaged fields are:
\begin{gather*}
  \begin{cases}
    \mathbf{u}(\x, t) =40 \text{ms}^{-1},\quad & (\x, t)\in\mathring{\Omega}_h\times \{t=0\}\\
    p(\x, t) = 0 ,\quad & (\x, t)\in\mathring{\Omega}_h\times \{t=0\}
  \end{cases},\qquad
  \begin{cases}
    \mathbf{u}(\x, t) =40 \text{ms}^{-1},\quad &\x\in\Gamma_{\text{inflow}}\\
    \mathbf{n}\cdot\nabla p(\x, t) = 0 ,\quad &\x\in\Gamma_{\text{inflow}}
  \end{cases},\\ 
  \begin{cases}
    \mathbf{n}\cdot\nabla\mathbf{u}(\x, t) = 0,\quad &\x\in\Gamma_{\text{outflow}}\\
    p(\x, t) = 0 ,\quad &\x\in\Gamma_{\text{outflow}}
  \end{cases},\qquad
  \begin{cases}
    \mathbf{u}(\x, t) = 0,\quad &\x\in\Gamma_{\text{Ahmed}}\cup\Gamma_{\text{bottom}}\\
    \mathbf{n}\cdot\nabla p(\x, t) = 0 ,\quad &\x\in\Gamma_{\text{Ahmed}}\cup\Gamma_{\text{bottom}}
  \end{cases}, \qquad
  \begin{cases}
    \mathbf{n}\cdot\nabla\mathbf{u}(\x, t) = 0,\quad &\x\in\Gamma_{\text{sides}}\\
    \mathbf{n}\cdot\nabla p(\x, t) = 0 ,\quad &\x\in\Gamma_{\text{sides}}
  \end{cases}
\end{gather*}
where the computational domain and its boundaries $\Gamma_{\text{outflow}}, \Gamma_{\text{inflow}}, \Gamma_{\text{Ahmed}}, \Gamma_{\text{bottom}}$ and the remaining faces $\Gamma_{\text{sides}}$ are shown in Figure~\ref{fig:domain_ahmed}. The computational mesh considered has a fixed number of cells for each geometrical deformation equal to $M=\mathbf{198633}$ for a total of $d=\mathbf{993165}$ dofs considering the average velocity, pressure and Eddy viscosity fields.

The \texttt{OpenFoam} solver we employ is the transient PISO numerical scheme. The predictor-corrector scheme is similar to the one described for the sonicFoam solver in subsection~\ref{subsec:compr} and shown in Algorithm~\ref{alg:piso} in a simplified version for a comparison with the nonlinear manifold least-squares Petrov-Galerkin (NM-INS) method in Algorithm~\ref{alg:nmpiso}. As before, we consider only the $n$-th time instant and possible $i$-th intermediate optimization steps for the NM-INS method. First, the velocity field is obtained with an implicit predictor step in line $2$, where the time discretization is hidden inside the diagonal $A$ and over-diagonal $H[\mathbf{u}^*, \nu^n_t]$ parts of the finite volume discretization of the Reynolds averaged momentum equation~\eqref{eq:momINS}. Afterwards, the velocity is corrected to satisfy the continuity equation~\eqref{eq:contINS} through the kinematic pressure, obtained with a pressure-Poisson equation in line $5$. Finally, the new Eddy viscosity is obtained. 

\IncMargin{1em}
\begin{figure}[!thpb]
  \caption{\textbf{INS}. Comparison between the $n$-th time instant iterations of the full-order model numerical scheme PISO (\textbf{Left}) and the nonlinear manifold least-squares Petrov Galerkin (NM-INS)method (\textbf{Right}).}
  \begin{minipage}{0.47\linewidth}
  \begin{algorithm}[H]
    \caption{PISO $n$-th iteration}\label{alg:piso}
    Start with an initial pressure field $p^{n}$, velocity field $\mathbf{u}^{n}$ and kinematic eddy viscosity $\nu_t^{n}$ at the $n$-th time step.\\
    Momentum predictor step:
    $$ A\mathbf{u}^{*} = H[\mathbf{u}^{*}, \nu_t^{n}]-\nabla p^{n}.$$\\
    \While{PISO pressure-corrector loop}{
      \While{Non-orthogonal corrector loop}{
        Evaluate the pressure-corrector term $p^{*}$:
        $$ \nabla\cdot(A^{-1}\nabla p^{*}) = \nabla\cdot (A^{-1}H[\mathbf{u}^{*}, \nu_t^{n}]).$$\\
      }
    Correct velocity
    $$\mathbf{u^{*}} \leftarrow  \mathbf{u}^{*}-A^{-1}\nabla p^{*}$$\\
    }
    Solve for the kinematic Eddy viscosity.\\
  \end{algorithm}
  \end{minipage}
  \hspace{1cm}
  \begin{minipage}{0.47\linewidth}
    \begin{algorithm}[H]
      \caption{NM-INS $(i,n)$-th residual}\label{alg:nmpiso}
      Start with an initial pressure field $p^{n, i}$, velocity field $\mathbf{u}^{n, i}$ and kinematic eddy viscosity $\nu_t^{n, i}$ at the $n$-th time step.\\
      Momentum residual evaluation:
      $$ r_{\mathbf{u}}=A\mathbf{u}^{n, i} - H[\mathbf{u}^{n, i}, \nu_t^{n, i}]+\nabla p^{n, i}.$$\\
      Pressure-Poisson residual evaluation:
      $$ r_{p}=\nabla\cdot(A^{-1}\nabla p^{n, i}) - \nabla\cdot (A^{-1}H[\mathbf{u}^{n, i}, \nu_t^{n, i}]).$$\\
      Normalization of the residuals:
      $$r_{\mathbf{u}}\leftarrow \frac{r_{\mathbf{u}}}{\max_{i} \mathbf{u}_i},\quad r_{p}\leftarrow \frac{r_{p}}{\max_{i} p_i}$$\\
      \vspace{1.5cm}
    \end{algorithm}
    \end{minipage}
  \DecMargin{1em}
\end{figure}

The NM-INS method does not implement a predictor-correct strategy instead: only residual evaluations need to be computed and PISO and non-orthogonal corrector loops are omitted, as we search for the converged solutions corresponding to the last corrector steps. We remark that we don't consider the residual of the Eddy viscosity equation of the Spalart-Allmaras one-equation model in the NM-INS algorithm, but only the velocity and pressure ones. This choice and the difficulty of linearly approximating the Eddy viscosity, compared to the velocity and pressure fields, could be the reasons behind the worse accuracy in the predictions of the Eddy viscosity. For example, this can be observed in Figure~\ref{fig:trajsmallAhmed}.

As for the CNS test case, we consider a five times bigger time step for the nonlinear manifold method: instead of $\Delta t = 0.0001s$ employed for the full-order solutions, we set $\Delta t =0.0005$ for the NM-INS method.

\subsubsection{Small geometrical deformations (INS-1)}
\label{subsubec:smallAhmed}

\begin{figure}[bpht!]
  \centering
  \includegraphics[width=0.49\textwidth]{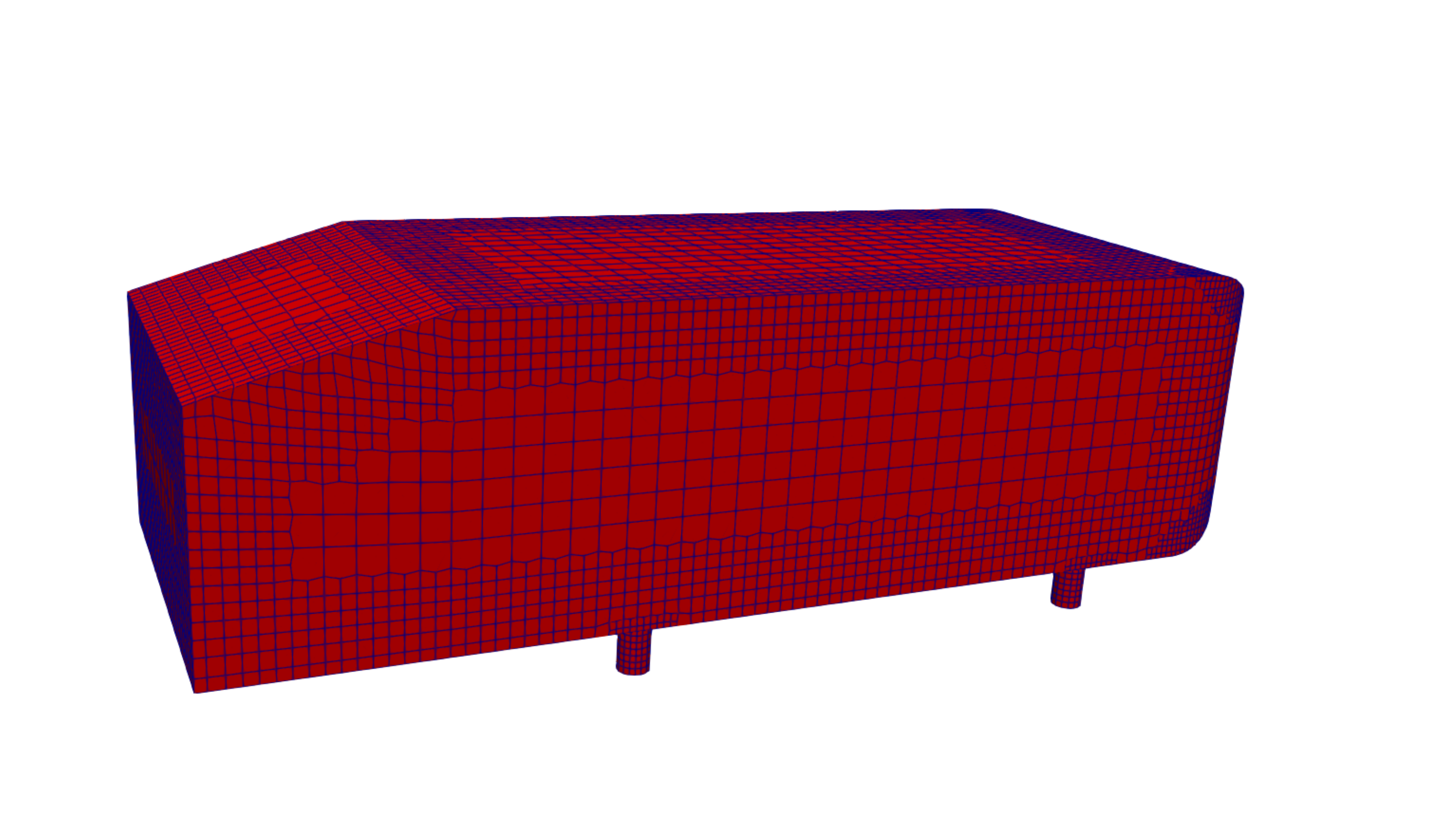}
  \includegraphics[width=0.49\textwidth]{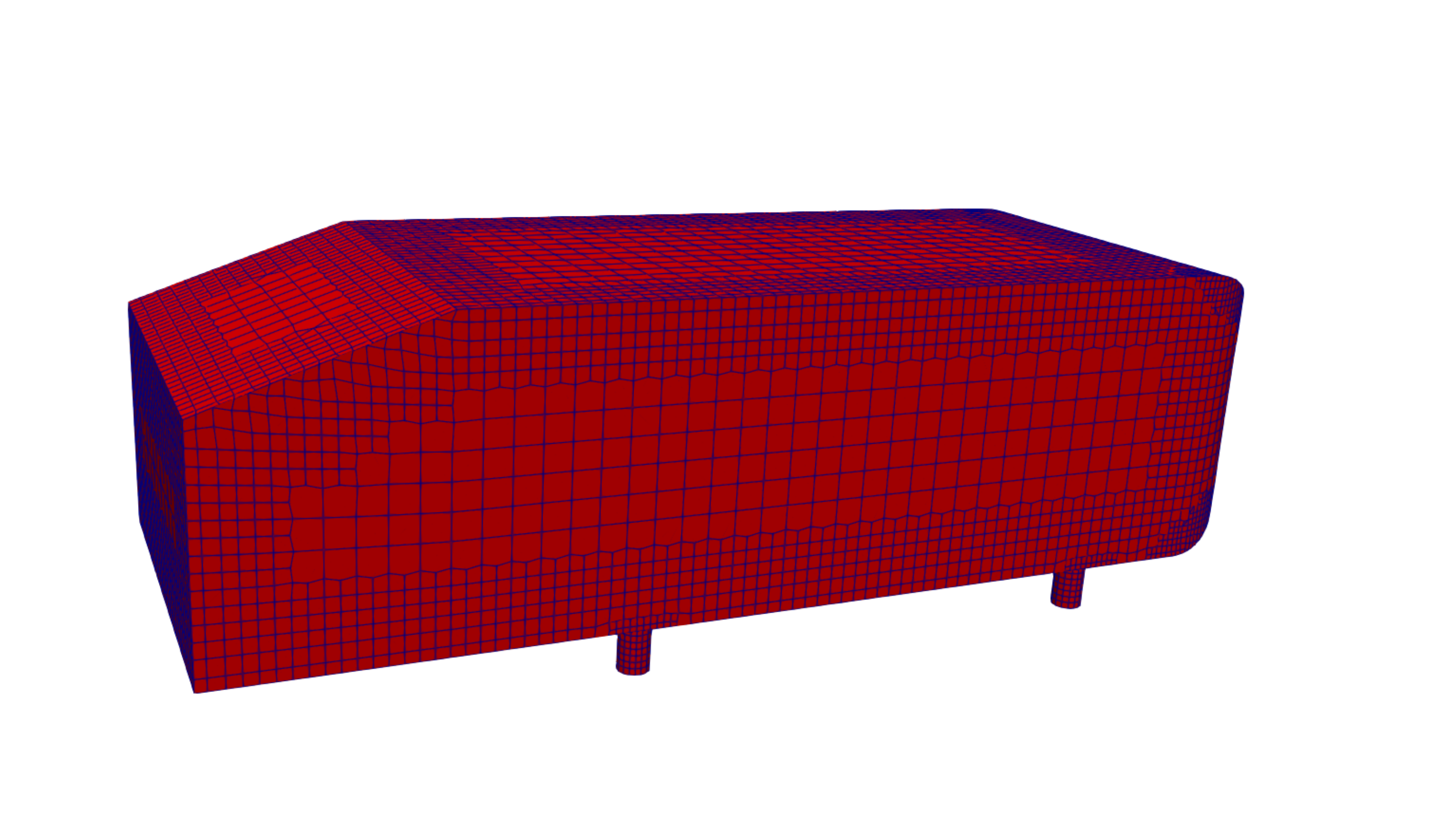}
  \caption{\textbf{INS-1}. \textbf{Left:} Ahmed body with slant angle $\theta=15.408^{\circ}$ from the training parameter set. \textbf{Right:} Ahmed body with slant angle $\theta=18.265^{\circ}$ from the training parameter set.}
  \label{fig:angles_ahmed}
\end{figure}

As can be understood from the small geometrical deformations in Figure~\ref{fig:angles_ahmed}, the difficulty resides in the approximation of the transient dynamics rather than in the influence of the geometrical parameter. We consider \textbf{5} training slant angles and \textbf{4} test slant angles inside the training range (no extrapolation):
\begin{subequations}
  \begin{align}
    \mathcal{P}_{\text{train}}=\{&15^{\circ}+ i\cdot \left((35^{\circ}-15^{\circ})/49\right)\vert i\in\{0,2,4,6,8\}\}\times V,\quad\vert\mathcal{P}_{\text{train}}\vert=5,\label{eq:trainINSFine}\\
    \mathcal{P}_{\text{test}}=\{&15^{\circ}+ i\cdot \left((35^{\circ}-15^{\circ})/49\right)\vert i\in\{1,3,5,7\}\}\times V,\quad\vert\mathcal{P}_{\text{test}}\vert=4,\label{eq:testINSFine}.
  \end{align}
\end{subequations}
In the order written in equations~\eqref{eq:trainINSFine} and~\eqref{eq:testINSFine}, we name the training and test parameters from $1$ to $5$ and from $1$ to $4$, respectively.

\begin{figure}[bpht!]
  \centering
  \includegraphics[width=1\textwidth, trim={0 0 0 0}, clip]{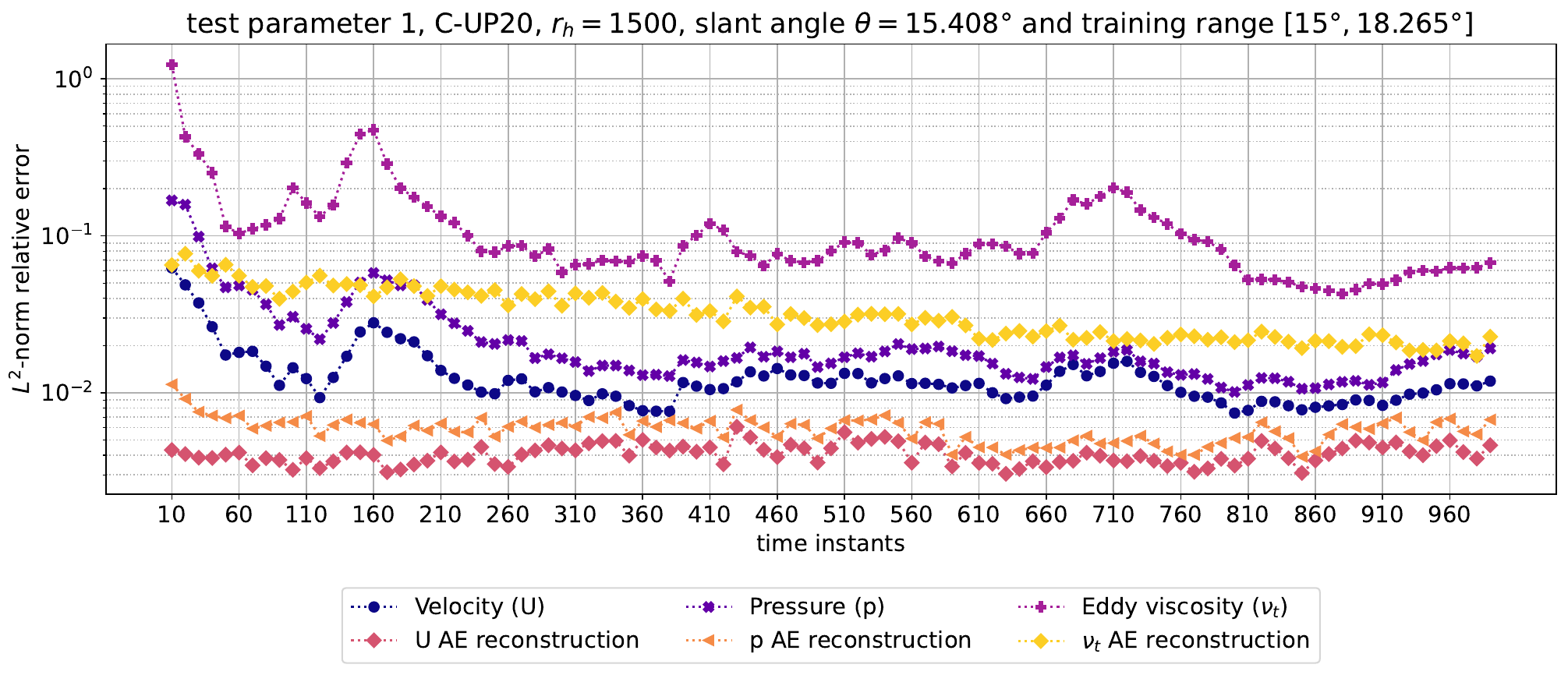}
  \caption{\textbf{INS-1}. Every $10$ time steps, the $L^2$ relative error is reported for the predicted physical fields of interest associated to test parameter \textbf{1}, corresponding to the slant angle $\theta=15.408^{\circ}$ and the \textbf{C-UP20} method with $\mathbf{r_h=1500}$ collocation nodes. The training parameter range is relatively small $[15^{\circ}, 18.265^{\circ}]$. The $L^2$ autoencoder relative reconstruction error represents an experimental lower bound to the prediction error. The accuracy is worse in the first time steps, a possible strategy to further reduce it is to employ local nonlinear manifold approximations as described in section~\ref{sec:loc}.}
  \label{fig:trajsmallAhmed}
\end{figure}

As anticipated we employ a five times bigger time step $\Delta t =0.0005\text{s}$ with respect to the full-order one $\Delta t=0.0001\text{s}$. We apply the C-UP-20 hyper-reduction with $r_h=1500$ collocation nodes. For the test parameter $1$, $\theta=15.4^{\circ}$, the mean $L_2$ relative errors corresponding to the $1000$ time instants from $0.0005\text{s}$ to $0.1\text{s}$ are shown in Figure~\ref{fig:trajsmallAhmed}. As can be seen, efficient and relatively accurate predictions can be obtained. The computational time spent is summarized in Table~\ref{tab: ahmed small timings}, reaching a speedup of around \textbf{26} with respect to the full-order model for this simple test case. Local nonlinear manifold approximations could be employed to achieve better predictions in the initial time instants, splitting the time interval.

\begin{figure}[bpht!]
  \centering
  \includegraphics[width=1\textwidth, trim={0 0 0 0}, clip]{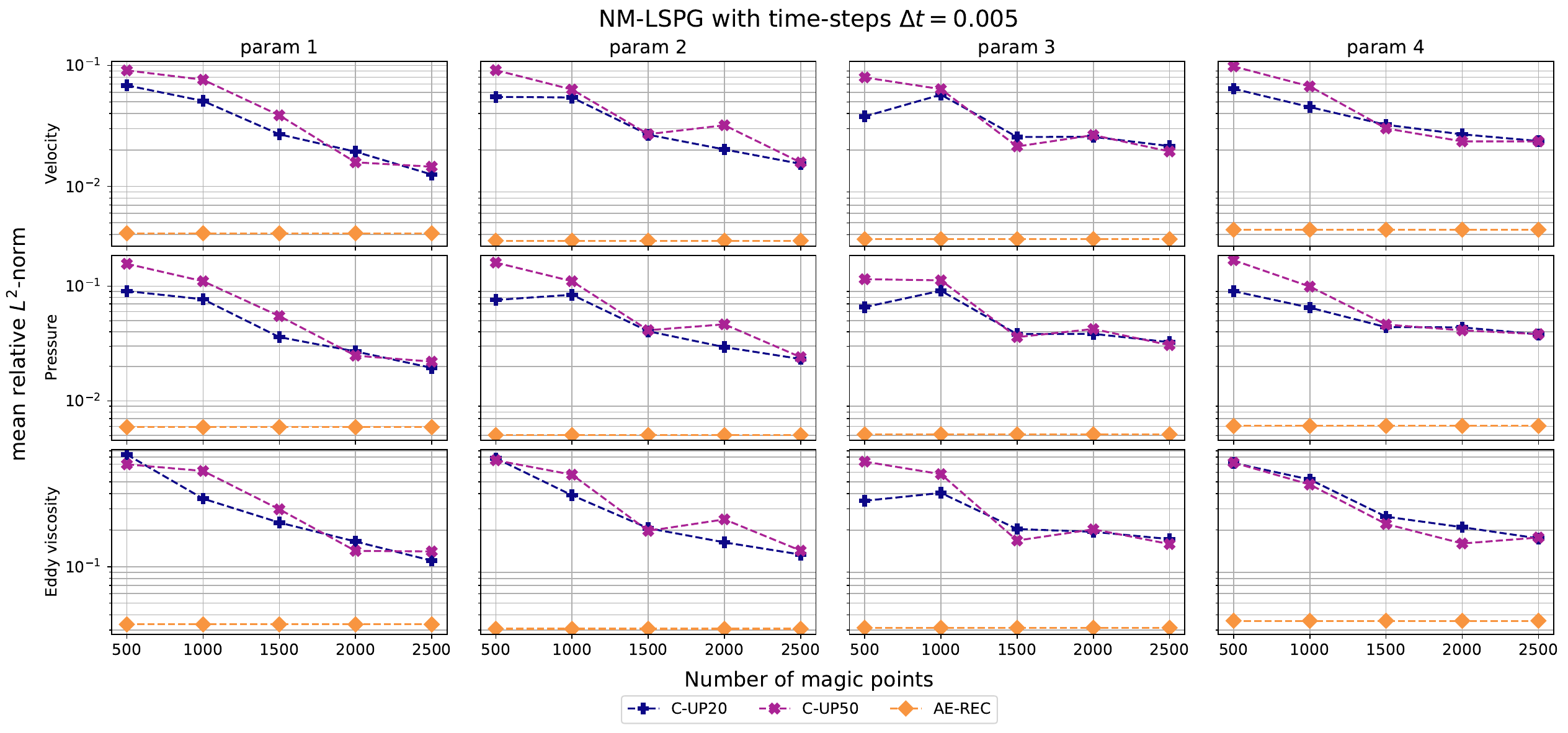}
  \caption{\textbf{INS-1}. The decay of the mean $L^2$ relative error is assessed for the \textbf{C-UP20}, and \textbf{C-UP50} reduced over-collocation methods with sampling strategies corresponding to the gradient-based adaptive one every $50$ and $20$ time steps, respectively. In this case we employ a reference time step of four times $(\Delta t)_{\text{ref}}=0.0005s$ the full-order model one $(\Delta t)_{\text{ref}}=0.0001s$. The baseline represented by the autoencoder mean $L^2$ reconstruction error \textbf{AE-REC} is reported.}
  \label{fig:hrAhmed}
\end{figure}

The convergence with respect to the number of collocation nodes is shown in Figure~\ref{fig:hrAhmed} for the hyper-reduction methods C-UP-20 and C-UP-50 and collocation nodes $r_h\in\{500, 1000, 1500, 2000, 2500\}$. The corresponding timings are reported in Table~\ref{tab: ahmed small timings} for a comparison with the full-order method. The computational cost of the submesh update can be reduced by restricting the collocation nodes that can be selected only to a neighborhood of the current submesh.

\begin{table}[htpb!]
  \centering
  \caption{\textbf{INS-1}. Timings of the INS for small geometrical deformations, subsection~\ref{subsubec:smallAhmed}. The timings of the reduced over-collocation methods \textbf{C-UP20} and \textbf{C-UP50} relative to the choice of $r_h=\{500, 1000, 1500, 2000, 2500\}$ collocation nodes is reported along the full-order model (\textbf{FOM}) timings of the PISO solver. The time step of \textbf{C-UP20} and \textbf{C-UP50} is $(\Delta t)_{\text{ref}}=0.0005$, five times bigger than the \textbf{FOM-1} with $(\Delta t)_{\text{ref}}=0.0001$ on \textbf{1} CPU core. The different reference time steps affect the average total time and not the mean time steps computational costs. The timings for the \textbf{FOM-8} run in parallel with \textbf{8} CPU cores are wallclock times.}
  \footnotesize
  \begin{tabular}{ l | c | c | c}
      \hline
      \hline
      MP  & mean time-step  & mean update every 20 & mean total time\\
      \hline
      \hline
      500  & 66.650 [ms]  & 660.136 [ms]   & 19.931 [s]     \\
      \hline
      1000  & 116.432 [ms]  & 749.300 [ms]   & 30.779 [s]     \\
      \hline
      1500  & 171.507 [ms]  & 809.385 [ms]   & 42.395 [s]     \\
      \hline
      2000  & 212.034 [ms]  & 777.112 [ms]   & 50.178 [s]     \\
      \hline
      2500  & 276.482 [ms]  & 933.621 [ms]   & 64.633 [s]     \\
      \hline
      \hline
      FOM-1 & 791.318 [ms]  & -   & 13.189 [min]    \\
      \hline
      FOM-8 & 365.355 [ms]  & -   & 3.589 [min]    \\
      \hline
      \hline
  \end{tabular}
  \begin{tabular}{ l | c | c}
      \hline
      \hline
     mean time-step  & mean update every 50 & mean total time\\
      \hline
      \hline
       49.171 [ms]  & 572.604 [ms]   & 12.125 [s]     \\
      \hline
      136.320 [ms]  & 778.478 [ms]   & 30.378 [s]     \\
      \hline
      160.974 [ms]  & 842.653 [ms]   & 35.565 [s]     \\
      \hline
      194.183 [ms]  & 723.026 [ms]   & 41.729 [s]     \\
      \hline
      312.016 [ms]  & 973.697 [ms]   & 66.298 [s]     \\
      \hline
      \hline
      791.318 [ms]  & -   & 13.189 [min]    \\
      \hline
      365.355 [ms]  & -   & 3.589 [min]    \\
      \hline
      \hline
  \end{tabular}
  \label{tab: ahmed small timings}
\end{table}

\subsubsection{Large geometrical deformations (INS-2)}
\label{subsubec:largeAhmed}

A situation where the methodology devised may fail is considered. One of the main problems of employing rSVD modes to linearly approximate solution manifolds with respect to nonlinear dimension reduction methods that employ neural networks, is a slow Kolmogorov n-width decay or, from a different point of view, the high generalization error on the test set. This is evident when increasing the dimension of the linear reduced space, the error on the training set decreases, but the error in the test set does not. 

\begin{figure}[bpht!]
  \centering
  \includegraphics[width=1\textwidth, trim={0 0 0 0}, clip]{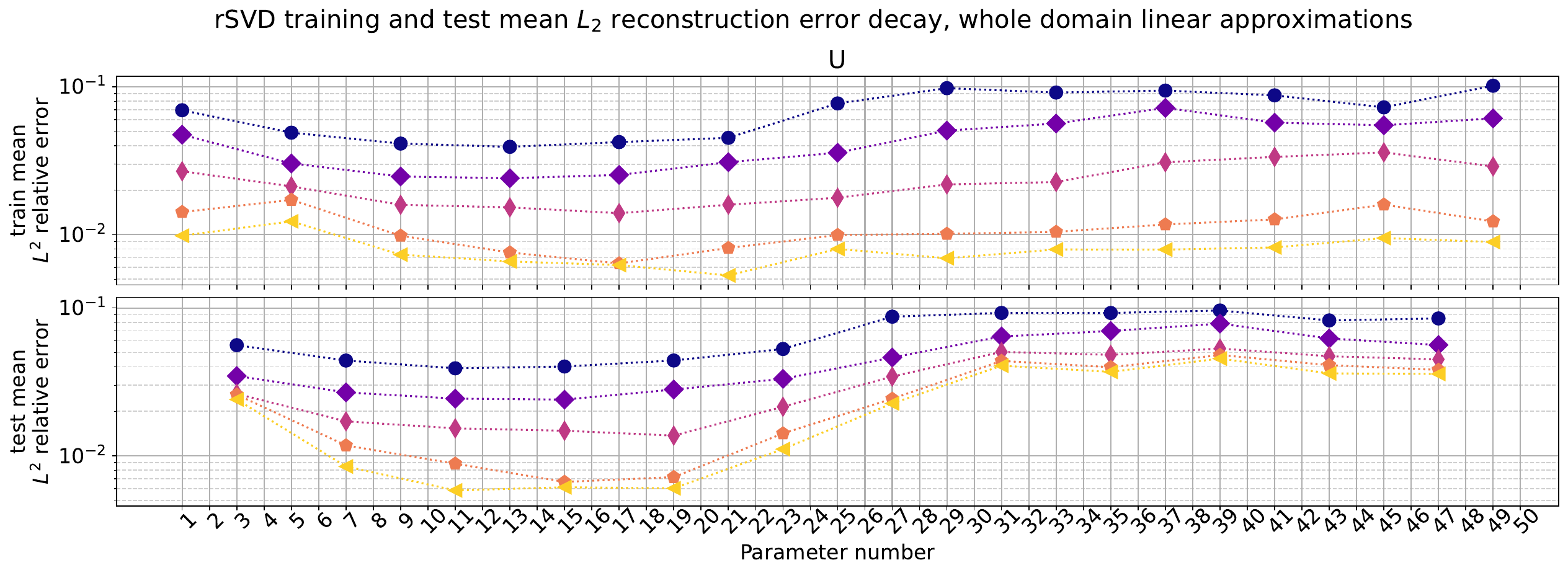}
  \includegraphics[width=1\textwidth, trim={0 0 0 0}, clip]{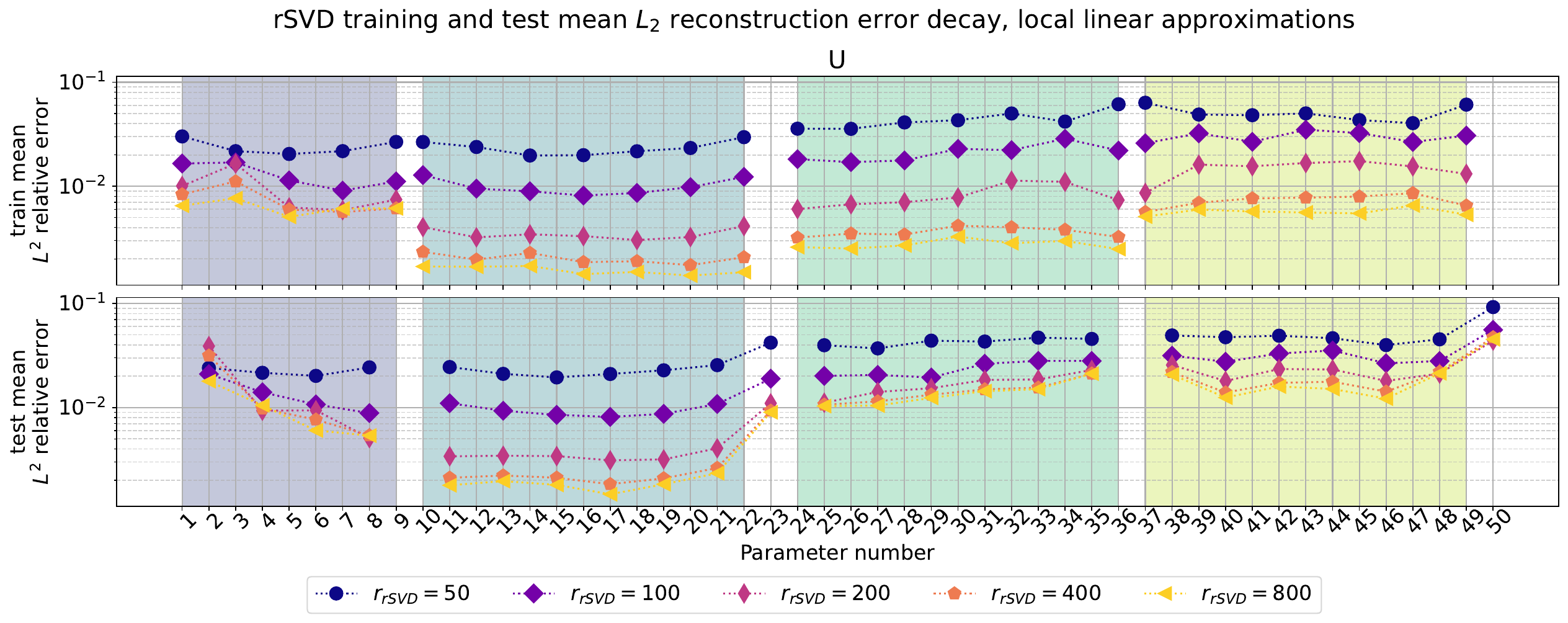}
  \caption{\textbf{INS-2}. Studies of the decay of the mean $L^2$ reconstruction errors defined trough equations~\eqref{eq:rec_rsvd} for the $50$ geometrical deformations in~\eqref{eq:INSLarge}. For brevity, we report only the errors on the velocity field. \textbf{Above:} approximation of the whole solution manifold with $13$ training time series. It is evident the major drawback in employing linear approximants: the generalization error on the test set remains high even if the reduced dimension is increased from $50$ to $800$. \textbf{Below:} approximation of the solution manifold with $4$ local linear approximants. They are highlighted by the shaded regions in background. The leftmost corresponds to the solution manifold studied in the previous section~\ref{subsubec:smallAhmed}.}
  \label{fig:decayLargeAhmed}
\end{figure}

An example of this behavior is shown in Figure~\ref{fig:decayLargeAhmed}. This time we want to approximate the solution manifold corresponding to the parameter range for the slant angle $\theta\in[15^\circ,35^{\circ}]$. We sample uniformly $50$ slant angles:
\begin{align}
  \mathcal{P}=\{&15^{\circ}+ i\cdot \left((35^{\circ}-15^{\circ})/49\right)\vert i\in\{0,\dots,49\}\}\times V,\quad\vert\mathcal{P}\vert=50\cdot 1000,\label{eq:INSLarge}
\end{align}
and, as usual, we number the test parameters in order from $1$ to $50$. The training and test parameters of the previous section~\ref{subsubec:smallAhmed}, correspond to the indices $\{1, 3, 5, 7 ,9\}$ for the training slant angles and $\{2, 4, 6, 8\}$ for the test slant angles.

\begin{figure}[bpht!]
  \centering
  \includegraphics[width=0.49\textwidth]{figures/ahmedDomain/1.pdf}
  \includegraphics[width=0.49\textwidth]{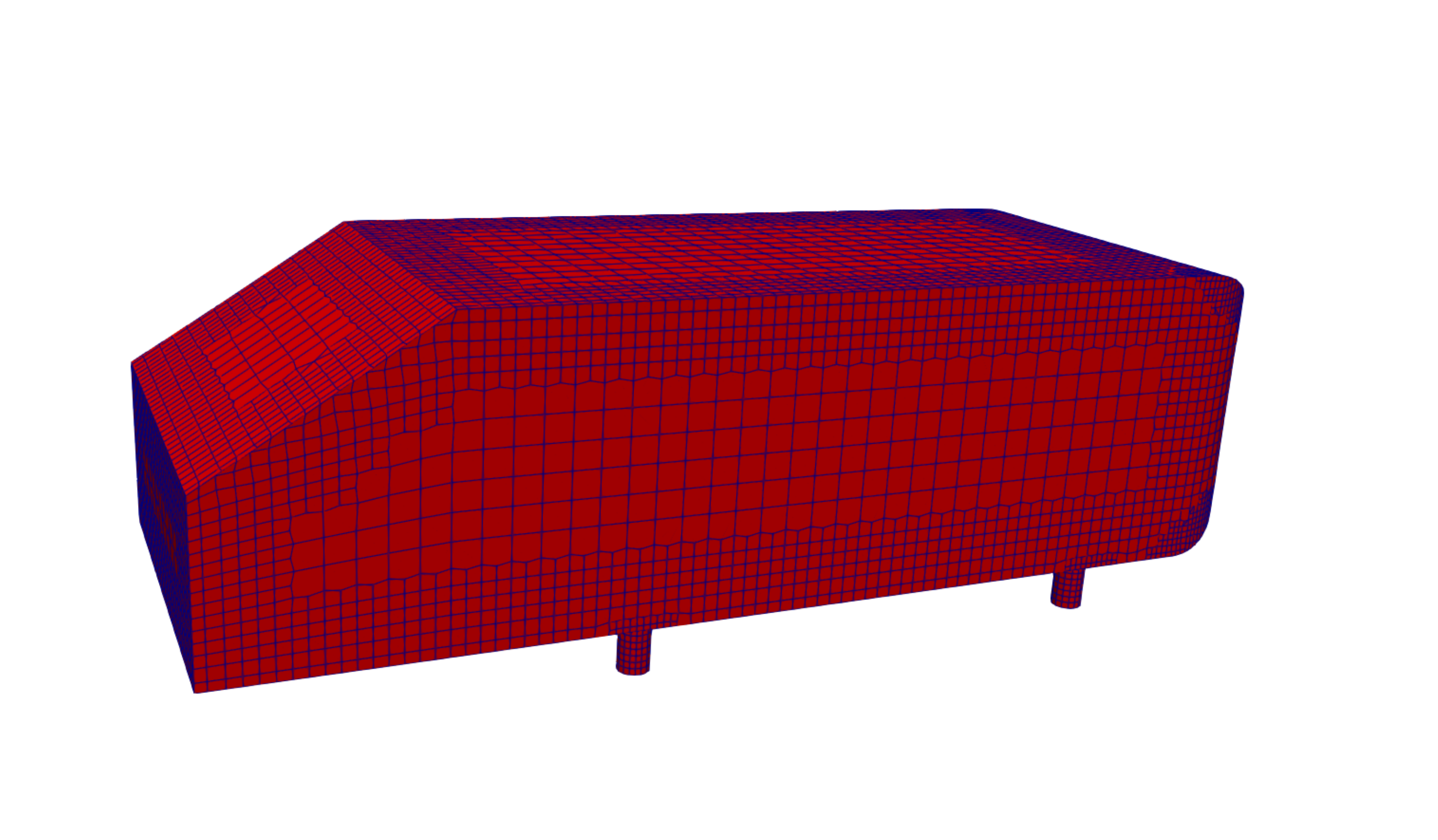}
  \caption{\textbf{INS-2}. \textbf{Left:} Ahmed body with slant angle $\theta=15^{\circ}$ from the training parameter set. \textbf{Right:} Ahmed body with slant angle $\theta=35^{\circ}$ from the training parameter set.}
  \label{fig:angles_ahmed_large}
\end{figure}

With reference to Figure~\ref{fig:angles_ahmed_large} above, it is clear that a linear approximant of the whole solution manifold cannot be used, if the computational budget is limited to only \textbf{13} training time series. Increasing the computational budget, local linear reduced basis achieve substantial improvements. The results in terms of mean $L^2$ reconstruction error are reported in Figure~\ref{fig:angles_ahmed_large} below. Each local linear solution manifold approximant is highlighted by shaded backgrounds of different colors. The leftmost corresponds to the parameters' range of the previous section~\ref{subsubec:largeAhmed}.

Neural networks are known to overcome this problem with truly nonlinear dimension reduction algorithms, that is, with respect to the methodology introduced in this work, without the direct involvement of linear rSVD basis. Possible hyper-reduction strategies that can be applied to a general truly nonlinear neural network architecture are studied in~\cite{romor2023non}. Nevertheless, the methodology presented in this work could effectively be applied for each of the four subdomains in Figure~\ref{fig:decayLargeAhmed} with satisfactory accuracy in terms of test reconstruction error.

\section{Discussion}
\label{sec:discussions}
\label{sec:discussion}
Some crucial considerations not yet underlined, are the subject of this section:
\begin{itemize}
  \item choice of latent dimension. The latent dimensionality equal to $\mathbf{r=4}$ and the number of rSVD modes equal to $\mathbf{r_{\text{rSVD}}=150}$ or $\mathbf{r_{\text{rSVD}}=300}$ are not chosen after parametric studies and therefore more optimal values for the problems at hand can be generally found. As long as we obtain a satisfactory approximation of the solution manifold for some values of the latent dimension of the autoencoder and of the reduced dimension of the linear projections, we did not change them, thus exploiting the possibility to employ the same neural networks architectures shown in Tables~\ref{tab: cae 150} and~\ref{tab: cae 300} for our numerical investigations. In fact, our focus is in the hyper-reduction methodologies applied after a relatively accurate nonlinear approximant of the solution manifold is obtained. Fixed a tolerance for the $L^2$ relative reconstruction error, the best value for $r_{\text{rSVD}}$ can be efficiently found observing the decay of the reconstruction error as done in section~\ref{subsubec:largeAhmed} for the INS-2 test case. Also, the latent dimensionality of the autoencoder cannot be established \textit{a priori}, for some theoretical bounds see~\cite{franco2023deep}. Nevertheless, many empirical studies can be performed to determine the latent dimensionality~\cite{denti2022generalized}.
  \item non collocated hyper-reduction methods. Due to our choice of monolithic normalization and reduction of the physical fields of interest as described in section~\ref{subsec:rSVD}, the non collocated versions of the hyper-reduction methods introduced in section~\ref{sec:hr} do not perform well in terms of accuracy as shown in Figure~\ref{fig:coarseHrConvergence}. The same gappy DEIM implementation performs well when applied in synergy with teacher-student training of a reduced decoder for a 2d nonlinear conservation law parametric model in~\cite{romor2023non}, the only physical field considered is velocity, so normalization is not needed. Anyway, in terms of efficiency, collocation methods are faster since they do not perform additional matrix-vector multiplications (for DEIM and DEIM-SOPT) or scalar products (ECSW), after the evaluation of the residuals at the magic points. Moreover, employing adaptive hyper-reduction strategies in the online stage has higher computational costs for non collocated approaches since psuedo-inverse matrices need to be evaluated (for DEIM).
  \item stability issues. Reduced over-collocation methods may bring unstable numerical schemes, especially if the solution manifold is not approximated with enough accuracy by the nonlinear approximant. Possible solutions involve the training of NN architectures with a more smooth latent space and more regular maps such as those that should be guaranteed by variational autoencoders~\cite{kingma2013auto} or other machine learning architectures. On this line of thought, many additional inductive biases can be imposed. However, from the point of view of numerical analysis, stabilization strategies for ROMs have frequently been employed also for the classical projection-based methotodologies as they often suffer from stability issues, especially when hyper-reduction is applied. Possible solutions include structure-preserving and/or regularized versions~\cite{klein2023structure, herkert2023dictionary} of hyper-reduction methods.
  \item inductive biases and autoencoder regularity. The training and the NN architecture itself can be enriched by inductive biases. Among them, there are first principles (e.g. conservation laws), geometrical simmetries (group equivariant filters~\cite{bronstein2017geometric}), latent operators/numerical schemes (e.g. operator inference~\cite{peherstorfer2016data, uy2023operator}) latent regularity (supposedly imposed by variational autoencoders and other architectures), structure-preservation~\cite{buchfink2023symplectic}, other numerical and mathematical properties (e.g. positiveness).
\end{itemize}

\section{Conclusions}
\label{sec:conclusions}
We developed and tested on challenging benchmarks a new method that exploits nonlinear solution manifold approximants. It relies on convolutional autoencoders and linear transforms/filter maps (specifically parallel randomized singular value decomposition) to approximate solution manifolds even when affected by a moderately slow Kolmogorov n-width decay. The main novelty resides in the implementation of new efficient collocated and adaptive gradient-based hyper-reduction strategies specifically tailored for our choice of nonlinear approximants. Local solution manifold approximations and efficient ways to perform the change of basis are also taken into consideration. We managed to achieve significant speedups while keeping a satisfactory accuracy even when we considered small benchmarks in terms of degrees of freedom, implemented in \texttt{OpenFoam}. Our test cases model complex physics such as the compressible and incompressible turbulent Navier-Stokes equations and suffer from a moderately slow Kolmogorov n-width decay that would need in the order of hundreds of linear basis to be well-approximated by classical projection-based ROMs.  

The crucial objective that we want to reach with our methodology is the development of numerically and physically explicable model order reduction methods while exploiting machine learning architectures. The majority of scientific machine learning strategies employ neural networks in the predictive/online stage as black box surrogate models, without exploiting the underlying physics of the models embedded in the numerical schemes of the full-order models' solvers. Differently, our approach efficiently exploits the numerical schemes also in the predictive phase, with the possibility to characterize the latent solutions found as minimizers of the residuals of conservation laws. The interpretability of the results is evidently increased.

The main disadvantages regard the employment of linear basis: parametric models that suffer for a slow Kolmogorov n-width decay, such as the incompressible flow around the Ahmed body with geometrical deformations studied in section~\ref{subsubec:largeAhmed}, may require too many computational resources to obtain local linear approximations of the solution manifold. The use of more generic truly nonlinear neural networks architectures should bring lower generalization errors with less training data. Some implementations of hyper-reduction for generic NN architectures involve teacher-student training and are presented in~\cite{romor2023non}.

Other aspects that can be substantially improved are stabilization issues: the development of stabilization mechanisms that would aid the nonlinear least-squares optimizers in the search for a physically accurate latent solution would greatly improve our methodology. On the same subject, structure-preserving and other regularizing frameworks, without mentioning additional useful inductive biases, would also help to achieve more accurate solutions.

\section*{Acknowledgements}
This work was partially funded by European Union Funding for Research and
Innovation --- Horizon 2020 Program --- in the framework of European Research
Council Executive Agency: H2020 ERC CoG 2015 AROMA-CFD project 681447 ``Advanced
Reduced Order Methods with Applications in Computational Fluid Dynamics'' P.I.
Professor Gianluigi Rozza. We also acknowledge the PRIN 2017 “Numerical Analysis for Full and Reduced Order Methods for the efficient and accurate solution of complex
systems governed by Partial Differential Equations” (NA-FROM-PDEs)

\appendix
\renewcommand\theequation{\arabic{equation}}
\section{Architectures}
\label{sec:arch}
\begin{table}[htp!]
  \caption{The 1d Convolutional Autoencoder for the CNS1 on the coarse mesh in subsection~\ref{subsec:coarseAirfoil}.}
  \centering
  \footnotesize
  \begin{tabular}{ l | c | c | c | c}
      \hline
      \hline
      Encoder & Activation & Weights & Padding & Stride \\
      \hline
      \hline
      Conv1d & ELU   & [1, 8, 4, 4] & 1 & 2 \\
      \hline
      Conv1d & ELU   & [8, 16, 4, 4] & 1 & 2 \\
      \hline
      Conv1d & ELU   & [16, 32, 4, 4] & 1 & 2 \\
      \hline
      Conv1d & ELU   & [32, 64, 4, 4] & 1 & 2 \\
      \hline
      Conv1d & ELU   & [64, 128, 4, 4] & 1 & 2 \\
      \hline
      Linear & ELU   & [512, 4] & - \\
      \hline
      \hline
  \end{tabular}
  \hspace{1mm}
  \begin{tabular}{ l | c | c | c | c }
      \hline
      \hline
      Decoder & Activation & Weights & Padding & Stride \\
      \hline
      \hline
      Linear & ELU   & [4, 512] & - \\
      \hline
      ConvTr1d & ELU   & [128, 64, 5, 5] & 1 & 2 \\
      \hline
      ConvTr1d & ELU   & [64, 32, 4, 4] & 1 & 2 \\
      \hline
      ConvTr1d & ELU   & [32, 16, 5, 5] & 1 & 2 \\
      \hline
      ConvTr1d & ELU   & [16, 8, 5, 5] & 1 & 2 \\
      \hline
      ConvTr1d & ReLU   & [8, 1, 4, 4] & 1 & 2 \\
      \hline
      \hline
  \end{tabular}
  \label{tab: cae 150}
\end{table}

\begin{table}[htp!]
  \centering
  \caption{The 1d Convolutional Autoencoder for the CNS2 on the finer mesh in subsections~\ref{subsubsec:sudomains_finer_airfoil} and for the INS1 in subsection~\ref{subsubec:smallAhmed}.}
  \footnotesize
  \begin{tabular}{ l | c | c | c | c}
      \hline
      \hline
      Encoder & Activation & Weights & Padding & Stride \\
      \hline
      \hline
      Conv1d & ELU   & [1, 8, 4, 4] & 1 & 2 \\
      \hline
      Conv1d & ELU   & [8, 16, 4, 4] & 1 & 2 \\
      \hline
      Conv1d & ELU   & [16, 32, 4, 4] & 1 & 2 \\
      \hline
      Conv1d & ELU   & [32, 64, 4, 4] & 1 & 2 \\
      \hline
      Conv1d & ELU   & [64, 128, 4, 4] & 1 & 2 \\
      \hline
      Linear & ELU   & [512, 4] & - \\
      \hline
      \hline
  \end{tabular}
  \hspace{1mm}
  \begin{tabular}{ l | c | c | c | c }
      \hline
      \hline
      Decoder & Activation & Weights & Padding & Stride \\
      \hline
      \hline
      Linear & ELU   & [4, 512] & - \\
      \hline
      ConvTr1d & ELU   & [128, 64, 5, 5] & 1 & 2 \\
      \hline
      ConvTr1d & ELU   & [64, 32, 4, 4] & 1 & 2 \\
      \hline
      ConvTr1d & ELU   & [32, 16, 5, 5] & 1 & 2 \\
      \hline
      ConvTr1d & ELU   & [16, 8, 5, 5] & 1 & 2 \\
      \hline
      ConvTr1d & ReLU   & [8, 1, 4, 4] & 1 & 2 \\
      \hline
      \hline
  \end{tabular}
  \label{tab: cae 300}
\end{table}

\section{Predicted snapshots}
\label{sec:snaps}
\begin{figure}[ht!]
  \centering
  \includegraphics[width=0.325\textwidth, trim={0 0 500 0}, clip]{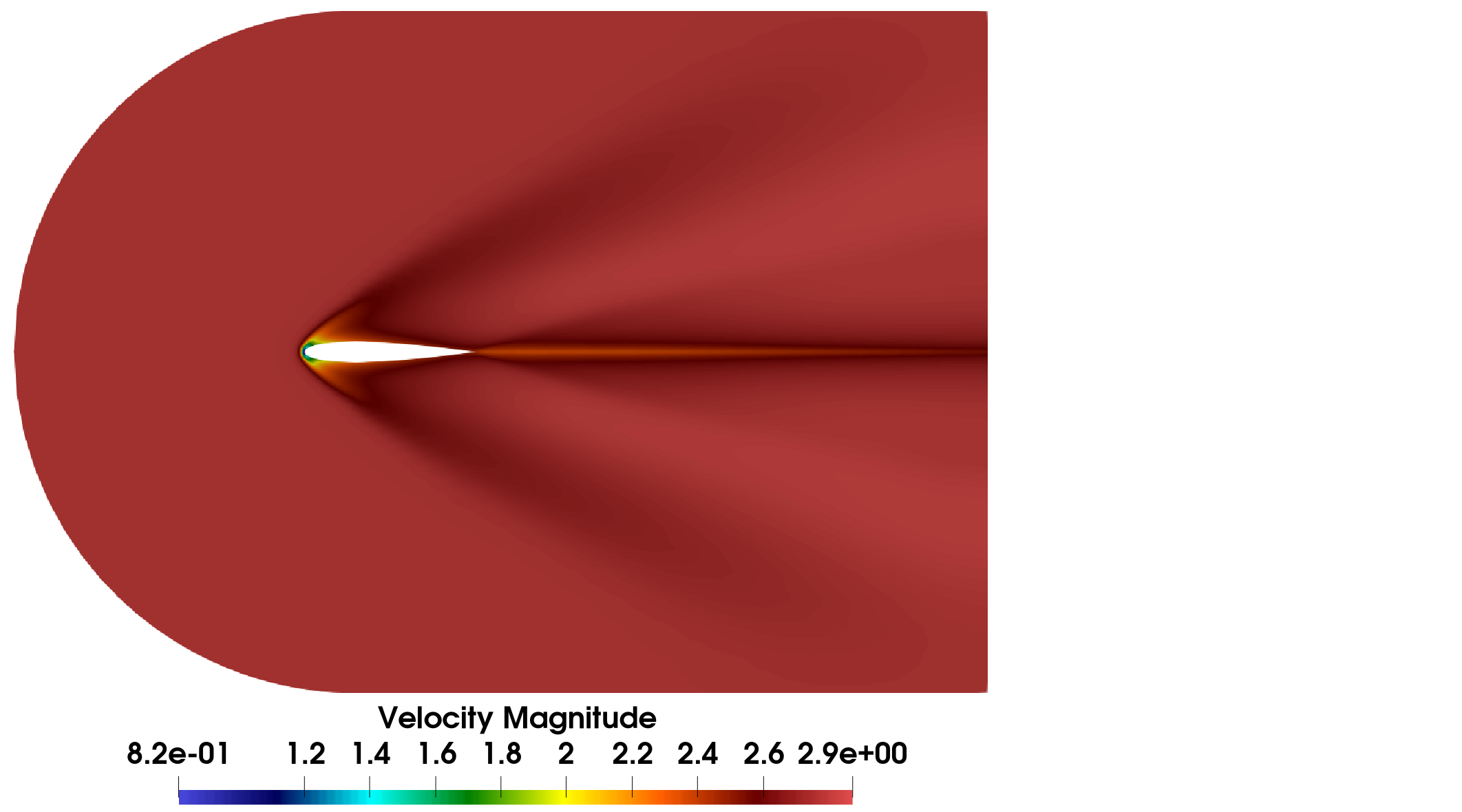}
  \includegraphics[width=0.325\textwidth, trim={0 0 500 0}, clip]{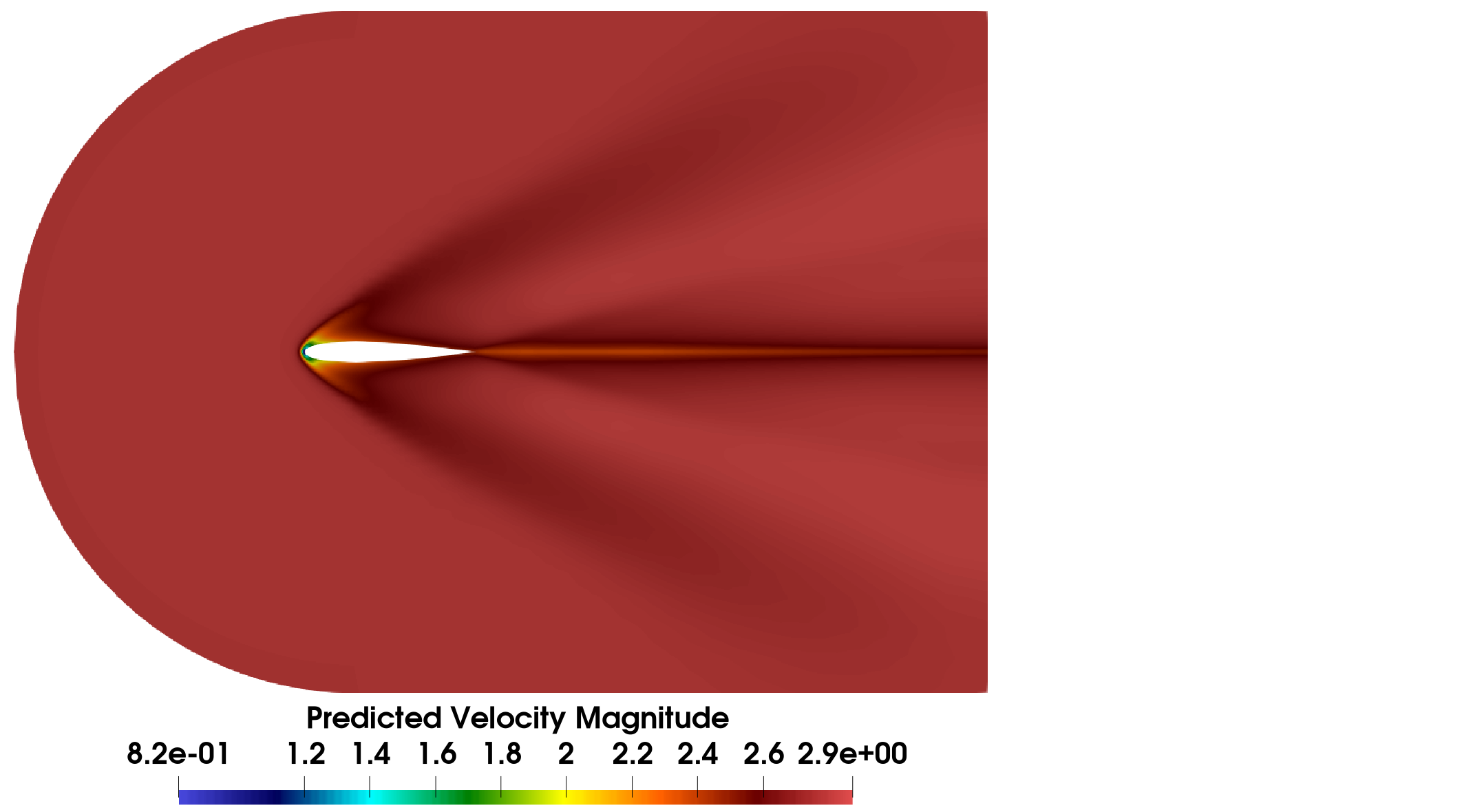}
  \includegraphics[width=0.325\textwidth, trim={0 0 500 0}, clip]{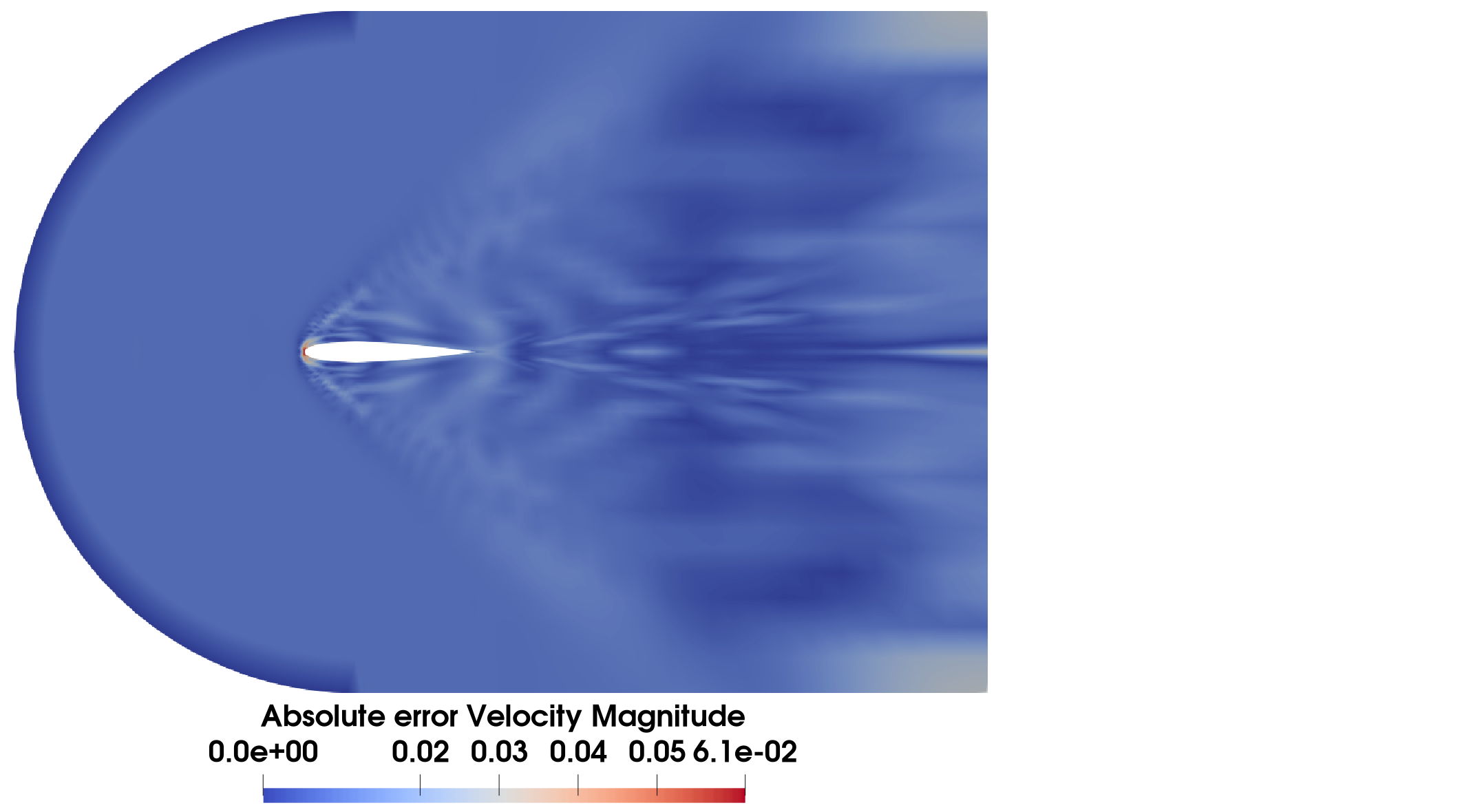}\\
  \includegraphics[width=0.325\textwidth, trim={0 0 500 0}, clip]{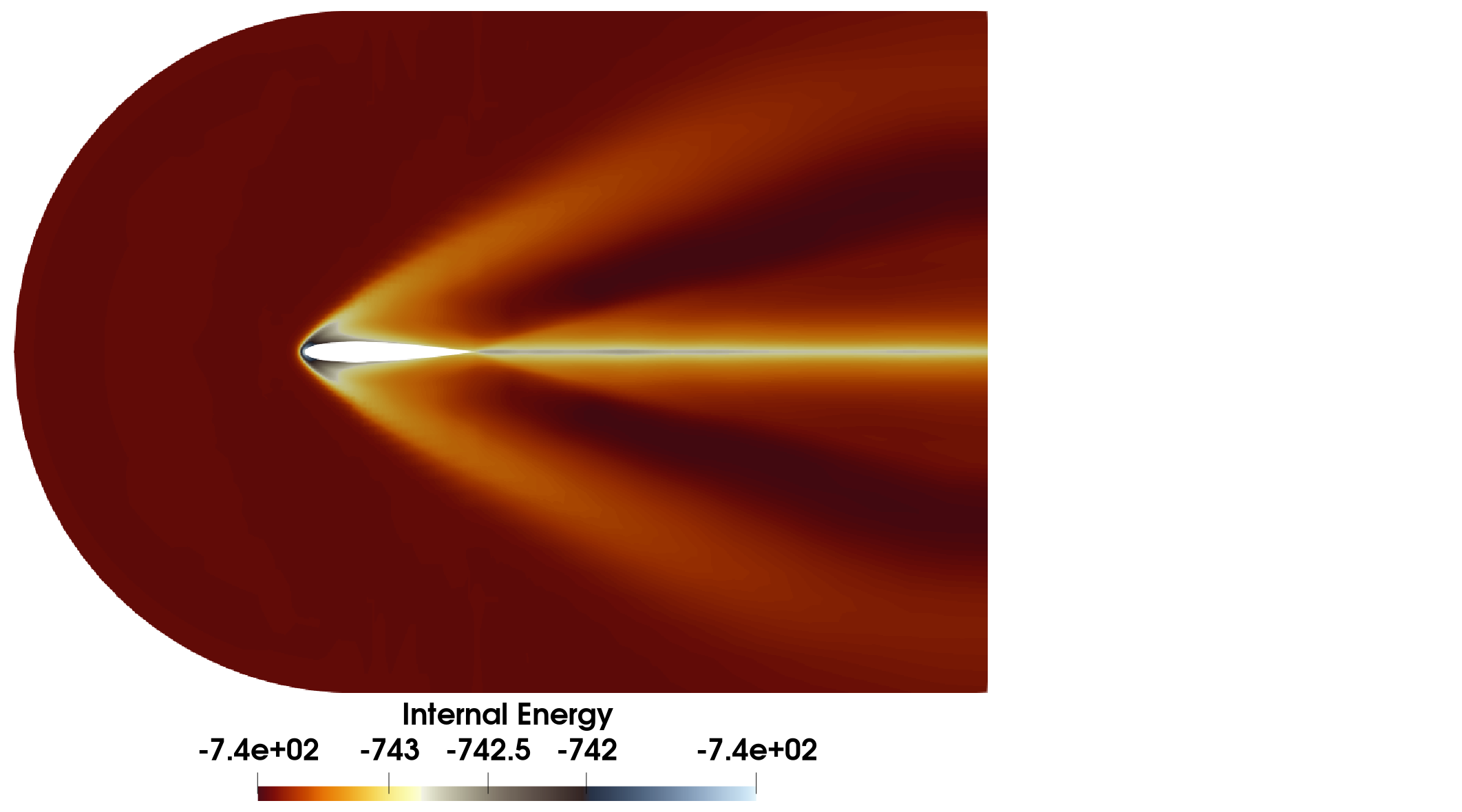}
  \includegraphics[width=0.325\textwidth, trim={0 0 500 0}, clip]{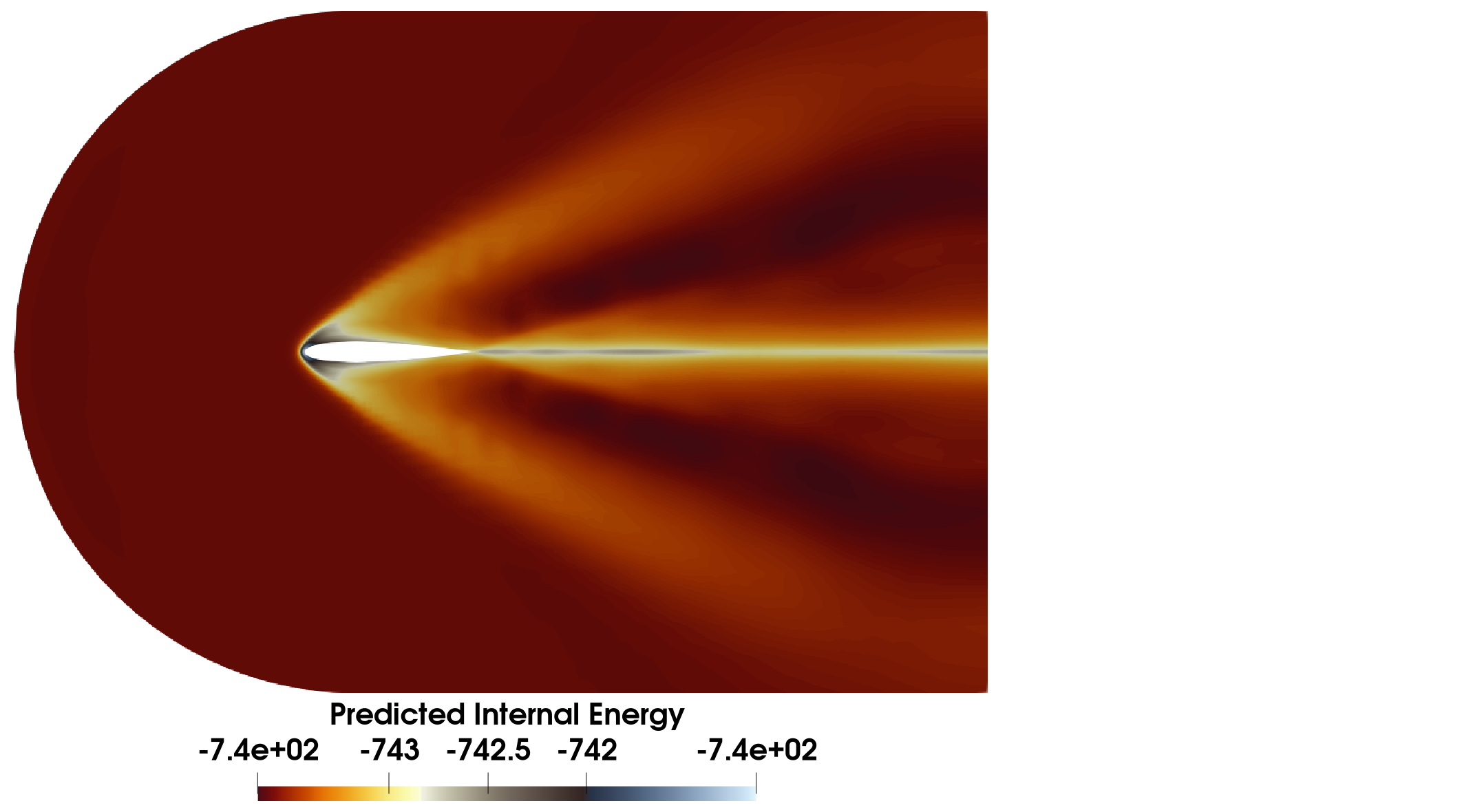}
  \includegraphics[width=0.325\textwidth, trim={0 0 500 0}, clip]{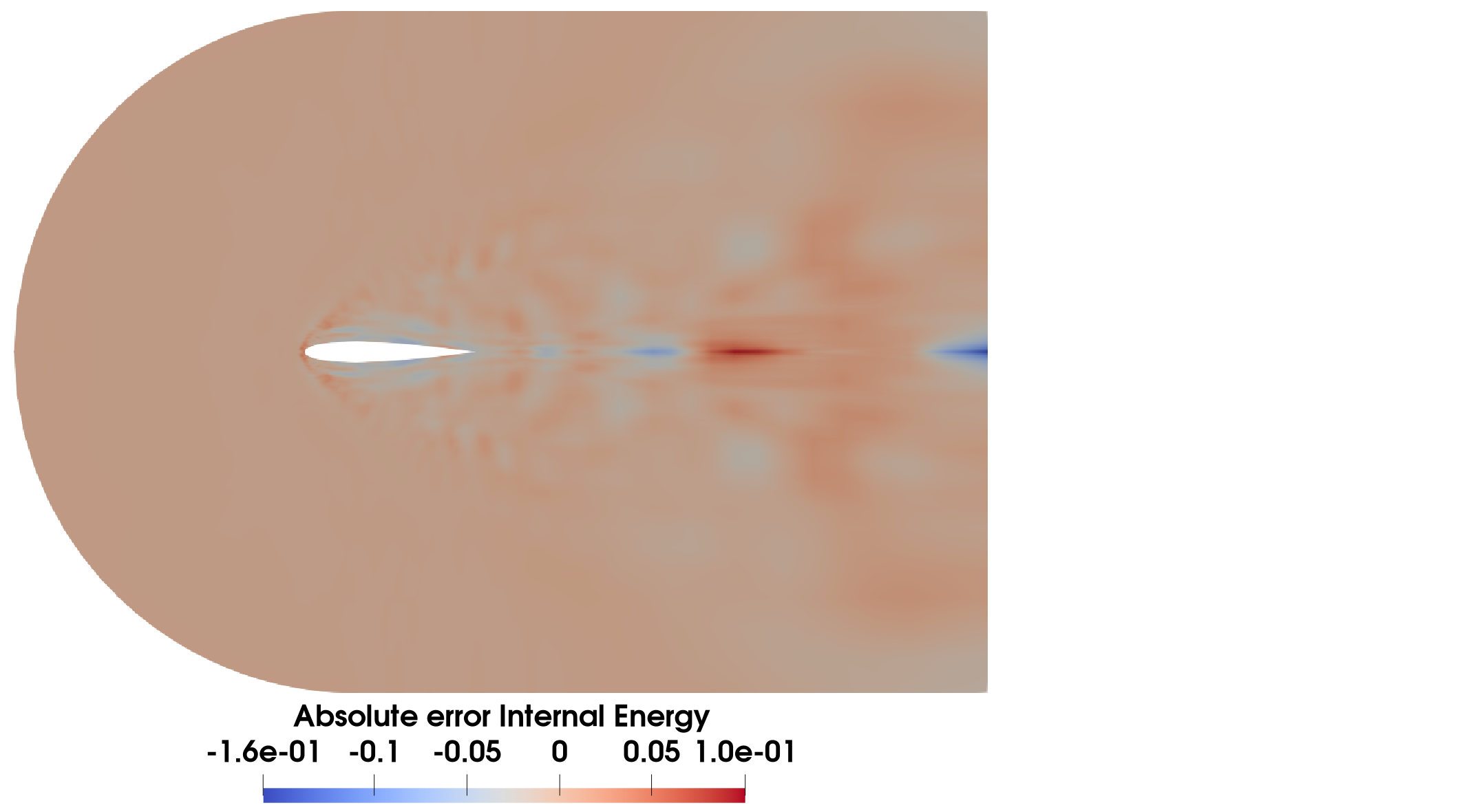}\\
  \includegraphics[width=0.325\textwidth, trim={0 0 500 0}, clip]{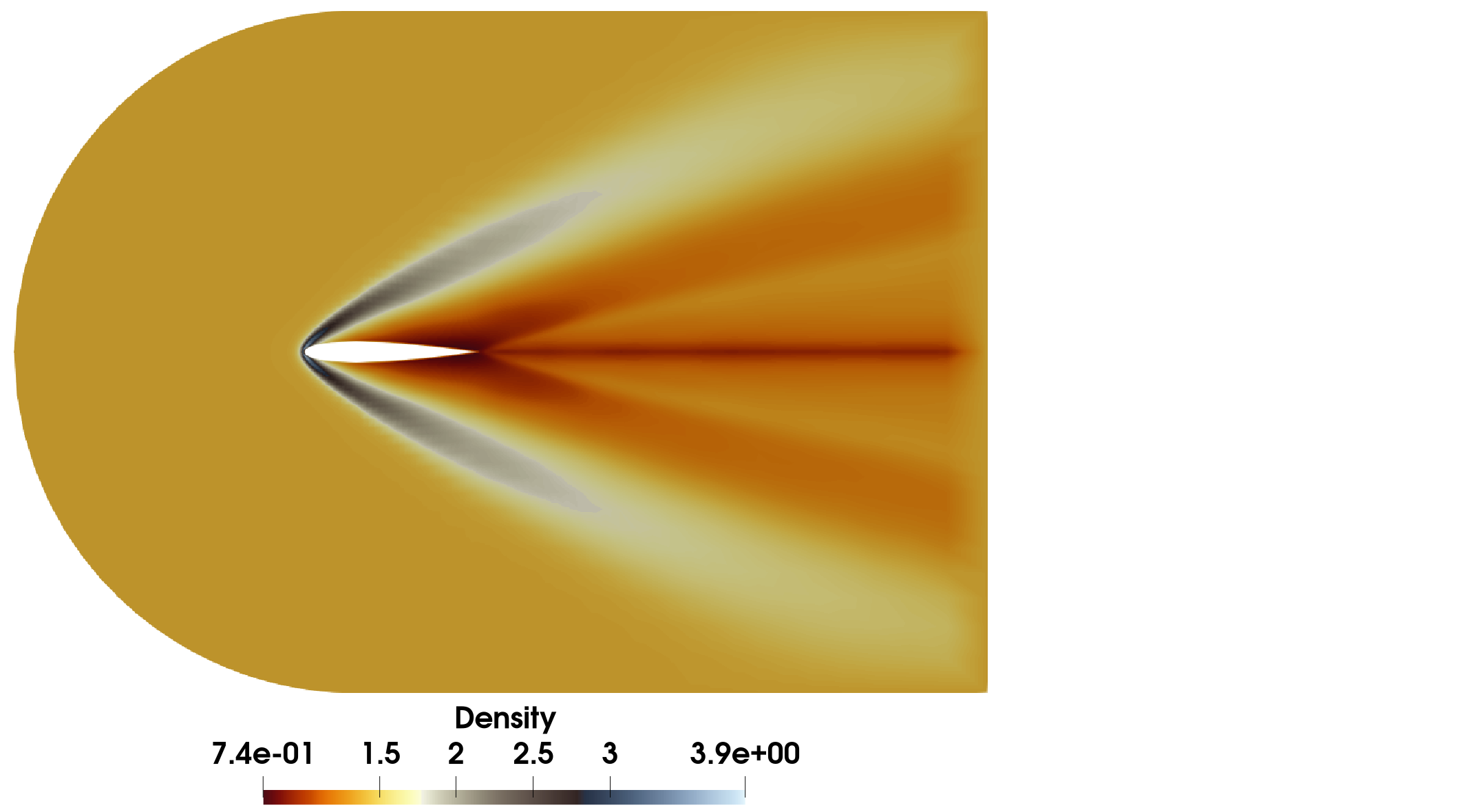}
  \includegraphics[width=0.325\textwidth, trim={0 0 500 0}, clip]{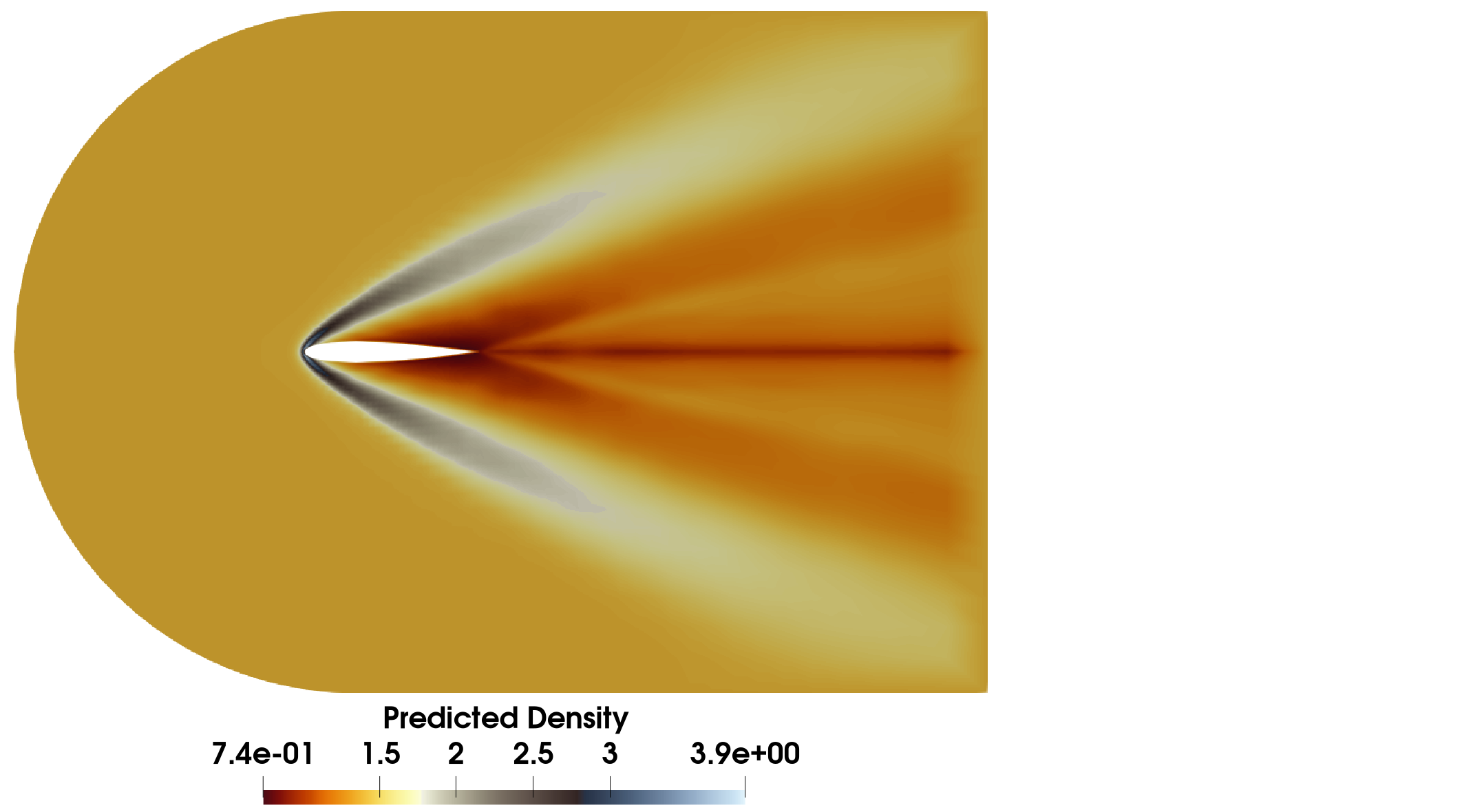}
  \includegraphics[width=0.325\textwidth, trim={0 0 500 0}, clip]{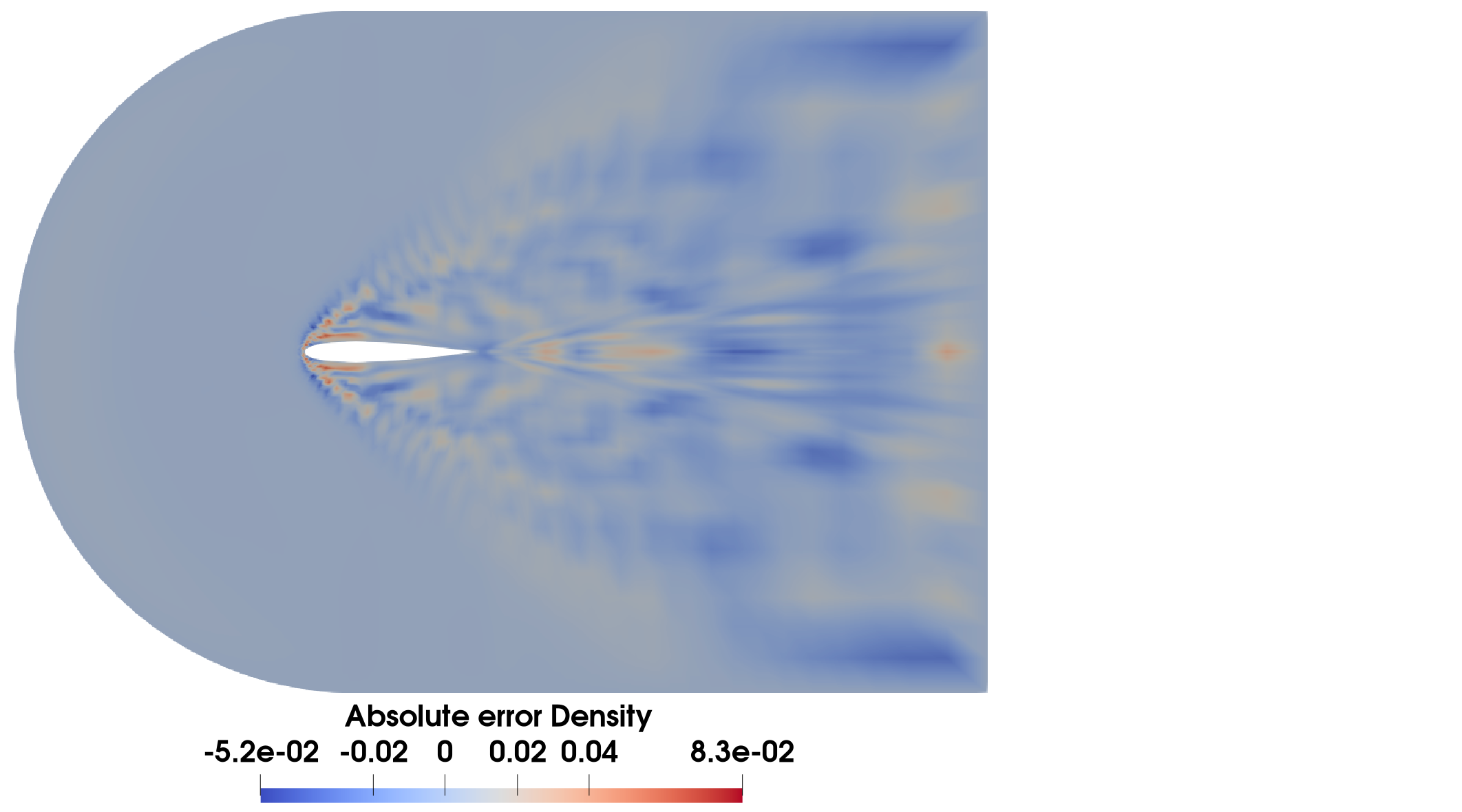}\\
  \includegraphics[width=0.325\textwidth, trim={0 0 500 0}, clip]{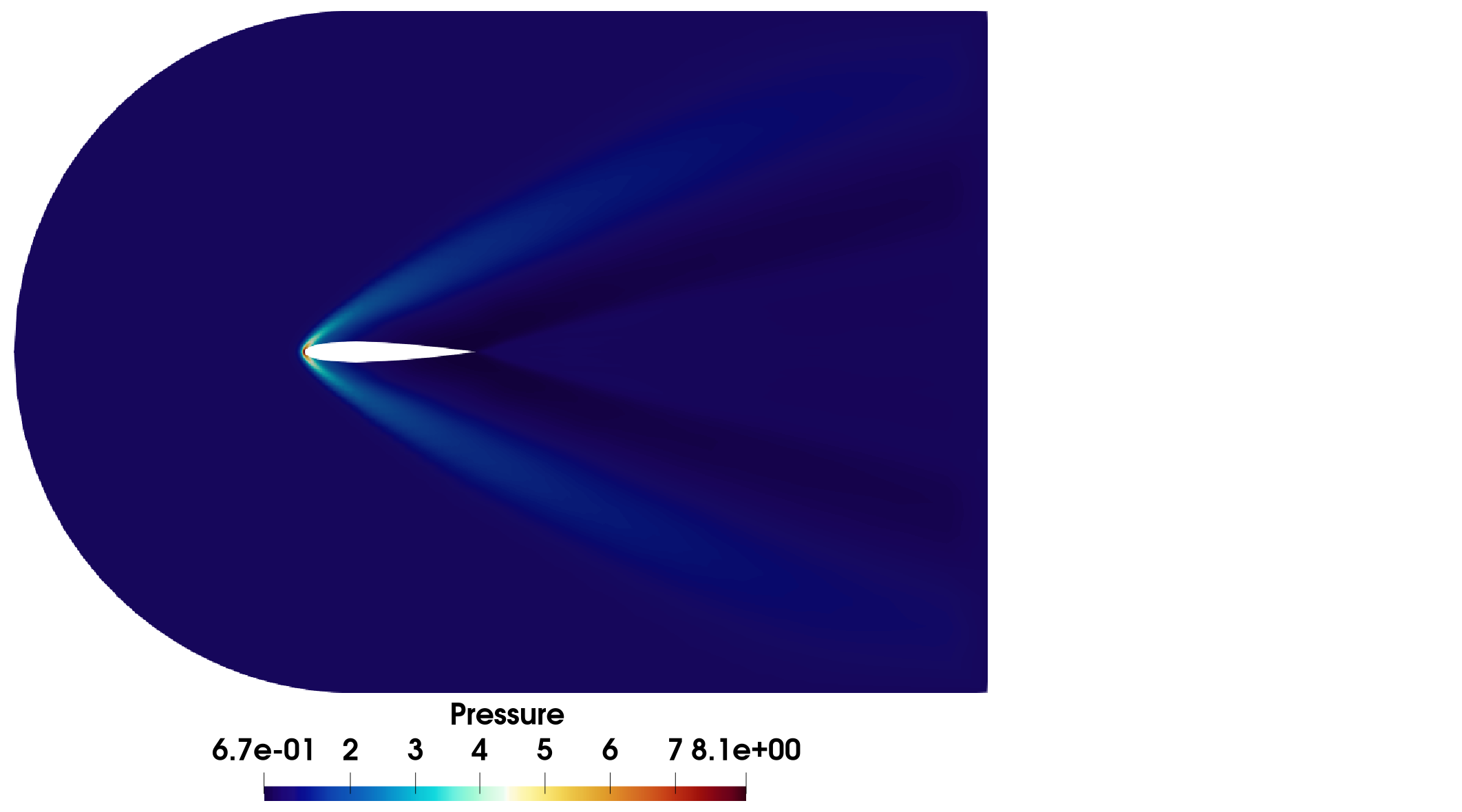}
  \includegraphics[width=0.325\textwidth, trim={0 0 500 0}, clip]{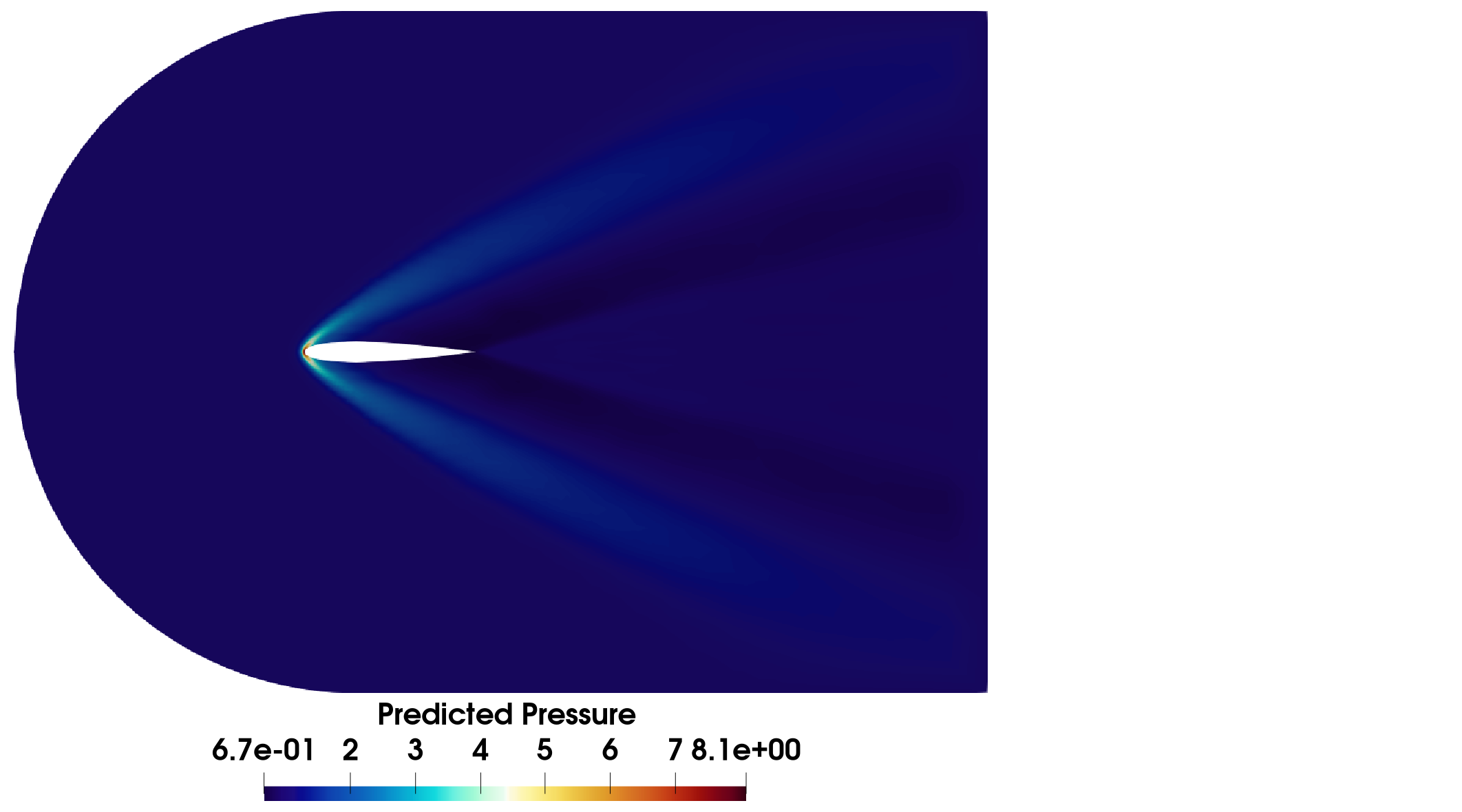}
  \includegraphics[width=0.325\textwidth, trim={0 0 500 0}, clip]{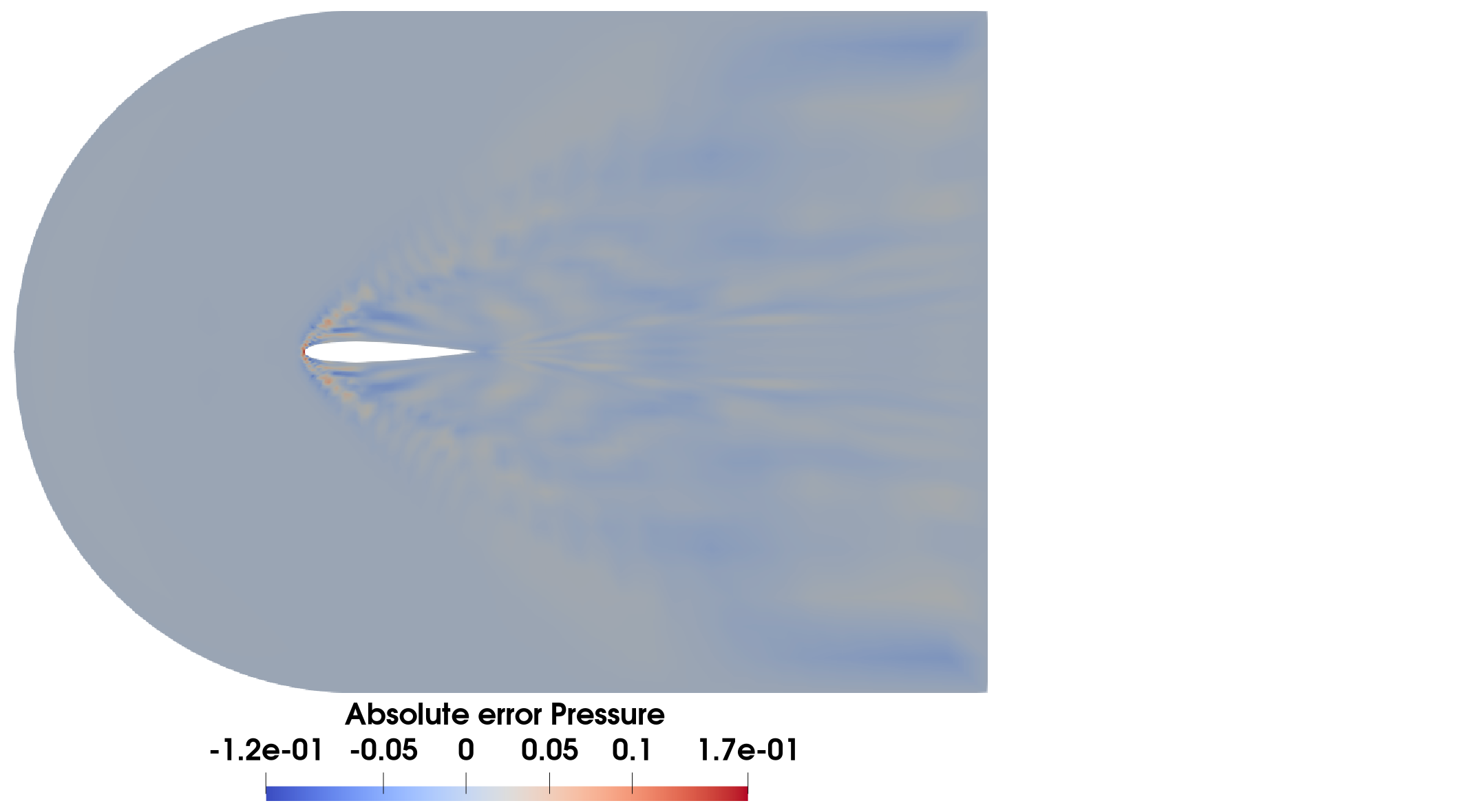}\\
  \caption{\textbf{CNS1.} Predicted velocity, density, internal energy and pressure fields for the test case described in~\cref{subsec:coarseAirfoil} with $\text{Ma}=3.314$ corresponding to test parameter number \textbf{3} at the final time instant $t=2.5\text{s}$. The number of cells is \textbf{4500} the total number of degrees of freedom is \textbf{27000}. The number of collocation nodes is $r_h=\textbf{500}$}
  \label{fig:snapsCoarseAirfoil}
\end{figure}

\begin{figure}[ht!]
  \centering
  \includegraphics[width=0.325\textwidth, trim={0 0 500 0}, clip]{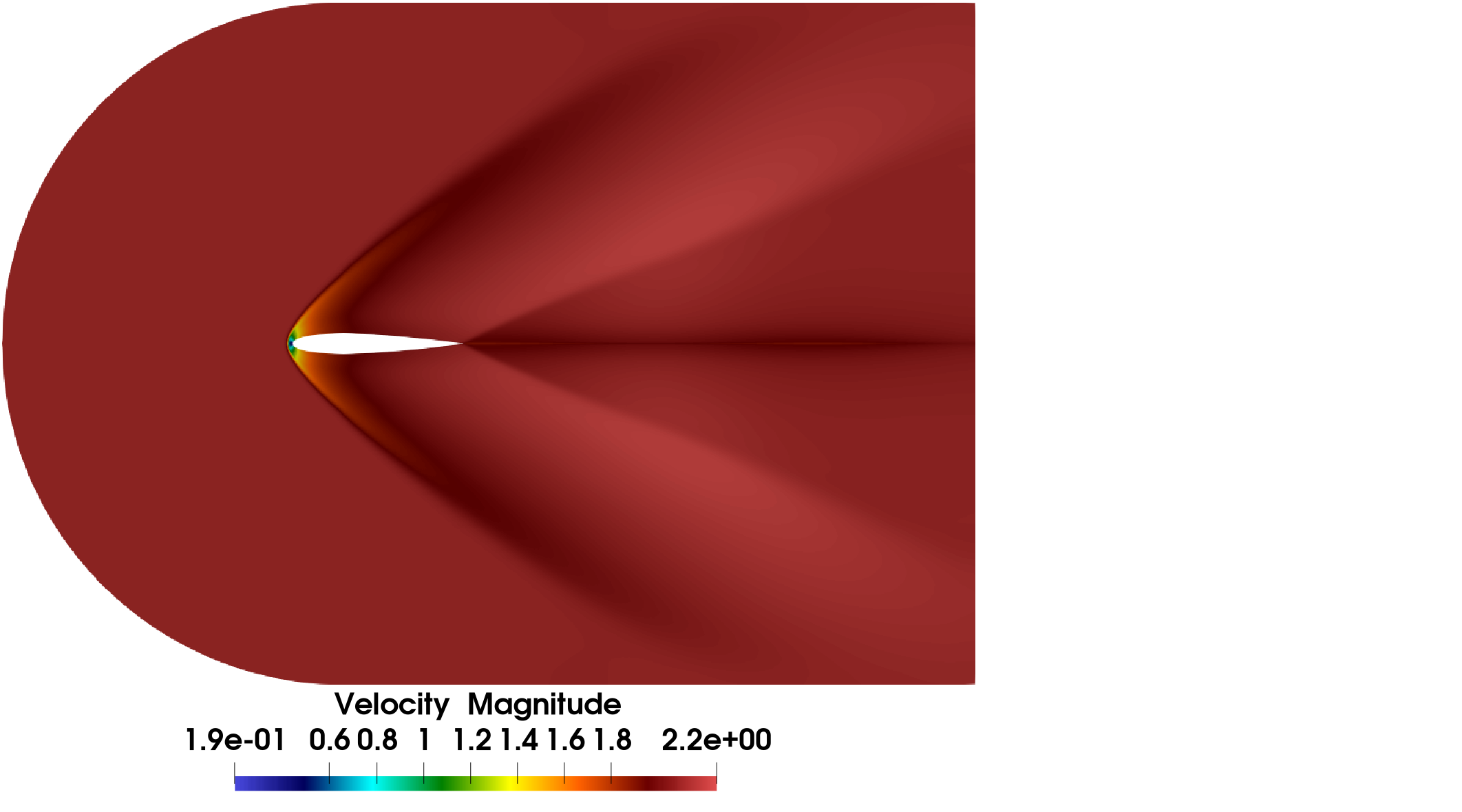}
  \includegraphics[width=0.325\textwidth, trim={0 0 500 0}, clip]{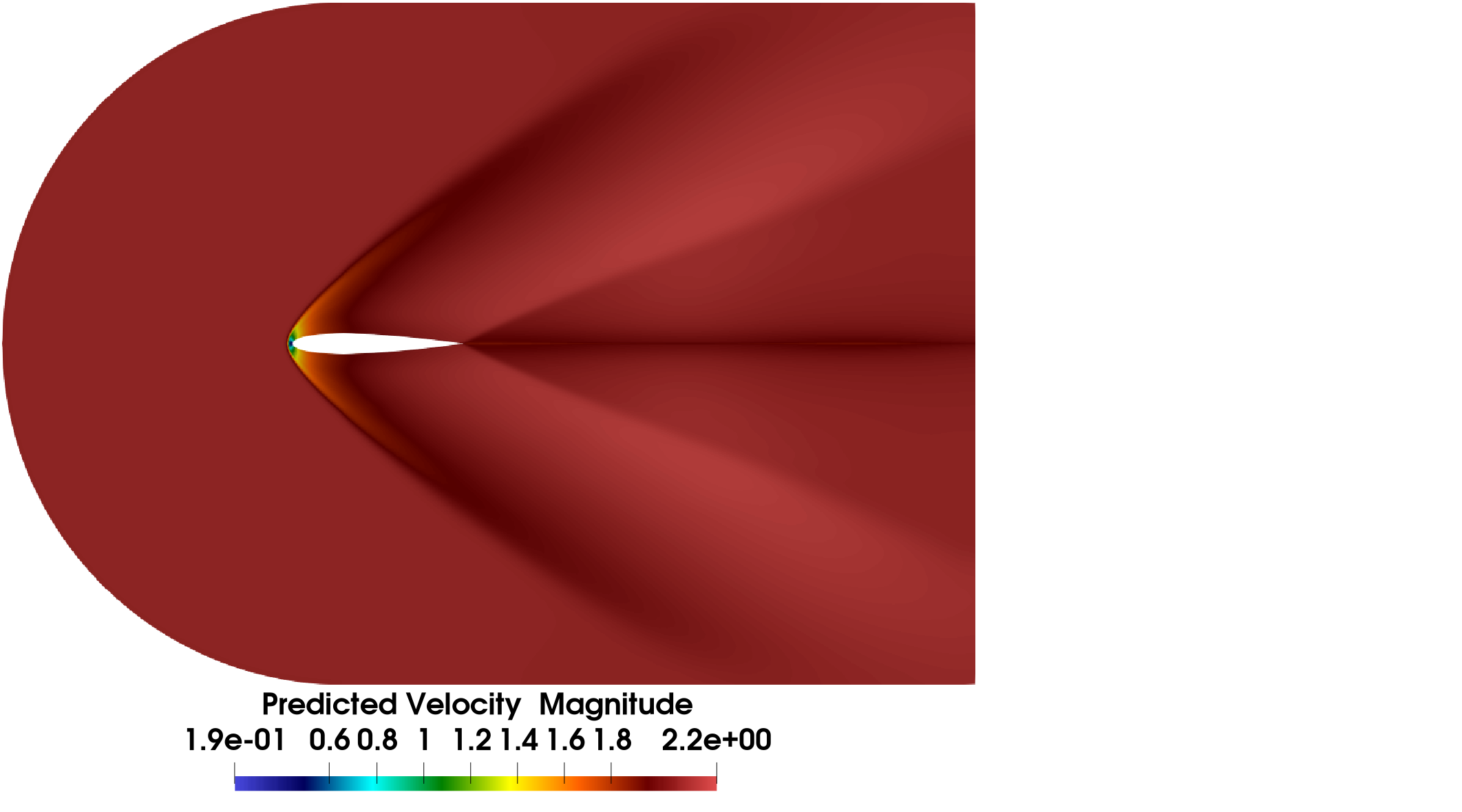}
  \includegraphics[width=0.325\textwidth, trim={0 0 500 0}, clip]{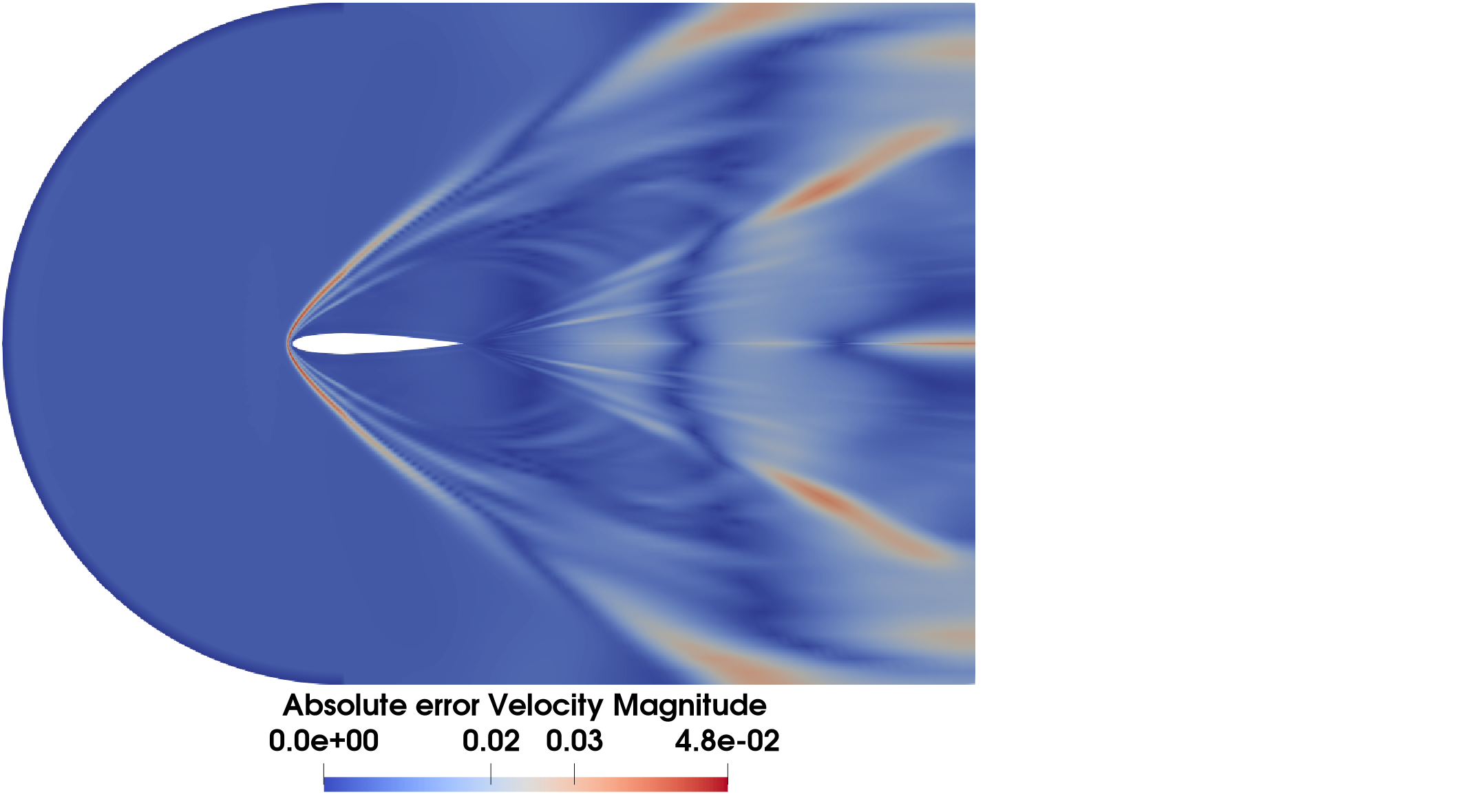}\\
  \includegraphics[width=0.325\textwidth, trim={0 0 500 0}, clip]{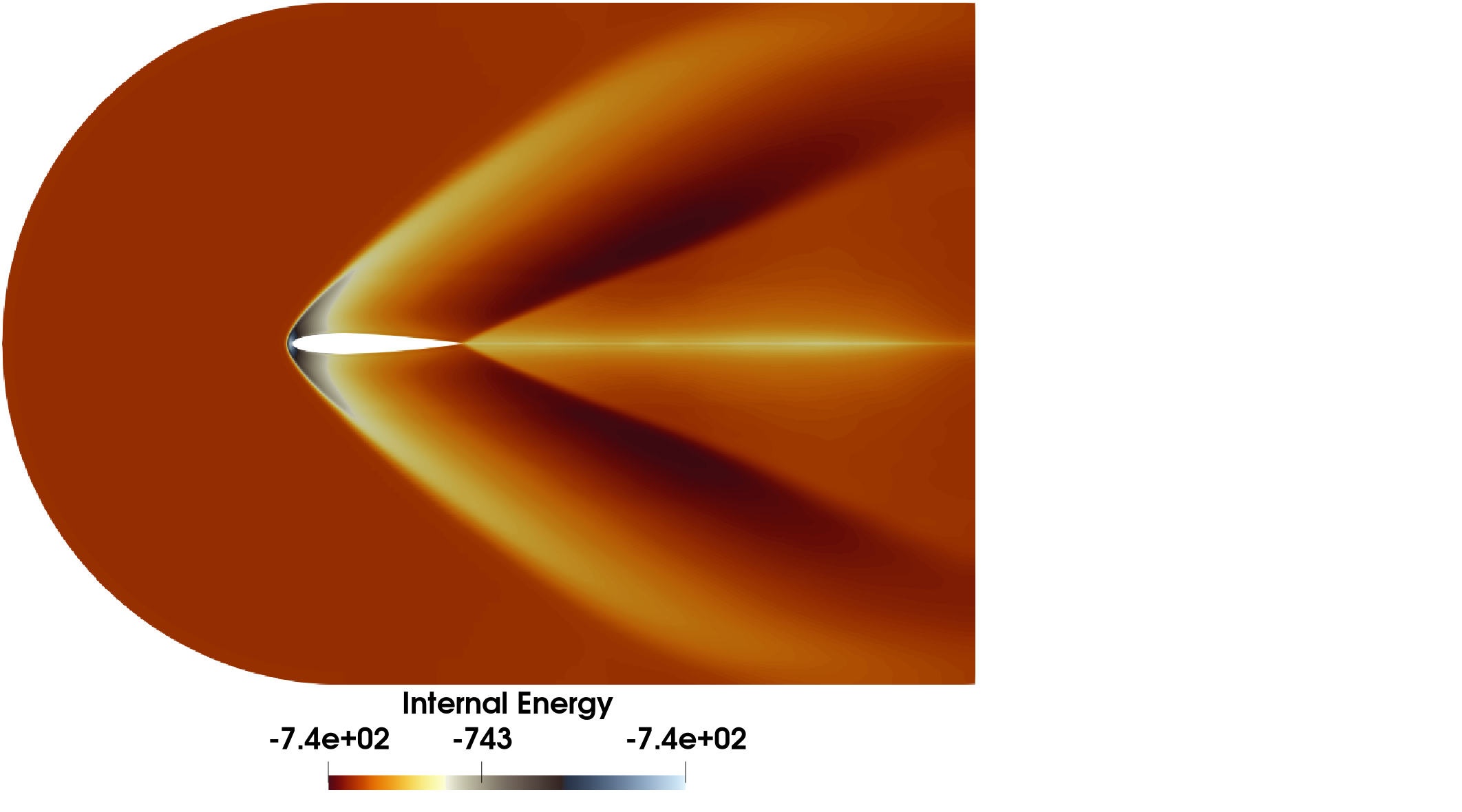}
  \includegraphics[width=0.325\textwidth, trim={0 0 500 0}, clip]{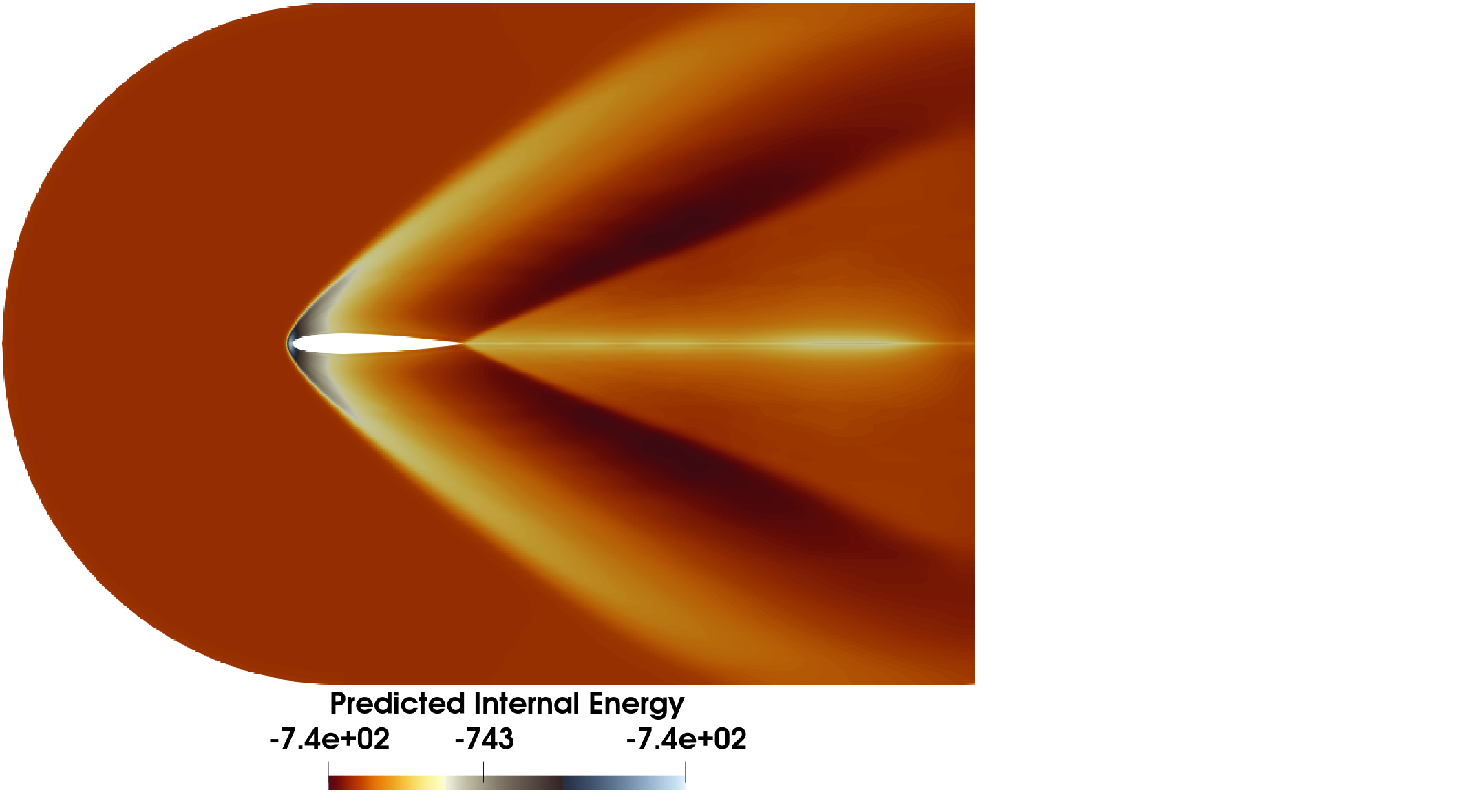}
  \includegraphics[width=0.325\textwidth, trim={0 0 500 0}, clip]{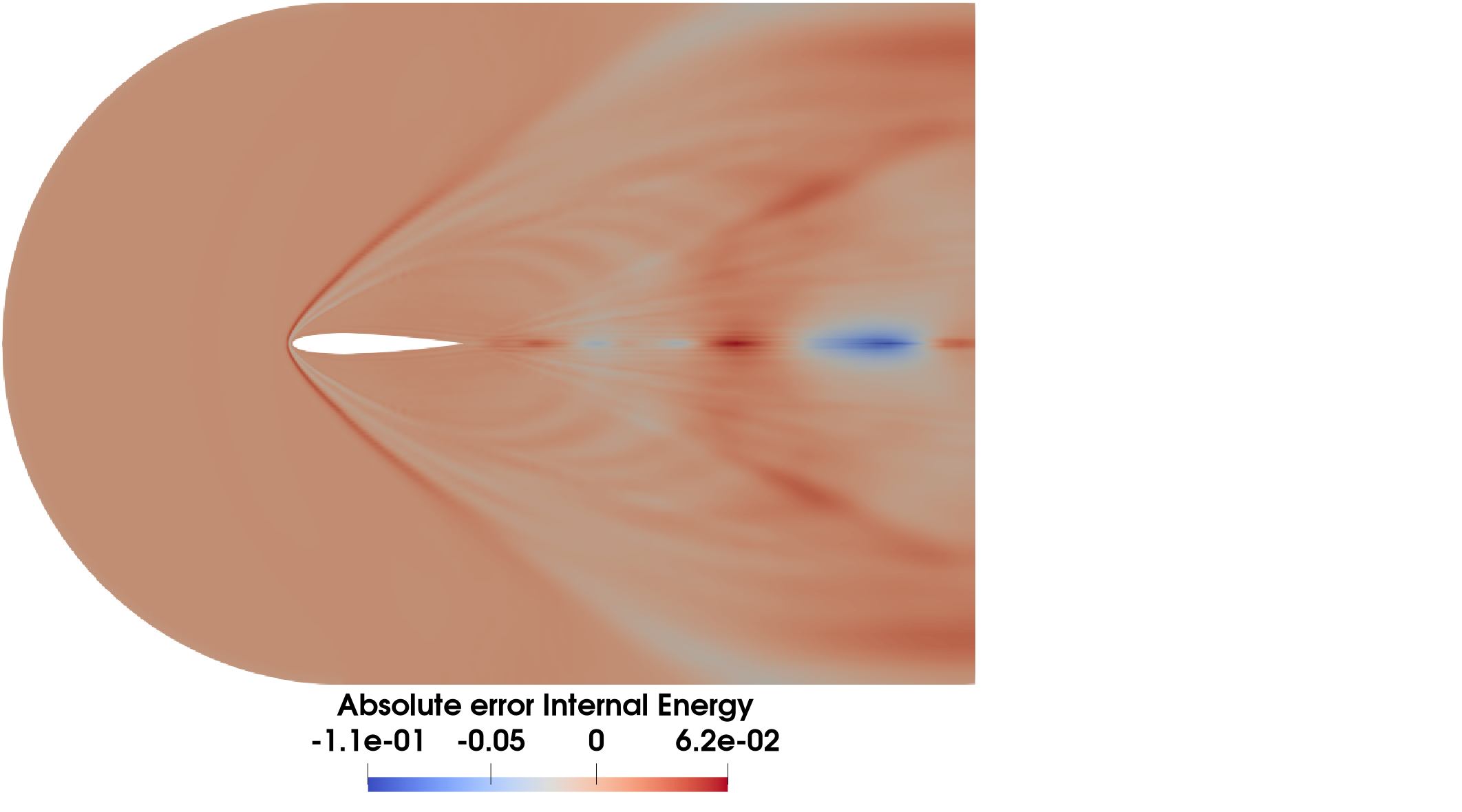}\\
  \includegraphics[width=0.325\textwidth, trim={0 0 500 0}, clip]{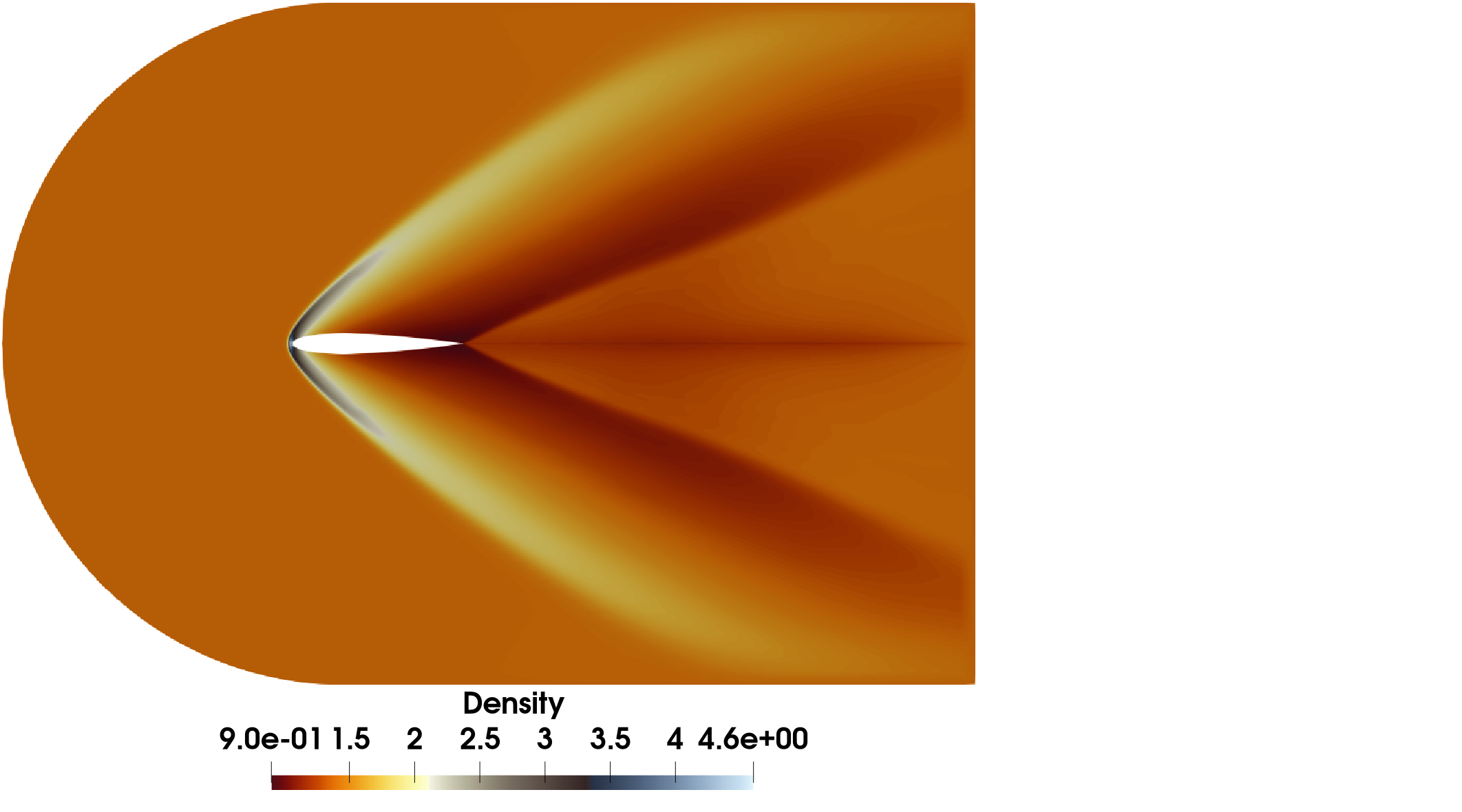}
  \includegraphics[width=0.325\textwidth, trim={0 0 500 0}, clip]{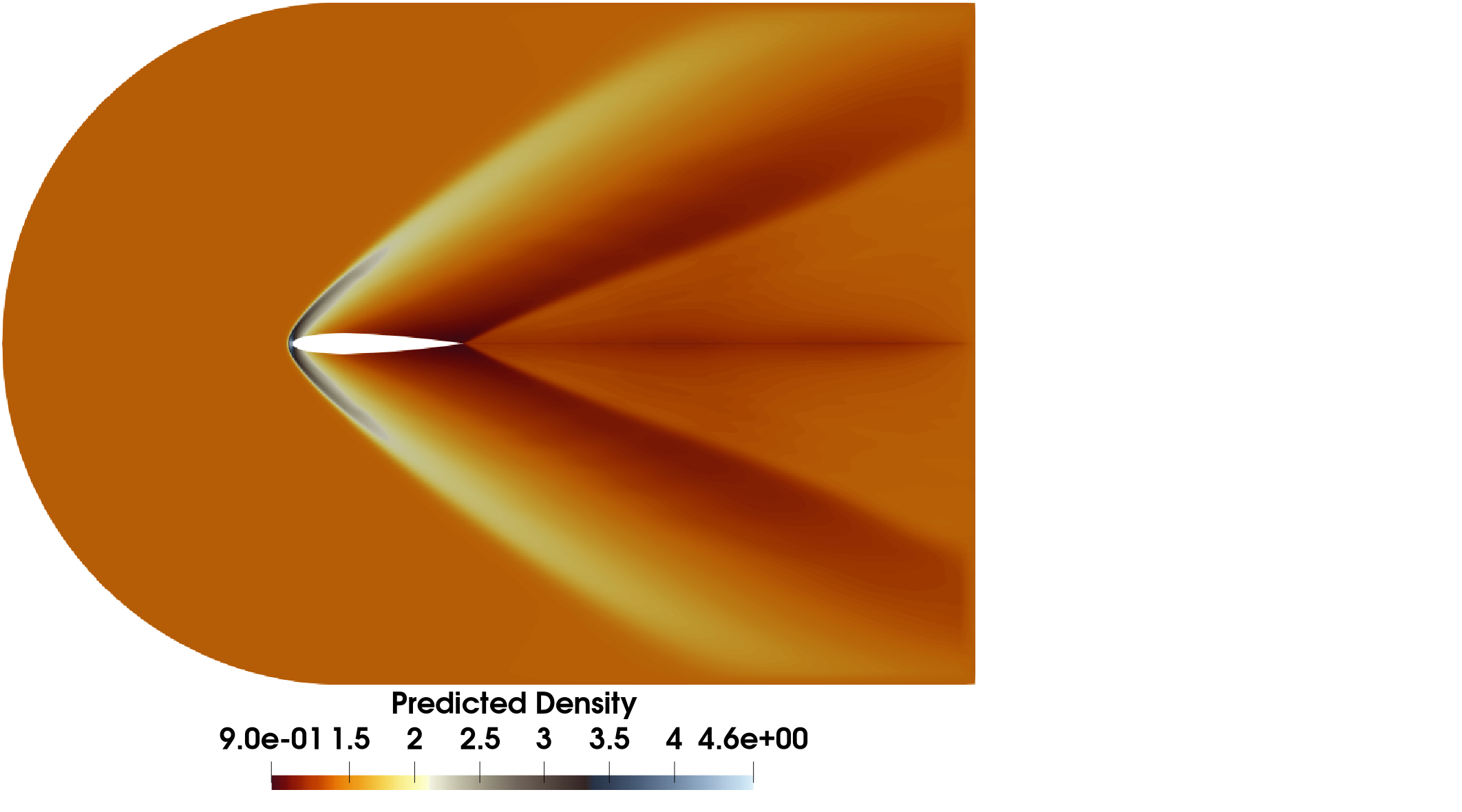}
  \includegraphics[width=0.325\textwidth, trim={0 0 500 0}, clip]{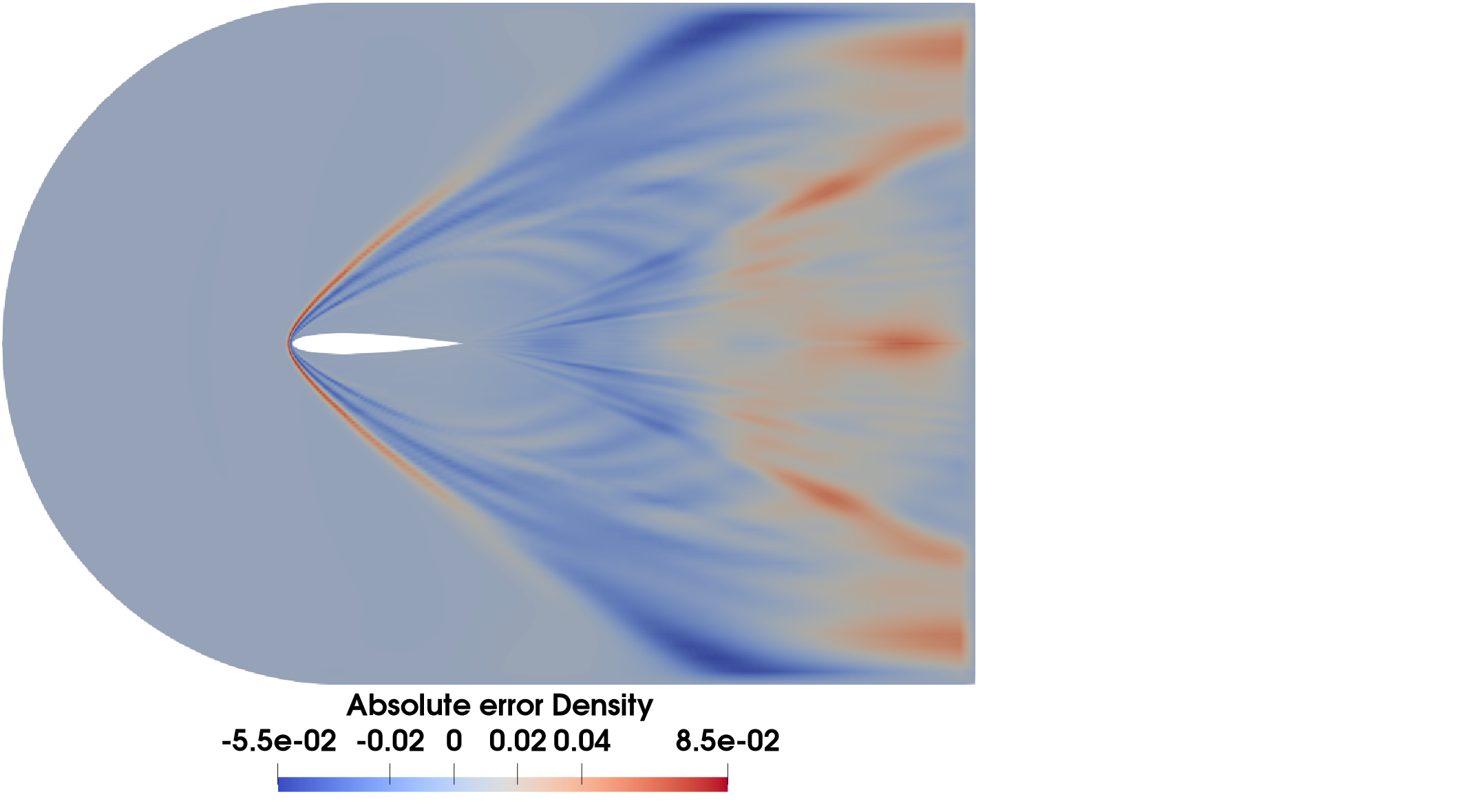}\\
  \includegraphics[width=0.325\textwidth, trim={0 0 500 0}, clip]{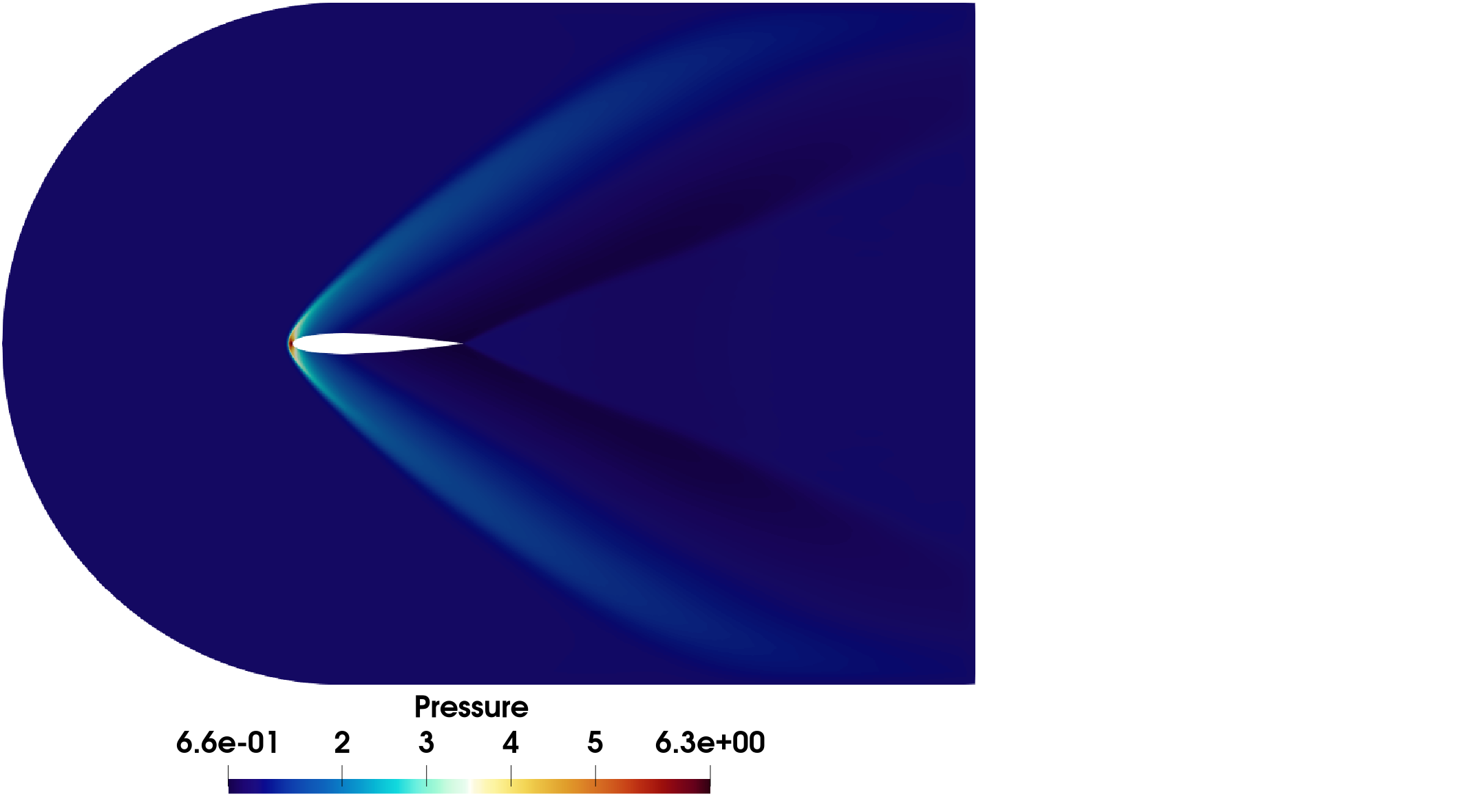}
  \includegraphics[width=0.325\textwidth, trim={0 0 500 0}, clip]{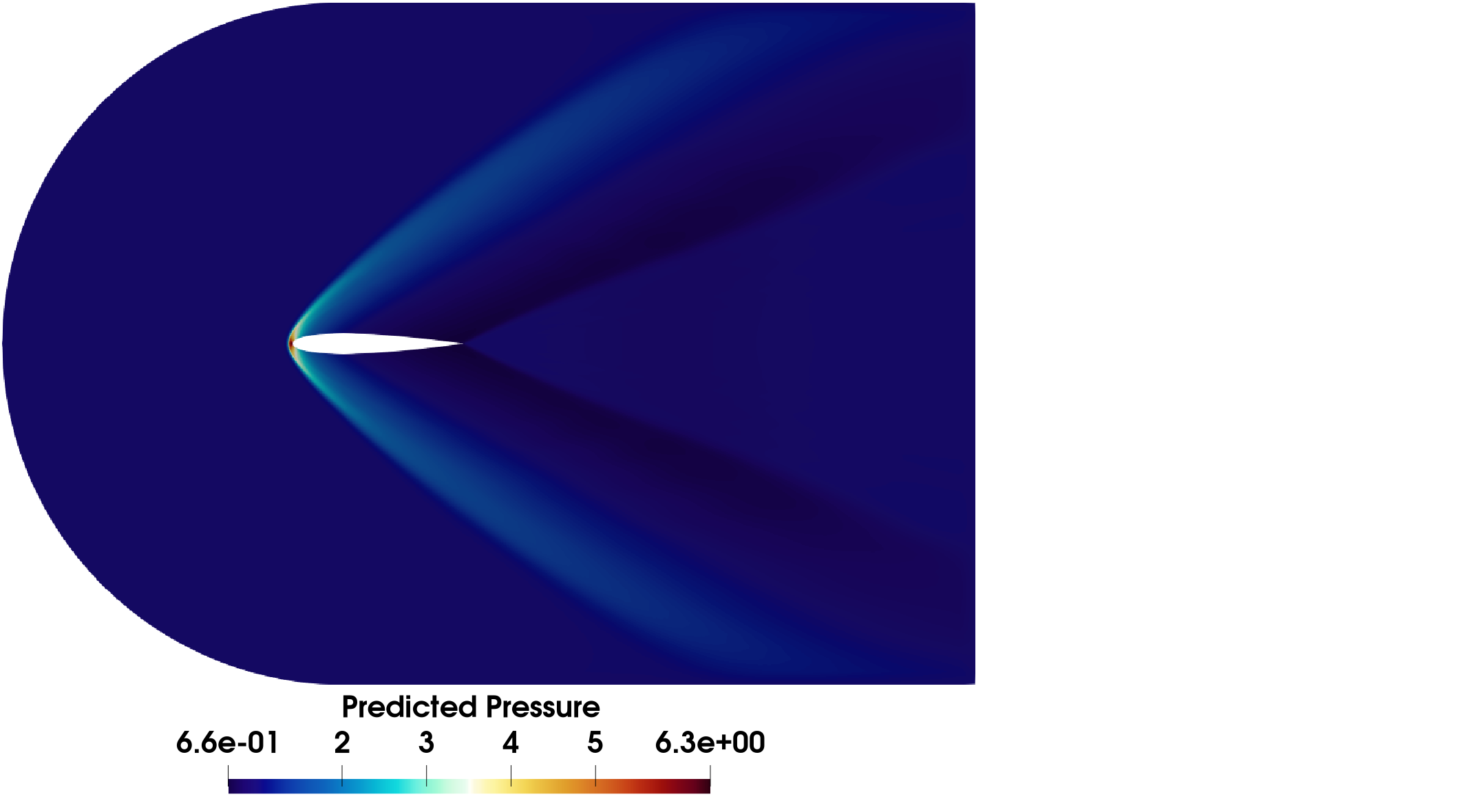}
  \includegraphics[width=0.325\textwidth, trim={0 0 500 0}, clip]{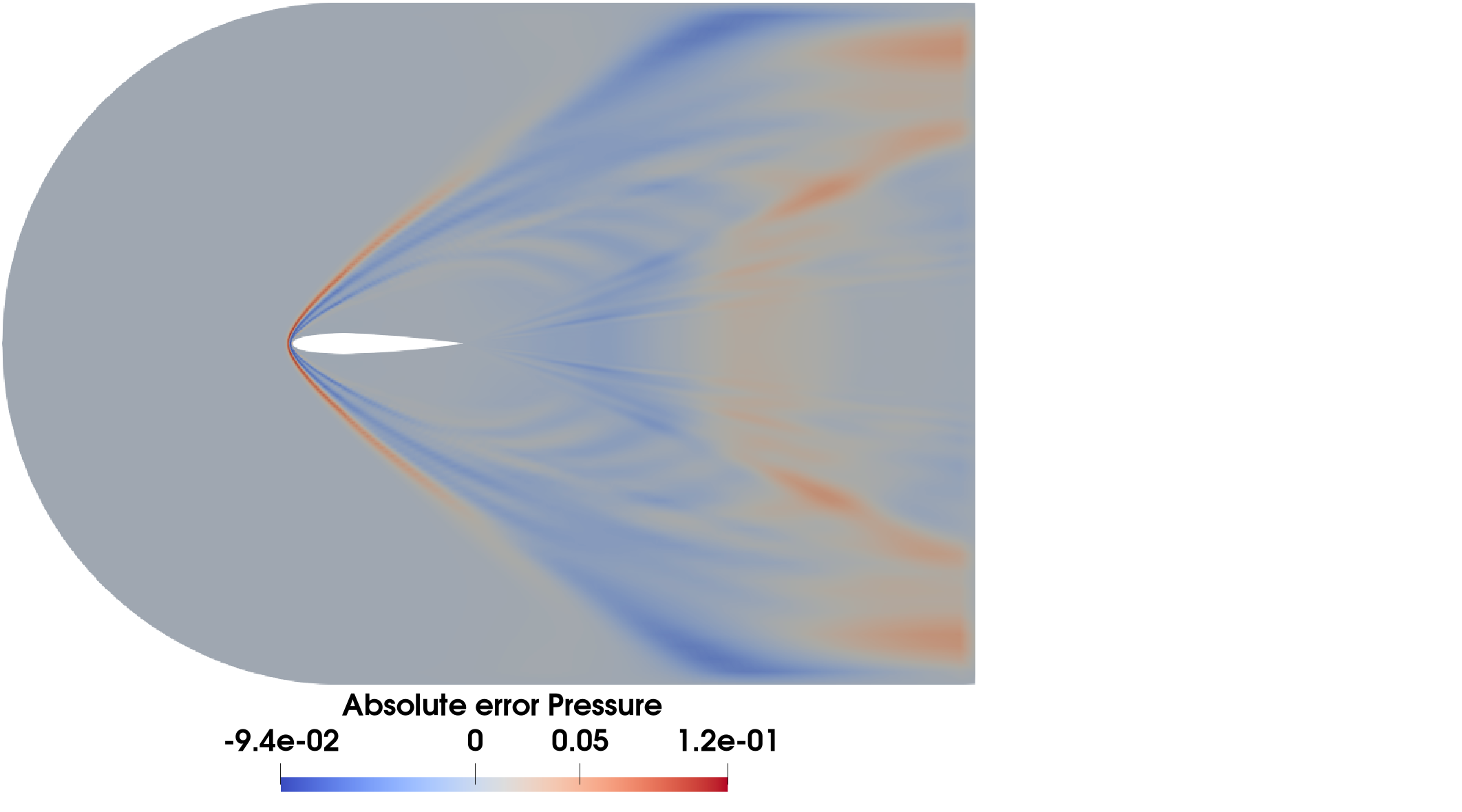}\\
  \caption{\textbf{CNS2.} Predicted velocity, density, internal energy and pressure fields for the test case described in~\cref{subsec:coarseAirfoil} with $\text{Ma}=2.622$ corresponding to test parameter number \textbf{3} at the final time instant $t=2.5\text{s}$. The number of cells is \textbf{32160} the total number of degrees of freedom is \textbf{192960}. The adaptive magic points are shown in Figure~\ref{fig:adaptive}. The number of collocation nodes is $r_h=\textbf{700}$}
  \label{fig:snapsFinerAirfoil}
\end{figure}

\begin{figure}[ht!]
  \centering
  \includegraphics[width=0.325\textwidth, trim={0 300 0 0}, clip]{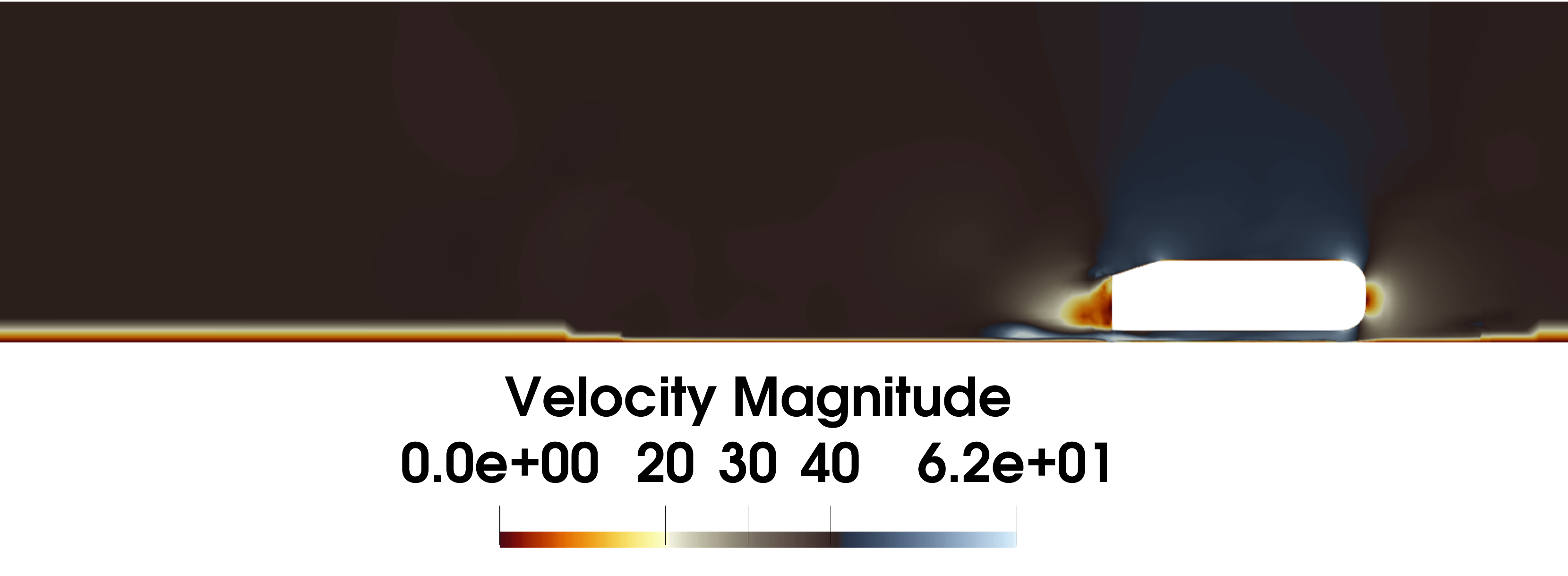}
  \includegraphics[width=0.325\textwidth, trim={0 300 0 0}, clip]{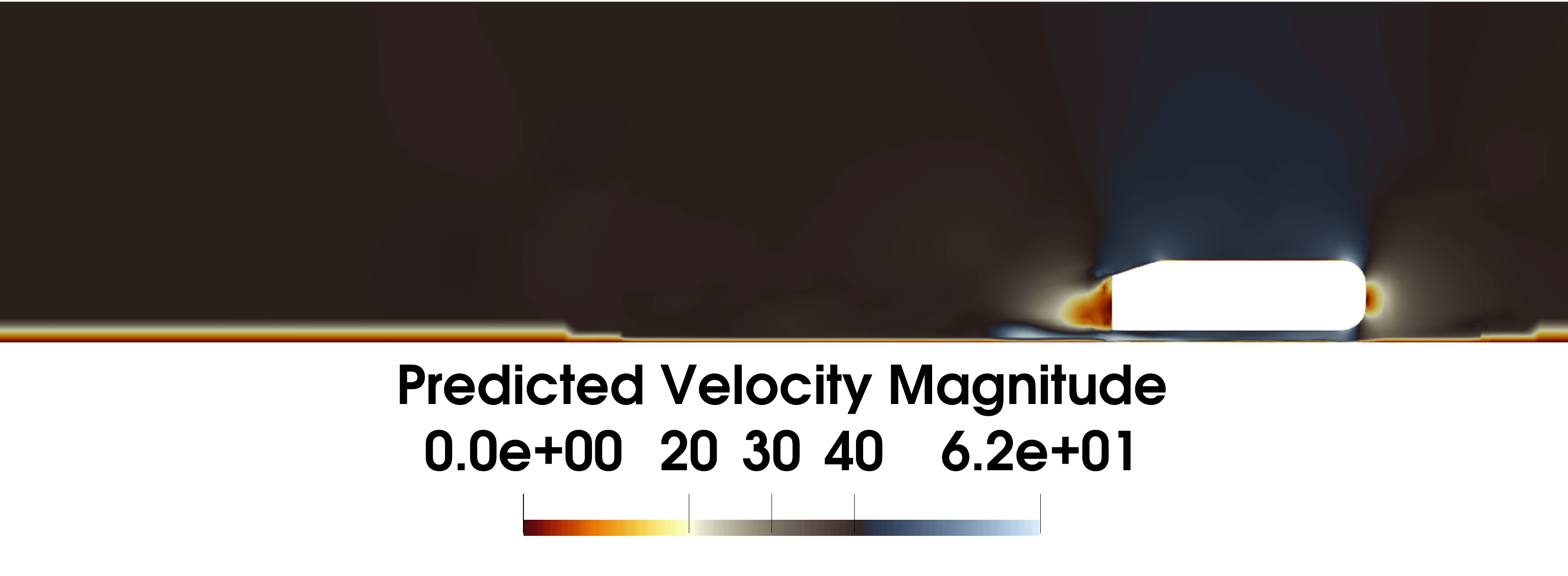}
  \includegraphics[width=0.325\textwidth, trim={0 300 0 0}, clip]{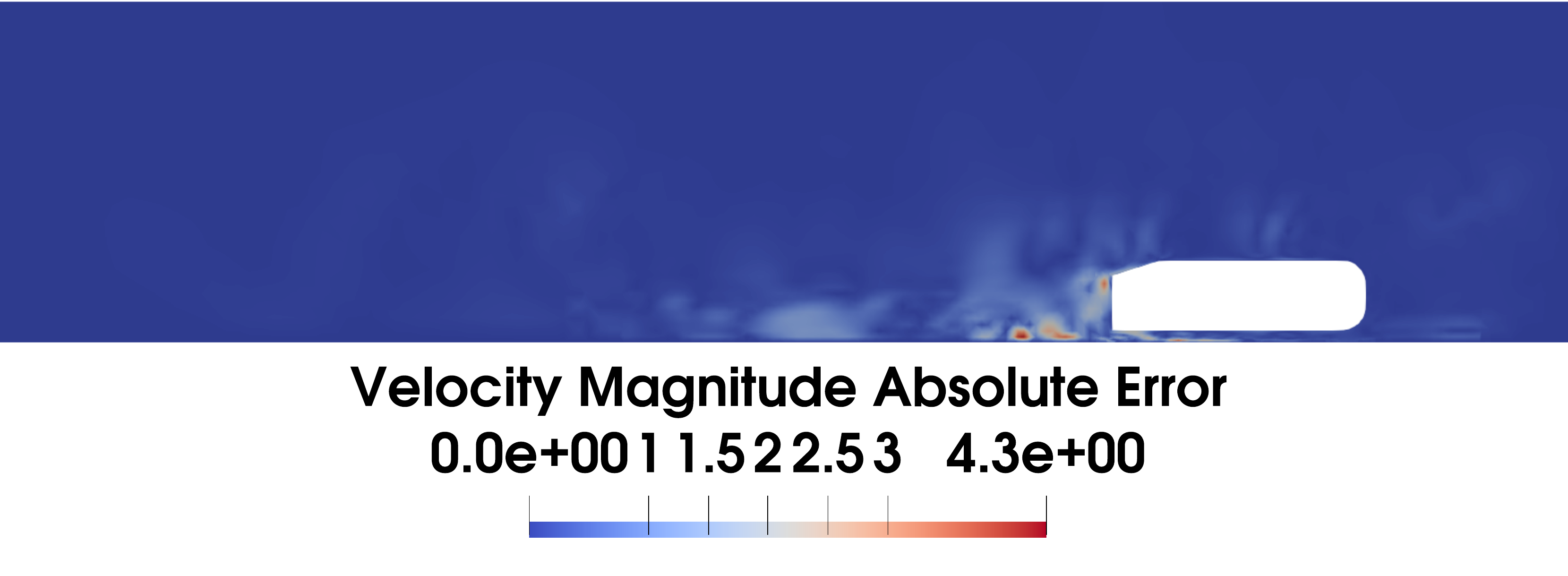}\\
  \includegraphics[width=0.325\textwidth, trim={0 300 0 0}, clip]{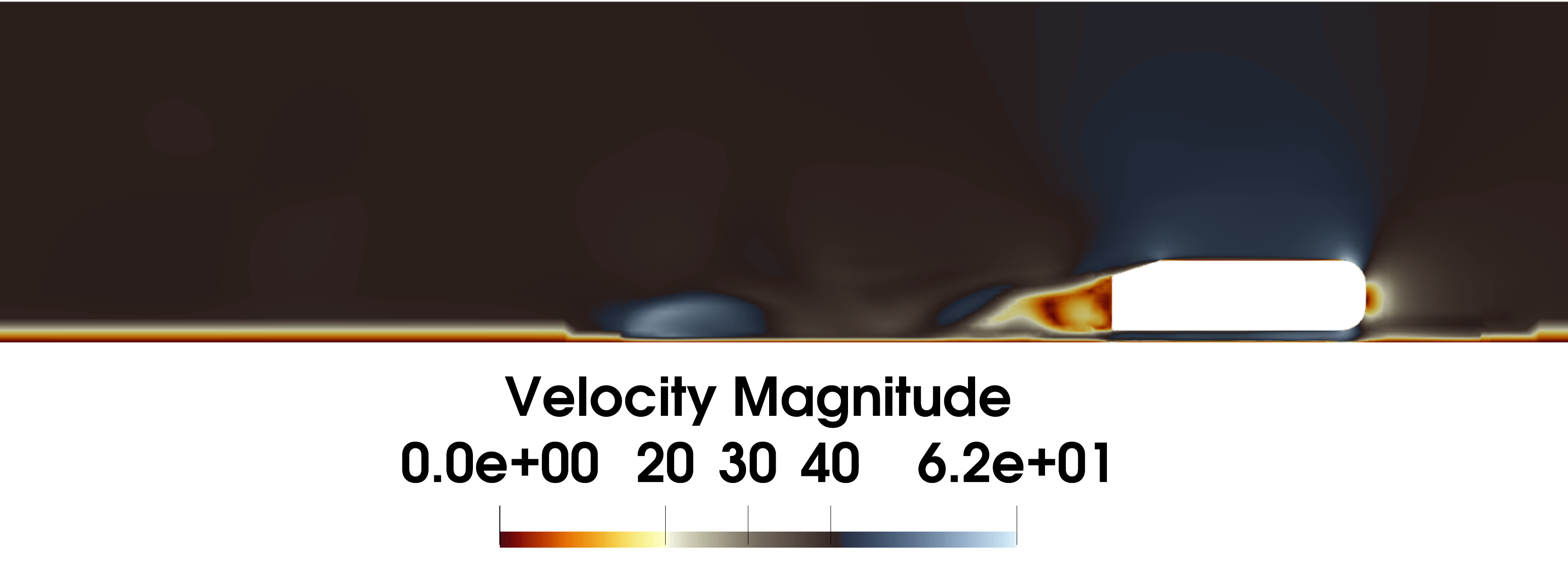}
  \includegraphics[width=0.325\textwidth, trim={0 300 0 0}, clip]{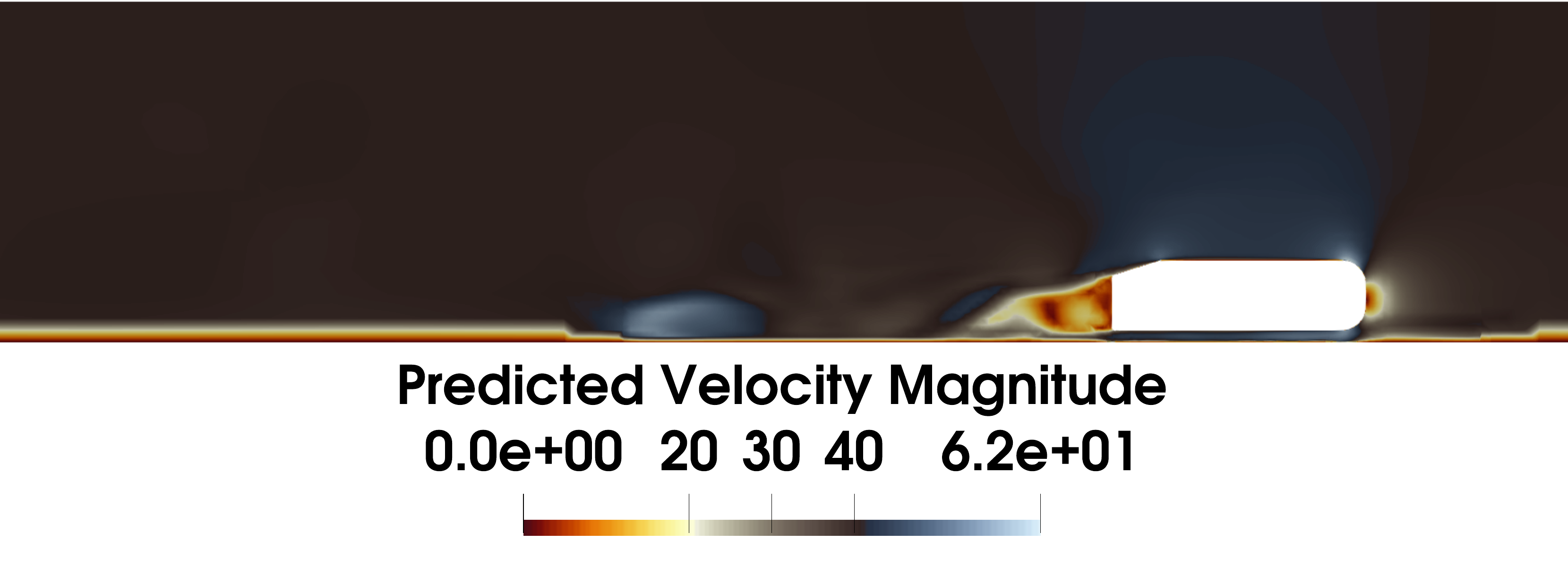}
  \includegraphics[width=0.325\textwidth, trim={0 300 0 0}, clip]{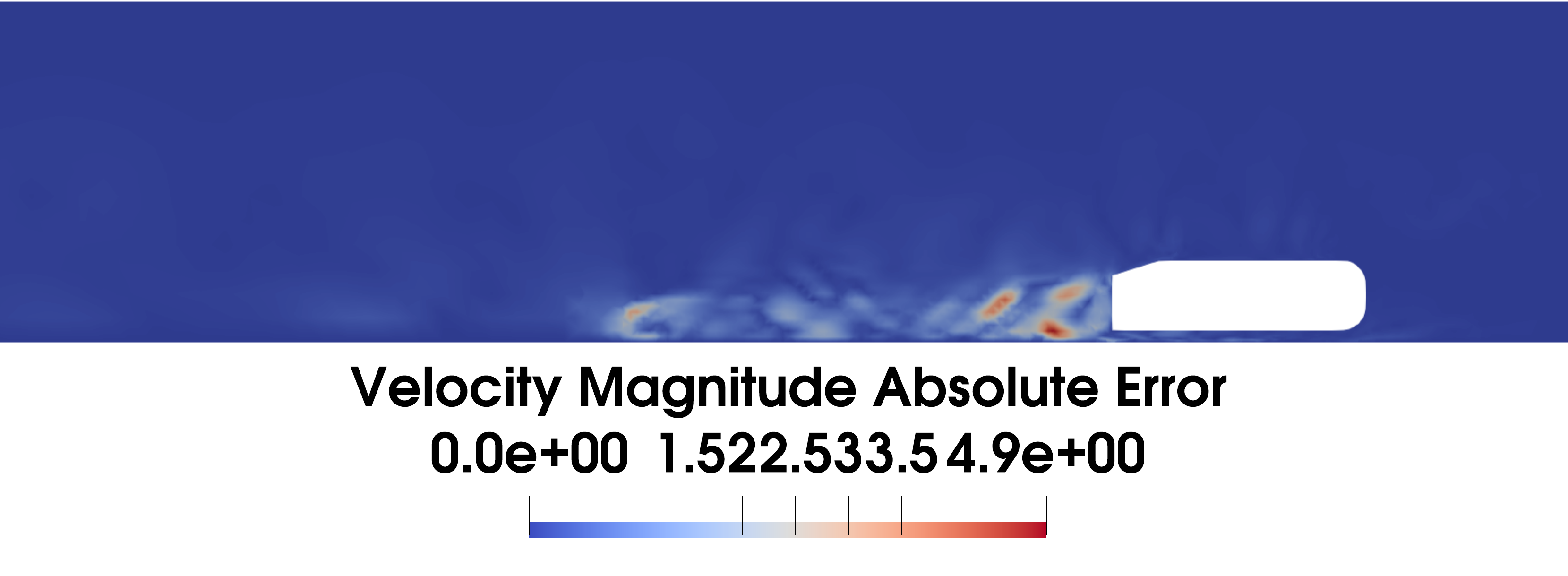}\\
  \includegraphics[width=0.325\textwidth, trim={0 300 0 0}, clip]{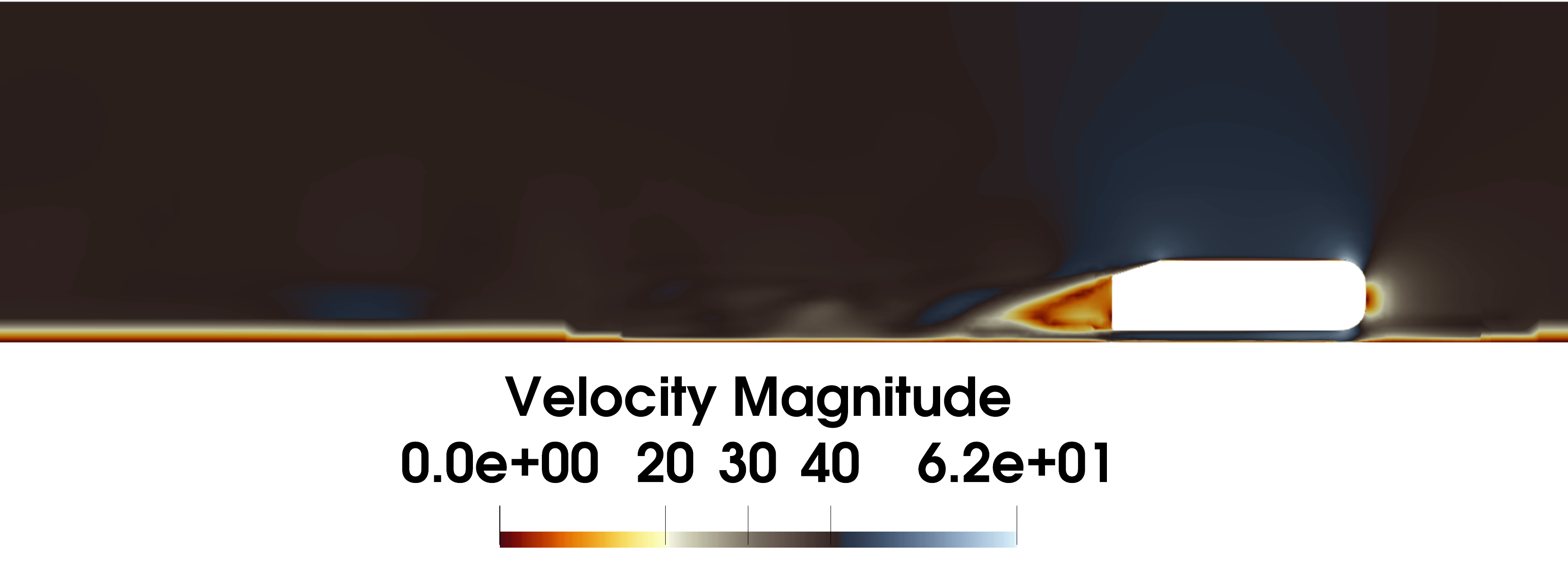}
  \includegraphics[width=0.325\textwidth, trim={0 300 0 0}, clip]{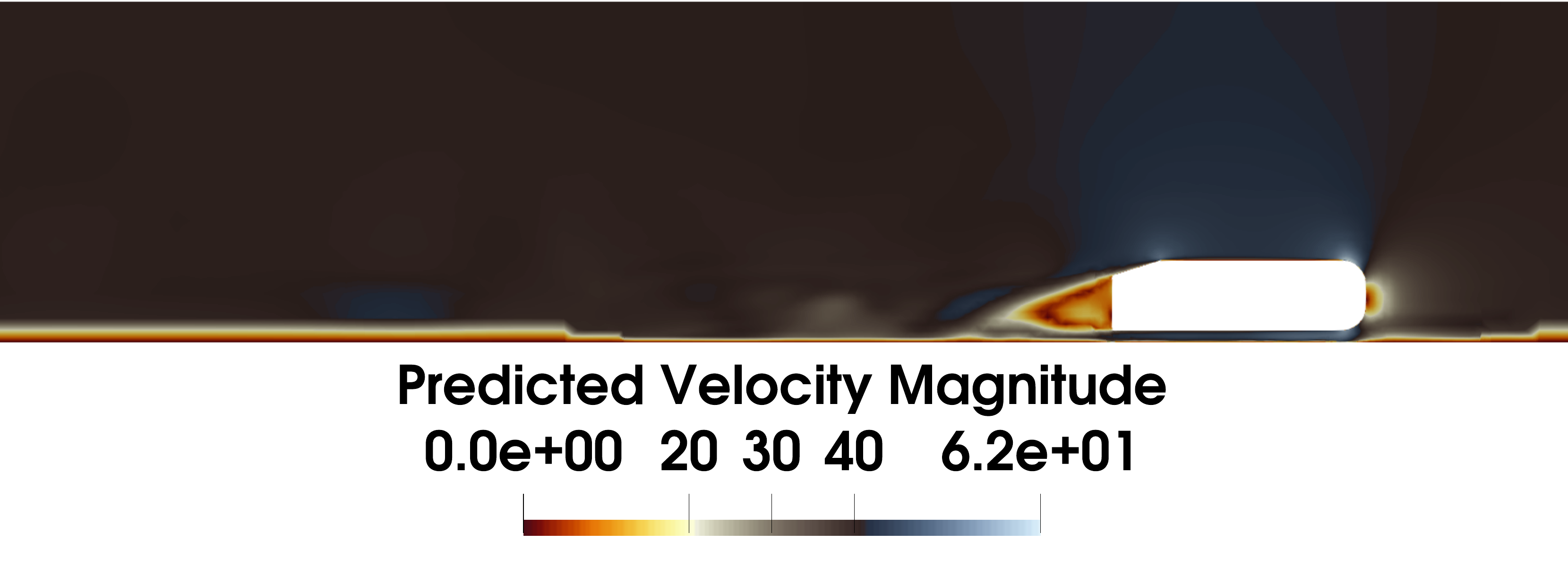}
  \includegraphics[width=0.325\textwidth, trim={0 300 0 0}, clip]{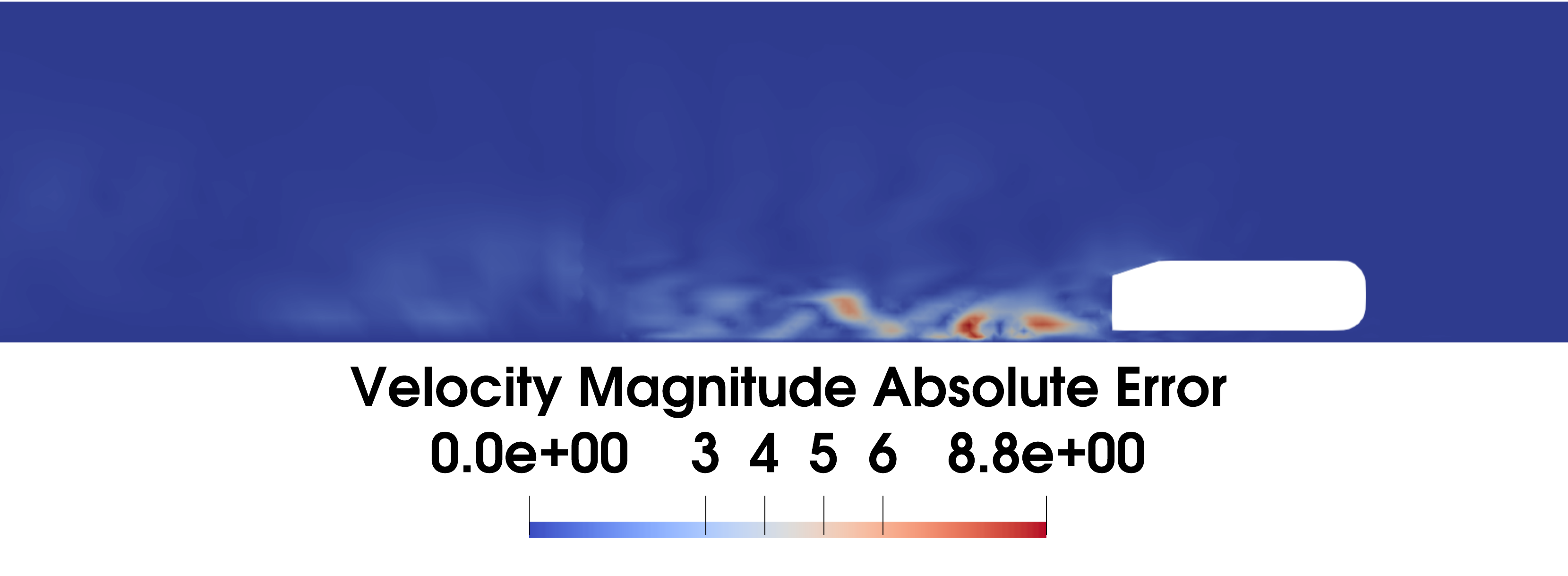}\\
  \includegraphics[width=0.325\textwidth, trim={0 300 0 0}, clip]{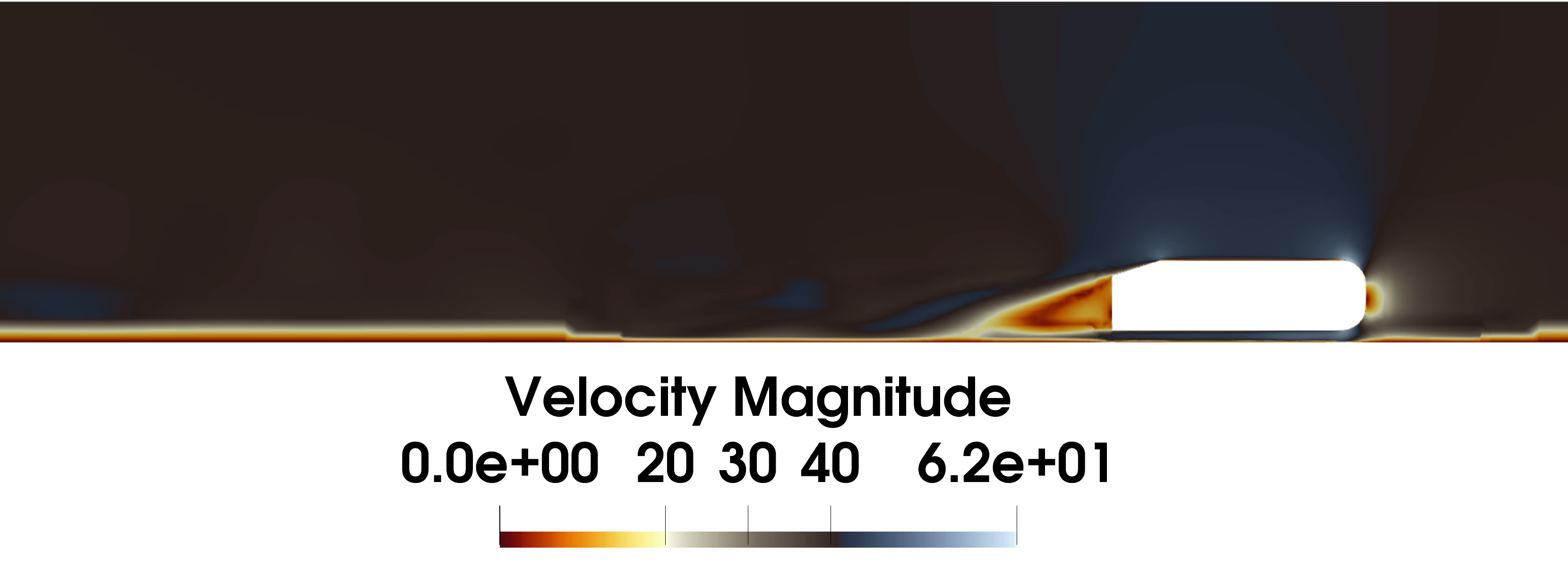}
  \includegraphics[width=0.325\textwidth, trim={0 300 0 0}, clip]{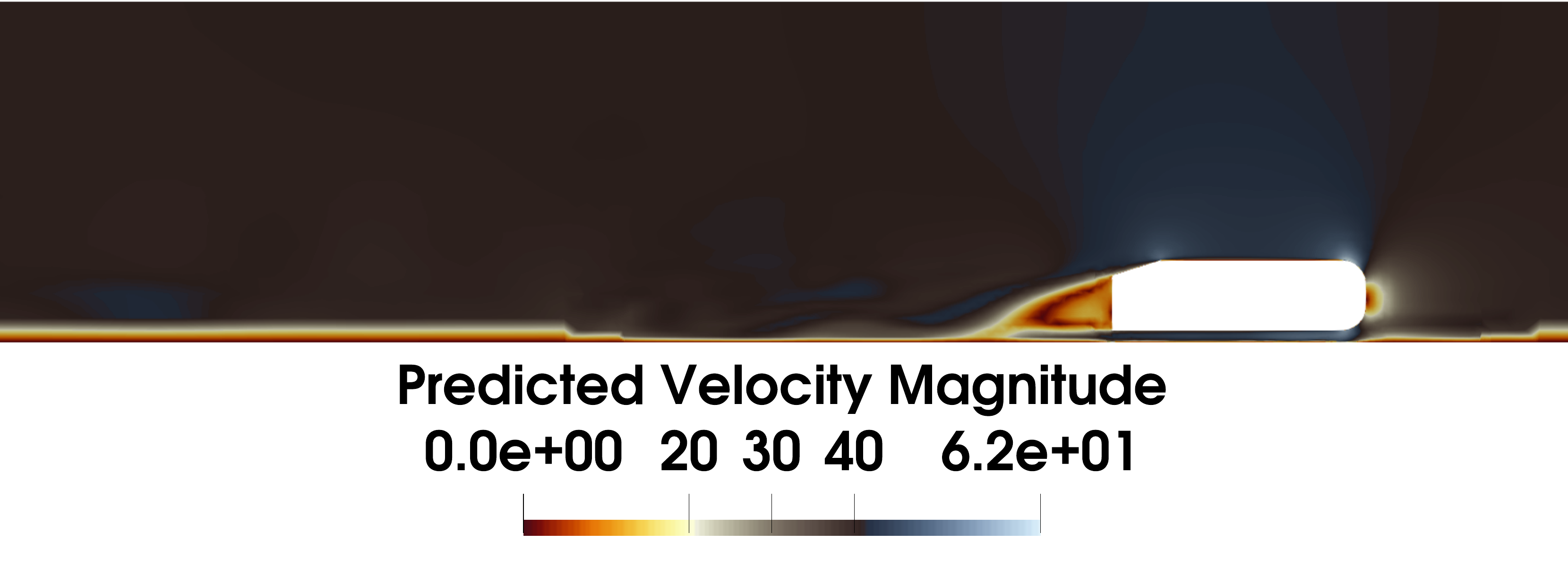}
  \includegraphics[width=0.325\textwidth, trim={0 300 0 0}, clip]{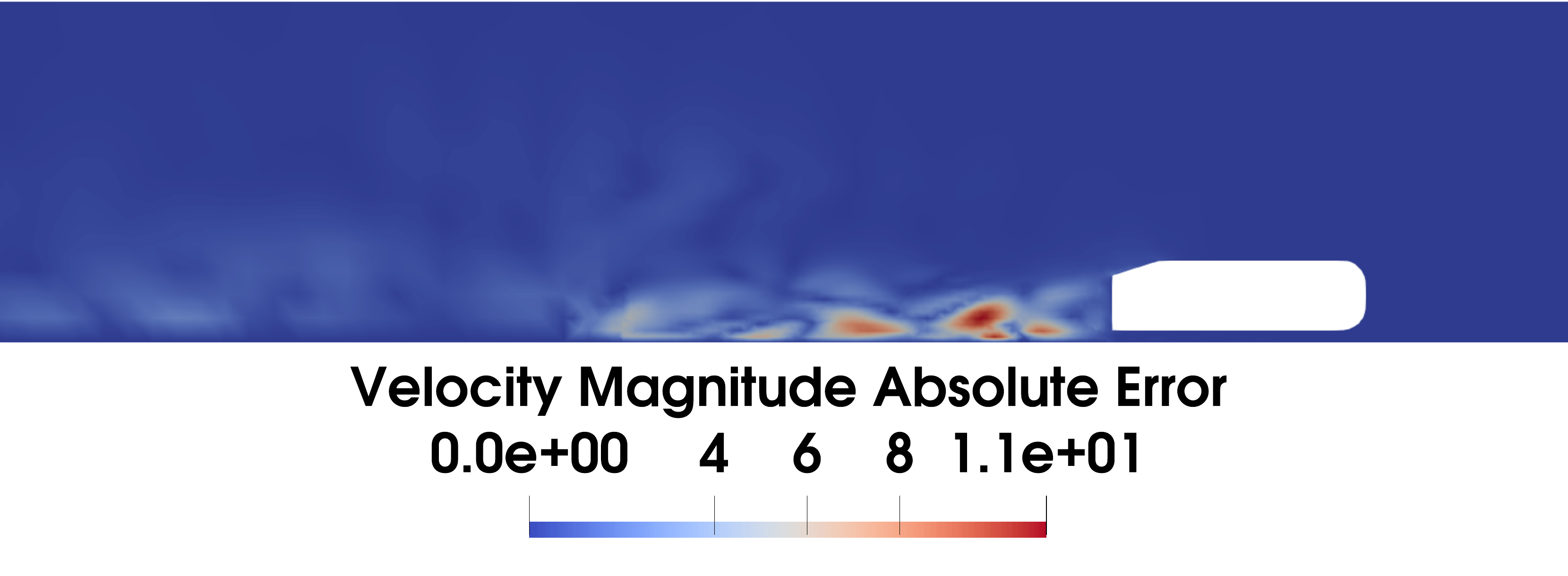}\\
  \includegraphics[width=0.325\textwidth, trim={0 300 0 0}, clip]{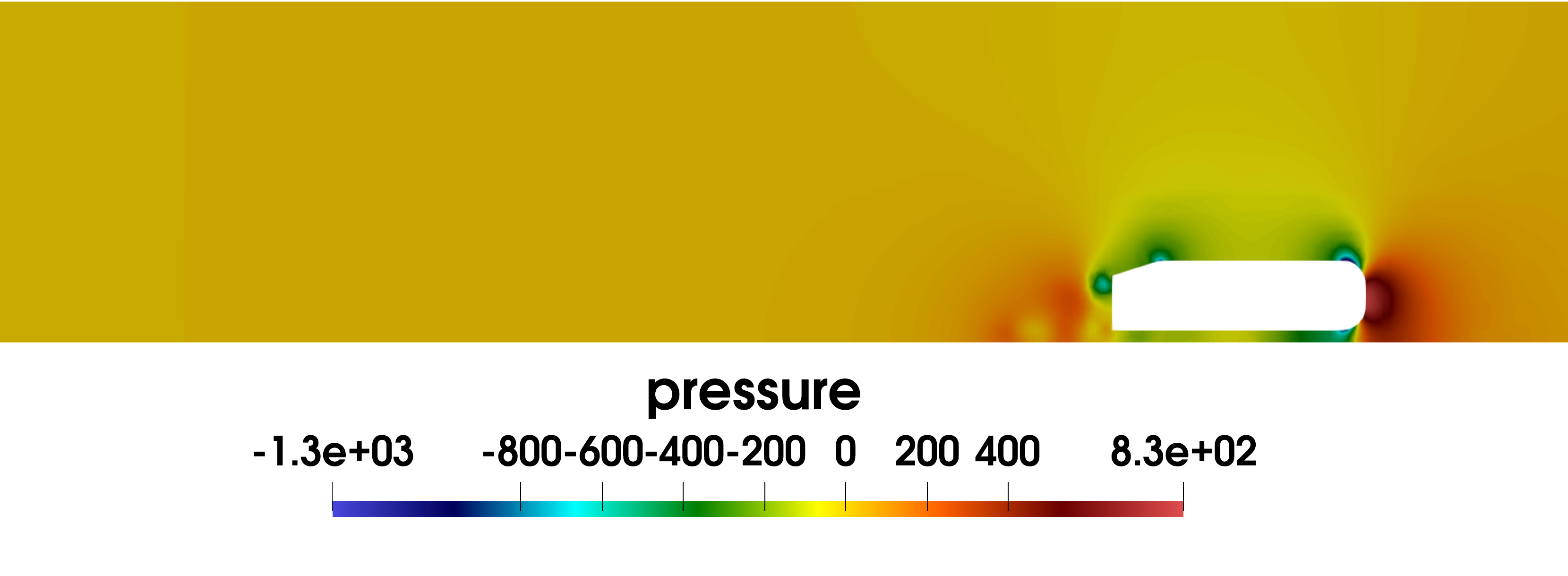}
  \includegraphics[width=0.325\textwidth, trim={0 300 0 0}, clip]{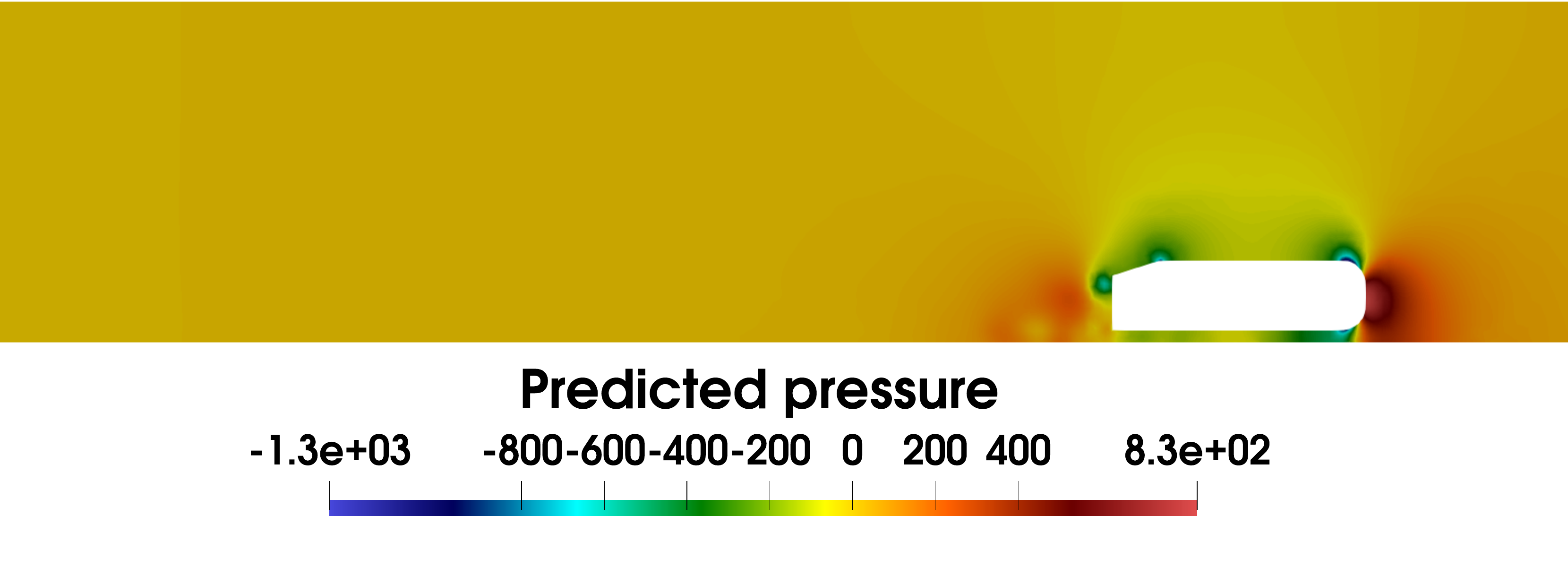}
  \includegraphics[width=0.325\textwidth, trim={0 300 0 0}, clip]{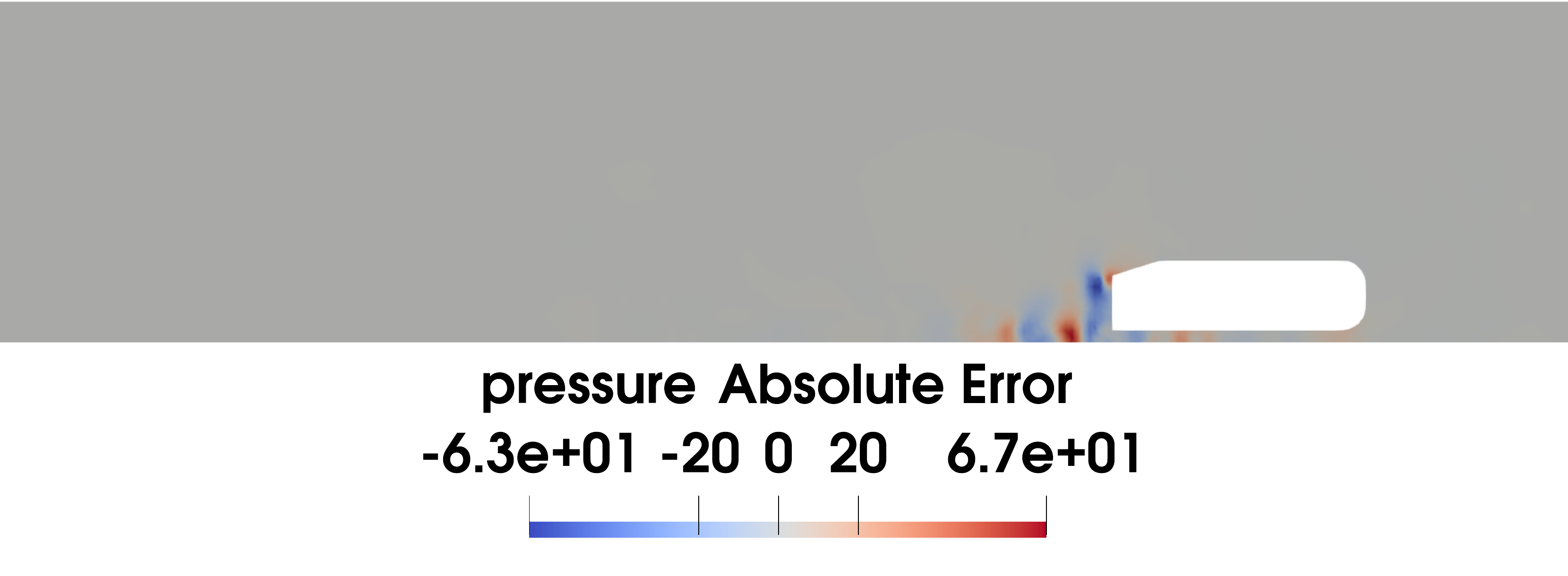}\\
  \includegraphics[width=0.325\textwidth, trim={0 300 0 0}, clip]{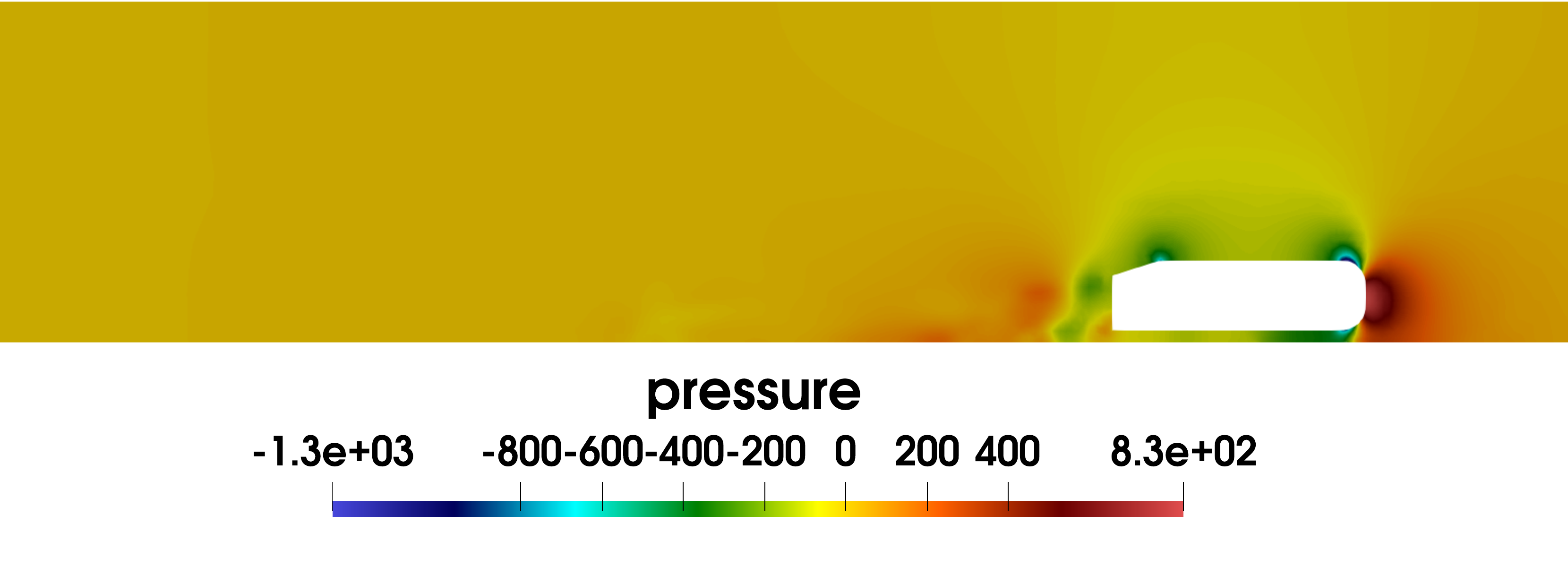}
  \includegraphics[width=0.325\textwidth, trim={0 300 0 0}, clip]{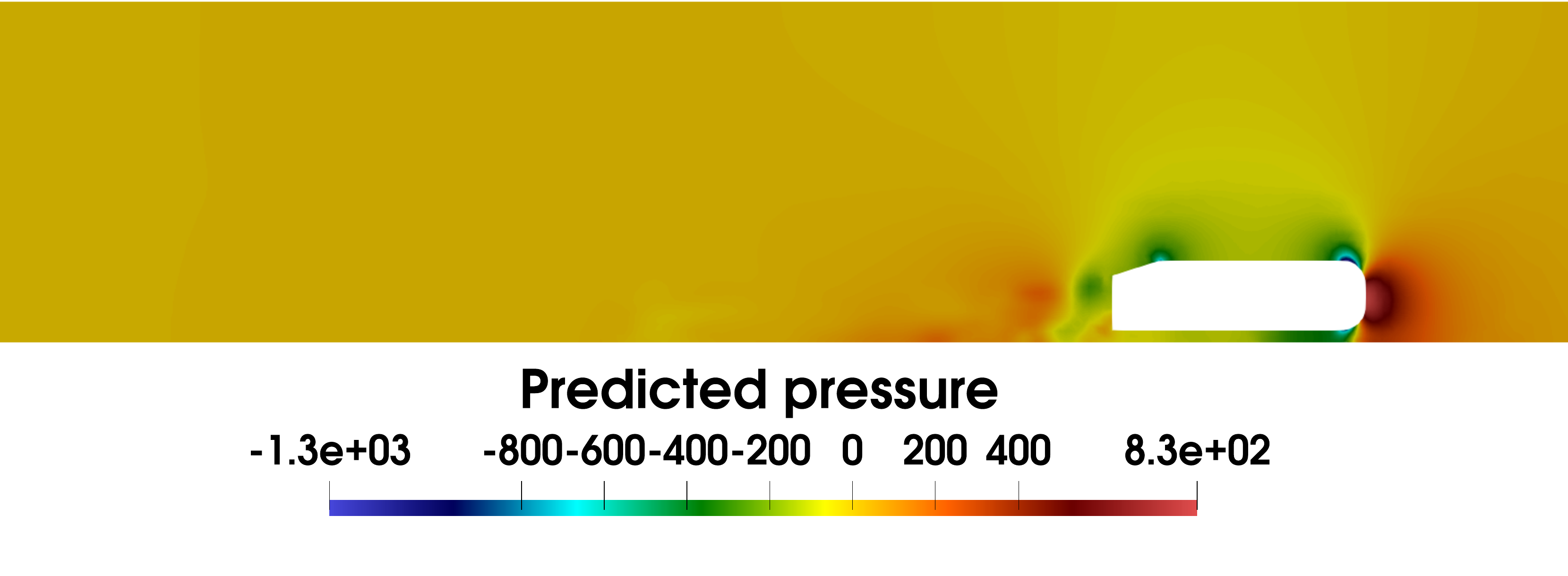}
  \includegraphics[width=0.325\textwidth, trim={0 300 0 0}, clip]{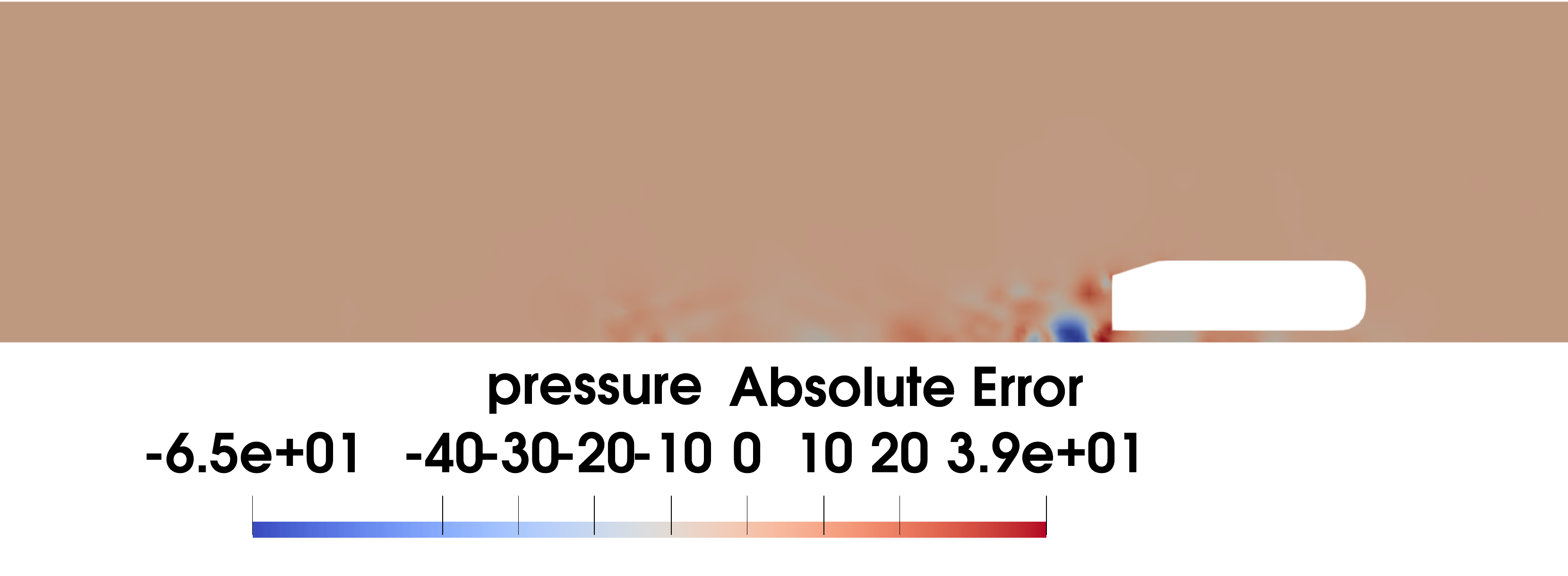}\\
  \includegraphics[width=0.325\textwidth, trim={0 300 0 0}, clip]{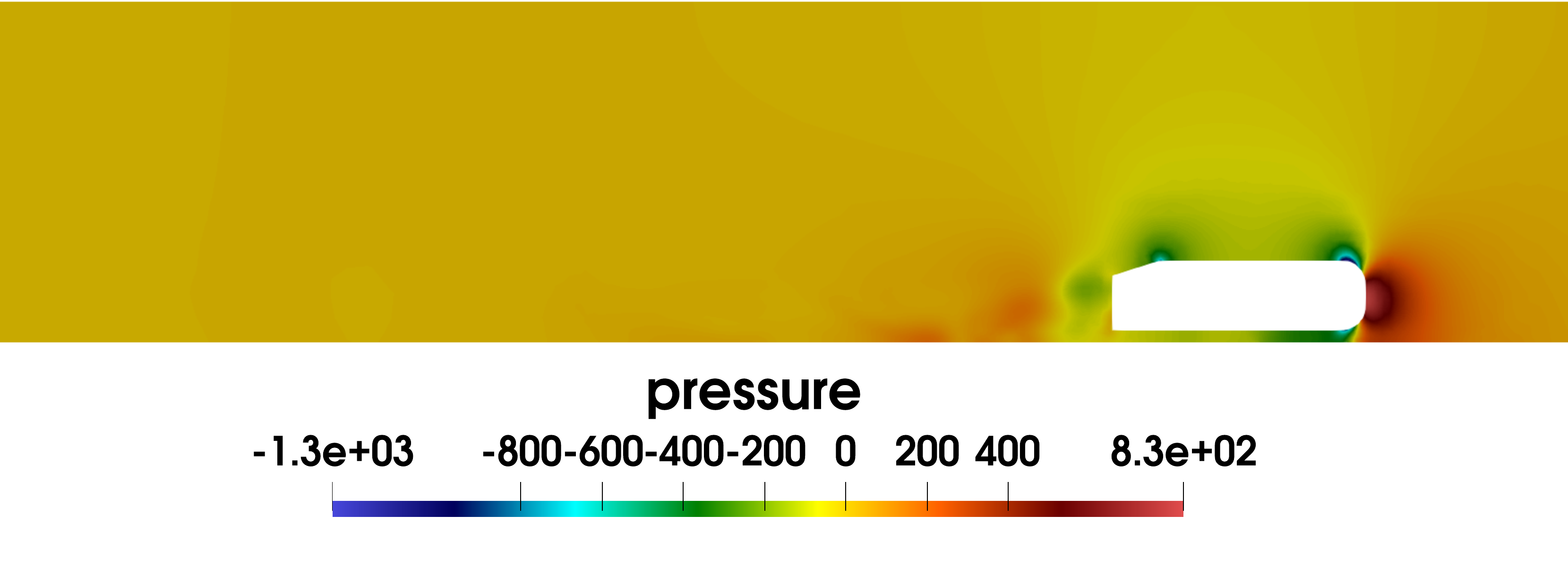}
  \includegraphics[width=0.325\textwidth, trim={0 300 0 0}, clip]{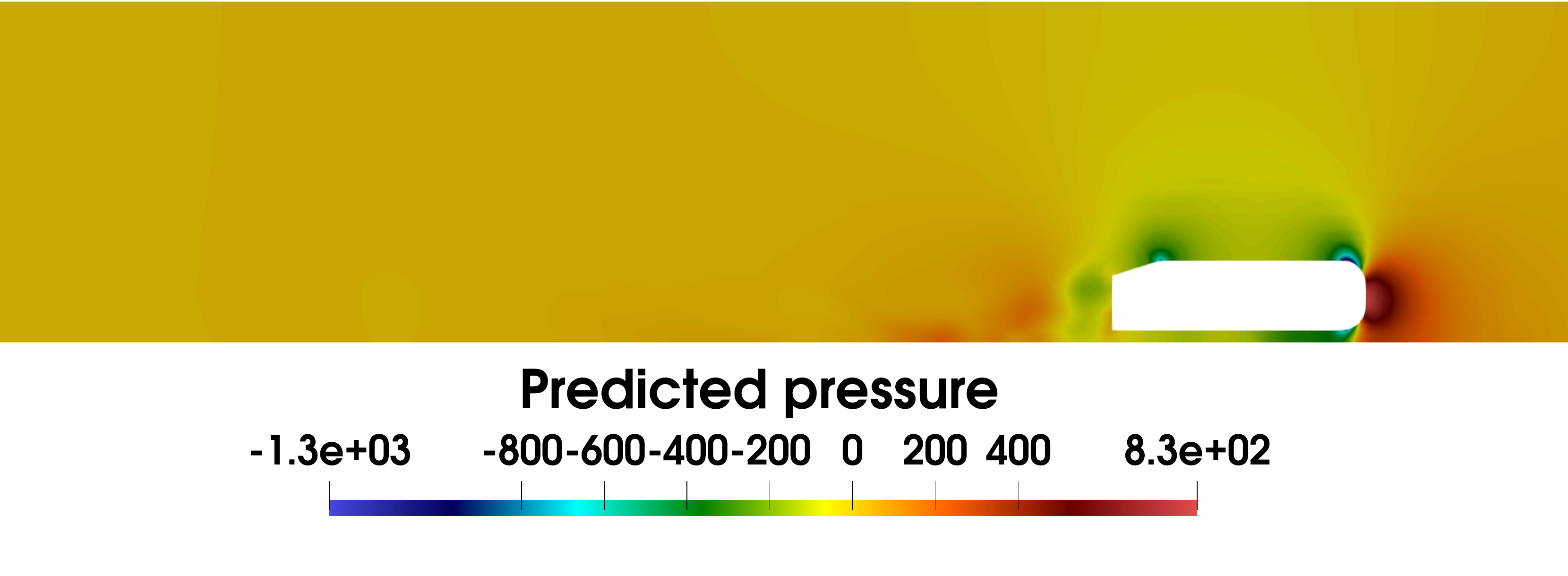}
  \includegraphics[width=0.325\textwidth, trim={0 300 0 0}, clip]{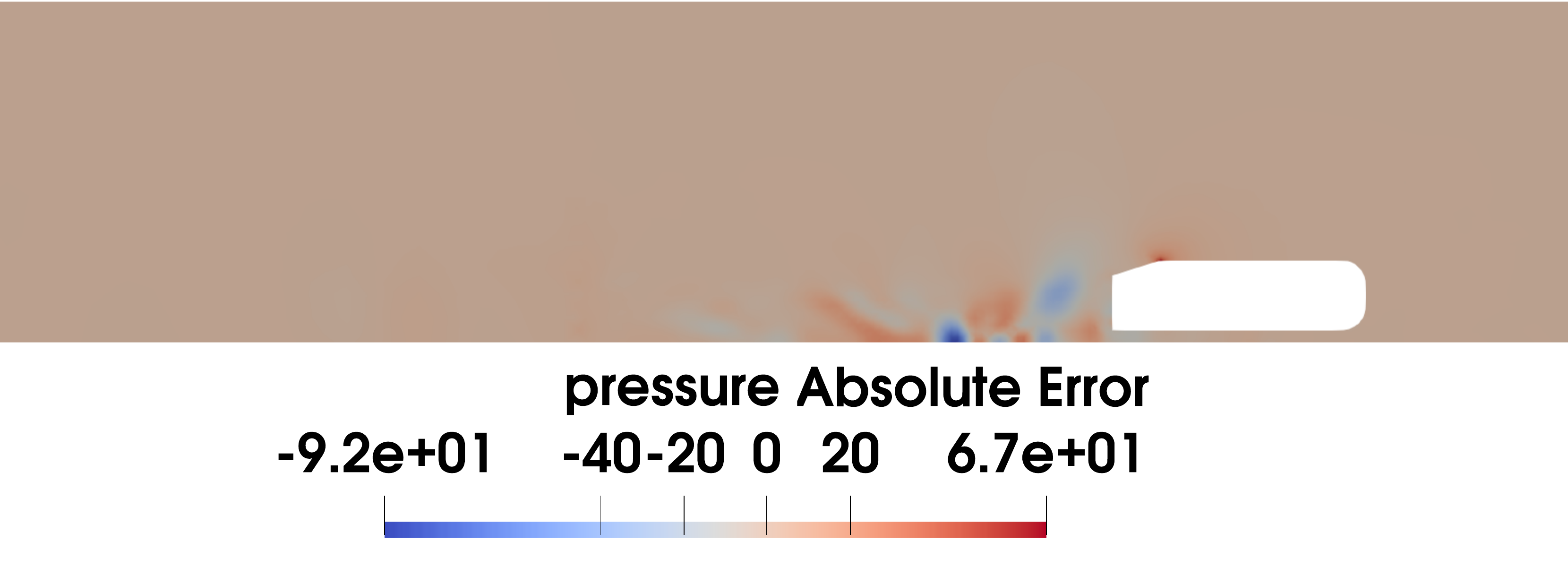}\\
  \includegraphics[width=0.325\textwidth, trim={0 300 0 0}, clip]{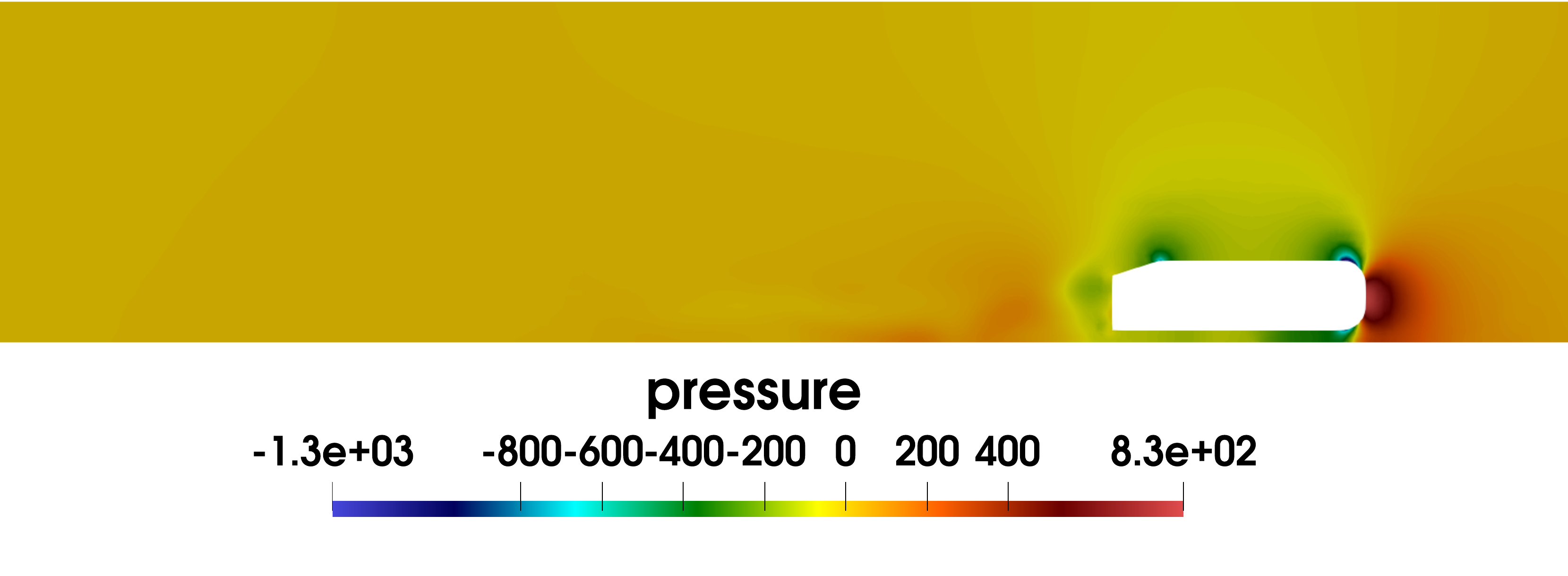}
  \includegraphics[width=0.325\textwidth, trim={0 300 0 0}, clip]{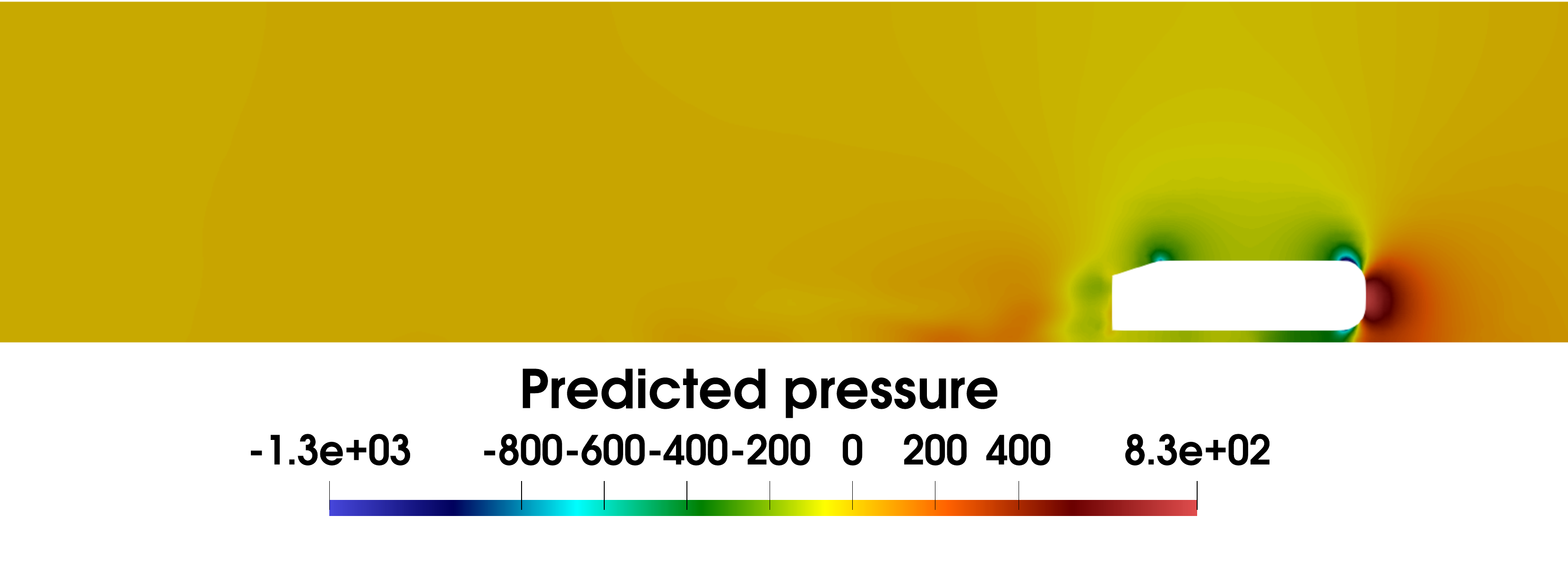}
  \includegraphics[width=0.325\textwidth, trim={0 300 0 0}, clip]{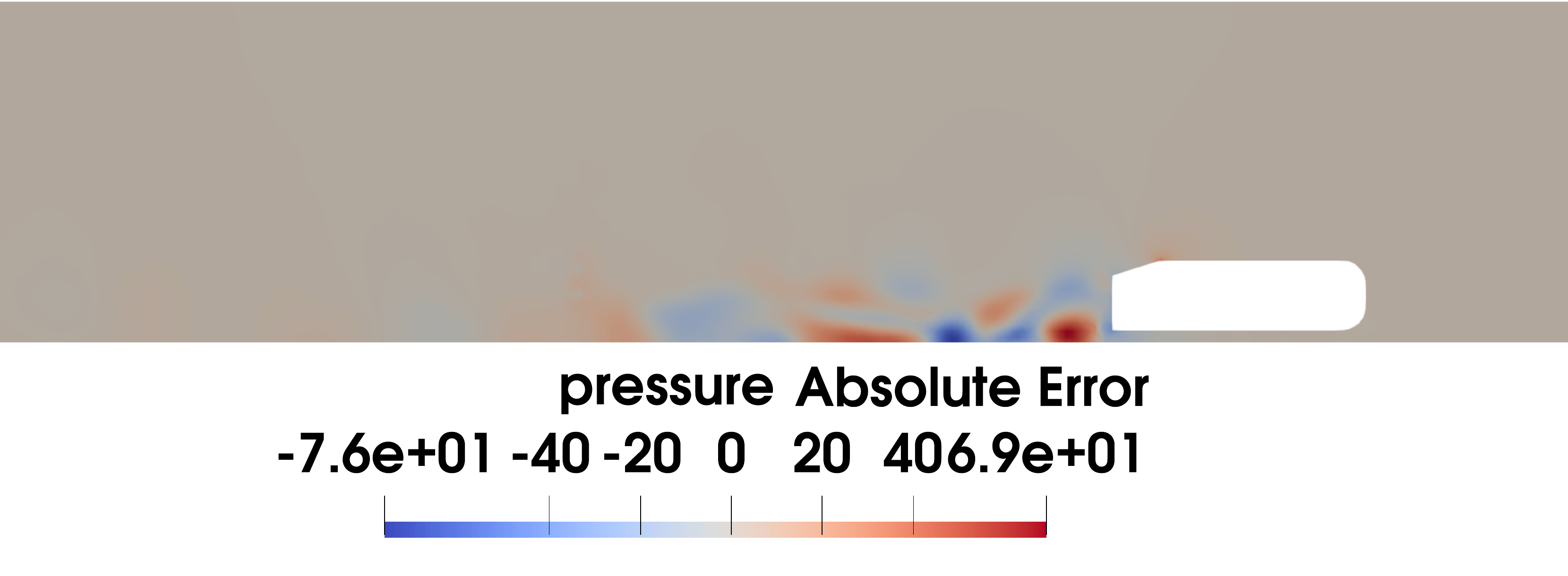}\\
  \caption{\textbf{INS1.} Predicted velocity and pressure fields for the test case described in~\cref{subsubec:smallAhmed} with corresponding to test parameter number \textbf{1} at the time instant s $t\in\{0.025\text{s},0.05\text{s}, 0.075\text{s}, 0.1\text{s}\}$. The number of cells is \textbf{32160} the total number of degrees of freedom is \textbf{192960}. The number of collocation nodes is $r_h=\textbf{2000}$. The adaptive magic points are shown in Figure~\ref{fig:adaptiveSmall}.}
  \label{fig:ahmedSmallSnaps}
\end{figure}

\begin{figure}[ht!]
  \centering
  \includegraphics[width=0.49\textwidth, trim={0 200 0 200}, clip]{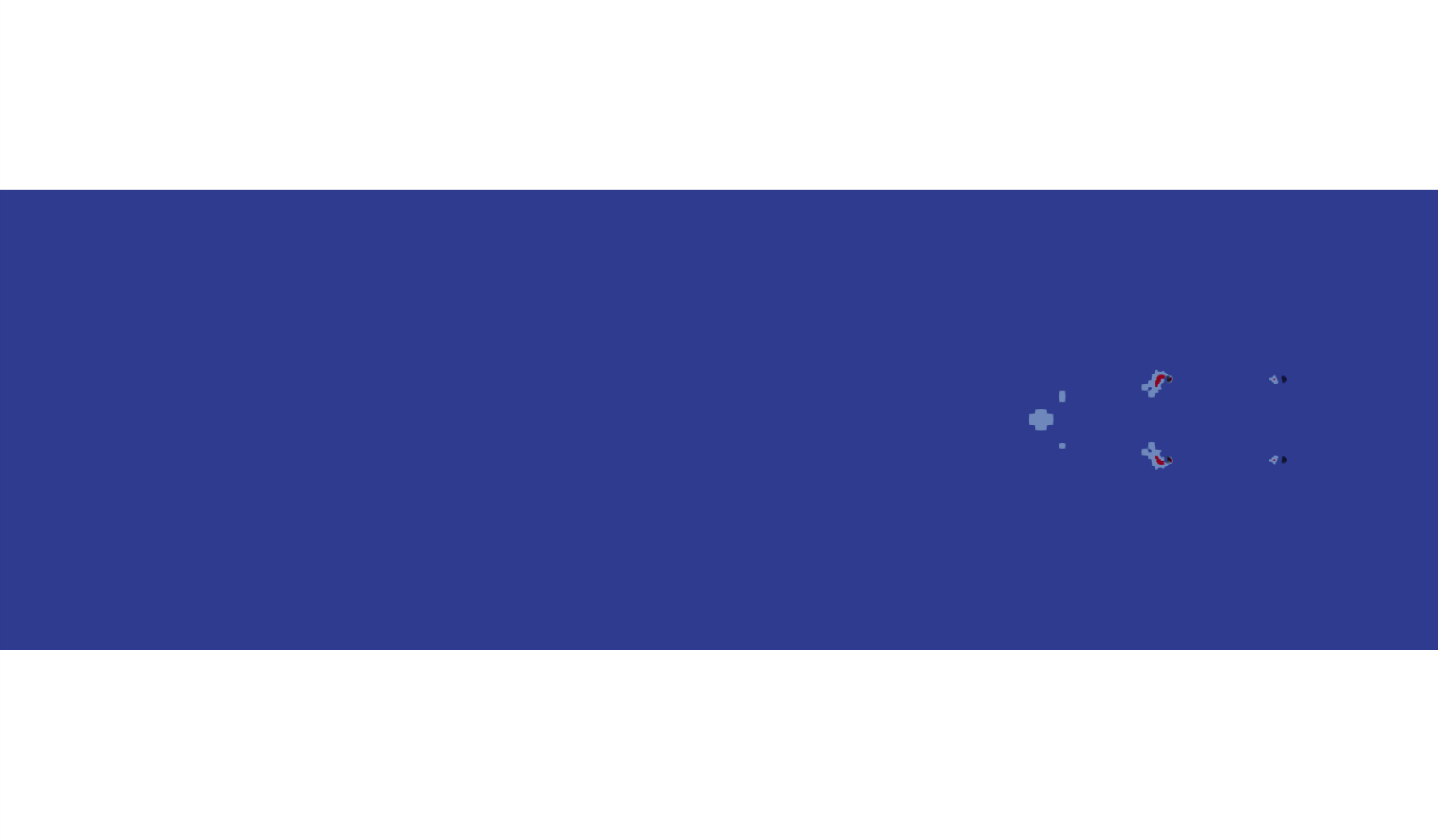}
  \includegraphics[width=0.49\textwidth, trim={0 200 0 200}, clip]{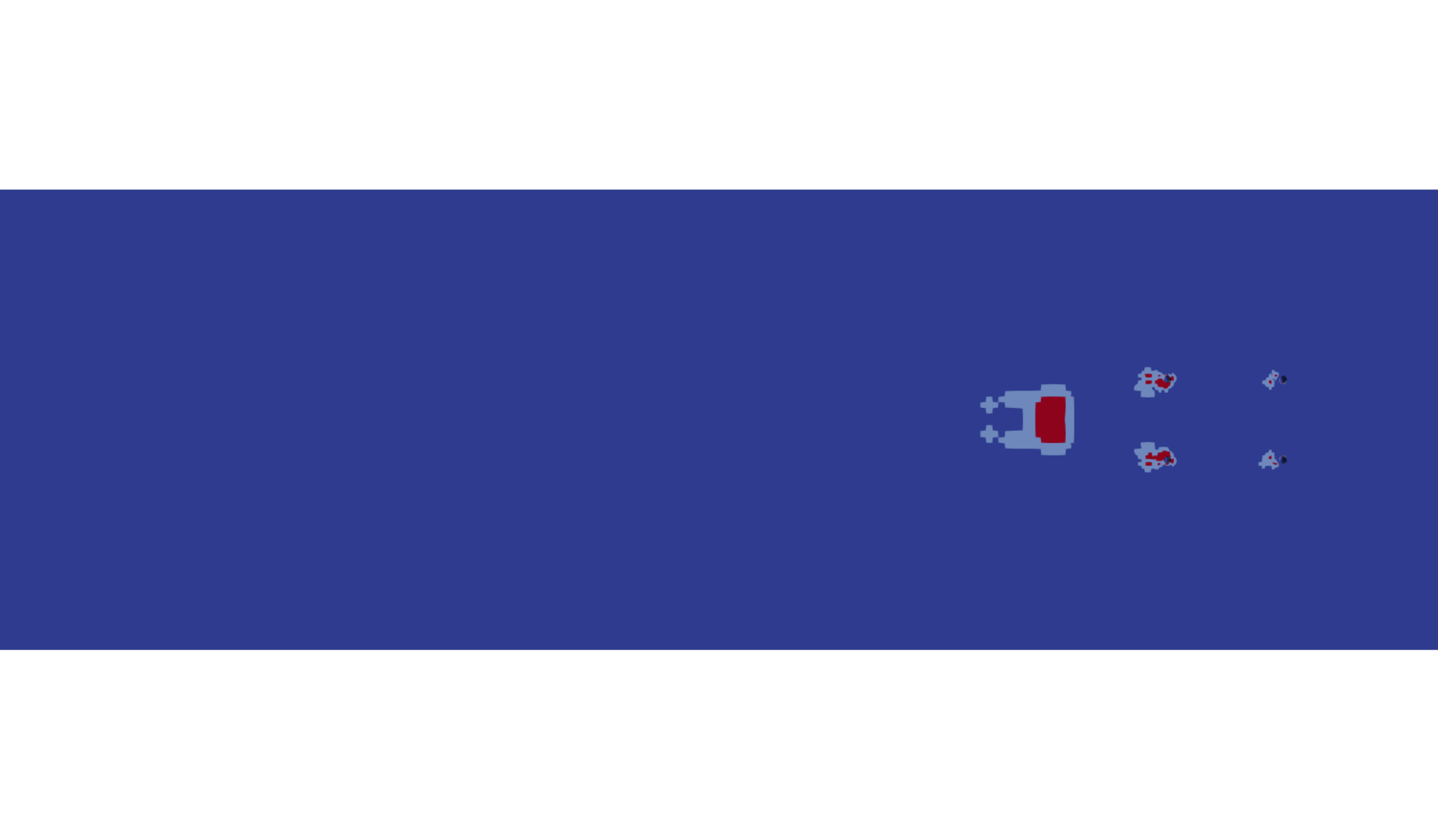}\\
  \includegraphics[width=0.49\textwidth, trim={0 200 0 200}, clip]{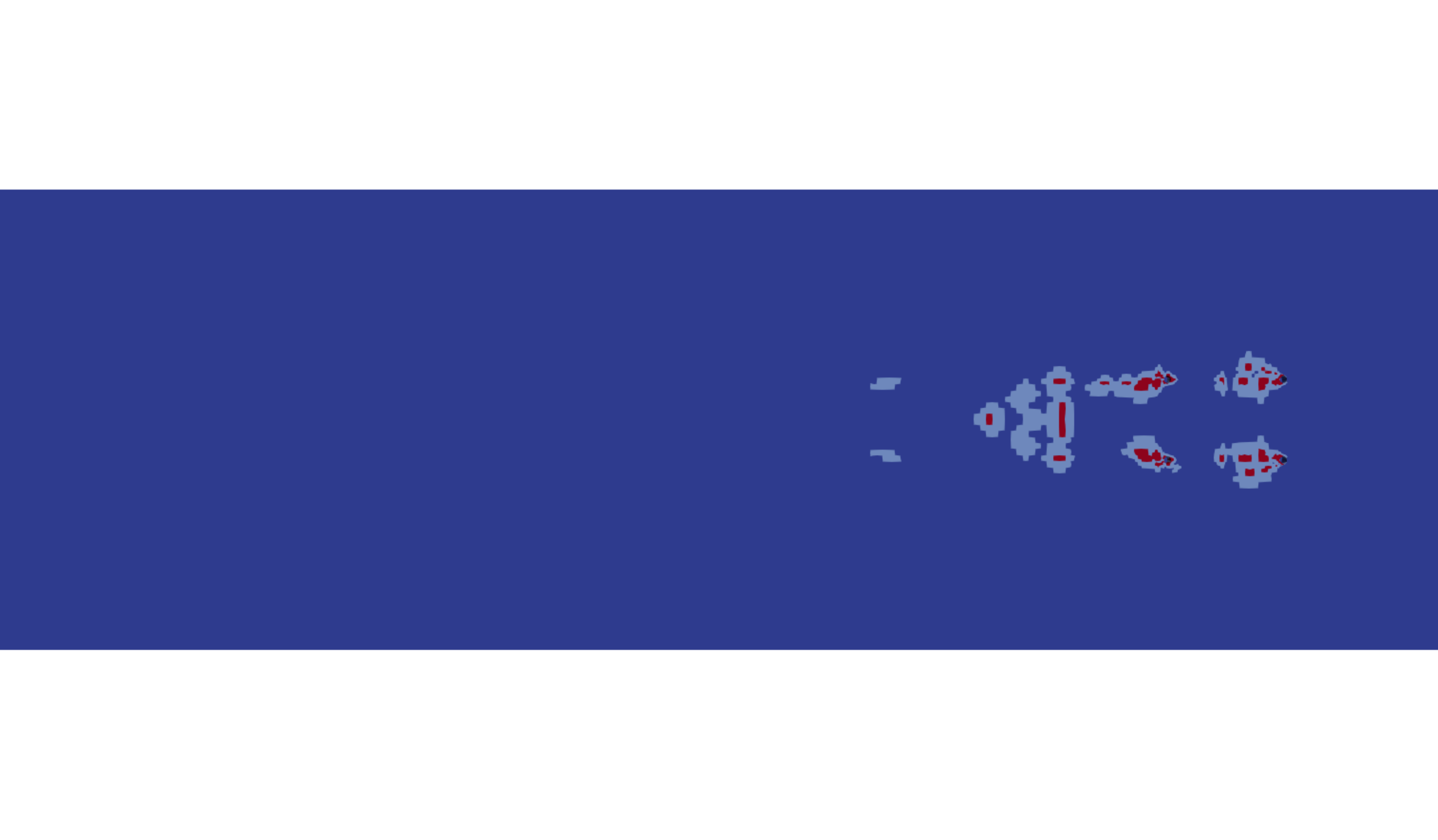}
  \includegraphics[width=0.49\textwidth, trim={0 200 0 200}, clip]{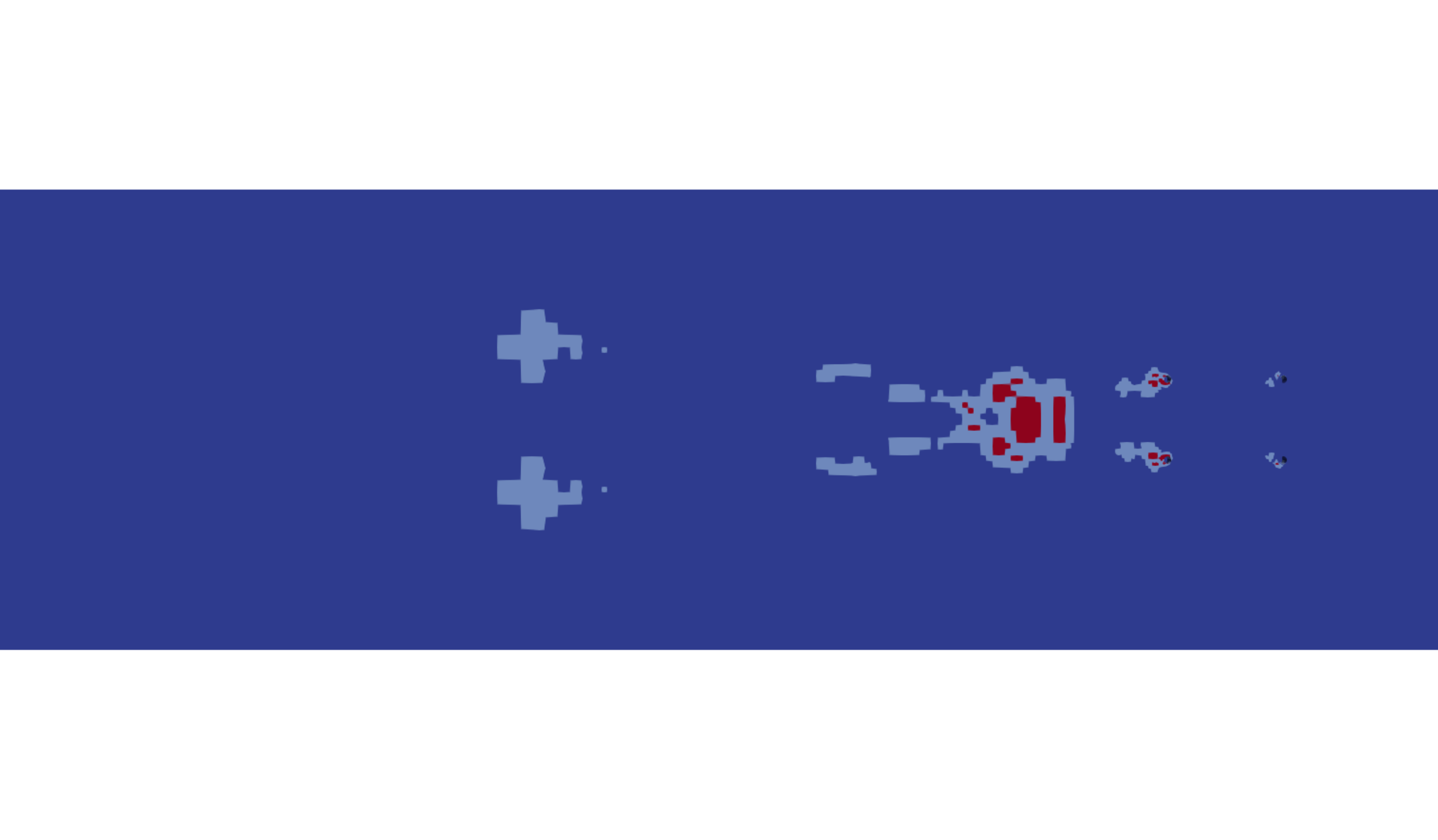}
  \includegraphics[width=0.49\textwidth, trim={0 270 0 200}, clip]{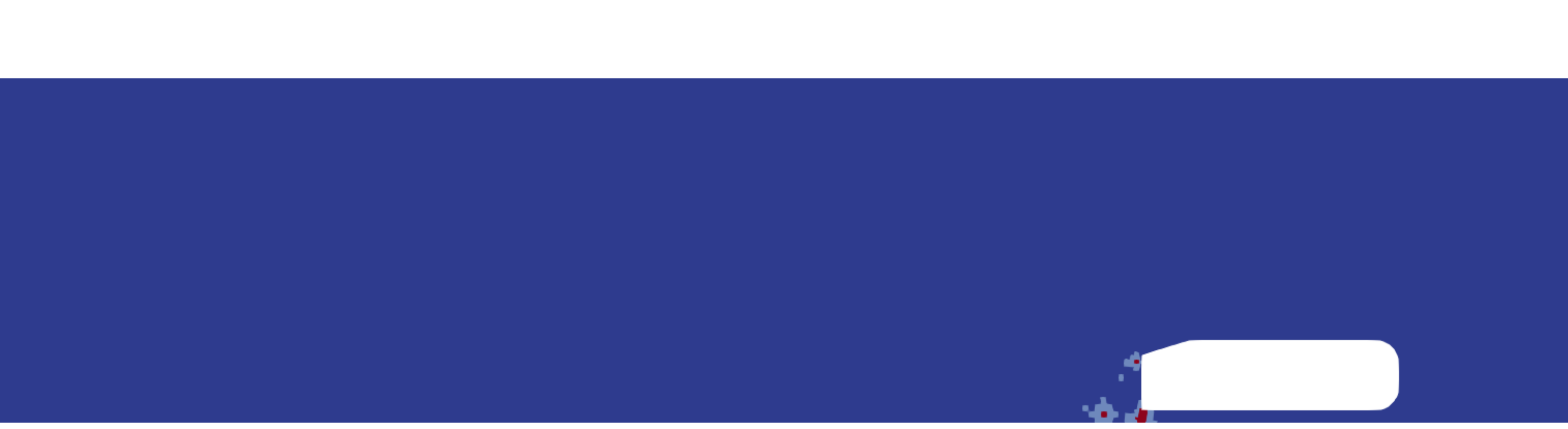}
  \includegraphics[width=0.49\textwidth, trim={0 270 0 200}, clip]{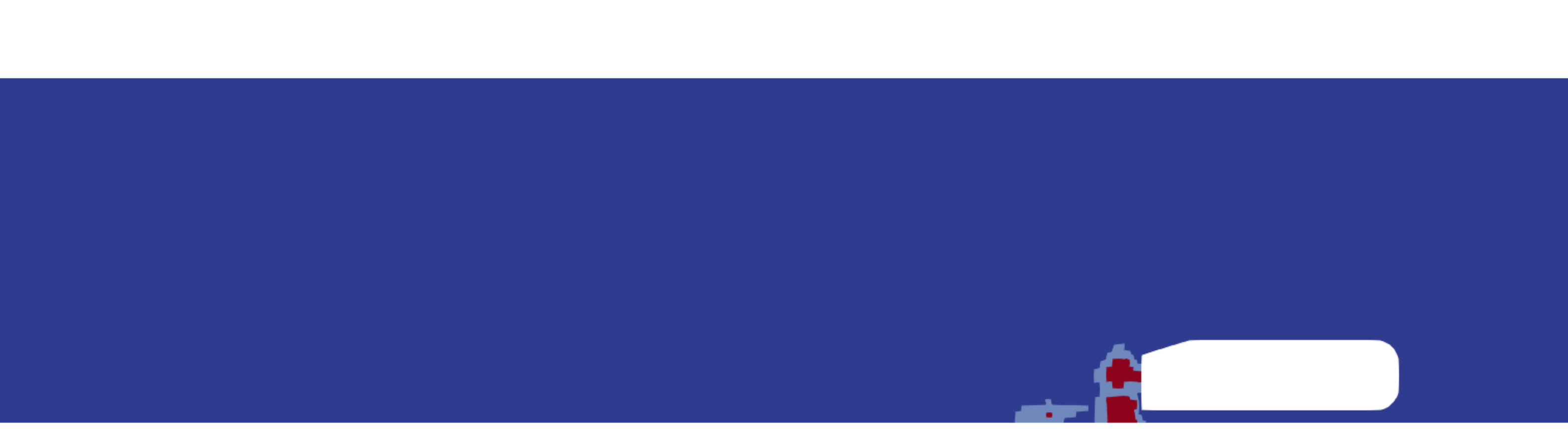}\\
  \includegraphics[width=0.49\textwidth, trim={0 250 0 300}, clip]{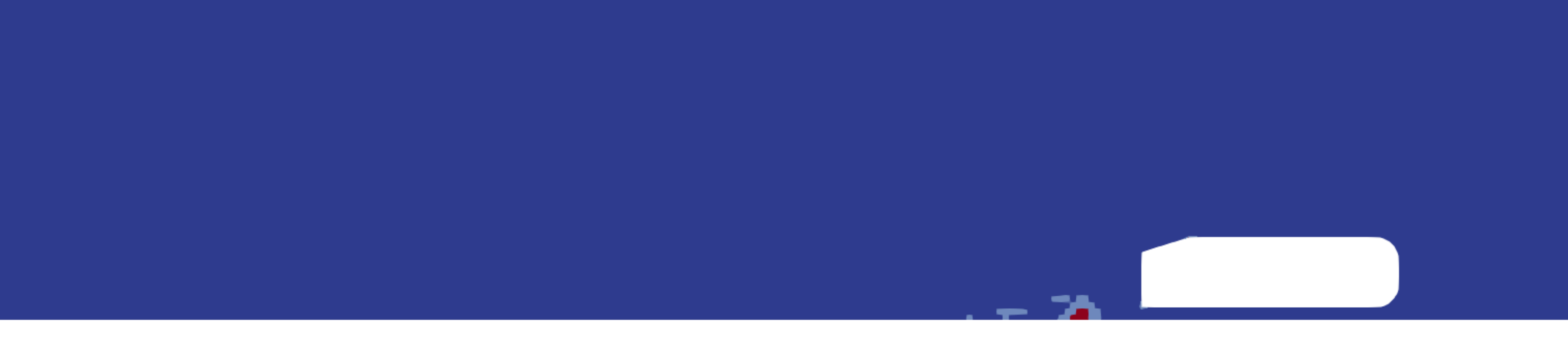}
  \includegraphics[width=0.49\textwidth, trim={0 250 0 300}, clip]{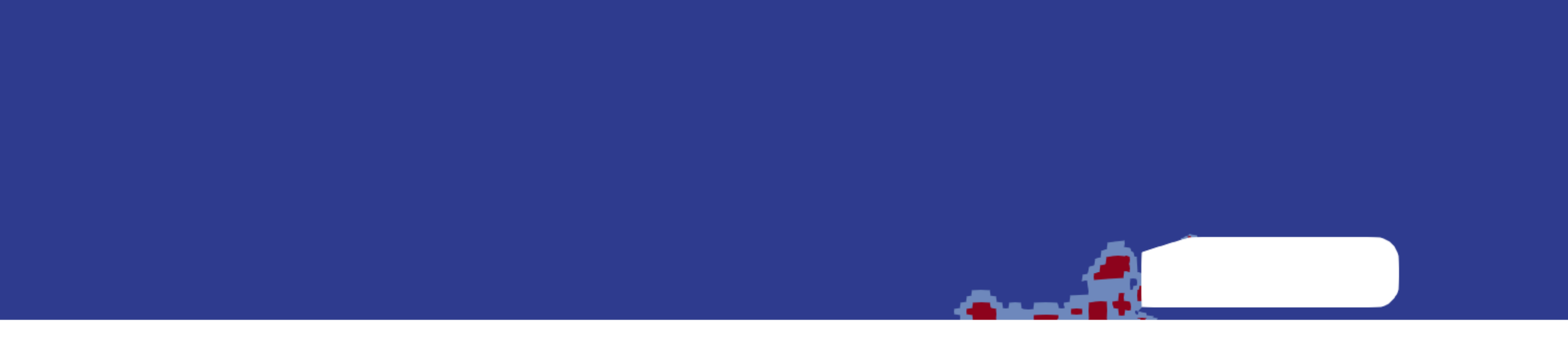}
  \caption{\textbf{INS1.} Adaptive collocation nodes for the test case described in~\cref{subsubec:smallAhmed} with corresponding test parameter number \textbf{1} at the time instants $t\in\{0.025\text{s},0.05\text{s}, 0.075\text{s}, 0.1\text{s}\}$. The number of collocation nodes is $r_h=\textbf{2000}$.The number of cells is \textbf{198633} the total number of degrees of freedom is \textbf{993165}. \textbf{Above 4 slices:} bottom view of the Ahmed body.  \textbf{Below 4 slices:} lateral view of the Ahmed body.}
  \label{fig:adaptiveSmall}
\end{figure}

Some snapshots associated to the test cases studedied are reported: Figure~\ref{fig:snapsCoarseAirfoil} refers to the \textbf{CNS1} test case in~\cref{subsec:coarseAirfoil}, Figure~\ref{fig:snapsFinerAirfoil} refers to the \textbf{CNS2} test case in~\cref{subsubsec:sudomains_finer_airfoil} and finally Figures~\ref{fig:ahmedSmallSnaps} and~\ref{fig:adaptiveSmall} refers to the test case \textbf{INS1} in~\cref{subsubec:smallAhmed}.

\clearpage

\bibliographystyle{abbrv}
\bibliography{biblio}


\end{document}